\documentclass[12pt]{article}

\usepackage{indentfirst}
\usepackage{amsfonts}
\usepackage{amsmath}
\usepackage{amssymb}
\usepackage{amsbsy, amsthm}
\usepackage[margin=0.8in]{geometry}
\usepackage{imakeidx,xcolor}

\newtheorem{dfn}{Definition}[section]
\newtheorem{theorem}[dfn]{Theorem}
\newtheorem{lemma}[dfn]{Lemma}

\newtheorem{corollary}[dfn]{Corollary}
\newtheorem{question}[dfn]{Question}
\newtheorem{conjecture}[dfn]{Conjecture}
\newenvironment{pf}{\noindent{\bf Proof.}}
{\enspace\vrule height5pt depth0pt width5pt} 

\newcommand{\defn}[1]{{\emph{#1}}{\index{#1}}}

\def\F {{\mathcal F}}

\def\P {{\mathcal P}}
\def\X {{\mathcal X}}
\def\W {{\mathcal W}}
\def\L {{\mathcal L}}
\def\U {{\mathcal U}}

\def\ad {{\rm asdim}}
\def\dist {{\rm dist}}
\def\leng {{\rm leng}}

\makeindex

\begin{document}

\title{Assouad-Nagata dimension of minor-closed metrics}
\author{Chun-Hung Liu\thanks{chliu@tamu.edu. Partially supported by NSF under award DMS-1954054 and CAREER award DMS-2144042.} \\
\small Department of Mathematics, \\
\small Texas A\&M University,\\
\small College Station, TX 77843-3368, USA \\
}
%\author{Chun-Hung Liu}
%\address{Department of Mathematics, Texas A\&M University, College Station, TX 77843-3368, USA.}
%\thanks{Email: chliu@tamu.edu. Partially supported by NSF under award DMS-1954054 and CAREER award DMS-2144042.}

\maketitle

\begin{abstract}
Assouad-Nagata dimension addresses both large and small scale behaviors of metric spaces and is a refinement of Gromov's asymptotic dimension.
A metric space $M$ is a minor-closed metric if there exists an (edge-)weighted graph $G$ satisfying a fixed minor-closed property such that the underlying space of $M$ is the vertex-set of $G$, and the metric of $M$ is the distance function in $G$. 
Minor-closed metrics naturally arise when removing redundant edges of the underlying graphs by using edge-deletion and edge-contraction.
In this paper, we determine the Assouad-Nagata dimension of every minor-closed metric.
Our main theorem simultaneously generalizes known results about the asymptotic dimension of $H$-minor free unweighted graphs and about the Assouad-Nagata dimension of complete Riemannian surfaces with finite Euler genus [Bonamy et al., {\it Asymptotic dimension of minor-closed families and Assouad–Nagata dimension of surfaces}, JEMS (2024)].
\end{abstract}


\section{Introduction} \label{sec:intro}

Asymptotic dimension of metric spaces is an important notion in geometric group theory introduced by Gromov \cite{g}.
(See \cite{bd} for a survey.)
For example, it is known that every metric space of finite asymptotic dimension admits a coarse embedding into a Hilbert space (see \cite{bd,r_lectures}); Yu \cite{y_novikov} proved that finitely generated groups with the word metric that admit coarse embeddings into Hilbert spaces satisfy the Novikov higher signature conjecture and hence proved this conjecture for manifolds whose fundamental group has finite asymptotic dimension.

There are a number of equivalent definitions for asymptotic dimension.
Here we use the following version.
For a nonnegative integer $n$, the \defn{$n$-dimensional control function of a pseudometric space} $(X,d)$ is a function $f: {\mathbb R}^+ \rightarrow {\mathbb R}^+$ such that for every positive real number $r$, there exist collections $\U_1,\U_2,...,\U_{n+1}$ of subsets of $X$ such that 
	\begin{itemize}
		\item $\bigcup_{i=1}^{n+1}\bigcup_{U \in \U_i}U \supseteq X$,
		\item for each $i \in [n+1]$, if $U,U'$ are distinct elements of $\U_i$, then $d(x,x')>r$ for every $x \in U$ and $x' \in U'$, and
		\item for each $i \in [n+1]$ and $U \in \U_i$, if $x,x' \in U$, then $d(x,x') \leq f(r)$.
	\end{itemize}
The \defn{asymptotic dimension of a pseudometric space} $(X,d)$ is the infimum integer $n$ such that there exists an $n$-dimensional control function.
Asymptotic dimension is an invariant under coarse equivalence\footnote{See Section \ref{subsubsec:continuous_space}.}.

Assouad-Nagata dimension is a refinement of asymptotic dimension, introduced by Assouad \cite{a_nagata0} under the name Nagata dimension.
We say that a function $f: {\mathbb R}^+ \rightarrow {\mathbb R}^+$ is a \defn{dilation} if there exists a constant $c>0$ such that $f(x) \leq cx$ for every $x \in {\mathbb R}^+$.
The \defn{Assouad-Nagata dimension of a pseudometric space} $(X,d)$ is the infimum integer $n$ such that some dilation is an $n$-dimensional control function of $(X,d)$.
Assouad-Nagata dimension is an invariant under bi-Lipschitz equivalence and is upper bounded by asymptotic dimension.

Assouad-Nagata dimension not only looks for slow growth rate (i.e.\ linear) of the control function to give information for the large scale, but also requires that the control function $f(x)$ goes to 0 when $x$ tends to 0 to give information for the small scale.
See, for example, \cite{bdhm_nagata,ls_nagata} for properties of Assouad-Nagata dimension. 

The metric spaces studied in this paper are defined by weighted graphs.
A \defn{weighted graph} is a pair $(G,\phi)$ such that $G$ is a graph\footnote{In this paper, graphs are loopless and for any two distinct vertices, there are only finitely many parallel edges between them. That is, every graph $G$ consists of a set $V(G)$ of vertices and a multiset $E(G)$ of 2-element subsets of $V(G)$ such that every 2-element subset appears zero or finitely many times in $E(G)$.} and $\phi: E(G) \rightarrow {\mathbb R}^+$ is a function.
For every edge $e$ of $G$, we call $\phi(e)$ the \defn{weight} or the \defn{length} of $e$.
Note that every (unweighted) graph can be viewed as a weighted graph where each edge has weight 1.

Given a weighted graph $(G,\phi)$ and vertices $x,y \in V(G)$, we define the \defn{distance in $(G,\phi)$ between $x$ and $y$}, denoted by \defn{$\dist_{(G,\phi)}(x,y)$}, to be the infimum\footnote{Note that we do not assume that the infimum is attained by a path.} $\sum_{e \in E(P)}\phi(e)$ over all paths $P$ in $G$ between $x$ and $y$. 
Note that $(V(G),\dist_{(G,\phi)})$ is a pseudometric space, and we call it \defn{the pseudometric space generated by $(G,\phi)$}; when $G$ is locally finite, $(V(G),\dist_{(G,\phi)})$ is a metric space.
We remark that every metric space can be generated by a weighted complete graph.

For a weighted graph $(G,\phi)$, we define the \defn{asymptotic dimension} and the \defn{Assouad-Nagata dimension} of $(G,\phi)$ to be the asymptotic dimension and the Assouad-Nagata dimension of the pseudometric space generated by $(G,\phi)$, respectively.

In fact, asymptotic dimension and Assouad-Nagata dimension of weighted and unweighted graphs have been extensively studied.
The pseudometric space generated by the Cayley graph of a group $\Gamma$ with a symmetric generating set is exactly the word metric.
Hence the asymptotic dimension of a finitely generated group (with word metric) is exactly the asymptotic dimension of its Cayley graph.
In addition, an algorithmic version of the Assouad-Nagata dimension is called a $(\sigma,\tau)$-weak sparse partition scheme studied in computer science (see \cite{ap_weak_part_scheme,jlnrs_weak_part_scheme}), where $\sigma$ and $\tau$ are real numbers.
It was shown in \cite{jlnrs_weak_part_scheme} that if a weighted graph on $n$ vertices admits a $(\sigma,\tau)$-weak sparse partition scheme, then one can compute certain universal structures, such as a tour for the traveling salesman problem with stretch $O(\sigma\tau\log_\sigma n)$ and a universal spanning tree with stretch $O(\sigma\tau\log_\sigma n)$, in polynomial time.
We omit the formal definitions of the above notions, but we remark that by a result of Gromov \cite{g}, if a weighted graph has a $d$-dimensional control function which is a dilation whose coefficient is independent from $\lvert V(G) \rvert$, such that the corresponding collections $\U_1,\U_2,...,\U_{d+1}$ can be computed in polynomial time for every positive real number $r$, then this weighted graph admits an $(N,d+1)$-weak sparse partition scheme, where $N$ is a constant independent from the weighted graph.

As we will see later in this section, it is more beneficial to consider the more general setting for the asymptotic dimension or Assouad-Nagata dimension of a class of pseudometric spaces.
The \defn{asymptotic dimension of a class} $\F$ of pseudometric spaces or weighted graphs, denoted by $\ad(\F)$, is the infimum integer $n$ such that there exists a common $n$-dimensional control function for all pseudometric spaces or weighted graphs in $\F$.
The \defn{Assouad-Nagata dimension of a class} $\F$ of pseudometric spaces or weighted graphs is the infimum integer $n$ such that there exists a dilation that is a common $n$-dimensional control function for all pseudometric spaces or weighted graphs in $\F$.

\subsection{Minor-closed metrics}

Recall that every metric space can be generated by a weighted complete graph.
It is therefore more interesting if the weighted graph is required to satisfy certain properties or admits certain structural restriction. 
The restriction considered in this paper is motivated by two elementary graph operations that can be used to simplify the graphs that generate the metric spaces. 

If some edge $e$ in a weighted graph $(G,\phi)$ has very large weight so that it does not belong to any path that defines the distance between two vertices, then the pseudometric space defined by the weighted graph remains the same after we delete $e$. 
So we expect that the properties that our weighted graphs are required to satisfy are preserved under edge-deletion.

In addition, one may consider putting weight 0 on edges instead of putting positive weights on locally finite weighted graphs.
If we allow edges with weight 0, then the distance between some vertices can be 0, so the locally finite weighted graph only defines a pseudometric space.
But if we contract all edges with weight 0 and the edges of weight 0 do not form an infinite component, then the resulting weighted graph is locally finite and defines a metric space that is essentially the same as the original pseudometric space. 
Moreover, if there exists a connected subgraph whose diameter is tiny, then we may contract the connected subgraph by repeatedly contracting edges to derive a new weighted graph, and the pseudometric space generated by the new weighted graph is quasi-isometric to the pseudometric space generated by the original weighted graph.
So we also expect that the properties that our weighted graphs are required to satisfy are preserved under edge-contraction.

The above observation motivates our study of minor-closed metrics.
We say that a graph $H$ is a \defn{minor} of a graph $G$ if $H$ is isomorphic to a graph that can be obtained from $G$ by repeatedly deleting vertices and edges and contracting edges.
We also say that a weighted graph $H$ is a \defn{minor} of a weighted graph $G$ if the underlying graph of $H$ is a minor of the underlying graph of $G$.
In this case, we also say that \defn{$G$ contains $H$ as a minor}.
A class $\F$ of (weighted or unweighted) graphs is \defn{minor-closed} if any minor of a member of $\F$ is also a member of $\F$.
Typical examples of minor-closed families include the class of planar graphs and the more general classes such as the class of graphs embeddable in a fixed surface, the class of linkless embeddable graphs, and the class of knotless embeddable graphs.

In this paper we consider the Assouad-Nagata dimension of a family of pseuometric spaces that are defined by weighted graphs belonging to a fixed minor-closed family $\F$.
The case that the minor-closed family $\F$ contains all finite graphs is not interesting.
In this case we can find weighted graphs $G$ in $\F$ containing arbitrarily large graphs as minors; in particular, they contain arbitrarily large expanders as minors.
So we can create metric spaces that are isometric to or quasi-isometric to (unweighted) expanders by assigning weight 0 or a very small weight to each edge that we would like to contract, assigning a very large weight to each edge that we would like to delete, and assigning weight 1 to the remaining edges.
Every infinite family of bounded degree expanders has infinite asymptotic dimension \cite{g} and hence has infinite Assouad-Nagata dimension.

Note that for every graph property that is preserved under vertex-deletion, edge-deletion and edge-contraction, there exists a minor-closed family consisting of the graphs satisfying this property.
In this paper we consider the following question.

\begin{question} \label{question_AN_minor}
Given a minor-closed family $\F$ of weighted graphs, what is the Assouad-Nagata dimension of $\F$? 
\end{question}

The special case of Question \ref{question_AN_minor} for unweighted graphs and asymptotic dimension was proposed by Fujiwara and Papasoglu \cite[Question 5.2]{fp}, solved in \cite{l_unweighted} and later published in its journal version \cite{bbeglps_merged}. 
This special case proved in \cite{bbeglps_merged,l_unweighted} not only improves an earlier result of Ostrovskii and Rosenthal \cite{or} who gave a bound depending on the number of vertices of the minimal graphs not in $\F$, but also answers a number of questions of Bonamy, Bousquet, Esperet, Groenland, Pirot and Scott \cite{bbegps} and a question of Ostrovskii and Rosenthal \cite{or}.
In addition, this special case and other related results in \cite{bbeglps_merged,l_unweighted} on asymptotic dimension imply many results about clustered coloring of graphs of bounded maximum degree in the literature \cite{adov,lo,lw,lw_simple}.

Another special case of Question \ref{question_AN_minor} is for the class of weighted $K_{3,p}$-minor free graphs for any fixed integer $p$, which was solved in \cite{bbegps} and later published in \cite{bbeglps_merged}. 
This result was used in \cite{bbeglps_merged} to prove that the Assouad-Nagata dimension of any complete Riemannian surface\footnote{A \defn{surface} is a 2-dimensional connected manifold without boundary.} of finite Euler genus is at most 2.

Some natural lower bounds for Question \ref{question_AN_minor} are known.
For example, the class of 2-dimensional grids has Assouad-Nagata dimension at least 2 \cite{g} and the class of paths has Assouad-Nagata dimension at least 1.

The main result of this paper is the complete solution to Question \ref{question_AN_minor}, showing that the above lower bounds are tight and being a common strengthening of the aforementioned results in \cite{bbeglps_merged,l_unweighted,or}.

\begin{theorem} \label{minor_AN_chara_intro}
Let $\F$ be a minor-closed family of weighted graphs.
Then the following statements hold.
	\begin{enumerate}
		\item The Assouad-Nagata dimension of $\F$ is infinite if and only if $\F$ contains all finite graphs.
		\item The Assouad-Nagata dimension of $\F$ is at most 2 if and only if some finite graph is not in $\F$. 
		\item The Assouad-Nagata dimension of $\F$ is at most 1 if and only if some finite planar graph is not in $\F$. 
		\item The Assouad-Nagata dimension of $\F$ is 0 if and only if some finite path is not in $\F$. 
	\end{enumerate}
\end{theorem}

Statement 4 and the ``if'' part of Statement 1 in \ref{minor_AN_chara_intro} are in fact trivial.
The other part of Statement 1 follows from Statement 2.
So Theorem \ref{minor_AN_chara_intro} is equivalent to the following two theorems, by combining the lower bound given by 2-dimensional grids mentioned above. 

\begin{theorem} \label{minor_AN_intro}
For any finite graph $H$, the Assouad-Nagata dimension of the class of $H$-minor free weighted graphs is at most 2.
\end{theorem}

\begin{theorem} \label{minor_planar_AN_intro}
For any finite planar graph $H$, the Assouad-Nagata dimension of the class of $H$-minor free weighted graphs is at most 1.
\end{theorem}

\subsection{Applications}

\subsubsection{Groups}
We can apply our results to study Assouad-Nagata dimension of groups.
Given a group $\Gamma$ and a symmetric generating set $S$, the \defn{Cayley graph} ${\rm Cay}(\Gamma,S)$ is the graph with vertex-set $\Gamma$ such that two vertices $u,v \in \Gamma$ are adjacent if and only if $u=vs$ for some $s \in S$.
Gromov \cite{g} observed that the asymptotic dimension of ${\rm Cay}(\Gamma,S)$ is independent from the choice of the finite generating set $S$ when $\Gamma$ is finitely generated.
So the asymptotic dimension of a finitely generated group is defined to be the asymptotic dimension of a Cayley graph and hence is a group invariant.

One particular question in this direction was asked by Ostrovskii and Rosenthal \cite{or} about the asymptotic dimension of minor-excluded groups.
This question was solved in \cite{bbeglps_merged,bbegps,l_unweighted}.
Theorem \ref{minor_AN_intro} implies a stronger answer that the same result actually holds for Assouad-Nagata dimension.

\begin{corollary}
For every group $\Gamma$ with a symmetric (not necessarily finite) generating set $S$, if there exists a finite graph $H$ such that $H$ is not a minor of the Cayley graph ${\rm Cay}(\Gamma,S)$, then the Assouad-Nagata dimension of ${\rm Cay}(\Gamma,S)$ is at most 2.
\end{corollary}

\subsubsection{Continuous spaces} \label{subsubsec:continuous_space}

Though Theorems \ref{minor_AN_intro} and \ref{minor_planar_AN_intro} are useful for discrete metric spaces, they lead to easy applications to continuous metric spaces and show that they are significantly more powerful than the version for unweighted graphs in \cite{bbeglps_merged,l_unweighted}.

We say that a metric space $(X,d_X)$ is \defn{coarsely embeddable} in a metric space $(Y,d_Y)$ if there exist a function $\iota: X \rightarrow Y$ and increasing functions $f: {\mathbb R}^+ \rightarrow {\mathbb R}^+$ and $g: {\mathbb R}^+ \rightarrow {\mathbb R}^+$ with $\lim_{x \rightarrow \infty}f(x)=\infty$ such that for any $x_1,x_2 \in X$, $f(d_X(x_1,x_2)) \leq d_Y(\iota(x_1),\iota(x_2)) \leq g(d_X(x_1,x_2))$; if in addition there exists a positive real number $C$ such that for every $y \in Y$, there exists $x_y \in X$ such that $d_Y(y,\iota(x_y)) \leq C$, then we say that $(X,d_X)$ is \defn{coarsely equivalent} to $(Y,d_Y)$.
It is not hard to see that if $(X,d_X)$ is coarsely equivalent to $(Y,d_Y)$, then $(Y,d_Y)$ is coarsely equivalent to $(X,d_X)$.
It is well-known that the asymptotic dimension is invariant under coarsely equivalence. 
Similarly, we say that a class $\F_1$ of metric spaces is \defn{coarsely equivalent} to a class $\F_2$ of metric spaces if each member of $\F_1$ is coarsely equivalent to a member of $\F_2$ by using the same functions $f$ and $g$ and constant $C$ mentioned above.
It is straightforward to show that if $\F_1$ is coarsely equivalent to $\F_2$, then the asymptotic dimension of $\F_1$ equals the asymptotic dimension of $\F_2$.
Moreover, a simple result (\cite[Lemma 7.5]{bbeglps_merged}) states that the Assouad-Nagata dimension of any scaling-closed metric space equals its asymptotic dimension, where a class $\F$ of metric spaces is \defn{scaling closed} if $(X,d_X) \in \F$ implies that $(X,c \cdot d_X) \in \F$ for every positive real number $c$.
Hence we immediately obtain the following corollary.

\begin{corollary} \label{cor_coarse_equiv}
Let $H$ be a finite graph, and let $\F_H$ be the class of metric spaces generated by an $H$-minor free weighted graphs.
Let $\F$ be a class of metric spaces coarsely equivalent to $\F_H$.
Then the following statements hold.
	\begin{enumerate}
		\item The asymptotic dimension of $\F$ is at most 2.
		\item If $H$ is planar, then the asymptotic dimension of $\F$ is at most 1.
		\item If $\F$ is scaling-closed, then the Assouad-Nagata dimension of $\F$ is at most 2.
		\item If $\F$ is scaling-closed and $H$ is planar, then the Assouad-Nagata dimension of $\F$ is at most 1.
	\end{enumerate}
\end{corollary}

Note that a particular use of Corollary \ref{cor_coarse_equiv} is stated in \cite{bbeglps_merged}.
For each surface $\Sigma$ of finite Euler genus, the class $\F$ of all complete Riemannian surfaces whose underlying space is $\Sigma$ is scaling-closed and is coarsely equivalent to the class of the metric spaces generated by weighted graphs embeddable in $\Sigma$, so $\F$  has Assouad-Nagata dimension at most 2.
We remark that the class of graphs embeddable in a fixed surface is a very special case of minor-closed families.
So Corollary \ref{cor_coarse_equiv} is expected to have applications to other spaces with certain geometric or topological feature.
For example, if $H=K_6$ (or $H=K_7$, respectively), then the class $\F_H$ mentioned in Corollary \ref{cor_coarse_equiv} contains all linkless embeddable graphs (or all knotless embeddable graphs, respectively).

Corollary \ref{cor_coarse_equiv} can also be applied to metric graphs, which allow us to consider points on the edges.
For a graph $G$, each edge of $G$ can be viewed as a topological space that is homeomorphic to $[0,1]$, so $G$ can be viewed as a 1-dimensional complex which is the union of all edges.
That is, for each edge $e=xy$ of $G$, if a homeomorphism $f_e$ from $e$ to $[0,1]$ with $f_e(x)=0$ and $f_e(y)=1$ is given, then for any point $z$ in the edge $e$, there exists $t_z \in [0,1]$ such that $f_e(z)=t_z$. 
As long as the weight $\phi(e)$ of $e$ is given, we can define the distance between $x$ and $z$ to be $f_e(z) \cdot \phi(e)$ and the distance between $y$ and $z$ is $(1-f_e(z)) \cdot \phi(e)$.
This allows us to define a distance function on a space that is the union of all edges of a weighted graph to obtain a metric space.
Such a metric space is called a \defn{metric graph} and is the underlying space of a quantum graph studied in the quantum graph theory (in the sense in \cite{bk_quantum}).

\subsubsection{Algorithmic aspects}

Note that our proof of Theorems \ref{minor_AN_intro} and \ref{minor_planar_AN_intro} is constructive and, as long as a finite $H$-minor free weighted graph $(G,\phi)$ and a real number $r$ are given, provides an algorithm for finding the collections $\U_1,\U_2,...,\U_{n+1}$ of subsets of $V(G)$ witnessing the Assouad-Nagata dimension of the class of $H$-minor free weighted graphs in time polynomial in the size of the graph\footnote{One might think that the running time of the algorithm given by our proof also depends on $\max_{e \in E(G)}\phi(e)$. But as long as $r$ is given, we can first delete edges with weight much bigger than $r$ without affecting the correctness of those covers $\U_1,\U_2,...,\U_{n+1}$, so every remaining edges has bounded weight.}, where the exponent of the polynomial is independent from $G$.
Hence our proof gives a $(c_H,3)$-weak sparse partition scheme for $H$-minor free weighted graphs, and a $(c_H,2)$-weak partition scheme for such graphs if $H$ is planar, where $c_H$ is a constant only depending on $H$.

\subsubsection{Other formulations and other graph classes}

By a simple compactness argument (or explicitly, by an immediate consequence of the combination of \cite[Theorem A.2]{bbeglps_merged} and \cite[Lemma 7.5]{bbeglps_merged}), we only have to consider classes of finite weighted graphs.
When restricted to finite graphs, Theorem \ref{minor_planar_AN_intro} can be stated in terms of tree-width.

A \defn{tree-decomposition} of a finite graph $G$ is a pair $(T,\X)$ such that $T$ is a tree and $\X$ is a collection $(X_t: t \in V(T))$ of subsets of $V(G)$, called the \defn{bags}, such that
	\begin{itemize}
		\item $\bigcup_{t \in V(T)}X_t = V(G)$,
		\item for every $e \in E(G)$, there exists $t \in V(T)$ such that $X_t$ contains the ends of $e$, and
		\item for every $v \in V(G)$, the set $\{t \in V(T): v \in X_t\}$ induces a connected subgraph of $T$.
	\end{itemize}
For a tree-decomposition $(T,\X)$, the \defn{adhesion} of $(T,\X)$ is $\max_{tt' \in E(T)}\lvert X_t \cap X_{t'} \rvert$, and the \defn{width} of $(T,\X)$ is $\max_{t \in V(T)}\lvert X_t \rvert-1$.
The \defn{tree-width} of $G$ is the minimum width of a tree-decomposition of $G$.
We define the \defn{tree-decomposition} and the \defn{tree-width} of a finite weighted graph to be the tree-decomposition and tree-width of its underlying graph, respectively.

By the Grid Minor Theorem \cite{rs_V}, excluding any planar finite graph as a minor from a finite graph is equivalent to having bounded tree-width.
Hence the following Theorem \ref{tw_AN_intro} is an equivalent form of Theorem \ref{minor_planar_AN_intro}.

\begin{theorem} \label{tw_AN_intro}
For every positive integer $w$, the Assouad-Nagata dimension of the class of finite weighted graphs of tree-width at most $w$ is at most 1.
\end{theorem}

As we mentioned earlier, the special case of Theorem \ref{minor_AN_intro} for asymptotic dimension of unweighted graphs was proved in \cite{l_unweighted}.
A key step in the proof of this special case in \cite{l_unweighted} is proving a result on the class of finite graphs of bounded layered tree-width.
Classes of graphs of bounded layered tree-width are common generalizations of classes of graphs of bounded tree-width and classes of graphs of bounded Euler genus, and are incomparable with minor-closed families.

A \defn{layering} of a graph $G$ is an ordered partition $(V_1,V_2,...)$ of $V(G)$ into (possibly empty) subsets such that for every edge $e$ of $G$, there exists $i_e$ such that $V_{i_e} \cup V_{i_e+1}$ contains both ends of $e$.
We call each $V_i$ in $(V_1,V_2,...)$ a \defn{layer}.
The \defn{layered tree-width} of a finite graph $G$ is the minimum $w$ such that there exist a tree-decomposition of $G$ and a layering of $G$ such that the size of the intersection of any bag and any layer is at most $w$.
Clearly, every graph of tree-width at most $w$ has layered tree-width at most $w+1$.

A number of classes of graphs with some geometric properties have bounded layered tree-width.
For example, \defn{$(g,k)$-planar graphs}, which are the graphs that can be drawn in a surface of Euler genus at most $g$ with at most $k$ crossings on each edge, have layered tree-width at most $(4g+6)(k+1)$ \cite{dew}.
Note that $(g,0)$-planar graphs are exactly the graphs of Euler genus at most $g$.
In addition, it is well-known that $(0,1)$-planar graphs can contain arbitrary graph as a minor, so the class of $(g,k)$-planar graphs is not minor-closed when $k>0$.
We refer readers to \cite{dew,djmnw} for other examples of graphs of bounded layered tree-width. 

It is known that the class of finite graphs of layered tree-width at most $w$, for any fixed integer $w$, has asymptotic dimension at most 2 \cite{bbeglps_merged,l_unweighted} but has infinite Assouad-Nagata dimension even when $w=1$ \cite{bbeglps_merged}.
Recall that the result for bounded layered tree-width graphs is a key step for proving the special case of Theorem \ref{minor_AN_intro} for unweighted graphs and asymptotic dimension.
The fact that the class of weighted graphs of bounded layered tree-width has infinite Assouad-Nagata dimension shows that significantly new ideas besides the ones in \cite{l_unweighted} are required to prove Theorem \ref{minor_AN_intro}.
In fact, our proof of Theorem \ref{minor_AN_intro} does not require anything related to layered tree-width.
Nonetheless, a byproduct (Theorem \ref{layered_tw_AN_intro}) of our proof gives a strengthening of the result about unweighted graphs of layered tree-width in \cite{bbeglps_merged,l_unweighted}.

For a subset $I$ of ${\mathbb R}^+$, an \defn{$I$-bounded weighted graph} is a weighted graph such that the weight of every edge belongs to $I$.
Note that $\{1\}$-bounded weighted graphs are exactly unweighted graphs, and $(0,\infty)$-bounded weighted graphs are exactly weighted graphs.

\begin{theorem} \label{layered_tw_AN_intro}
For every positive integer $w$ and positive real number $\epsilon$, the asymptotic dimension of the class of finite $[\epsilon,\infty)$-bounded weighted graphs of layered tree-width at most $w$ is at most 2.
\end{theorem}

Since 2-dimensional grids have layered tree-width at most 2 and have asymptotic dimension 2, the bound in Theorem \ref{layered_tw_AN_intro} is optimal.

\subsection{Remarks}

As we mentioned earlier, by a simple compactness argument, we only have to consider classes of finite weighted graphs.
So we assume that all graphs are finite in the rest of the paper.

The main body of this draft was written during the pandemic and was sent to a number of researchers at the time when it was written.
A few days before the first version of this paper appeared on arXiv, Distel \cite{d} announced an independent proof of Theorems \ref{minor_AN_intro} and \ref{minor_planar_AN_intro} on arXiv.
His proof and the proof in this paper are similar in the sense that both of them use the strategies in \cite{l_unweighted} (or its journal version \cite{bbeglps_merged}) by developing a machinery for handling tree-decompositions of weighted graphs and for handling near embeddings on surfaces without using layered tree-width (see Section \ref{sec:sketch} for details).
His proof is shorter because he uses results in \cite{bbeglps_merged,bbegps} and handles near embeddings in a way that is different from ours so that his result about tree-decompositions, which looks similar to but seems weaker than the one in this paper, could be sufficient.

Our machinery for tree-decompositions might be applicable to other classes of (weighted) graphs and of independent interests.
In particular, our machinery works for hereditary graph classes, while the one in \cite{d} seems not.
In addition, his proof in \cite{d} uses the result in \cite{bbeglps_merged} about the Assouad-Nagata dimension of weighted graphs embedded in surfaces as a block box.
That result in \cite{bbeglps_merged} was proved by using fat minors, which is different from the proof in \cite{bbeglps_merged,l_unweighted} of the result about asymptotic dimension of unweighted $H$-minor free graphs and relies on \cite{fp}.
Our paper does not require it or any result in \cite{bbeglps_merged,bbegps,fp,l_unweighted} (except a simple observation in \cite{bbeglps_merged} for scaling-closed families and an explicit statement for a standard compactness argument).
Hence, the present paper gives a new proof for the result about the complete Riemannian surfaces without relying on \cite{bbeglps_merged} (except a fairly standard argument about quasi-isometry used in \cite{bbeglps_merged}) and gives a direct and unified proof for all such minor-closed metrics.
Our result about layered tree-width is simply a by-product of the machinery developed along the way toward the proof of Theorems \ref{minor_AN_intro} and \ref{minor_planar_AN_intro} and is not needed to proving Theorems \ref{minor_AN_intro} and \ref{minor_planar_AN_intro}.

\subsection{Proof sketch} \label{sec:sketch}

We sketch the proof of Theorems \ref{minor_AN_chara_intro} in this subsection.
As explained earlier, it suffices to prove Theorems \ref{minor_AN_intro} and \ref{minor_planar_AN_intro}.

The first observation is that the classes mentioned in Theorems \ref{minor_AN_intro} and \ref{minor_planar_AN_intro} are scaling-closed.
(Recall that a class $\F$ of metric spaces is \defn{scaling-closed} if $(X,d_X) \in \F$ implies that $(X,c \cdot d_X) \in \F$ for every positive real number $c$.)
So by the following simple result, it suffices to prove that the asymptotic dimension of $H$-minor free weighted graphs is at most 2 (and at most 1, if $H$ is planar, respectively).

\begin{lemma}[{\cite[Lemma 7.5]{bbeglps_merged}}] \label{scaling_closed_equiv}
Let $n$ be a positive integer.
Let $\F$ be a scaling-closed classes of metric spaces.
If $\F$ has asymptotic dimension at most $n$, then $\F$ has Assouad-Nagata dimension at most $n$.
\end{lemma}

The very high level overview of the proof for bounding asymptotic dimension of weighted $H$-minor free graphs is similar to the version for unweighted graphs in \cite{bbeglps_merged,l_unweighted}.
However, the weighted graph case contains extra difficulties in each step and requires significantly more effort as well as new ideas.

To prove our bound for asymptotic dimension, we prove that the existence of $d$-dimensional control functions is equivalent to the existence of $(d+1)$-colorings of graphs with bounded weak diameter in powers of graphs.
This equivalence for unweighted graph is straightforward from the definition.
For the weighted graphs dealt with in this paper, the equivalence is still quite straightforward, except that we have to deal with powers of weighted graphs a little more carefully.
This equivalence is shown in Section \ref{sec:ad_wd_color}.

Hence the task is reduced to the one for finding, for every positive real number $\ell$, a $(d+1)$-coloring of an arbitrary $H$-minor free weighted graph $(G,\phi)$, where $H$ is a fixed graph and $d=1$ or $2$ depending on whether $H$ is planar or not, such that every monochromatic component in $(G,\phi)^\ell$ has weak diameter at most a real number $N$ that only depends on $H$ and $\ell$.
Colorings and weak diameter are defined in Section \ref{sec:ad_wd_color}.

When $H$ is planar, the Grid Minor Theorem \cite{rs_V} implies that $G$ has bounded tree-width; that is, $G$ has a tree-decomposition such that every bag has bounded size.
When $H$ is non-planar, the structure theorem proved by Robertson and Seymour \cite{rs_XVI} implies that $G$ has a tree-decomposition with bounded adhesion such that every ``torso'', which is the graph that can be slightly modified from the subgraph induced by a bag by adding edges, is ``nearly embeddable'' in a surface of bounded Euler genus.
Here a ``nearly embeddable graph'' is a graph that can be obtained from a graph drawn in a surface (with possible crossings) by adding a bounded number of vertices such that the crossings can be covered by a bounded number of disks in the surface, and the subgraph contained in each disk is ``well-structured''.
Torsos and nearly embeddable graphs are formally defined in Section \ref{sec:near_embedding_geo}.

Hence the remaining proof is divided into two main steps.
The first step is to show that if $G$ has a tree-decomposition with bounded adhesion such that the torsos have desired colorings, then $G$ has a desired coloring.
(A desired coloring mentioned here refers to a coloring with small weak diameter using the correct number of colors.)
The second step is to show that graphs with bounded size have a desired 2-coloring and every nearly embeddable graph has a desired 3-coloring.

\subsubsection{Handling tree-decompositions}
We first explain how to deal with the first step mentioned above.
Note that torsos are not necessarily subgraphs of the original graph $G$ and usually do not belong to the same graph class that $G$ belongs to.
So even though developing tools to solve this step by reducing problems on $G$ to problems on its torsos is sufficient as long as we can solve the second step, the generality of these tools are somehow limited and probably not easily exploitable for future applications. 

The actual first step we prove in this paper is that if $G$ has a tree-decomposition with bounded adhesion such that each ``bag'' has a desired coloring, then $G$ has a desired coloring. 
(See Lemma \ref{weighted_tree_extension_clean} for a formal description.)
Note that our tools for reducing problems to ``bags'' instead of to torsos provide more generality for potential future applications on graphs in hereditary classes, which are more general than minor-closed families.
But we should point out that the ``bags'' mentioned here are not just subgraphs induced by bags and we will explain the details soon.

Section \ref{sec:bdd_adhesion} is dedicated to a proof of Lemma \ref{weighted_tree_extension_clean}.
To prove Lemma \ref{weighted_tree_extension_clean}, we prove a stronger statement that any coloring $c_Z$ of a subset of vertices $Z$ consisting of vertices that are not far away from the root bag can be extended to a desired coloring for the whole graph $G$ (see Lemma \ref{weighted_tree_extension}).
This stronger statement gives us stronger inductive hypothesis, where its necessity will be clear later.
Note that we can decompose the whole graph $G$ into two parts: one part is the central part, consisting of the bags intersecting $Z$; the other part is the peripheral part, consisting of the bags disjoint from $Z$.
We will show that each part has a desired coloring such that those colorings can be combined into a desired coloring of the whole graph $G$.

To prove that the central part has a desired coloring, we delete $Z$ to obtain a graph $G_1$ with a tree-decomposition with smaller adhesion so that we can obtain a desired coloring $c_1$ for $G_1$ by the induction hypothesis, and then show that the union of $c_1$ and $c_Z$ is a desired coloring of the central part.
Showing that $c_1 \cup c_Z$ is a desired coloring is simple and is proved in Section \ref{sec:center}.
The main defect of this argument about coloring the central part is that even though we obtain a desired coloring for the central part, it is not necessarily extendable to a desired coloring for the whole graph $G$ because the metric of the central part can be very different from the metric on those vertices in the whole graph.
Therefore, in order to make the coloring of the central part extendable to the whole graph, we will add ``gadgets'' to the central part to approximate the original metric when we color the central part.
Hence, in Lemma \ref{weighted_tree_extension}, what we really need to assume and prove is that the central part together with the gadgets have a desired coloring, and the ``bags'' mentioned above refers to the union of this central part and the gadgets.

Once the central part with gadgets is colored, there is a natural way to color the ``boundary'' of the peripheral part according to the colors on the gadgets so that any further coloring on the peripheral part does not create a monochromatic component in the whole graph intersecting the central part beyond the control the original coloring of the central part.
We then apply induction to the peripheral part, where the coloring on the boundary of the peripheral part serves the coloring $c_Z$ of the precolored set $Z$ for the peripheral part.
It is the main motivation for considering the stronger statement that allows a precolored set for the whole graph at the beginning, and we can further obtain a desired coloring of the peripheral part by applying the induction hypothesis.
Note that the way we precolor the boundary corresponds to the coloring of gadgets for the central part, so it can be shown that the union of the colorings of the central part and peripheral part gives a desired coloring of the whole graph with small weak diameter.

We remark that the strategy for proving Lemma \ref{weighted_tree_extension_clean} described above is similar to the proof of a similar result in \cite{bbeglps_merged,l_unweighted} for unweighted graphs.
The main difference and challenge here is that the objects we deal with in this paper are weighted graphs, so the corresponding gadgets are more complicated than the gadgets used for unweighted graphs in \cite{bbeglps_merged,l_unweighted}.
The gadgets that we will use are descried in Section \ref{sec:hierarchy}, and the central part with gadgets is dealt with in Section \ref{sec:condesation}.

\subsubsection{Handling near embeddings}
Now we explain how to handle the second step.
That is, we should show that weighted graphs with bounded size have a 2-coloring with bounded weak diameter and nearly embeddable weighted graphs have a 3-coloring with bounded weak diameter.
In fact, in order to combine with the tool proved in the first step, what we really have to show is that the aforementioned graphs together with the gadgets have a desired coloring.
Fortunately, as the gadgets we use are ``small'' and ``well-structured'', it is not hard to show that the desired colorings of weighted graphs with bounded size and of nearly embeddable weighted graphs can be extended to include those gadgets.
So we only sketch how to handle weighted graphs with bounded size and nearly embeddable weighted graphs.

Note that coloring weighted graphs with bounded size is trivial because every coloring is desired.
So Theorem \ref{minor_planar_AN_intro} easily follows from the tool developed in the first step.
Then our result for graphs of bounded layered tree-width (Theorem \ref{layered_tw_AN_intro}) easily follows from a combination of Theorem \ref{minor_planar_AN_intro} and a Hurewicz type result (Theorem \ref{large_scale_metric_projection}) proved in \cite{bdlm}.

So the difficulty lies in coloring nearly embeddable weighted graphs.
This difficulty does not exist in the unweighted graph version in \cite{bbeglps_merged,l_unweighted}: nearly embeddable unweighted graphs have bounded layered tree-width, so they can be colored by the result for unweighted graphs with bounded tree-width and the aforementioned Hurewicz type result.
In contrast, as shown in \cite{bbeglps_merged}, weighted graphs with bounded layered tree-width have infinite asymptotic dimension.
So new ideas for dealing with nearly embeddable weighted graphs are required.

We will show that every nearly embeddable weighted graph has a tree-decomposition such that every bag is not far away from a union of a bounded number of geodesics in Section \ref{sec:near_embedding_geo}.
Then by a Hurewicz type result proved in \cite{bdlm} (Theorem \ref{large_scale_metric_projection}), it suffices to show that every weighted graph that has a tree-decomposition for which every bag is a union of a bounded number sets with bounded radius has a 2-coloring with bounded weak diameter.
This 2-coloring result (Lemma \ref{strong_weighted_tree_extension_control_clean}) is proved based on a similar strategy to the proof of Lemma \ref{weighted_tree_extension_clean}, but there are significant differences.
One reason for this is that we have to deal with tree-decompositions with unbounded adhesion here, so we cannot do induction on the adhesion like in Lemma \ref{weighted_tree_extension_clean}; the other reason is that we cannot afford to delete vertices like in the proof of Lemma \ref{weighted_tree_extension_clean} because deleting vertices changes the metric dramatically.
Instead, the way we prove Lemma \ref{strong_weighted_tree_extension_control_clean} is to consider a subset $R$ of vertices that are not required to be colored, and this set $R$ plays the role like the set of deleted vertices.
Then Lemma \ref{strong_weighted_tree_extension_control_clean} is proved by considering a technical form (Lemma \ref{strong_weighted_tree_extension_control}) that allows us to apply stronger induction hypotheses consistent with the setup for the set $R$.

Note that Lemma \ref{strong_weighted_tree_extension_control} also implies Theorem \ref{minor_planar_AN_intro}.
So in fact, this paper includes two proofs for Theorem \ref{minor_planar_AN_intro} (and hence for Theorem \ref{layered_tw_AN_intro}).

Moreover, graphs of bounded Euler genus are a very special case of nearly embeddable graphs.
So our proof for coloring nearly embeddable weighted graphs gives a new proof for coloring weighted graphs of bounded Euler genus, which was proved in \cite{bbeglps_merged,bbegps} via ``fat minors'' based on work in \cite{fp}, without requiring fat minors.
Note that the coloring result for weighted graphs of bounded Euler genus implies that the Assouad-Nagata dimension of every complete Riemmanian surface of finite Euler genus is at most 2 by a fairly standard argument about quasi-isometry as shown in \cite{bbeglps_merged}.
We include some remarks about fat minors in Section \ref{sec:concluding_remarks}.

\subsection{Organization of this paper}

In Section \ref{sec:ad_wd_color}, we define coloring, weak diameter and powers of weighted graphs and show the equivalence between $d$-dimension control functions and $(d+1)$-colorings with bounded weak diameter in powers of weighted graphs.
In Section \ref{sec:center}, we show that deleting a set of vertices that are not far away from a set of bounded size does not affect the existence of colorings with bound weak diameter, which will be repeatedly used in Section \ref{sec:tree_weighted} when handling tree-decompositions.

In Section \ref{sec:tree_weighted}, we prove results about handling tree-decompositions mentioned in the above proof sketch.
We describe the gadgets that we will use in Section \ref{sec:hierarchy}.
In Section \ref{sec:condesation}, we prove that every desired coloring of the central part with gadgets can be extended to a coloring of the whole graph such that the coloring of the whole graph is desired as long as the coloring restricted on the peripheral part is desired.
(``Desired coloring'' refers to the one mentioned in the previous proof sketch.)
We use the result in Section \ref{sec:condesation} to prove the two main results about handling tree-decompositions.
In Section \ref{sec:bdd_adhesion}, we handle tree-decompositions with bounded adhesion.
In Section \ref{sec:weak_control}, we handle tree-decompositions whose bags are under control.
Then we prove Theorems \ref{minor_planar_AN_intro} and \ref{layered_tw_AN_intro} in Section \ref{sec:app_ad} by combining results in Sections \ref{sec:center} and \ref{sec:bdd_adhesion}.

In Section \ref{sec:near_embedding_geo}, we prove that every nearly embedding (with gadgets) has a tree-decomposition whose bags consist of vertices not far away from a bounded number of geodesics.
We use this result in Section \ref{sec:weighted_minor} to prove Theorem \ref{minor_AN_intro} by combining the results in Section \ref{sec:weak_control}.

We include some concluding remarks in Section \ref{sec:concluding_remarks}.

\section{Basic terminologies}

In the rest of the paper, graphs are finite unless otherwise specified.

For a function $f$ and a subset $S$ of its domain, we define \defn{the restriction of $f$ on $S$}, denoted by \defn{$f|_S$}, to be the function $g$ with domain $S$ such that $g(x)=f(x)$ for every $x \in S$.
For every integer $k$, we define $[k]$ to be the set $\{i \in {\mathbb N}: 1 \leq i \leq k\}$.

Let $(G,\phi)$ be a weighted graph.
Let $S$ be a set. 
We denote the subgraph of $G$ induced on $V(G) \cap S$ by \defn{$G[S]$}.
(Note that we do not assume that $S$ is a subset of $V(G)$.)
The \defn{subgraph of $(G,\phi)$ induced on $S$}, denoted by \defn{$(G,\phi)[S]$}, is the weighted graph $(G[S],\phi|_{E(G[S])})$.
We define \defn{$(G,\phi)-S$} to be $(G,\phi)[V(G)-S]$.

Let $(G,\phi)$ be a weighted graph.
That is, $G$ is a graph, and $\phi:E(G) \rightarrow {\mathbb R}^+$, where ${\mathbb R}^+$ is the set of all positive real numbers.
The \defn{length in $(G,\phi)$} of a path $P$ in $G$, denoted by \defn{$\leng_{(G,\phi)}(P)$}, is $\sum_{e \in E(P)}\phi(e)$.
Recall that the distance in $(G,\phi)$ between two vertices $u$ and $v$, denoted by $\dist_{(G,\phi)}(u,v)$, is defined to be the infimum of the length in $(G,\phi)$ of a path in $G$ between $u$ and $v$.
Since $G$ is finite, there exists a path attaining this infimum.
The \defn{weak diameter in $(G,\phi)$ of a set $S \subseteq V(G)$} is $\sup_{u,v \in S}\dist_{(G,\phi)}(u,v)$; the \defn{weak diameter in $(G,\phi)$ of a subgraph $H$ of $G$} is the weak diameter in $(G,\phi)$ of $V(H)$.
Let $k$ be a positive integer.
A \defn{$k$-coloring} of $(G,\phi)$ is a function $c:V(G) \rightarrow [k]$.
For a $k$-coloring $c$ of $(G,\phi)$, a \defn{$c$-monochromatic component in $(G,\phi)$} is a component of the subgraph of $(G,\phi)$ induced by $\{v \in V(G): c(v)=i\}$ for some $i \in [k]$. 
For a subgraph $(H,\phi|_{E(H)})$ of $(G,\phi)$ and an integer $d$, a \defn{$k$-coloring $c$ of $(H,\phi|_{E(H)})$ has weak diameter in $(G,\phi)$} at most $d$ if every $c$-monochromatic component in $(H,\phi|_{E(H)})$ has weak diameter in $(G,\phi)$ at most $d$.

Let $G$ be a graph.
The \defn{distance in $G$} between any two vertices $u,v$ of $G$, denoted by $\dist_G(u,v)$, the \defn{weak diameter in $G$ of a subset} $S$ of $V(G)$, and the \defn{weak diameter in $G$ of a subgraph} $H$ of $G$, are defined to be $\dist_{(G,\phi)}(u,v)$, the weak diameter in $(G,\phi)$ of $S$, and the weak diameter in $(G,\phi)$ of $H$, respectively, where $\phi$ is the constant function that maps each edge of $G$ to 1.
Let $k$ be a positive integer.
A \defn{$k$-coloring} of $G$ is a function $c:V(G) \rightarrow [k]$.
For a $k$-coloring $c$ of $G$, a \defn{$c$-monochromatic component in $G$} is a component of the subgraph of $G$ induced by $\{v \in V(G): c(v)=i\}$ for some $i \in [k]$. 
For a subgraph $H$ of $G$ and an integer $d$, a $k$-coloring $c$ of $H$ has \defn{weak diameter in $G$} at most $d$ if every $c$-monochromatic component in $H$ has weak diameter in $G$ at most $d$.

\section{Asymptotic dimension and weak diameter coloring} \label{sec:ad_wd_color}

We reduce the problem on asymptotic dimension of weighted graphs to a problem on weak diameter coloring of (unweighted) graph in this section.
The idea of this reduction is similar to the one in \cite{bbeglps_merged,l_unweighted} for asymptotic dimension of unweighted graphs.
The unweighted case in \cite{bbeglps_merged,l_unweighted} is fairly straightforward, but weighted graphs are more subtle and require careful treatments.

Let $(G,\phi)$ be a weighted graph.
Let $r \in {\mathbb R}^+$.
Define $(G,\phi,r)$ to be the weighted graph $(H,\phi')$ such that
	\begin{itemize}
		\item $H$ is obtained from $G$ by for each $e \in E(G)$, say with ends $u$ and $v$, replacing $e$ by two internally disjoint paths $P_{e,u}$ and $P_{e,v}$ between $u$ and $v$ with $\lceil \frac{\phi(e)}{r} \rceil$ edges,
		\item the image of $\phi'$ is contained in $(0,r]$ such that for every edge $z$ of $H$, 
			\begin{itemize}
				\item if $z$ is incident with $x$, where $x$ is the vertex in $G$ such that $z$ is an edge of $P_{e,x}$ for some $e \in E(G)$ incident with $x$, then $\phi'(z)=\phi(e)-r \cdot (\lceil \frac{\phi(e)}{r} \rceil-1)$, and 
				\item otherwise, $\phi'(z)=r$.
			\end{itemize}
	\end{itemize}
(For example, if $G=K_2$ and $\phi$ maps the unique edge to $\frac{13}{3}$, then $(G,\phi,2)$ is the weighted graph that consists of a 6-cycle and a function that assigns weight $2,2,\frac{1}{3},2,2,\frac{1}{3}$ to the edges in the cyclic order.)
Note that for every $xy \in E(G)$ and $z \in \{x,y\}$, $\phi(xy)$ equals the length of $P_{xy,z}$ in $(G,\phi,r)$.
In addition, the underlying graph $H$ of $(G,\phi,r)$ can be obtained from $G$ by duplicating edges and subdividing edges.
Furthermore, if the image of $\phi$ is contained in $(0,r]$, then $(G,\phi,r)$ is obtained from $(G,\phi)$ by duplicating each edge and defining the weight of the copy equal to the weight of its original. 

\begin{lemma} \label{length_subdiv}
Let $(G,\phi)$ be a weighted graph.
Let $r>0$ and $k \geq 0$ be real numbers.
Let $u,v \in V(G)$.
Denote $(G,\phi,r)$ by $(H, \phi')$.
Then there exists a path $P_G$ in $G$ between $u$ and $v$ with length $k$ in $(G,\phi)$ if and only if there exists a path $P_H$ in $H$ between $u$ and $v$ with length $k$ in $(H,\phi')$. 
In particular, the distance in $(G,\phi)$ between $u$ and $v$ equals the distance in $(G,\phi,r)$ between $u$ and $v$.
\end{lemma}

\begin{pf}
Assume that $P_G$ is a path in $G$ between $u$ and $v$ with length $k$ in $(G,\phi)$.
For each $e \in E(P_G)$, let $v_e$ be an end of $e$.
Let $P_H = \bigcup_{e \in E(P_G)}P_{e,v_e}$.
Then the length of $P_H$ in $(H,\phi')$ is $\sum_{e \in E(P_G)}\sum_{e' \in E(P_{e,v_e})}\phi'(e) = \sum_{e \in E(P_G)}\phi(e)=k$.

Now we assume that $P_H$ is a path in $H$ between $u$ and $v$ with length $k$ in $(H,\phi')$.
Since $u,v \in V(G)$, and every vertex in $V(H)-V(G)$ has degree two in $H$, for every $xy \in E(G)$ and $z \in \{x,y\}$, $E(P_H) \cap E(P_{xy,z}) \neq \emptyset$ if and only if $P_{xy,z} \subseteq P_H$.
Let $P_G$ be the path in $G$ with $E(P_G)=\{xy \in E(G):P_{xy,z} \subseteq P_H$ for some $z \in \{x,y\}\}$.
So the length of $P_G$ in $(G,\phi)$ equals $\sum_{xy \in E(G),z \in \{x,y\},P_{xy,z} \subseteq P_H}\phi(xy) = \sum_{xy \in E(G),z \in \{x,y\},P_{xy,z} \subseteq P_H}\sum_{e \in E(P_{xy,z})}\phi'(e) = \sum_{e \in E(P_H)}\phi'(e)=k$.

Since there exists a path $P_G$ in $G$ between $u$ and $v$ with length $k$ in $(G,\phi)$ if and only if there exists a path $P_H$ in $H$ between $u$ and $v$ with length $k$ in $(H,\phi')$, the distance in $(G,\phi)$ between $u$ and $v$ equals the distance in $(G,\phi,r)$ between $u$ and $v$.
\end{pf}

\bigskip

Let $(G,\phi)$ be a weighted graph.
Let $r \in {\mathbb R}^+$.
Define \defn{$(G,\phi)^r$} to be the (unweighted) graph whose vertex-set is the vertex-set of $(G,\phi,r)$, and for any two vertices $x$ and $y$, $xy$ is an edge of $(G,\phi)^r$ if and only if $x$ and $y$ are distinct and $\dist_{(G,\phi,r)}(x,y) \leq r$.
Since the weight of every edge of $(G,\phi,r)$ is at most $r$, we know that any pair of adjacent vertices in $(G,\phi,r)$ is a pair of adjacent vertices in $(G,\phi)^r$. 
In addition, if $\phi$ is the constant function with image $\{1\}$, then $(G,\phi)^\ell=G^\ell$ (except for the parallel edges) for every positive integer $\ell$.

Note that a main reason why we define $(G,\phi)^r$ via $(G,\phi,r)$ is to make sure that any two vertices with bounded distance between them in $(G,\phi)$ remains having bounded distance in $(G,\phi)^r$.
For example, if we just define $(G,\phi)^r$ from the edgeless graph with vertex-set $V(G)$ by adding an edge for any pair of vertices $u,v$ with $\dist_{(G,\phi)}(u,v) \leq r$, as in the standard way for defining the $r$-th power of a unweighted graph, then it is possible that there exist adjacent vertices $x,y$ with $\dist_{(G,\phi)}(x,y)=r+1=\phi(xy)$ such that $x$ and $y$ lie in different components of the $r$-th power and hence have infinite distance in $(G,\phi)^r$ between them.

\begin{lemma} \label{weighted_wd_to_ad}
Let $(G,\phi)$ be a weighted graph.
Let $r \in {\mathbb R}^+$, and let $m,N$ be positive integers.
If $(G,\phi)^r$ is $m$-colorable with weak diameter in $(G,\phi)^r$ at most $N$, then there exist collections $\X_1,\X_2,...,\X_m$ of subsets of $V(G)$ such that 
	\begin{itemize}
		\item $\bigcup_{i=1}^m\bigcup_{X \in \X_i}X \supseteq V(G)$,
		\item for any $i \in [m]$ and $X \in \X_i$, the weak diameter of $X$ in $(G,\phi)$ is at most $rN$, and
		\item for any $i \in [m]$, distinct $X,X' \in \X_i$ and elements $x \in X$ and $x' \in X'$, $\dist_{(G,\phi)}(x,x')>r$.
	\end{itemize}
\end{lemma}

\begin{pf}
Let $c$ be an $m$-coloring of $(G,\phi)^r$ with weak diameter in $(G,\phi)^r$ at most $N$.
For each $i \in [m]$, define $\X_i=\{V(M) \cap V(G): M$ is a $c$-monochromatic component in $(G,\phi)^r$ with $c(M)=i\}$.
So $\bigcup_{i=1}^m\bigcup_{X \in \X_i}X \supseteq V(G)$.
Denote $(G,\phi,r)$ by $(H,\phi')$.

Let $i \in [m]$ and $X \in \X_i$.
Let $u,v \in X$.
Since $X$ is contained in some $c$-monochromatic component in $(G,\phi)^r$, the weak diameter in $(G,\phi)^r$ of $X$ is at most $N$.
So there exists a path $P$ in $(G,\phi)^r$ of length at most $N$ between $u$ and $v$.
Hence there exists a path in $H$ between $u$ and $v$ of length in $(G,\phi,r)$ at most $rN$.
By Lemma \ref{length_subdiv}, the distance in $(G,\phi)$ between $u,v$ is at most $rN$.

Therefore, for every $j \in [m]$ and $X \in \X_j$, the weak diameter of $X$ in $(G,\phi)$ is at most $rN$.

Let $j \in [m]$.
Let $X',X'' \in \X_j$ be distinct.
Let $x' \in X'$ and $x'' \in X''$.
So $x'$ and $x''$ belong to different $c$-monochromatic components in $(G,\phi)^r$ with the same color.
Hence $x'$ and $x''$ are not adjacent in $(G,\phi)^r$.
So $\dist_{(G,\phi,r)}(x',x'')>r$.
Since $x',x'' \in V(G)$, by Lemma \ref{length_subdiv}, $\dist_{(G,\phi)}(x',x'')=\dist_{(G,\phi,r)}(x',x'')>r$.
This proves the lemma.
\end{pf}

\begin{lemma} \label{weight_wd_ad_1}
Let $\F$ be a class of weighted graphs. 
Let $m$ be a positive integer.
If there exists a function $f: {\mathbb R}^+ \rightarrow {\mathbb R}^+$ such that for every $(G,\phi) \in \F$ and every $\ell \in {\mathbb R}^+$, $(G,\phi)^\ell$ is $m$-colorable with weak diameter in $(G,\phi)^\ell$ at most $f(\ell)$, then $\ad(\F) \leq m-1$.
\end{lemma}

\begin{pf}
Define $g: {\mathbb R}^+ \rightarrow {\mathbb R}^+$ such that for every $x \in {\mathbb R}^+$, $g(x)=xf(x)$.
We shall show that $g$ is an $(m-1)$-dimensional control function of $\F$.

Let $(G,\phi) \in \F$ and $\ell \in {\mathbb R}^+$.
By assumption, $(G,\phi)^\ell$ is $m$-colorable with weak diameter in $(G,\phi)^\ell$ at most $f(\ell)$.
By Lemma \ref{weighted_wd_to_ad}, there exist collections $\X_1,\X_2,...,\X_m$ of subsets of $V(G)$ such that 
	\begin{itemize}
		\item $\bigcup_{i=1}^m\bigcup_{X \in \X_i}X \supseteq V(G)$,
		\item for any $i \in [m]$ and $X \in \X_i$, the weak diameter of $X$ in $(G,\phi)$ is at most $\ell f(\ell)=g(\ell)$, and
		\item for any $i \in [m]$, distinct $X,X' \in \X_i$ and elements $x \in X$ and $x' \in X'$, $\dist_{(G,\phi)}(x,x')>\ell$.
	\end{itemize}
So $g$ is an $(m-1)$-dimensional control function of $\F$.
Hence $\ad(\F) \leq m-1$.
\end{pf}

\bigskip

Let $G$ be a graph.
Let $e$ be an edge of $G$.
Let $x,y$ be the ends of $e$.
Then \defn{subdividing $e$} is the operation that deletes $e$ and add a path with two edges between $x$ and $y$.
A \defn{subdivision} of $G$ is a graph that can be obtained from $G$ by repeatedly subdividing edges.

Let $\F$ be a class of graphs.
We say that $\F$ is \defn{closed under taking subdivision} if every subdivision of any member of $\F$ belongs to $\F$.
We say that $\F$ is \defn{closed under duplication} if for every $G \in \F$, the graph obtained from $G$ by duplicating an edge belongs to $\F$.

\begin{lemma} \label{weight_wd_ad_2}
Let $\F_0$ be a class of graphs closed under duplication and taking subdivision.
Let $\F$ be a class of weighted graphs such that for every $(G,\phi) \in \F$, $G \in \F_0$.
Let $m$ be a positive integer.
If there exists a function $f: {\mathbb R}^+ \rightarrow {\mathbb R}^+$ such that for every $\ell \in {\mathbb R}^+$, every $G \in \F_0$ and every function $\phi: E(G) \rightarrow (0,\ell]$, $(G,\phi)^\ell$ is $m$-colorable with weak diameter in $(G,\phi)^\ell$ at most $f(\ell)$, then $\ad(\F) \leq m-1$.
\end{lemma}

\begin{pf}
Let $(G,\phi) \in \F$ and $\ell \in {\mathbb R}^+$.
So $G \in \F_0$.
Denote $(G,\phi,\ell)$ by $(H,\phi')$.
Since $\F_0$ is closed under duplication and taking subdivision, $H$ belongs to $\F_0$.
Note that the image of $\phi'$ is contained in $(0,\ell]$.
By assumption, $(H,\phi')^\ell$ is $m$-colorable with weak diameter in $(H,\phi')^\ell$ at most $f(\ell)$.

Since the image of $\phi'$ is contained in $(0,\ell]$, $(H,\phi',\ell)$ is obtained from $(H,\phi')$ by duplicating all edges and defining weights of the copies of the edges as the weights of their original.
So $(H,\phi')^\ell$ can be obtained from $(G,\phi)^\ell$ by duplicating some edges and defining weights of the copies of the edges as the weights of their original.
Since $(H,\phi')^\ell$ is $m$-colorable with weak diameter in $(H,\phi')^\ell$ at most $f(\ell)$, $(G,\phi)^\ell$ is $m$-colorable with weak diameter in $(G,\phi)^\ell$ at most $f(\ell)$.

Therefore, by Lemma \ref{weight_wd_ad_1}, $\ad(\F) \leq m-1$.
\end{pf}

\begin{lemma} \label{dist_original_power}
Let $\ell,k \in {\mathbb R}^+$.
Let $(G,\phi)$ be a weighted graph.
If $x,y \in V(G)$, then $\dist_{(G,\phi)^\ell}(x,y) \leq \lceil \frac{2}{\ell} \cdot \dist_{(G,\phi)}(x,y) \rceil$.
\end{lemma}

\begin{pf}
By Lemma \ref{length_subdiv}, $\dist_{(G,\phi,\ell)}(x,y)=\dist_{(G,\phi)}(x,y)$, so there exists a path $P$ in $(G,\phi,\ell)$ between $x$ and $y$ with $\leng_{(G,\phi,\ell)}(P)=\dist_{(G,\phi)}(x,y)$.
Let $(H,\phi')=(G,\phi,\ell)$.
Since the image of $\phi'$ is contained in $(0,\ell]$, we can partition $P$ into edge-disjoint subpaths $P_1,P_2,...,P_r$ (for some positive integer $r$) such that each $P_j$ has length in $(H,\phi')$ at most $\ell$, and for each $i \in [r-1]$, $P_i \cup P_{i+1}$ has length in $(H,\phi')$ greater than $\ell$.
So $(r-1)\ell <\sum_{i \in [r-1]}\sum_{e \in E(P_i) \cup E(P_{i+1})}\phi'(e) \leq 2\sum_{e \in E(P)}\phi'(e) \leq 2\dist_{(G,\phi)}(x,y)$.
Hence $r-1 < \frac{2}{\ell} \cdot \dist_{(G,\phi)}(x,y)$.
Since $r$ is an integer, $r \leq \lceil \frac{2}{\ell} \cdot \dist_{(G,\phi)}(x,y) \rceil$.
Since each $P_i$ has length in $(H,\phi')$ at most $\ell$, there exists an edge of $(G,\phi)^\ell$ between the ends of $P_i$.
So $\dist_{(G,\phi)^\ell}(x,y) \leq r \leq \lceil \frac{2}{\ell} \cdot \dist_{(G,\phi)}(x,y) \rceil$.
\end{pf}

\begin{lemma} \label{weighted_ad_wd_converse}
Let $\F_0$ be a class of graphs closed under duplication and taking subdivision.
Let $\F$ be a class of weighted graphs such that for every $(G,\phi) \in \F$, $G \in \F_0$.
Let $m$ be a positive integer.
If $\ad(\F) \leq m-1$, then there exists a function $f: {\mathbb R}^+ \rightarrow {\mathbb R}^+$ such that for every $\ell \in {\mathbb R}^+$ and $(G,\phi) \in \F$ in which the image of $\phi$ is contained in $(0,\ell]$, $(G,\phi)^\ell$ is $m$-colorable with weak diameter in $(G,\phi)^\ell$ at most $f(\ell)$.
\end{lemma}

\begin{pf}
Since $\ad(\F) \leq m-1$, there exists a function $f_0: {\mathbb R}^+ \rightarrow {\mathbb R}^+$ such that for every $(G,\phi) \in \F$ and $r \in {\mathbb R}^+$, there exist $m$ collections $\X_1,\X_2,...,\X_m$ such that 
	\begin{itemize}
		\item $\bigcup_{i=1}^m\bigcup_{X \in \X_i}X \supseteq V(G)$,
		\item for any $i \in [m]$ and $X \in \X_i$, the weak diameter in $(G,\phi)$ of $X$ is at most $f_0(r)$, and
		\item for any $i \in [m]$, distinct $X,X' \in \X_i$ and elements $x \in X$ and $x' \in X'$, the distance in $(G,\phi)$ between $x$ and $x'$ is greater than $r$.
	\end{itemize}
Define $f: {\mathbb R}^+ \rightarrow {\mathbb R}^+$ such that for every $x \in {\mathbb R}^+$, $f(x) = \lceil \frac{2f_0(x)}{x} \rceil$.

Let $\ell \in {\mathbb R}^+$.
Let $(G,\phi) \in \F$, where the image of $\phi$ is contained in $(0,\ell]$.
Hence there exist $m$ collections $\X_1,\X_2,...,\X_m$ such that 
	\begin{itemize}
		\item $\bigcup_{i=1}^m\bigcup_{X \in \X_i}X \supseteq V(G)$,
		\item for any $i \in [m]$ and $X \in \X_i$, the weak diameter in $(G,\phi)$ of $X$ is at most $f_0(\ell)$, and
		\item for any $i \in [m]$, distinct $X,X' \in \X_i$ and elements $x \in X$ and $x' \in X'$, the distance in $(G,\phi)$ between $x$ and $x'$ is greater than $\ell$.
	\end{itemize}

Since the image of $\phi$ is contained in $(0,\ell]$, $(G,\phi,\ell)$ is obtained from $(G,\phi)$ by duplicating edges and defining weights on the copied edges equal to their original.
Denote $(G,\phi,\ell)$ by $(H,\phi')$.
So $V((G,\phi)^\ell)=V(H)=V(G)$.
Define $c: V((G,\phi)^\ell) \rightarrow [m]$ such that for every $v \in V((G,\phi)^\ell) = V(G)$, $c(v) = \min\{i \in [m]: v \in \bigcup_{X \in \X_i}X\}$.
So $c$ is an $m$-coloring of $(G,\phi)^\ell$.

Let $M$ be a $c$-monochromatic component in $(G,\phi)^\ell$.
Let $i$ be the color of $M$.
Let $xy \in E(M)$.
So the distance in $(H,\phi')$ between $x$ and $y$ is at most $\ell$.
By Lemma \ref{length_subdiv}, the distance in $(G,\phi)$ between $x$ and $y$ is at most $\ell$.
Since $c(x)=c(y)$, there exists $X \in \X_i$ such that $X$ contains both $x$ and $y$.
Since this is true for arbitrary edges of $M$, we know that there exists $X^* \in \X_i$ such that $V(M) \subseteq X^*$.

Let $u,v \in V(M)$.
Since the weak diameter in $(G,\phi)$ of $X^*$ is at most $f_0(\ell)$ and $V(M) \subseteq X^*$, the distance in $(G,\phi)$ between $u,v$ is at most $f_0(\ell)$.
By Lemma \ref{dist_original_power}, the distance in $(G,\phi)^\ell$ between $u$ and $v$ is at most $\lceil \frac{2f_0(\ell)}{\ell} \rceil = f(\ell)$.

Hence the weak diameter of $M$ in $(G,\phi)^\ell$ is at most $f(\ell)$.
Therefore, $c$ is an $m$-coloring of $(G,\phi)^\ell$ with weak diameter in $(G,\phi)$ at most $f(\ell)$.
\end{pf}

\begin{theorem} \label{weighted_ad_wd_equiv}
Let $\F_0$ be a class of graphs closed under duplication and taking subdivision.
Let $\F$ be a class of weighted graphs such that 
	\begin{itemize}
		\item for every $(G,\phi) \in \F$, $G \in \F_0$, and
		\item for every $G \in \F_0$ and $\phi:E(G) \rightarrow {\mathbb R}^+$, $(G,\phi) \in \F$.
	\end{itemize}
Let $m$ be a positive integer.
Then $\ad(\F) \leq m-1$ if and only if there exists a function $f: {\mathbb R}^+ \rightarrow {\mathbb R}^+$ such that for every $\ell \in {\mathbb R}^+$, every $G \in \F_0$ and every function $\phi: E(G) \rightarrow (0,\ell]$, $(G,\phi)^\ell$ is $m$-colorable with weak diameter in $(G,\phi)^\ell$ at most $f(\ell)$.
\end{theorem}

\begin{pf}
The if part follows from Lemma \ref{weight_wd_ad_2}.
The only if part follows from Lemma \ref{weighted_ad_wd_converse}.
\end{pf}

\section{Centered sets} \label{sec:center}

The goal of this section is to provide some tools that allow us to simplify the weak diameter coloring problem on complicated graphs to the one on simpler graphs.
Those tools will be extensively used in later sections. 

Let $r$ be a nonnegative real number.
Let $(G,\phi)$ be a weighted graph.
Let $S \subseteq V(G)$.
We define \defn{$N_{(G,\phi)}^{\leq r}[S]$} $= \{v \in V(G):$ there exists a path in $G$ from $v$ to $S$ with length in $(G,\phi)$ at most $r\}$.

For $i \in [2]$, let $f_i$ be a function with domain $S_i$.
If $f_1(x)=f_2(x)$ for every $x \in S_1 \cap S_2$, then we define \defn{$f_1 \cup f_2$} to be the function with domain $S_1 \cup S_2$ such that for every $x \in S_1 \cup S_2$, $(f_1 \cup f_2)(x) = f_{i_x}(x)$, where $i_x$ is an element in $[2]$ such that $x \in S_{i_x}$.

Let $(G,\phi)$ be a weighted graph.
Let $Z \subseteq V(G)$.
Recall that we define $(G,\phi)-Z$ to be the weighted graph $(G-Z, \phi|_{E(G-Z)})$ and $(G,\phi)[Z]$ to be the weighted graph $(G[Z],\phi|_{E(G[Z])})$.
Note that for every $\ell \in {\mathbb R}^+$, the weighted graphs $(G,\phi)-Z, ((G,\phi)-Z,\ell)$ and $((G,\phi)-Z)^\ell$ are subgraphs of $(G,\phi),(G,\phi,\ell)$ and $\allowbreak (G,\phi)^\ell$, respectively.

Recall that for a subset $I$ of ${\mathbb R}^+$, we say that a weighted graph $(G,\phi)$ is $I$-bounded if the image of $\phi$ is contained in $I$.

The following lemma is the main result of this section, which roughly states that if the vertices in $Z$ is not far from a small set $S$ of vertices, then any coloring on $Z$ can be combined with any coloring of $G-(Z \cup R)$ (for any set $R$) without increasing the weak diameter much.
It would allow us to first delete $Z$ from $G$ to reduce the complexity of $G$, then find a desired coloring for $G-Z$, and finally combine this coloring with a coloring on $Z$ to obtain a desired coloring for $G$.

\begin{lemma} \label{patching_centered_set_gen}
For any nonnegative integer $k$, nonnegative real number $r$, and positive real numbers $\ell,\nu$, there exists a real number $\nu^* \geq (k+1)\nu$ such that the following holds.
Let $(G,\phi)$ be a weighted graph.
Let $S \subseteq V(G)$ with $\lvert S \rvert \leq k$.
Let $Z \subseteq N_{(G,\phi)}^{\leq r}[S]$.
Let $R \subseteq V(G)$.
Let $c_Z: Z \rightarrow [m]$ for some positive integer $m$.
Let $c$ be an $m$-coloring of $((G,\phi)-Z)^\ell[V(G)]-R$ with weak diameter in $(G,\phi)^\ell$ at most $\nu$.
Then the $m$-coloring $c \cup c_Z|_{Z-R}$ of $(G,\phi)^\ell[V(G)]-R$ has weak diameter in $(G,\phi)^\ell$ at most $\nu^*$.
\end{lemma}

\begin{pf}
Let $k$ be a nonnegative integer, $r$ be a nonnegative real number, and $\ell,\nu$ be positive real numbers.
Let $f: ({\mathbb N} \cup \{0\}) \times {\mathbb R} \rightarrow {\mathbb R}$ be the function such that for every $y \in {\mathbb R}$, $f(0,y)=y$ and for every $x \in {\mathbb N}$, $f(x,y) = 2f(x-1,\lceil \frac{4}{\ell} \cdot (\ell+r+\ell y) \rceil+y) + 2 \cdot \lceil \frac{2(\ell+r)}{\ell} \rceil$.
Define $\nu^* = f(k,\nu)$.
Note that $\nu^* \geq (k+1)\nu$.

We shall prove this lemma by induction on $k$.
Let $(G,\phi),S,Z,m,R,c_Z,c$ be the ones as defined in the lemma.

When $k=0$, $S=Z=\emptyset$, so $c \cup c_Z|_{Z-R}=c$ is an $m$-coloring of $((G,\phi)-Z)^\ell[V(G)]-R=(G,\phi)^\ell[V(G)]-R$ with weak diameter in $(G,\phi)^\ell$ at most $\nu=f(0,\nu)=\nu^*$.
So we may assume $k \geq 1$ and the lemma holds when $k$ is smaller.
In particular, $|S|=k \geq 1$, for otherwise we are done.

Let $s \in S$.
Let $S'=S-\{s\}$.
Let $Z' = Z \cap N_{(G,\phi)}^{\leq r}[S']$.
Let $Z_s = Z-Z'$.
Note that $Z_s \subseteq N_{(G,\phi)}^{\leq r}[\{s\}]$.
Let $R' = R \cup Z_s$.
So $R' \cup Z' = R \cup Z$.

Note that $V(((G,\phi)-Z')^\ell[V(G)]-R') = V(((G,\phi)-Z)^\ell[V(G)]-R)$, and $((G,\phi)-Z)^\ell[V(G)]-R$ is a subgraph of $((G,\phi)-Z')^\ell[V(G)]-R'$.
So $c$ is a coloring for both graphs, but the $c$-monochromatic components in those two graphs can be different.

Let $\nu'=\lceil \frac{4}{\ell} \cdot (\ell+r+\ell \nu) \rceil+\nu$.

\medskip

\noindent{\bf Claim 1:} $c$ is an $m$-coloring of $((G,\phi)-Z')^\ell[V(G)]-R'$ with weak diameter in $(G,\phi)^\ell$ at most $\nu'$.

\noindent{\bf Proof of Claim 1:}
Let $W'$ be a $c$-monochromatic component in $((G,\phi)-Z')^\ell[V(G)]-R'$.
It suffices to show that $W'$ has weak diameter in $(G,\phi)^\ell$ at most $\nu'$.

Let $W$ be a $c$-monochromatic component in $((G,\phi)-Z)^\ell[V(G)]-R$ with $V(W) \subseteq V(W')$.
By the assumption of $c$, $W$ has weak diameter in $(G,\phi)^\ell$ at most $\nu$.
If $V(W)=V(W')$, then $W'$ also has weak diameter in $(G,\phi)^\ell$ at most $\nu \leq \nu'$.
So we may assume $V(W) \subset V(W')$.
Hence $V(W)$ is incident with an edge in $((G,\phi)-Z')^\ell[V(G)]-R'$ but not in $((G,\phi)-Z)^\ell[V(G)]-R$.
So $V(W)$ contains a vertex $w \in N_{(G,\phi)}^{\leq \ell}[Z-Z'] = N_{(G,\phi)}^{\leq \ell}[Z_s] \subseteq N_{(G,\phi)}^{\leq \ell+r}[\{s\}]$.
Since $W$ has weak diameter in $(G,\phi)^\ell$ at most $\nu$, $W$ has weak diameter in $(G,\phi)$ at most $\ell \nu$, so $V(W) \subseteq N_{(G,\phi)}^{\leq \ell+r+\ell \nu}[\{s\}]$.

Since $W$ is an arbitrary $c$-monochromatic component in $((G,\phi)-Z)^\ell[V(G)]-R$ with $V(W) \subseteq V(W')$, and $Z-Z'=Z_s \subseteq N_{(G,\phi)}^{\leq r}[\{s\}]$, we know that for any $u,v \in V(W')$, $\dist_{(G,\phi)}(u,v) \leq \dist_{(G,\phi)}(u,s)+ \dist_{(G,\phi)}(s,v) \leq 2(\ell+r+\ell \nu)$.
So by Lemma \ref{dist_original_power}, for any $u,v \in V(W')$, $\dist_{(G,\phi)^\ell}(u,v) \leq \lceil \frac{2}{\ell} \cdot \dist_{(G,\phi)}(u,v) \rceil \leq \lceil \frac{4}{\ell} \cdot (\ell+r+\ell \nu) \rceil \leq \nu'$.
Hence $W'$ has weak diameter in $(G,\phi)^\ell$ at most $\nu'$.
$\Box$

\medskip

Since $Z' \subseteq N_{(G,\phi)}^{\leq r}[S']$ and $|S'| = |S|-1=k-1$, by the induction hypothesis (with $(S,Z,R,c_Z,c)$ replaced by $(S',Z',R',c_Z|_{Z'},c))$, Claim 1 implies that $c \cup c_Z|_{Z'-R'}$ is an $m$-coloring of $(G,\phi)^\ell[V(G)]-R' = ((G,\phi)^\ell[V(G)]-Z_s)-R$ with weak diameter in $(G,\phi)^\ell$ at most $f(k-1,\nu')$.
Note that $c \cup c_Z|_{Z-R} = (c \cup c_Z|_{Z'-R'}) \cup c_Z|_{Z_s-R}$.

Let $M$ be a $(c \cup c_Z|_{Z-R})$-monochromatic component in $(G,\phi)^\ell[V(G)]-R$.
To prove the lemma, it suffices to show that $M$ has weak diameter in $(G,\phi)^\ell$ at most $f(k,\nu)$.

Let $Q$ be a component of $M-Z_s$.
So $Q$ is contained in a $(c \cup c_Z|_{Z'-R'})$-monochromatic component in $((G,\phi)^\ell[V(G)]-R)-Z_s = (G,\phi)^\ell[V(G)]-R'$.
Hence $Q$ has weak diameter in $(G,\phi)^\ell$ at most $f(k-1,\nu')$.
If $V(Q)$ is not adjacent in $(G,\phi)^\ell[V(G)]-R$ to a vertex in $Z_s$, then $M=Q$ has weak diameter in $(G,\phi)^\ell$ at most $f(k-1,\nu') \leq f(k,\nu)$.
So we may assume that $V(Q)$ is adjacent in $(G,\phi)^\ell[V(G)]-R$ to a vertex in $Z_s$.
Hence there exists $q \in V(Q) \cap N_{(G,\phi)}^{\leq \ell}[Z_s] \subseteq V(Q) \cap N_{(G,\phi)}^{\leq \ell+r}[\{s\}]$.
By Lemma \ref{dist_original_power}, $\dist_{(G,\phi)^\ell}(q,s) \leq \lceil \frac{2}{\ell} \cdot \dist_{(G,\phi)}(q,s) \rceil \leq \lceil \frac{2(\ell+r)}{\ell} \rceil$.
Since $Q$ has weak diameter in $(G,\phi)^\ell$ at most $f(k-1,\nu')$, for every vertex $v \in V(Q)$, $\dist_{(G,\phi)^\ell}(v,s) \leq \dist_{(G,\phi)^\ell}(v,q) + \dist_{(G,\phi)^\ell}(q,s) \leq f(k-1,\nu') + \lceil \frac{2(\ell+r)}{\ell} \rceil$.

Since $Q$ is an arbitrary component of $M-Z_s$ and $Z_s \subseteq N_{(G,\phi)}^{\leq r}[\{s\}] \subseteq N_{(G,\phi)^\ell}^{\leq \lceil 2r/\ell \rceil}[\{s\}]$ by Lemma \ref{dist_original_power}, we know that for any vertices $x,y \in V(M)$, $\dist_{(G,\phi)^\ell}(x,y) \leq \dist_{(G,\phi)^\ell}(x,s) + \dist_{(G,\phi)^\ell}(s,y) \leq 2 \cdot (f(k-1,\nu') + \lceil \frac{2(\ell+r)}{\ell} \rceil) = f(k,\nu)$.
This proves the lemma.
\end{pf}

\bigskip

We will extensively use corollaries of Lemma \ref{patching_centered_set_gen} in this paper.
The following is the first corollary, which will allow us to reduce the complexity of a weighted graph by deleting vertices.

\begin{lemma} \label{deleting_centered_set}
For any nonnegative integer $k$, nonnegative real number $r$, and positive real numbers $\ell,\nu$, there exists a real number $\nu^* \geq (k+1)\nu$ such that the following holds.
Let $(G,\phi)$ be a $(0,\ell]$-bounded weighted graph.
Let $S \subseteq V(G)$ with $\lvert S \rvert \leq k$.
Let $Z \subseteq N_{(G,\phi)}^{\leq r}[S]$.
Let $m$ be a positive integer.
Let $c_Z: Z \rightarrow [m]$.
Let $c$ be an $m$-coloring of $((G,\phi)-Z)^\ell$ with weak diameter in $((G,\phi)-Z)^\ell$ at most $\nu$.
Then the $m$-coloring $c \cup c_Z$ of $(G,\phi)^\ell$ has weak diameter in $(G,\phi)^\ell$ at most $\nu^*$.
\end{lemma}

\begin{pf}
Since $((G,\phi)-Z)^\ell$ is a subgraph of $(G,\phi)^\ell$, $c$ is an $m$-coloring of $((G,\phi)-Z)^\ell$ with weak diameter in $(G,\phi)^\ell$ at most $\nu$.
Since $(G,\phi)$ is $(0,\ell]$-bounded, $((G,\phi)-Z)^\ell = ((G,\phi)-Z)^\ell[V(G)]$.
Then the lemma follows from Lemma \ref{patching_centered_set_gen} with $R=\emptyset$. 
\end{pf}

\bigskip

Let $\ell,\nu$ be positive real numbers and $m$ be a positive integer.
We say a class $\F$ of $(0,\ell]$-bounded weighted graphs is \defn{$(m,\ell,\nu)$-nice} if for every $(G,\phi) \in \F$, $(G,\phi)^\ell$ is $m$-colorable with weak diameter in $(G,\phi)^\ell$ at most $\nu$.

For a nonnegative integer $n$, a subset $I$ of ${\mathbb R}^+$, and a class $\F$ of weighted graph, \defn{$\F^{+n,I}$} is the class consisting of all the $I$-bounded weighted graphs $(G,\phi)$ in which there exists $Z \subseteq V(G)$ with $\lvert Z \rvert \leq n$ such that $(G,\phi)-Z \in \F$.

The following is a more concise version of Lemma \ref{deleting_centered_set}. 

\begin{lemma} \label{apex_extension}
For any positive real numbers $\ell,\nu$, and nonnegative integer $n$, there exists a real number $\nu^*$ such that if $m$ is a positive integer and $\F$ is an $(m,\ell,\nu)$-nice class of $(0,\ell]$-bounded weighted graphs, then $\F^{+n,(0,\ell]}$ is an $(m,\ell,\nu^*)$-nice class of $(0,\ell]$-bounded weighted graphs. 
\end{lemma}

\begin{pf}
Let $\ell,\nu$ be positive real numbers, and let $n$ be a nonnegative integer.
Define $\nu^*$ to be the real number $\nu^*$ given by Lemma \ref{deleting_centered_set} by taking $(k,r,\ell,\nu)=(n,0,\ell,\nu)$.

Let $(G,\phi) \in \F^{+n,(0,\ell]}$.
Then there exists $Z \subseteq V(G)$ with $\lvert Z \rvert \leq n$ such that $(G,\phi)-Z \in \F$.
Since $\F$ is an $(m,\ell,\nu)$-nice $(0,\ell]$-bounded class, there exists an $m$-coloring $c$ of $((G,\phi)-Z)^\ell$ with weak diameter in $((G,\phi)-Z)^\ell$ at most $\nu$. 
Let $c_Z$ be an $m$-coloring of $Z$.
Since $Z \subseteq N_{(G,\phi)}^{\leq 0}[Z]$ and $|Z| \leq n$, by Lemma \ref{deleting_centered_set}, $c \cup c_Z$ is an $m$-coloring of $(G,\phi)^\ell$ with weak diameter in $(G,\phi)^\ell$ at most $\nu^*$.

Therefore, $\F^{+n,(0,\ell]}$ is an $(m,\ell,\nu^*)$-nice class of $(0,\ell]$-bounded weighted graphs.
\end{pf}

\bigskip

Let $r$ and $k$ be nonnegative real numbers.
Let $(G,\phi)$ be a graph.
We say that a subset $Z$ of $V(G)$ is \defn{$(k,r)$-centered in $(G,\phi)$} if there exists $S \subseteq V(G)$ with $\lvert S \rvert \leq k$ such that $Z \subseteq N_{(G,\phi)}^{\leq r}[S]$.

The following two lemmas show that certain kinds of graphs can be trivially colored, and we will use them as building blocks to color more complicated graphs.

\begin{lemma} \label{all_centered}
For any nonnegative integer $k$, nonnegative real number $r$, and positive real number $\ell$, there exists a positive real number $\nu$ such that the following holds.
Let $(G,\phi)$ be a weighted graph.
Let $m$ be a positive integer.
Let $R \subseteq V(G)$.
If $V(G)-R$ is $(k,r)$-centered in $(G,\phi)$, then any $m$-coloring of $(G,\phi)^\ell[V(G)]-R$ has weak diameter in $(G,\phi)^\ell$ at most $\nu$. 
\end{lemma}

\begin{pf}
Let $k$ be a nonnegative integer, $r$ be a nonnegative real number, and $\ell$ be a positive real number.
Define $\nu$ to be the real number $\nu^*$ given by Lemma \ref{patching_centered_set_gen} by taking $(k,r,\ell,\nu)=(k,r,\ell,1)$.

Let $(G,\phi)$ be a weighted graph.
Let $m$ be a positive integer.
Let $R \subseteq V(G)$.
Since $V(G)-R$ is $(k,r)$-centered in $(G,\phi)$, there exists $S \subseteq V(G)$ with $|S| \leq k$ such that $V(G)-R \subseteq N_{(G,\phi)}^{\leq r}[S]$.
Let $c_Z$ be an arbitrary $m$-coloring of $(G,\phi)^\ell[V(G)]-R$.
It suffices to show that $c_Z$ has weak diameter in $(G,\phi)^\ell$ at most $\nu$. 
Note that the domain of $c_Z$ is $V(G)-R$.
Then this lemma follows from Lemma \ref{patching_centered_set_gen} (with $\nu=1$, $Z=V(G)-R$ and $c$ an arbitrary coloring on the empty graph).
\end{pf}

\bigskip

For a positive integer $w$, a \defn{$w$-vertex-cover} of a graph $G$ is a subset $S$ of $V(G)$ such that every component of $G-S$ has at most $w$ vertices.

\begin{lemma} \label{weighted_vc_color}
For any nonnegative integer $k$, positive real number $\ell$ and positive integer $w$, there exists a positive real number $\nu$ such that if $m$ is a positive integer and $(G,\phi)$ is a $(0,\ell]$-bounded weighted graph such that $G$ has a $w$-vertex-cover of size at most $k$, then any $m$-coloring of $(G,\phi)^\ell$ has weak diameter in $(G,\phi)^\ell$ at most $\nu$.
\end{lemma}

\begin{pf}
Let $k$ be a nonnegative integer, $\ell$ be a positive real number, and $w$ be a positive integer.
Let $\nu_1$ be the real number $\nu$ given by Lemma \ref{all_centered} by taking $(k,r,\ell)=(w,0,\ell)$.
Define $\nu$ to be the real number $\nu^*$ given by Lemma \ref{deleting_centered_set} by taking $(k,r,\ell,\nu)=(k,0,\ell,\nu_1)$.

Let $m$ be a positive integer.
Let $(G,\phi)$ be a $(0,\ell]$-bounded weighted graph such that $G$ has a $w$-vertex-cover $S$ of size at most $k$.
Let $c$ be an $m$-coloring of $(G,\phi)^\ell$.
It suffices to show that $c$ has weak diameter in $(G,\phi)^\ell$ at most $\nu$.

For every component $C$ of $(G,\phi)-S$, $|V(C)| \leq w$, so $V(C)$ is $(w,0)$-centered in $C$, and hence $c|_{V(C)}$ is an $m$-coloring of $C^\ell$ with weak diameter in $C^\ell$ at most $\nu_1$ by Lemma \ref{all_centered}.
Since $((G,\phi)-S)^\ell = \bigcup_C C^\ell$, where the union is over all components $C$ of $(G,\phi)-S$, $c|_{V(G)-S}$ is an $m$-coloring of $((G,\phi)-S)^\ell$ with weak diameter in $((G,\phi)-S)^\ell$ at most $\nu_1$.
Since $S$ is $(k,0)$-centered in $(G,\phi)$, by Lemma \ref{deleting_centered_set}, $c$ is an $m$-coloring of $(G,\phi)^\ell$ with weak diameter in $(G,\phi)^\ell$ at most $\nu$.
\end{pf}

\bigskip

We also need the following variation of Lemma \ref{deleting_centered_set}, which concerns the weak diameter in the whole graph and allows a set $R$ of vertices that are not required to be colored.
It will be used in later sections when deleting vertices would destroy the property of having small weak diameter.

\begin{lemma} \label{patching_centered_set}
For any nonnegative integer $k$, nonnegative real number $r$, and positive real numbers $\ell,\nu$, there exists a real number $\nu^* \geq \nu$ such that the following holds.
Let $(G,\phi)$ be a weighted graph.
Let $S \subseteq V(G)$ with $\lvert S \rvert \leq k$.
Let $Z \subseteq N_{(G,\phi)}^{\leq r}[S]$.
Let $m$ be a positive integer.
Let $R \subseteq V(G)$.
Let $c_Z: Z \rightarrow [m]$.
Let $c$ be an $m$-coloring of $(G,\phi)^\ell[V(G)]-(R \cup Z)$ with weak diameter in $(G,\phi)^\ell$ at most $\nu$.
Then the $m$-coloring $c \cup c_Z|_{Z-R}$ of $(G,\phi)^\ell[V(G)]-R$ has weak diameter in $(G,\phi)^\ell$ at most $\nu^*$.
\end{lemma}

\begin{pf}
Let $k$ be a nonnegative integer, $r$ be a nonnegative real number, and $\ell,\nu$ be positive real numbers.
Define $\nu^*$ to be the real number $\nu^*$ given by Lemma \ref{patching_centered_set_gen} by taking $(k,r,\ell,\nu)=(k,r,\ell,\nu)$.

Let $(G,\phi),S,Z,m,R,c_Z,c$ be as defined in the lemma.
Since $((G,\phi)-Z)^\ell[V(G)]-R$ is a spanning subgraph of $(G,\phi)^\ell[V(G)]-(R \cup Z)$, every $c$-monochromatic component in $((G,\phi)-Z)^\ell[V(G)]-R$ is contained in some $c$-monochromatic component in $(G,\phi)^\ell[V(G)]-(R \cup Z)$.
Since $c$ is an $m$-coloring of $(G,\phi)^\ell[V(G)]-(R \cup Z)$ with weak diameter in $(G,\phi)^\ell$ at most $\nu$, $c$ is also an $m$-coloring of $((G,\phi)-Z)^\ell[V(G)]-R$ with weak diameter in $(G,\phi)^\ell$ at most $\nu$.
So this lemma follows from Lemma \ref{patching_centered_set_gen}.
\end{pf}

\section{Gluing weighted graphs along trees} \label{sec:tree_weighted}

We will develop the main tools related to tree-decompositions in this section.
The goals are to prove Lemmas \ref{weighted_tree_extension_clean} and \ref{strong_weighted_tree_extension_control_clean}, which roughly state that if a graph admits a tree-decomposition such that every ``bag'' has a ``nice'' coloring or a ``nice'' property, then the whole graph has a ``nice'' coloring.
The strategies used to prove them are similar: we shall first divide the graphs into the ``central part'' and the ``peripheral parts'', and then apply the induction hypothesis to each of the central part and the peripheral parts to obtain a coloring, and finally combine the colorings of different parts to obtain a coloring of the whole graph.
In order to do the last step about combining colorings, we should add some ``gadgets'' to the central part to include information of the peripheral parts; the condensation defined in Section \ref{sec:condesation} is the formal form of the central part with the added gadgets.
Section \ref{sec:hierarchy} introduces the gadgets that we will use in Section \ref{sec:condesation}.
In Section \ref{sec:condesation}, besides formally defining condensations, we study properties of condensations and show how the coloring of the condensation allows us to combine the colorings.
Then we prove Lemma \ref{weighted_tree_extension_clean} in Section \ref{sec:bdd_adhesion} and Lemma \ref{strong_weighted_tree_extension_control_clean} in Section \ref{sec:weak_control}. 

We define some terminologies before we proceed.

Let $G$ be a graph, and let $(T,\X)$ be a tree-decomposition of $G$, where $\X=(X_t: t \in V(T))$.
For every $S \subseteq V(T)$, we define \defn{$X_S$} $= \bigcup_{t \in S}X_t$; for every subgraph $H$ of $T$, we define \defn{$X_H$} $= X_{V(H)}$.

A \defn{rooted tree} is a directed graph whose underlying graph is a tree such that there exists a unique vertex $r$ with in-degree 0.
We call $r$ the \defn{root} of this rooted tree.
A \defn{rooted tree-decomposition} of a graph $G$ is a tree-decomposition $(T,\X)$ of $G$ such that $T$ is a rooted tree.

We say that two edges $e$ and $e'$ of a rooted tree $T$ are \defn{incomparable} if no directed path in $T$ from the root of $T$ contains both $e$ and $e'$.

\subsection{Hierarchy} \label{sec:hierarchy}

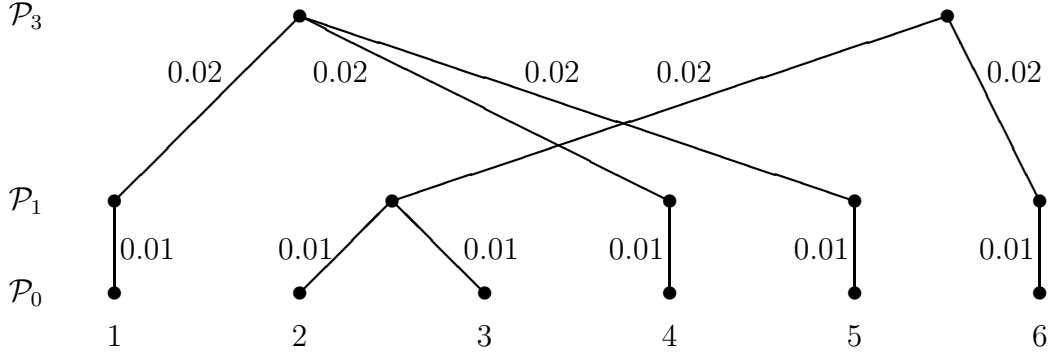
\begin{figure} 
	\begin{picture}(100,200) (-5,-20)

		\thicklines
		
		\multiput(80,20)(70,0){6}{\circle*{5}}
		\put(40,17){{$\P_0$}}
		\put(82,33){{$0.01$}}
		\put(142,33){{$0.01$}}
		\put(212,33){{$0.01$}}
		\put(267,33){{$0.01$}}
		\put(337,33){{$0.01$}}
		\put(407,33){{$0.01$}}
		\put(77,0){{$1$}}
		\put(147,0){{$2$}}
		\put(217,0){{$3$}}
		\put(287,0){{$4$}}
		\put(357,0){{$5$}}
		\put(427,0){{$6$}}

		\put(80,55){\circle*{5}}
		\put(185,55){\circle*{5}}
		\put(290,55){\circle*{5}}
		\put(360,55){\circle*{5}}
		\put(430,55){\circle*{5}}
		\put(40,52){{$\P_1$}}
		\put(80,20){\line(0,1){35}}
		\put(150,20){\line(1,1){35}}
		\put(220,20){\line(-1,1){35}}
		\put(290,20){\line(0,1){35}}
		\put(360,20){\line(0,1){35}}
		\put(430,20){\line(0,1){35}}

		\put(150,125){\circle*{5}}
		\put(395,125){\circle*{5}}
		\put(40,122){{$\P_3$}}
		\put(80,55){\line(1,1){70}}
		\put(290,55){\line(-2,1){140}}
		\put(360,55){\line(-3,1){210}}
		\put(185,55){\line(3,1){210}}
		\put(430,55){\line(-1,2){35}}
		\put(100,100){{$0.02$}}
		\put(155,100){{$0.02$}}
		\put(235,100){{$0.02$}}
		\put(285,100){{$0.02$}}
		\put(410,100){{$0.02$}}

	\end{picture}
	\caption{An example of the $(8,\frac{1}{4},3,1)$-heirachy $(H,\phi_H,I,B)$ of $\P=(\P_i: i \in {\mathbb N})$, where $X=[6]$, $\P_0=\{\{x\}: x \in [6]\}$, $\P_1=\P_2 = \{\{1\},\{2,3\},\{4\},\{5\},\{6\}\}$, $\P_i = \{\{1,4,5\},\{2,3,6\}\}$ for every $i \geq 3$. Note that $I=\{0,1,3\}$, $H$ is the graph in the picture, and $B$ is the set of the six vertices at the bottom.} \label{fig_hierachy}
\end{figure}

Let $X$ be a finite set.
Let $\P_0$ be the partition $\{\{x\}: x \in X\}$.
For each $i \in {\mathbb N}$, let $\P_i$ be a partition of $X$ such that $\P_{i-1}$ is a refinement of $\P_i$.
That is, every member of $\P_{i-1}$ is a subset of a member of $\P_i$.
Let $\ell, \epsilon$ be positive real numbers with $\epsilon \leq \ell$, and let $\mu$ be a nonnegative real number.
Let $\theta$ be a positive integer\footnote{The use of the parameters $\ell,\epsilon,\mu,\theta$ will be clear later. Roughly speaking, recall that the hierarchy that we will define serves the gadgets that will be used to encode information about the vertices not far from the central part. The parameters $\ell$ and $\mu$ determine how far the vertices that will be concerned. The parameter $\epsilon$ will usually denote the minimum length of an edge of the graph, which determines the number of ``levels'' of the gadgets. The parameter $\theta$ is used for normalizing the weight of the edges of the gadgets.}.
We define the following.
	\begin{itemize}
		\item $I = \{0,1\} \cup \{i \in [\lceil \frac{3\ell+\mu}{\epsilon} \rceil]: \P_{i} \neq \P_{i-1}\}$. 
		\item Let $H$ be the graph with 
			\begin{itemize}
				\item $V(H) = \{v_Y: Y \in \P_{i}$ for some $i \in I\}$\footnote{For simplicity of notation, even though a set $Y$ can be a member of both $\P_i$ and $\P_{i'}$ for some $i \neq i'$, we consider the set $Y$ in different partitions $\P_i$ and $\P_{i'}$ as distinct copies when we define $V(H)$, so the vertices $v_Y$ given by $\P_i$ and $\P_{i'}$ are actually distinct vertices in $H$. In the rest of the paper, whenever we use $v_Y$, it will be clear that this $Y$ is contained in which partition $\P_i$. See Figure \ref{fig_hierachy} for an example.}, and 
				\item $uv \in E(H)$ if and only if there exist $i,j \in I$ with $i<j$ such that 
					\begin{itemize}
						\item there exists no $k \in I$ with $i<k<j$, and 
						\item there exist $Y_i \in \P_{i}$ and $Y_j \in \P_{j}$ with $Y_i \subseteq Y_j$ such that either $u=v_{Y_i}$ and $v=v_{Y_j}$, or $u=v_{Y_j}$ and $v=v_{Y_i}$. 
					\end{itemize}
			\end{itemize}
			(In other words, $H$ is the underlying graph of the Hasse diagram of the following poset: the ground set is $\{(Y,i): i \in I, Y \in \P_i\}$, and for any two elements $(Y_1,i_1)$ and $(Y_2,i_2)$, $(Y_1,i_1)$ is equal to or strictly smaller than $(Y_2,i_2)$ in partial ordering if and only if $Y_1 \subseteq Y_2$ and $i \leq j$. The vertex $v_Y$ in the definition of $V(H)$ above refers to the point $(Y,i)$. See Figure \ref{fig_hierachy} for an example.) 
		\item Let $\phi_{H}: E(H) \rightarrow {\mathbb R}^+$ such that for every $v_Yv_{Y'} \in E(H)$, $\phi_{H_e}(v_Yv_{Y'})=(j-i)\cdot\frac{\epsilon}{8(\theta+\frac{\mu}{\ell})}$, where $i,j$ are the elements in $I$ such that $i<j$, $Y \in \P_{i}$ and $Y' \in \P_{j}$. 
		\item Let $B = \{v_Y: Y \in \P_{0}\}$.
	\end{itemize}
We call $(H,\phi_H,I,B)$ the \defn{$(\ell,\epsilon,\theta,\mu)$-hierarchy} of $\P$, where $\P=(\P_0,\P_1,...)$.

\begin{lemma} \label{hierachy_distance}
Let $X$ be a finite set.
Let $\P_0$ be the partition $\{\{x\}: x \in X\}$.
For each $i \in {\mathbb N}$, let $\P_i$ be a partition of $X$ such that $\P_{i-1}$ is a refinement of $\P_i$.
Let $\P=(\P_0,\P_1,...)$.
Let $\ell, \epsilon$ be positive real numbers with $\epsilon \leq \ell$.
Let $\mu$ be a nonnegative real number.
Let $\theta$ be a positive integer.
Let $(H,\phi_H,I,B)$ be the $(\ell,\epsilon,\theta,\mu)$-hierarchy of $\P$. 
Then
	\begin{enumerate}
		\item the image of $\phi_{H}$ is contained in $(0,\frac{\ell+\mu}{2(\theta+\frac{\mu}{\ell})}] \subseteq (0,\frac{\ell}{2}]$, 
		\item $V(H) \subseteq N_{(H,\phi_{H})}^{\leq \frac{\ell}{2}}[B]$, 
		\item for any two vertices $x,y$ in the same component of $H$, $\dist_{(H,\phi_H)}(x,y) \leq \ell$, and 
		\item $\lvert V(H)-B \rvert \leq \lvert B \rvert^2$.
	\end{enumerate}
\end{lemma}

\begin{pf}
Since $\epsilon \in (0,\ell]$ and $\max I \leq \lceil \frac{3\ell+\mu}{\epsilon} \rceil$, we have $\max I \cdot \frac{\epsilon}{8(\theta+\frac{\mu}{\ell})} \leq \lceil \frac{3\ell+\mu}{\epsilon} \rceil \cdot \frac{\epsilon}{8(\theta+\frac{\mu}{\ell})} \leq (\frac{3\ell+\mu}{\epsilon}+1) \frac{\epsilon}{8(\theta+\frac{\mu}{\ell})} \leq \frac{3\ell+\mu+\epsilon}{8(\theta+\frac{\mu}{\ell})} \leq \frac{4\ell+\mu}{8(\theta+\frac{\mu}{\ell})} \leq \frac{\ell+\mu}{2(\theta+\frac{\mu}{\ell})}$.
So the image of $\phi_{H}$ is contained in $(0,\frac{\ell+\mu}{2(\theta+\frac{\mu}{\ell})}] \subseteq (0,\frac{\ell}{2}]$ and $V(H) \subseteq N_{(H,\phi_H)}^{\leq \frac{\ell+\mu}{2(\theta+\frac{\mu}{\ell})}}[B] \subseteq N_{(H,\phi_H)}^{\leq \frac{\ell}{2}}[B]$.

Let $x,y$ be two vertices in the same component of $H$.
So there exist $i_x,i_y \in I \subseteq [\lceil \frac{3\ell+\mu}{\epsilon} \rceil] \cup \{0\}$, $Y_x \in \P_{i_x}$ and $Y_y \in \P_{i_y}$ such that $x=v_{Y_x}$ and $y=v_{Y_y}$.
Since for every $i \in {\mathbb N}$, $\P_{i-1}$ is a refinement of $\P_{i}$, there exists $Y \in \P_{\max I}$ with $Y_x \cup Y_y \subseteq Y$ such that for each $z \in \{x,y\}$, $\dist_{(H,\phi_H)}(z,v_Y) = (\max I - i_z) \cdot \frac{\epsilon}{8(\theta+\frac{\mu}{\ell})} \leq \frac{\ell}{2}$. 
So $\dist_{(H,\phi_H)}(x,y) \leq \dist_{(H,\phi_H)}(x,v_Y) + \dist_{(H,\phi_H)}(v_Y,y) \leq \ell$. 

If $X=\emptyset$, then $|V(H)-B|=0= |B|=|B|^2$.
So we may assume $|X| \geq 1$.
For each $i \in {\mathbb N}$, $\P_{i-1}$ is a refinement of $\P_{i}$, so $|\P_{i-1}| \geq |\P_{i}|$, and $|\P_{i-1}|>|\P_i|$ unless $\P_i=\P_{i-1}$.
Denote the elements of $I$ by $i_0<i_1<i_2<...<i_{|I|-1}$.
Note that $i_0=0$ and $i_1=1$.
By the definition of $I$, we know $|\P_1| \leq |\P_0| = |B|$, and for every $j \in [|I|-1]$, $1 \leq |\P_{i_j}| \leq |\P_1|-j+1 \leq |B|-j+1$.
Hence $|I|-1 \leq |B|$, and $\lvert V(H)-B \rvert = \sum_{j=1}^{|I|-1}|\P_{i_j}| \leq \sum_{j=1}^{|B|}(|B|-j+1) \leq \lvert B \rvert^2$.
\end{pf}

\subsection{Condensation} \label{sec:condesation}

Let $\ell$ be a positive real number.
Let $(G,\phi)$ be a weighted graph. 
Let $(T,\X)$ be a rooted tree-decomposition of $G$. 
Denote the root of $T$ by $t^*$, and denote $\X=(X_t: t \in V(T))$.
For each $e=tt' \in E(T)$, define $T_e$ to be the component of $T-e$ disjoint from $t^*$, and define $X_e=X_t \cap X_{t'}$.
We assume that $(G,\phi)$ is \defn{$(T,\X,\ell)$-bounded}, that is, for every edge $xy \in E(G)$ with $\phi(xy)>\ell$, there exists $t \in V(T)-\{t^*\}$ such that
	\begin{itemize}
		\item $t$ has no child in $T$,
		\item $\{x,y\} \subseteq X_t-X_e$ and $X_t \subseteq N_{(G,\phi)}^{\leq \ell}[X_e]$, where $e$ is the unique edge of $T$ incident with $t$.
	\end{itemize}

Let $\epsilon = \min_{e \in E(G)}\phi(e)$.
Let $U_E \subseteq E(T)$ such that edges in $U_E$ are pairwise incomparable.
For each $e \in U_E$ and $i \in {\mathbb N}$, define $\preceq_{e,i}$ to be the binary relation on $X_e$ such that for any $x,y \in X_e$, $x \preceq_{e,i} y$ if and only if there exists a path $P$ in $(G,\phi)[X_{T_e}]$ from $x$ to $y$ such that 
	\begin{itemize}
		\item $V(P) \subseteq N_{(G,\phi)[X_{T_e}]}^{\leq i\epsilon}[X_e]$, and 
		\item for every $e' \in E(P)$, $\phi(e') \leq i\epsilon$. 
	\end{itemize}
Clearly, for each $e \in U_E$ and $i \in {\mathbb N}$, $\preceq_{e,i}$ is an equivalence relation, so there exists a partition $\P_{e,i}$ of $X_e$ such that each member of $\P_{e,i}$ is an equivalence class of $\preceq_{e,i}$. 
Note that for each $e \in U_E$, $\P_{e,0}=\{\{v\}: v\in X_e\}$.

Note that for each $e \in U_E$ and $i \in {\mathbb N} \cup \{0\}$, $\P_{e,i}$ is a refinement of $\P_{e,i+1}$, and $\lvert \P_{e,i} \rvert \leq \lvert X_e \rvert$.
Let $\mu$ be a nonnegative real number.
Let $\theta$ be a positive integer.
For each $e \in U_E$, let $(H_e,\phi_{H_e},I_e,B_e)$ be the $(\ell,\epsilon,\theta,\mu)$-hierarchy of $\P_e$, where $\P_e=(\P_{e,0},\P_{e,1},...)$.  

\begin{figure} 
	\begin{picture}(100,200) (10,-20)

		%T
		\put(70,170){{$T$}}
		\put(40,150){\line(1,0){70}}
		\put(40,150){\line(0,-1){70}}
		\put(110,150){\line(0,-1){70}}
		\put(40,80){\line(1,0){70}}
		\put(70,120){{$T_0$}}

		\put(40,50){\line(1,0){30}}
		\put(70,50){\line(0,-1){70}}
		\put(40,50){\line(0,-1){70}}
		\put(40,-20){\line(1,0){30}}
		\put(50,10){{$T_{e_1}$}}
		\put(55,50){\line(0,1){30}}
		\put(45,62){{$e_1$}}

		\put(80,50){\line(1,0){30}}
		\put(110,50){\line(0,-1){70}}
		\put(80,50){\line(0,-1){70}}
		\put(80,-20){\line(1,0){30}}
		\put(90,10){{$T_{e_2}$}}
		\put(95,50){\line(0,1){30}}
		\put(84,62){{$e_2$}}

		%G
		\put(220,170){{$G$}}
		\put(170,150){\line(1,0){110}}
		\put(170,150){\line(0,-1){90}}
		\put(280,150){\line(0,-1){90}}
		\put(170,60){\line(1,0){110}}
		\put(220,120){{$X_{T_0}$}}

		\put(180,70){\line(1,0){50}}
		\put(180,70){\line(0,-1){90}}
		\put(230,70){\line(0,-1){90}}
		\put(180,-20){\line(1,0){50}}
		\put(195,10){{$X_{T_{e_1}}$}}
		\multiput(190,65)(15,0){3}{\circle*{5}}
		
		\put(240,70){\line(1,0){30}}
		\put(240,70){\line(0,-1){90}}
		\put(270,70){\line(0,-1){20}}
		\put(240,-20){\line(1,0){83}}
		\put(270,50){\line(1,0){53}}
		\put(323,50){\line(0,-1){70}}
		\put(265,10){{$X_{T_{e_2}}$}}
		\multiput(247,65)(15,0){2}{\circle*{5}}

		%condensation
		\put(370,170){{condensation}}
		\put(345,150){\line(1,0){110}}
		\put(345,150){\line(0,-1){90}}
		\put(455,150){\line(0,-1){90}}
		\put(345,60){\line(1,0){110}}
		\put(395,120){{$X_{T_0}$}}

		\put(355,70){\line(1,0){50}}
		\put(355,70){\line(0,-1){80}}
		\put(405,70){\line(0,-1){80}}
		\put(355,-10){\line(1,0){50}}
		\put(357,10){{hierarchy}}
		\put(357,-2){{for $\P_{e_1}$}}
		\multiput(365,65)(15,0){3}{\circle*{5}}
		\put(372.5,50){\circle*{5}}
		\put(395,50){\circle*{5}}
		\thicklines
		\put(372.5,50){\line(-1,2){7.5}}
		\put(372.5,50){\line(1,2){7.5}}
		\put(395,50){\line(0,1){15}}
		\put(383.75,35){\circle*{5}}
		\put(383.75,35){\line(-3,4){11}}
		\put(383.75,35){\line(3,4){11}}

		\thinlines
		\put(415,70){\line(1,0){30}}
		\put(415,70){\line(0,-1){63}}
		\put(445,70){\line(0,-1){20}}
		\put(415,7){\line(1,0){83}}
		\put(445,50){\line(1,0){53}}
		\put(498,50){\line(0,-1){43}}
		\put(416,30){{$N_{(G,\phi)[X_{T_{e_2}]}}^{\leq \ell}[X_{e_2}]$}}
		\put(416,12){{$+$ some edges}}
		\multiput(422,65)(15,0){2}{\circle*{5}}
		
	\end{picture}
	\caption{An example of the $(U_E,U_E',\ell,\theta,\mu)$-condensation of $(G,\phi,T,\X)$, where $U_E=\{e_1,e_2\}$ and $U_E'=\{e_1\}$.} \label{fig_condensation}
\end{figure}
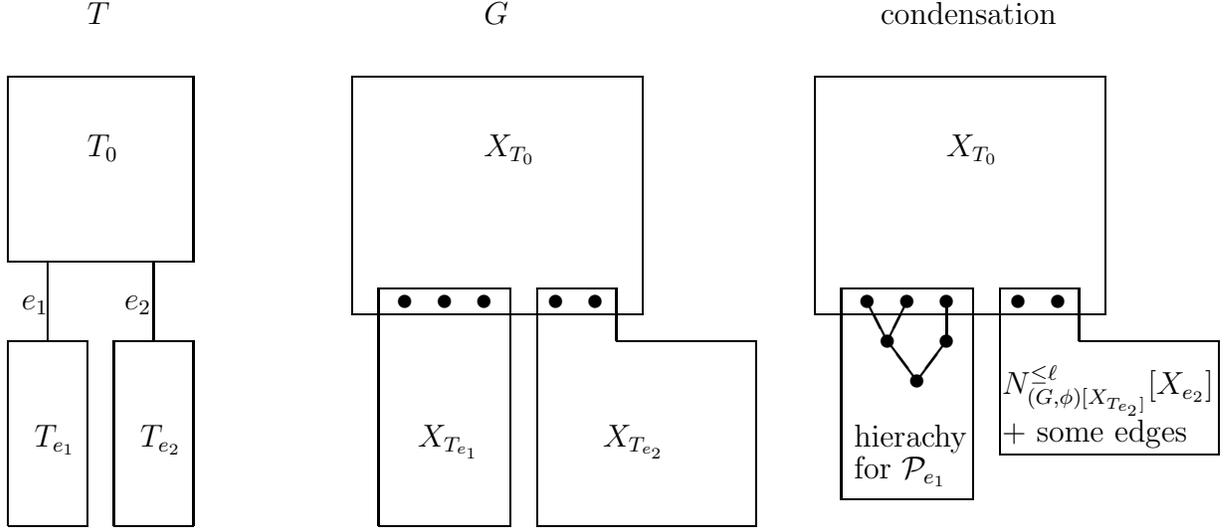

Let $U_E'$ be a subset of $U_E$ such that $\lvert X_e \rvert \leq \theta$ for every $e \in U_E'$.
Let $T_0$ be the component of $T-U_E$ containing $t^*$.
The \defn{$(U_E,U_E',\ell,\theta,\mu)$-condensation of $(G,\phi,T,\X)$} is the weighted graph $(G_0,\phi_0)$ obtained from 
$$(G,\phi)[X_{T_0} \cup \allowbreak \bigcup_{e \in U_E-U_E'}N_{(G,\phi)[X_{T_e}]}^{\leq \ell}[X_e]]$$ 
by 
	\begin{itemize}
		\item for each $e \in U_E'$, creating a copy of $(H_e,\phi_{H_e})$ and then for each $x \in X_e$, identifying $x$ with $v_{\{x\}} \in B_e \subseteq V(H_e)$,  
		\item for each $e \in U_E-U_E'$, for each pair of vertices $u,v \in N_{(G,\phi)[X_{T_e}]}^{\leq \ell}[X_e]-X_e$ with distance in $(G,\phi)[X_{T_e}]$ at most $3\ell+\mu$, adding an edge $e'$ between $u$ and $v$ and defining $\phi_0(e') = \dist_{(G,\phi)[X_{T_e}]}(u,v)$. 
	\end{itemize}
(See Figure \ref{fig_condensation} for an example.)
Note that if $(G,\phi)$ is $(0,\ell]$-bounded and $U_E=U_E'$, then the $(U_E,U_E',\ell,\theta,\mu)$-condensation of $(G,\phi,T,\X)$ is a $(0,\ell]$-bounded weighted graph by Lemma \ref{hierachy_distance}. 

\begin{lemma} \label{con_qu_iso}
Let $\ell$ be a positive real number.
Let $\theta$ be a positive integer.
Let $\mu$ be a nonnegative real number.
Let $(G,\phi)$ be a weighted graph. 
Let $(T,\X)$ be a rooted tree-decomposition such that $(G,\phi)$ is $(T,\X,\ell)$-bounded. 
Let $U_E$ be a subset of $E(T)$ such that edges in $U_E$ are pairwise incomparable.
Let $U_E'$ be a subset of $U_E$ such that $\lvert X_e \rvert \leq \theta$ for every $e \in U_E'$.
Let $(G_0,\phi_0)$ be the $(U_E,U_E',\ell,\theta,\mu)$-condensation of $(G,\phi,T,\X)$.
Then the following statements hold.
	\begin{enumerate}
		\item For every path $P$ in $G$ between two distinct vertices in $V(G_0) \cap V(G)$ with $\leng_{(G,\phi)}(P) \leq 3\ell+\mu$, there exists a path $\overline{P}$ in $G_0$ having the same ends as $P$ with $\leng_{(G_0,\phi_0)}(\overline{P}) \leq \leng_{(G,\phi)}(P)$. 
		\item For every path $P$ in $G_0$ between two distinct vertices $x,y \in V(G_0) \cap V(G)$, there exists a path $\widehat{P}$ in $G$ between $x$ and $y$ such that 
			\begin{enumerate}
				\item $\leng_{(G,\phi)}(\widehat{P}) \leq 4(3\theta+1)(\theta+\frac{\mu}{\ell}) \cdot \leng_{(G_0,\phi_0)}(P)$, and
				\item if $U_E'=\emptyset$, then $\leng_{(G,\phi)}(\widehat{P}) \leq \leng_{(G_0,\phi_0)}(P)$.	
			\end{enumerate}	
	\end{enumerate}
\end{lemma}

\begin{pf}
Let $\epsilon = \min_{e \in E(G)}\phi(e)$.
We define $T_0,X_e,\preceq_{e,i},\P_{e,i},T_e,(H_e,\phi_{H_e},I_e,B_e)$ as mentioned in the definition of the $(U_E,U_E',\ell,\theta,\mu)$-condensation of $(G,\phi,T,\X)$.

We first prove Statement 1.
We need the following claim.

\medskip

\noindent{\bf Claim 1:} For any $e \in U_E$, distinct vertices $x,y \in X_e$ and path $Q$ in $G[X_{T_e}]$ between $x$ and $y$ internally disjoint from $X_e$ with $\leng_{(G,\phi)}(Q) \leq 3\ell+\mu$, there exists a path $Q'$ in $G_0$ between $x$ and $y$ with $\leng_{(G_0,\phi_0)}(Q') \leq \leng_{(G,\phi)}(Q)$.

\noindent{\bf Proof of Claim 1:}
Let $e \in U_E$.
Let $x,y$ be distinct vertices in $X_e$, and $Q$ be a path $Q$ in $G[X_{T_e}]$ between $x$ and $y$ internally disjoint from $X_e$ with $\leng_{(G,\phi)}(Q) \leq 3\ell+\mu$. 

We first assume that $e \in U_E'$. 
Note that $V(Q) \subseteq N_{(G,\phi)[X_{T_e}]}^{\leq \leng_{(G,\phi)}(Q)}[X_e] \subseteq N_{(G,\phi)[X_{T_e}]}^{\leq \lceil \frac{\leng_{(G,\phi)}(Q)}{\epsilon} \rceil \epsilon}[X_e]$.
Since the length in $(G,\phi)$ of each edge of $Q$ is at most $\leng_{(G,\phi)}(Q) \leq \lceil \frac{\leng_{(G,\phi)}(Q)}{\epsilon} \rceil \epsilon$, we have $x \preceq_{i,\lceil \frac{\leng_{(G,\phi)}(Q)}{\epsilon} \rceil} y$.
Hence $x$ and $y$ are contained in the same part of $\P_{e,\lceil \frac{\leng_{(G,\phi)}(Q)}{\epsilon} \rceil}$.
So there exists a path $Q'$ in $G_0$ between $x$ and $y$ with 
	\begin{align*}
		\leng_{(G_0,\phi_0)}(Q') \leq & 2 \cdot \lceil \frac{\leng_{(G,\phi)}(Q)}{\epsilon} \rceil \cdot \frac{\epsilon}{8(\theta+\frac{\mu}{\ell})} \\
		\leq & (\frac{\leng_{(G,\phi)}(Q)}{\epsilon}+1) \cdot \frac{\epsilon}{4\theta} = \frac{\leng_{(G,\phi)}(Q)+\epsilon}{4\theta} \leq \frac{\leng_{(G,\phi)}(Q)}{2\theta} \leq \leng_{(G,\phi)}(Q).
	\end{align*}

So we may assume that $e \in U_E-U_E'$. 
Since $x \neq y$, we know that $x$ has a neighbor $x'$ in $Q$, and $y$ has a neighbor $y'$ in $Q$.
Since $Q$ is internally disjoint from $X_e$, we know that $x'$ and $y'$ are not in $X_e$.
Since $(G,\phi)$ is $(T,\X,\ell)$-bounded, $\phi(xx') \leq \ell$ and $\phi(yy') \leq \ell$.
So $\{x',y'\} \subseteq N_{(G,\phi)[X_{T_e}]}^{\leq \ell}[X_e]-X_e$.
Note that the subpath $P$ of $Q$ between $x'$ and $y'$ is contained in $(G,\phi)[X_{T_e}]$, and $\leng_{(G,\phi)[X_{T_e}]}(P) \leq \leng_{(G,\phi)}(Q) \leq 3\ell+\mu$.
So there exists an edge $x'y'$ of $G_0$ with $\phi_0(x'y') \leq \dist_{(G,\phi)[X_{T_e}]}(x',y') \leq \leng_{(G,\phi)}(P)$.
Therefore, $xx'y'y$ is a path in $G_0$ between $x$ and $y$ of length in $(G_0,\phi_0)$ at most $\phi(xx')+\leng_{(G,\phi)}(P)+\phi(yy') = \leng_{(G,\phi)}(Q)$.
$\Box$

\medskip

For every path $P$ in $G$ between two distinct vertices in $V(G_0) \cap V(G)$ with $\leng_{(G,\phi)}(P) \leq 3\ell+\mu$, we can replace each maximal subpath $Q$ of $P$ between $X_e$ with $\lvert V(Q) \rvert \geq 3$ whose all internal vertices are in $X_{T_e}-X_e$ for some $e \in U_E'$ by the path $Q'$ mentioned in Claim 1 to obtain a walk $\overline{P}$ in $G_0$ with the same ends as $P$ such that $\leng_{(G_0,\phi_0)}(\overline{P}) \leq \leng_{(G,\phi)}(P)$. 
This proves Statement 1.

\medskip

\noindent{\bf Claim 2:} For any $e \in U_E'$, if $a$ and $b$ are two vertices in $X_e$, and $Q$ is a path in $G_0[V(H_e)]$ between $a$ and $b$ internally disjoint from $X_e$ with $\lvert V(Q) \rvert \geq 3$, then there exists a path $\widehat{Q}$ in $(G,\phi)[X_{T_e}]$ between $a$ and $b$ with $\leng_{(G,\phi)}(\widehat{Q}) \leq 4(3\theta+1)(\theta+\frac{\mu}{\ell}) \cdot \leng_{(G_0,\phi_0)}(Q)$. 

\noindent{\bf Proof of Claim 2:}
Let $e \in U_E'$. 
Let $a,b \in X_e$.
Let $Q$ be a path in $G_0[V(H_e)]$ from $a$ to $b$ internally disjoint from $X_e$ with at least one internal vertex.
Recall that $H_e$ is a forest.
So there exist $i_Q \in I_e$ and $Y_Q \in \P_{e,i_Q}$ such that $Q$ is the union of a path in $H_e$ from $a$ to $v_{Y_Q}$ and a path in $H_e$ from $v_{Y_Q}$ to $b$.
Hence $a,b \in Y_Q$, and $\leng_{(G_0,\phi_0)}(Q) = 2i_Q \cdot \frac{\epsilon}{8(\theta+\frac{\mu}{\ell})}$.
So $i_Q\epsilon = 4(\theta+\frac{\mu}{\ell}) \cdot \leng_{(G_0,\phi_0)}(Q)$.

Since $\{a,b\} \subseteq Y_Q \in \P_{e,i_Q}$, there exists a path $Q'$ in $(G,\phi)[X_{T_e}]$ between $a$ and $b$ such that $V(Q') \subseteq N_{(G,\phi)[X_{T_e}]}^{\leq i_Q\epsilon}[X_e]$, and the length in $(G,\phi)$ of every edge of $Q'$ is at most $i_Q\epsilon$.
We choose such $Q'$ such that $\leng_{(G,\phi)}(Q')$ is as small as possible.

Since $V(Q') \subseteq N_{(G,\phi)[X_{T_e}]}^{\leq i_Q\epsilon}[X_e]$, we know that for every vertex $z \in V(Q')$, there exists $v_z \in X_e$ with $z \in N_{(G,\phi)[X_{T_e}]}^{\leq i_Q\epsilon}[\{v_z\}]$, so there exists a path $P_z$ in $(G,\phi)[X_{T_e}]$ between $v_z$ and $z$ with $V(P_z) \subseteq N_{(G,\phi)[X_{T_e}]}^{\leq i_Q\epsilon}[\{v_z\}]$ and with $\leng_{(G,\phi)}(P_z) \leq i_Q\epsilon$.
For any two vertices $z_1,z_2$ of $Q'$ with $v_{z_1}=v_{z_2}$, since $P_{z_1} \cup P_{z_2}$ has length in $(G,\phi)$ at most $2i_Q\epsilon$, the minimality of the length of $Q'$ implies that the subpath of $Q'$ between $z_1$ and $z_2$ has length in $(G,\phi)$ at most $2i_Q\epsilon$.

Since $e \in U_E'$, we have $\lvert X_e \rvert \leq \theta$.
So the vertices of $Q'$ can be covered by at most $\theta$ subpaths of $Q'$ of length in $(G,\phi)$ at most $2i_Q\epsilon$.
Since each edge of $Q'$ has length in $(G,\phi)$ at most $i_Q\epsilon$, the length in $(G,\phi)$ of $Q'$ is at most $\theta \cdot 2i_Q\epsilon + (\theta+1) \cdot i_Q\epsilon = (3\theta+1)i_Q\epsilon= (3\theta+1) \cdot 4(\theta+\frac{\mu}{\ell}) \cdot \leng_{(G_0,\phi_0)}(Q)$.
So $Q'$ is a desired path $\widehat{Q}$ for this claim.
$\Box$

\medskip

\noindent{\bf Claim 3:} For each $e \in U_E-U_E'$, if $a$ and $b$ are two vertices in $X_{T_e}$, and $Q$ is a path in $G_0[N_{(G,\phi)[X_{T_e}]}^{\leq \ell}[X_e]]$ between $a$ and $b$, then there exists a path $\widehat{Q}$ in $(G,\phi)[X_{T_e}]$ between $a$ and $b$ with $\leng_{(G,\phi)}(\widehat{Q}) \leq \leng_{(G_0,\phi_0)}(Q)$. 

\noindent{\bf Proof of Claim 3:}
Let $a,b \in X_{T_e}$.
Let $Q$ be a path in $G_0[N_{(G,\phi)[X_{T_e}]}^{\leq \ell}[X_e]]$ between $a$ and $b$.
Note that for each edge $e'$ of $Q-E(G)$, both ends of $e'$ are in $N_{(G,\phi)[X_{T_e}]}^{\leq \ell}[X_e]-X_e$, and there exists a path in $(G,\phi)[X_{T_e}]$ connecting the ends of $e'$ with length in $(G,\phi)[X_{T_e}]$ equal to $\phi_0(e')$. 
Hence there exists a path $\widehat{Q}$ in $(G,\phi)[X_{T_e}]$ between $a$ and $b$ of length in $(G,\phi)$ at most $\leng_{(G_0,\phi_0)}(Q)$. 
$\Box$

\medskip

Now we prove Statement 2.
Let $x,y \in V(G) \cap V(G_0)$.
Let $P$ be a path in $G_0$ between $x$ and $y$.
We say that a subpath $P'$ of $P$ is ``special'' if there exists $e_{P'} \in U_E$ such that either 
	\begin{itemize}
		\item $e_{P'} \in U_E'$ and $P'$ is a path between two vertices in $X_{e_{P'}}$ internally disjoint from $X_{e_{P'}}$ such that it has at least one internal vertex and all its internal vertices are not in $X_{T_0}$, or 
		\item $e_{p'} \in U_E-U_E'$ such that $P' \subseteq G_0[N_{(G,\phi)[X_{T_e}]}^{\leq \ell}[X_e]]$.
	\end{itemize}
Note that for each special subpath $P'$ of $P$, the path $\widehat{P'}$ in $(G,\phi)$ is defined in Claims 2 or 3.
Then by replacing each maximal special subpath $P'$ of $P$ by $\widehat{P'}$, we obtain a walk in $(G,\phi)$ between $x$ and $y$ of length in $(G,\phi)$ at most $\max\{4(3\theta+1)(\theta+\frac{\mu}{\ell}), 1\} \cdot \leng_{(G_0,\phi_0)}(P) = 4(3\theta+1)(\theta+\frac{\mu}{\ell}) \cdot \leng_{(G_0,\phi_0)}(P)$; moreover, if $U_E'=\emptyset$, then this walk in $(G,\phi)$ has length in $(G,\phi)$ at most $\leng_{(G_0,\phi_0)}(P)$.
\end{pf}

\bigskip

Let $c$ and $c'$ be functions.
We say that $c'$ \defn{can be obtained by extending} $c$ if the domain of $c'$ contains the domain of $c$ and $c'(x)=c(x)$ for every $x$ in the domain of $c$.

The following lemma (Lemma \ref{con_color}) is the main result of this subsection.
It roughly says that if we have a weighted graph with a rooted tree-decomposition and define the ``central part'' of the graph to be roughly the subgraph induced by the bags not contained in a descendant of an edge of a set $U_E$, then any good coloring of the condensation can give us a good coloring of the ``central part''.

Let us first sketch the proof of Lemma \ref{con_color}. 
We first show that the condensation is a good ``approximation'' of the central part (i.e. $(G,\phi)^\ell[X_{T_0} \cup \bigcup_{e \in U_E}N_{(G,\phi)[X_{T_e}]}^{\leq 3\ell}[X_e]]-R$).
Claims 1 and 2 are for this purpose and will be used later as well.
Then we use the coloring $c_0$ of the condensation to obtain a coloring $c$ of the central part.
The remaining work is to prove that $c$ and any extension $c'$ of $c$ have small weak diameter.
We prove that every monochromatic component that is far from $X_{T_0}$ has small weak diameter in Claim 3.
Then we consider a monochromatic component $M$ that is not far from $X_{T_0}$.
Claims 4-6 show the ``position'' of $M$.
Then we ``project'' $M$ to the condensation to obtain a subgraph $M'$ in the condensation.
We show that $M'$ has small weak diameter in Claim 8 by using Claim 7, which helps us translate the distance in $G$ to the distance in the condensation. 
Finally we show how to use the bound for the weak diameter of $M'$ to bound the weak diameter of $M$ according to how we construct the coloring $c$ based on the coloring $c_0$ of the condensation in the rest of the proof.

\begin{lemma} \label{con_color}
For any positive real numbers $\ell,\ell',\nu,\mu$, nonnegative real number $\lambda$ and positive integers $m,\theta$ with $m \geq 2$, there exists a real number $\nu^* \geq \nu$ such that the following holds.

Let $(G,\phi)$ be a $(0,\ell']$-bounded weighted graph.
Let $(T,\X)$ be a rooted tree-decomposition of $G$ such that $(G,\phi)$ is $(T,\X,\ell)$-bounded.
Denote $\X$ by $(X_t: t \in V(T))$, and for each $tt' \in E(T)$, let $X_{tt'}=X_t \cap X_{t'}$.
Let $U_E \subseteq E(T)$ be a subset of pairwise incomparable edges.
Let $U_E'$ be a subset of $U_E$ such that $\lvert X_e \rvert \leq \theta$ for every $e \in U_E'$.
Let $(G_0,\phi_0)$ be the $(U_E,U_E',\ell,\theta,\lambda)$-condensation of $(G,\phi,T,\X)$.
Let $R \subseteq V(G)$.
Assume that for every $e \in U_E$ with $\lvert X_e \rvert >\theta$, there exists $A_e \subseteq X_e$ with $\lvert A_e \rvert \leq \theta$ such that $X_e \subseteq N_{(G,\phi)}^{\leq \mu}[A_e]$.

If there exists an $m$-coloring $c_0$ of $(G_0,\phi_0)^\ell[V(G_0)]-R$ with weak diameter in $(G_0,\phi_0)^\ell$ at most $\nu$, then there exists an $m$-coloring $c$ of $(G,\phi)^\ell[X_{T_0} \cup \bigcup_{e \in U_E}N_{(G,\phi)[X_{T_e}]}^{\leq 3\ell}[X_e]]-R$ with weak diameter in $(G,\phi)^\ell$ at most $\nu^*$ such that  
	\begin{enumerate}
		\item $c(v)=c_0(v)$ for every $v \in V(G_0) \cap V(G)-R$, and
		\item for every $m$-coloring $c'$ of $(G,\phi)^\ell[V(G)]-R$ that can be obtained by extending $c$, and for every $c'$-monochromatic component $M$ in $(G,\phi)^\ell[V(G)]-R$, if $V(M) \cap (X_{T_0} \cup \allowbreak \bigcup_{e \in U_E}N_{(G,\phi)[X_{T_e}]}^{\leq \ell}[X_e]) \allowbreak \neq \emptyset$, then $M$ has weak diameter in $(G,\phi)^\ell$ at most $\nu^*$.
	\end{enumerate}
\end{lemma}

\begin{pf}
Let $\ell,\ell',\nu,\mu$ be positive real numbers, $\lambda$ be a nonnegative real number, and $m,\theta$ be positive integers with $m \geq 2$. 
Let $\nu_0$ be the real number $\nu^*$ given by Lemma \ref{patching_centered_set} by taking $(k,r,\ell,\nu)=(\theta,3\ell+\mu,\ell',1)$. 
Let $\nu_1 = \lceil \frac{2\ell' \nu_0}{\ell} \rceil$.
Define $\nu^* = \lceil \nu_1 + (24+\frac{8\mu+4\ell'}{\ell})\theta + (8(\nu+3)(3\theta+1)(\theta+\frac{\lambda}{\ell})+4) \rceil$.

Let $(G,\phi),(T,\X),U_E,U_E',(G_0,\phi_0),R,c_0$ be as defined in the lemma.
Let $\epsilon= \min_{e \in E(G)}\phi(e)$.
For every $e \in U_E$ and $i \in {\mathbb N}$, let $\preceq_{e,i}$ be the relation, $\P_{e,i}$ be the partition, and $(H_e,\phi_{H_e},I_e,B_e)$ be the $(\ell,\epsilon,\theta,\lambda)$-hierarchy of $\P_e=(\P_{e,0},\P_{e,1},...)$ as stated in the definition of the $(U_E,U_E',\ell,\theta,\lambda)$-condensation. 

Let $t^*$ be the root of $T$.
For every $e \in E(T)$, let $T_e$ be the component of $T-e$ disjoint from $t^*$.

For $i \in [3]$, let $Z_i = \bigcup_{e \in U_E}N_{(G,\phi)[X_{T_e}]}^{\leq i\ell}[X_e]$.
Note that $N_{(G,\phi)[X_{T_e}]}^{\leq i\ell}[X_e] \subseteq X_{T_e}$ for each $e \in U_E$.

\medskip

\noindent{\bf Claim 1:} For every $e \in U_E$ and every path $P$ in $(G,\phi)$ with $\leng_{(G,\phi)}(P) \leq \ell$ between $u \in X_{T_0}$ and $v \in X_{T_e}-X_e$, we have $v \in Z_1 \cap X_{T_e}-X_e$. 

\noindent{\bf Proof of Claim 1:}
Since $u \in X_{T_0}$, there exists $t \in V(T_0)$ such that $u \in X_t$.
Since $P$ is between $u \in X_t$ and $v \in X_{T_e}-X_e$, $P$ contains a subpath $P'$ between $X_e$ and $v$ internally disjoint from $X_e$.
So $v \in N_{(G,\phi)[X_{T_e}]}^{\leq \ell}[X_e]-X_e \subseteq Z_1 \cap X_{T_e}-X_e$.
$\Box$

\medskip

For every $v \in Z_3-X_{T_0}$, 
	\begin{itemize}
		\item we define $e_v$ to be the element $e$ in $U_E$ with $v \in X_{T_e}-X_e$, and
		\item if $e_v \in U_E'$, then 
			\begin{itemize}
				\item define $p_v = \min\{j \in [\lceil \frac{3\ell}{\epsilon} \rceil]:v \in N_{(G,\phi)[X_{T_{e_v}}]}^{\leq j\epsilon}[X_{e_v}]\}$, and  
				\item define $j_v$ to be the largest element in $I_{e_v}$ with $j_v \leq p_v$;
			\end{itemize}
	\end{itemize}
Note that $p_v$ and $j_v$ are undefined when $e_v \not \in U_E'$.
For every $v \in Z_3-X_{T_0}$ with $e_v \in U_E'$, since $Z_3 = \bigcup_{e \in U_E}N_{(G,\phi)[X_{T_e}]}^{\leq 3\ell}[X_e]$ and $\epsilon=\min_{e \in E(G)}\phi(e)>0$, we know that the set mentioned in the definition of $p_v$ is nonempty, so $p_v \geq 1$; this implies that $j_v$ is a well-defined positive integer since $1 \in I_{e_v}$; since $1 \leq j_v \leq p_v \leq \lceil \frac{3\ell}{\epsilon} \rceil \leq \lceil \frac{3\ell+\lambda}{\epsilon} \rceil$, by the definition of $I_{e_v}$, we known that $\P_{e_v,\beta}$ are identical for all $\beta$ with $j_v \leq \beta \leq p_v$, so $\P_{e_v,p_v}=\P_{e_v,j_v}$.

\medskip

\noindent{\bf Claim 2:} For every $v \in Z_3-X_{T_0}$ with $e_v \in U_E'$, there exists $Y_v \in \P_{e_v,p_v}=\P_{e_v,j_v}$ such that 
	\begin{itemize}
		\item $v \in N_{(G,\phi)[X_{T_{e_v}}]}^{\leq p_v\epsilon}[Y_v]$, 
		\item for every $Y' \in \P_{e_v,j_v}-\{Y_v\}$, we have $v \not \in N_{(G,\phi)[X_{T_{e_v}}]}^{\leq p_v\epsilon}[Y']$, and
		\item there exists a path $P$ in $(G,\phi)[X_{T_{e_v}}]$ from $v$ to $Y_v$ such that $V(P) \subseteq N_{(G,\phi)[X_{T_{e_v}}]}^{\leq p_v\epsilon}[Y_v]$, $\leng_{(G,\phi)}(P) \leq p_v\epsilon$, and every edge of $P$ has length in $(G,\phi)$ at most $p_v\epsilon$.
	\end{itemize}

\noindent{\bf Proof of Claim 2:}
By the definition of $p_v$, we have $v \in N_{(G,\phi)[X_{T_e}]}^{\leq p_v\epsilon}[X_{e_v}] \cap X_{T_{e_v}}-X_{e_v}$, so there exists a path $P$ in $G[X_{T_e}]$ from $v$ to $X_{e_v}$ internally disjoint from $X_{e_v}$ of length in $(G,\phi)[X_{T_e}]$ at most $p_v\epsilon$.
Let $y$ be the vertex in $V(P) \cap X_{e_v}$.
Let $Y_v$ be the member of $\P_{e_v,p_v}$ containing $y$.
So $v \in V(P) \subseteq N_{(G,\phi)[X_{T_{e_v}}]}^{\leq p_v\epsilon}[Y_v]$.
Since $\leng_{(G,\phi)}(P) \leq p_v\epsilon$, the length in $(G,\phi)$ of every edge of $P$ is at most $p_v\epsilon$.

Let $Y'$ be a member of $\P_{e_v,p_v}-\{Y_v\}$.
Suppose that $v \in N_{(G,\phi)[X_{T_{e_v}}]}^{\leq p_v\epsilon}[Y']$.
Then there exists a path $P'$ in $(G,\phi)[X_{T_{e_v}}]$ such that $V(P') \subseteq N_{(G,\phi)[X_{T_{e_v}}]}^{\leq p_v\epsilon}[Y']$, and the length in $(G,\phi)$ of every edge of $P'$ is at most $p_v\epsilon$.
So there exists a path $P''$ in $(G,\phi)[X_{T_{e_v}}]$ from $Y_v$ to $Y'$ internally disjoint from $X_{e_v}$ with $V(P'') \subseteq N_{(G,\phi)[X_{T_{e_v}}]}^{\leq p_v\epsilon}[X_e]$, and the length in $(G,\phi)$ of every edge of $P''$ is at most $p_v\epsilon$.
So $Y_v=Y'$ by the definition of $\P_{e_v,p_v}$, a contradiction.
Hence $v \not \in N_{(G,\phi)[X_{T_{e_v}}]}^{\leq p_v\epsilon}[Y']$.
$\Box$

\medskip

For every $x \in Z_3-V(G_0)$, if $e_x \in U_E'$, then by Claim 2, there exists a unique member $Y_x$ of $\P_{e_x,j_x}$ such that $x \in N_{(G,\phi)[X_{T_{e_x}}]}^{\leq p_x\epsilon}[Y_x]$.
Note that for every $x \in Z_3-V(G_0)$, if $e_x \in U_E'$, then $v_{Y_x} \in V(H_{e_x})-V(G) \subseteq V(G_0)-V(G)$ since $j_x \in I_{e_x}$ and $j_x \geq 1$, so $c_0(v_{Y_x})$ is defined.
For every $x \in Z_1-X_{T_0}$ with $e_x \in U_E-U_E'$, we know $x \in V(G_0) \cap V(G)$ and we define $v_{Y_x}=x$ and $Y_x=\{x\}$.

Let $c_1: (V(G_0) \cap V(G)-R) \cup Z_1 \rightarrow [m]$ be the function such that 
	\begin{itemize}
		\item $c_1(u)=c_0(u)$ for every $u \in (X_{T_0} \cup (Z_1 \cap V(G_0)))-R$, and 
		\item $c_1(u) = c_0(v_{Y_u})$ for every $u \in Z_1-V(G_0)$.
	\end{itemize}
For each $i \in \{2,3\}$, let $c_i: (V(G_0) \cap V(G)-R) \cup Z_i \rightarrow [m]$ be the function such that 
	\begin{itemize}
		\item $c_i(v)=c_{i-1}(v)$ for every $v \in (V(G_0) \cap V(G)-R) \cup Z_{i-1}$, and 
		\item $c_i(v)=i-1$ for every $v \in Z_i-Z_{i-1}$.
	\end{itemize}
Note that the domain of $c_3$ is $(X_{T_0}-R) \cup (Z_3-X_{T_0})=(X_{T_0}-R) \cup (\bigcup_{e \in U_E}N_{(G,\phi)[X_{T_e}]}^{\leq 3\ell}[X_e]-X_{T_0})$. 
Let $c=c_3|_{(X_{T_0} \cup Z_3)-R}$. 

We shall prove that $c$ satisfies the conclusion of this lemma.
Clearly, for every $v \in V(G_0) \cap V(G)-R$, $c(v)=c_0(v)$.

\medskip

\noindent{\bf Claim 3:} For every $c$-monochromatic component $M$ in $(G,\phi)^\ell[X_{T_0} \cup \bigcup_{e \in U_E}N_{(G,\phi)[X_{T_e}]}^{\leq 3\ell}[X_e]]-R$ disjoint from $X_{T_0} \cup \bigcup_{e \in U_E}N_{(G,\phi)[X_{T_e}]}^{\leq \ell}[X_e]$, $M$ has weak diameter in $(G,\phi)^\ell$ at most $\nu^*$.

\noindent{\bf Proof of Claim 3:}
Note that $V(M) \subseteq X_{T_0} \cup \bigcup_{e \in U_E}N_{(G,\phi)[X_{T_e}]}^{\leq 3\ell}[X_e]-R$, so $V(M) \subseteq V(G)$.
Since $M$ is disjoint from $X_{T_0} \cup \bigcup_{e \in U_E}N_{(G,\phi)[X_{T_e}]}^{\leq \ell}[X_e]$, we know $V(M) \subseteq Z_3-Z_1$.
In particular, there exists $e_M \in U_E$ such that $V(M) \cap N_{(G,\phi)[X_{T_{e_M}}]}^{\leq 3\ell}[X_{e_M}]-X_{e_M} \neq \emptyset$. 
Since $(G,\phi)$ is $(T,\X,\ell)$-bounded, there exists no edge of $G$ with length in $(G,\phi)$ greater than $\ell$ between $X_{T_{e_M}}-X_{e_M}$ and $V(G)-(X_{T_{e_M}}-X_{e_M})$.
This together with the fact $V(M) \cap Z_1 = \emptyset$ implies that $V(M) \subseteq N_{(G,\phi)[X_{T_{e_M}}]}^{\leq 3\ell}[X_{e_M}]$.

By the assumption of this lemma, for every $e \in U_E$ with $|X_e| > \theta$, there exists $A_e \subseteq X_e$ with $|A_e| \leq \theta$ such that $X_e \subseteq N_{(G,\phi)}^{\leq \mu}[A_e]$.
So for every $e \in U_E$, there exists $A_e \subseteq X_e$ with $|A_e| \leq \theta$ such that $X_e \subseteq N_{(G,\phi)}^{\leq \mu}[A_e]$.
Hence $V(M) \subseteq N_{(G,\phi)}^{\leq 3\ell}[X_{e_M}] \subseteq N_{(G,\phi)}^{\leq 3\ell+\mu}[A_{e_M}]$.
By applying Lemma \ref{patching_centered_set} by taking $(k,r,\ell,\nu,(G,\phi),S,Z,m,R,c_Z,c) = (\theta,3\ell+\mu,\ell',1,(G,\phi),A_{e_M},V(M),m,V(G)-V(M),c|_{V(M)},\emptyset)$, where the last entry $\emptyset$ denotes the coloring with domain $V(G)-((V(G)-V(M)) \cup V(M))=\emptyset$, we know that $c|_{V(M)}$ is a coloring of $(G,\phi)^{\ell'}-(V(G)-V(M))=(G,\phi)^{\ell'}[V(M)]$ with weak diameter in $(G,\phi)^{\ell'}$ at most $\nu_0$.
So $V(M)$ has weak diameter in $(G,\phi)$ at most $\ell' \nu_0$.
By Lemma \ref{dist_original_power}, $V(M)$ has weak diameter in $(G,\phi)^\ell$ at most $\lceil \frac{2\ell' \nu_0}{\ell} \rceil =\nu_1 \leq \nu^*$.
So $M$ has weak diameter in $(G,\phi)^\ell$ at most $\nu^*$.
$\Box$

\medskip

By Claim 3, to prove this lemma, it suffices to prove that for every $m$-coloring $c'$ of \linebreak $(G,\phi)^\ell[V(G)] -R$ that can be obtained by extending $c$, and for every $c'$-monochromatic component $M$ in $(G,\phi)^\ell[V(G)]-R$, if $V(M) \cap (X_{T_0} \cup \bigcup_{e \in U_E}N_{(G,\phi)[X_{T_e}]}^{\leq \ell}[X_e]) \neq \emptyset$, then $M$ has weak diameter in $(G,\phi)^\ell$ at most $\nu^*$.

Suppose to the contrary that there exists an $m$-coloring $c'$ of $(G,\phi)^\ell[V(G)]-R$ that can be obtained by extending $c$, and there exists a $c'$-monochromatic component $M$ in $(G,\phi)^\ell[V(G)]-R$ with $V(M) \cap (X_{T_0} \cup \bigcup_{e \in U_E}N_{(G,\phi)[X_{T_e}]}^{\leq \ell}[X_e]) \neq \emptyset$ such that $M$ has weak diameter in $(G,\phi)^\ell$ greater than $\nu^*$.

\medskip

\noindent{\bf Claim 4:} There exists no $e \in U_E$ with $V(M)$ intersecting both $(Z_1-X_e) \cap X_{T_e}$ and $(Z_3-Z_2) \cap X_{T_e}$. 

\noindent{\bf Proof of Claim 4:}
Suppose to the contrary that there exists $e \in U_E$ such that $V(M) \cap (Z_1-X_e) \cap X_{T_e} \neq \emptyset$ and $V(M) \cap (Z_3-Z_2) \cap X_{T_e} \neq \emptyset$.
Then there exists a path $P$ in $M \subseteq (G,\phi)^\ell[V(G)]$ from $V(M) \cap (Z_1-X_e) \cap X_{T_e}$ to $V(M) \cap (Z_3-Z_2) \cap X_{T_e}$. 
Hence $P$ contains an edge $e_1$ between $V(M) \cap (Z_1-X_e) \cap X_{T_e}$ and $V(M) \cap (Z_2-Z_1) \cap X_{T_e}$, and $P$ contains an edge $e_2$ between $V(M) \cap (Z_2-Z_1) \cap X_{T_e}$ and $V(M) \cap (Z_3-Z_2) \cap X_{T_e}$.
In particular, $M$ contains a vertex $v_1$ in $(Z_2-Z_1) \cap X_{T_e}$ and a vertex $v_2$ in $(Z_3-Z_2) \cap X_{T_e}$.
By the definition of the colorings $c',c_2$ and $c_3$, we know that $c'(v_1)=c_2(v_1)=1$ and $c'(v_2)=c_3(v_2)=2$.
So $M$ is not $c'$-monochromatic, a contradiction.
$\Box$

\medskip

\noindent{\bf Claim 5:} For every $e \in U_E$, if $V(M) \cap X_{T_e}-X_e \neq \emptyset$, then $V(M) \cap Z_1 \cap X_{T_e} -X_e \neq \emptyset$.

\noindent{\bf Proof of Claim 5:}
Let $e \in U_E$ with $V(M) \cap X_{T_e}-X_e \neq \emptyset$.
If $V(M) \cap N_{(G,\phi)[X_{T_e}]}^{\leq \ell}[X_e]-X_e \neq \emptyset$, then we are done.
So we may assume that $V(M) \cap N_{(G,\phi)[X_{T_e}]}^{\leq \ell}[X_e]-X_e = \emptyset$.
It implies that $V(M) \cap X_{T_0} \neq \emptyset$, since $V(M) \cap X_{T_e}-X_e \neq \emptyset$ and $V(M) \cap (X_{T_0} \cup \bigcup_{e \in U_E}N_{(G,\phi)[X_{T_e}]}^{\leq \ell}[X_e]) \neq \emptyset$.
Hence there exists a path $P$ in $M \subseteq (G,\phi)^\ell[V(G)]$ from $V(M) \cap X_{T_0}$ to $V(M) \cap X_{T_e}-X_e$.
So there exists a path $P'$ in $(G,\phi)$ of length in $(G,\phi)$ at most $\ell$ from $X_{T_0}$ to $V(M) \cap X_{T_e}-X_e$.
By Claim 1, $V(P') \cap V(M) \cap Z_1 \cap X_{T_e}-X_e \neq \emptyset$.
$\Box$

\medskip

Let $B = \{e \in U_E: V(M) \cap X_{T_e}-X_e \neq \emptyset\}$.

\medskip

\noindent{\bf Claim 6:} $V(M) \subseteq X_{T_0} \cup \bigcup_{e \in B}(X_{T_e} \cap Z_2)$.

\noindent{\bf Proof of Claim 6:}
For every $e \in B$, by Claim 5, $V(M) \cap Z_1 \cap X_{T_e} -X_e \neq \emptyset$, so by Claim 4, $V(M) \cap (Z_3-Z_2) \cap X_{T_e}=\emptyset$.
Hence for every $e \in B$, since $M \subseteq (G,\phi)^\ell[V(G)]$, we have $V(M) \cap X_{T_e} \subseteq X_{T_e} \cap Z_2$.
Therefore, $V(M) \subseteq X_{T_0} \cup \bigcup_{e \in B}(X_{T_e} \cap Z_2)$.
$\Box$

\medskip

\noindent{\bf Claim 7:} For every path $P$ in $(G,\phi)$ with $\leng_{(G,\phi)}(P) \leq \ell$ between a vertex $a$ in $X_e$ for some $e \in U_E'$ and a vertex $b$ in $X_{T_e}-X_e$ internally disjoint from $X_e$, there exists a path $\overline{P}$ in $(G_0,\phi_0)$ between $a$ and $v_{Y_b}$ with $\leng_{(G_0,\phi_0)}(\overline{P}) \leq \leng_{(G,\phi)}(P)$.

\noindent{\bf Proof of Claim 7:}
Let $P$ be a path in $(G,\phi)$ with $\leng_{(G,\phi)}(P) \leq \ell$ between a vertex $a$ in $X_e$ for some $e \in U_E'$ and a vertex $b$ in $X_{T_e}-X_e$ internally disjoint from $X_e$.
Note that $b \in N_{(G,\phi)[X_{T_e}]}^{\leq \leng_{(G,\phi)}(P)}[X_e]-X_e \subseteq N_{(G,\phi)[X_{T_e}]}^{\leq \ell}[X_e]-X_e \subseteq Z_3-X_{T_0}$.
So $p_b \leq \lceil \frac{\leng_{(G,\phi)}(P)}{\epsilon} \rceil$.
Hence $\dist_{(G_0,\phi_0)}(a,v_{Y_b}) \leq p_b \cdot \frac{\epsilon}{8(\theta+\frac{\lambda}{\ell})} \leq \lceil \frac{\leng_{(G,\phi)}(P)}{\epsilon} \rceil \cdot \frac{\epsilon}{8\theta}  \leq (\frac{\leng_{(G,\phi)}(P)}{\epsilon}+1) \cdot \frac{\epsilon}{8\theta} \leq \frac{\leng_{(G,\phi)}(P)+\epsilon}{8\theta} \leq \leng_{(G,\phi)}(P)$ since $\epsilon = \min_{e \in E(G)}\phi(e) \leq \leng_{(G,\phi)}(P)$.
So there exists a path $\overline{P}$ in $(G_0,\phi_0)$ between $a$ and $v_{Y_b}$ with $\leng_{(G_0,\phi_0)}(\overline{P}) \leq \leng_{(G,\phi)}(P)$.
$\Box$

\medskip

Next, we ``project'' $M$ to $G_0$ to obtain a graph $M'$, and then show that $M'$ has small weak diameter in Claim 8, and finally use it to bound the weak diameter of $M$ to complete the proof. 

Define $M'$ to be the graph as follows.
	\begin{itemize}
		\item Define $V(M') = (V(M) \cap V(G_0) \cap V(G)) \cup \{v_Y:$ there exists $x \in V(M) \cap Z_1-(V(G_0) \cap V(G))$ such that $Y_x=Y\}$.
        	\item Let $E_1(M')=\{xy \in E(M): x,y \in V(M) \cap V(M') \cap V(G_0) \cap V(G)\}$. 
        	\item Let $E_2(M') = \{xv_Y: x \in V(M) \cap V(G_0) \cap V(G),$ there exists $y \in V(M) \cap Z_1-(V(G_0) \cap V(G))$ such that $Y_y=Y$ and $xy \in E(M)\}$.
        	\item Let $E_3(M') = \{v_Yv_{Y'}:$ there exist $x,y \in V(M) \cap Z_1-(V(G_0) \cap V(G))$ with $xy \in E(M)$ such that $Y_x=Y$ and $Y_y=Y'\}$.
        	\item Let $E_4(M') = \{v_Yv_{Y'}:$ there exist $x,y \in V(M) \cap Z_1-(V(G_0) \cap V(G))$ with $e_x=e_y$, $Y_x=Y$, $Y_y=Y'$ and $Y \neq Y'$ such that there exists a path in $M$ between $x$ and $y$ whose all internal vertices are in $(Z_2-Z_1) \cap X_{T_e}-X_e\}$.
        	\item Define $E(M')= E_1(M') \cup E_2(M') \cup E_3(M') \cup E_4(M')$.
	\end{itemize}

\medskip

\noindent{\bf Claim 8:} $M'$ is contained in a $c_0$-monochromatic component in $(G_0,\phi_0)^\ell[V(G_0)]-R$. 

\noindent{\bf Proof of Claim 8:}
Clearly, $V(M') \subseteq V(G_0)-R$.
By the definition of $c_1$, since $M$ is $c'$-monochromatic, $M'$ is $c_0$-monochromatic. 
So to prove this claim, it suffices to prove $E(M') \subseteq E((G_0,\phi_0)^\ell)$.

For every path $P$ in $G$ of length in $(G,\phi)$ at most $\ell$ between two distinct vertices in $V(G_0) \cap V(G)$, we let $\overline{P}$ be the path in $G_0$ mentioned in Statement 1 of Lemma \ref{con_qu_iso}.
For every path $P$ in $G$ of length in $(G,\phi)$ at most $\ell$ between a vertex in $X_e$ for some $e \in U_E'$ and a vertex in $X_{T_e}-X_e$ internally disjoint from $X_e$, we let $\overline{P}$ be the path in $G_0$ mentioned in Claim 7. 

If $xy$ is an edge of $M'$ with $x,y \in V(G_0) \cap V(G)$, then by the definition of $E(M')$, we know $xy \in E_1(M') \subseteq E(M)$, and since $M \subseteq (G,\phi)^\ell[V(G)]$, there exists a path $P_{xy}$ in $G$ of length in $(G,\phi)$ at most $\ell$ between $x$ and $y$, so $\overline{P_{xy}}$ is a path in $G_0$ of length in $(G_0,\phi_0)$ at most $\ell$ between $x$ and $y$, so $xy \in E((G_0,\phi_0)^\ell)$.

Now we assume that there exists an edge $xy$ of $M'$ with $x \in V(G_0) \cap V(G)$ and $y \not \in V(G_0) \cap V(G)$.
So $xy \in E_2(M')$ by the definition of $E(M')$.
Since $y \not \in V(G_0) \cap V(G)$, there exists $y_0 \in V(M) \cap Z_1-(V(G_0) \cap V(G))$ such that $y=v_{Y_{y_0}}$ and $xy_0 \in E(M)$.
Since $M \subseteq (G,\phi)^\ell$, there exists a path $P_{xy_0}$ in $G$ of length in $(G,\phi)$ at most $\ell$ between $x$ and $y_0$.
Let $y'$ be the vertex in $V(P_{xy_0}) \cap X_{e_{y_0}}$ such that the subpath of $P_{xy_0}$ between $y_0$ and $y'$ is contained in $G[X_{T_{e_{y_0}}}]$ and internally disjoint from $X_{e_{y_0}}$.
So $x\overline{P_{xy'}}y'\overline{P_{y_0y'}}v_{Y_{y_0}}$ is a walk in $(G_0,\phi_0)$ of length at most $$\leng_{(G,\phi)}(P_{xy'})+\leng_{(G,\phi)}(P_{y_0y'}) \leq \leng_{(G,\phi)}(P_{xy_0}) \leq \ell,$$ where $P_{xy'}$ is the subpath of $P_{xy_0}$ between $x$ and $y'$, and $P_{y_0y'}$ is the subpath of $P_{xy_0}$ between $y'$ and $y_0$.
So $xy=xv_{Y_{y_0}} \in E((G_0,\phi_0)^\ell)$.

Hence every edge of $M'$ incident with a vertex of $V(G_0) \cap V(G)$ is an edge of $(G_0,\phi_0)^\ell$.
Note that every end of an edge of $M$ not in $V(G_0) \cap V(G)$ is in $X_{T_e}-X_e$ for some $e \in B \cap U_E'$.

Now we assume that there exist distinct edges $e_1,e_2 \in B$ and sets $Y_1 \in \bigcup_{i \in {\mathbb N}}\P_{e_1,i}$ and $Y_2 \in \bigcup_{i \in {\mathbb N}}\P_{e_2,i}$ such that $v_{Y_1}v_{Y_2} \in E(M')$.
So $v_{Y_1}v_{Y_2} \in E_3(M')$, and for each $i \in [2]$, there exists $x_i \in V(M) \cap Z_1 \cap X_{T_{e_i}}-X_{e_i}$ such that $Y_{x_i}=Y_i$, and $x_1x_2 \in E(M) \subseteq E((G,\phi)^\ell[V(G)])$.
Hence there exists a path $P_{x_1x_2}$ in $G$ of length in $(G,\phi)$ at most $\ell$ between $x_1$ and $x_2$.
For each $i \in [2]$, let $y_i$ be the vertex in $V(P_{x_1x_2}) \cap X_{e_i}$ such that the subpath $P_{x_iy_i}$ of $P_{x_1x_2}$ between $x_i$ and $y_i$ is contained in $G[X_{T_{e_i}}]$ and internally disjoint from $X_{e_{i}}$.
Then $v_{Y_1}\overline{P_{x_1y_1}}y_1\overline{P'_{x_1x_2}}y_2\overline{P_{x_2y_2}}v_{Y_2}$ is a walk in $G_0$ of length in $(G_0,\phi_0)$ at most $\leng_{(G,\phi)}(P_{x_1x_2}) \leq \ell$, where $P'_{x_1x_2}$ is the subpath of $P_{x_1x_2}$ between $y_1$ and $y_2$.
Therefore, $v_{Y_1}v_{Y_2} \in E((G_0,\phi_0)^\ell)$.

Finally, we assume that there exist $e \in B \cap U_E'$ and distinct $Y,Y' \in \bigcup_{i \in {\mathbb N}}\P_{e,i}$ such that $v_Yv_{Y'} \in E(M')$.
So $v_Yv_{Y'} \in E_3(M') \cup E_4(M')$, and there exist distinct $a,b \in V(M) \cap Z_1 \cap X_{T_e}-X_e$ such that $Y=Y_a$, $Y'=Y_b$, and there exists a path $Q$ in $M$ between $a$ and $b$ whose all internal vertices are contained in $(Z_2-Z_1) \cap X_{T_e}-X_e$.
Note that $Q$ can consist of an edge and have no internal vertex.
Since $M \subseteq (G,\phi)^\ell[V(G)]$, for each edge $z$ of $Q$, there exists a path $Q_z$ in $G$ of length in $(G,\phi)$ at most $\ell$ connecting the ends of $z$.

We first assume that $Q$ has exactly one edge $z_Q$ and $Q_{z_Q} \not \subseteq G[X_{T_e}]$.
Let $a'$ and $b'$ be the vertices in $X_e \cap V(Q_{z_Q})$ such that the subpaths of $Q_{z_Q}$ between $a$ and $a'$ and between $b$ and $b'$ are contained in $G[X_{T_e}]$ and internally disjoint from $X_e$, respectively.
Then $v_Y\overline{P_{aa'}}a'\overline{P_{a'b'}}b'\overline{P_{b'b}}v_{Y'}$ is a walk in $G_0$ of length in $(G_0,\phi_0)$ at most $\phi(z_Q) \leq \ell$, where $P_{aa'}$, $P_{a'b'}$ and $P_{b'b}$ are the subpaths of $Q_{z_Q}$ between $a$ and $a'$, between $a'$ and $b'$, and between $b'$ and $b$, respectively.
So $v_Yv_{Y'} \in E((G_0,\phi_0)^\ell)$.

Hence we may assume that either $Q$ has at least two edges, or $Q$ has only one edge $z_Q$ and $Q_{z_Q} \subseteq G[X_{T_e}]$.
For the former, for each edge $z$ of $Q$, since all internal vertices of $Q$ are contained in $(Z_2-Z_1) \cap X_{T_e}-X_e$, we know that $z$ is incident with a vertex in $Z_2-Z_1$, so $V(Q_z) \subseteq N_{(G,\phi)[X_{T_e}]}^{\leq 3\ell}[X_e]$.
For the latter, since $a,b \in Z_1$, we know $V(Q_z) \subseteq N_{(G,\phi)[X_{T_e}]}^{\leq 3\ell}[X_e]$.
So $V(Q_z) \subseteq N_{(G,\phi)[X_{T_e}]}^{\leq 3\ell}[X_e]$ in either case.
In addition, for every $z \in E(Q)$, the length in $(G,\phi)$ of any edge of $Q_z$ is at most $\ell \leq \lceil \frac{3\ell}{\epsilon} \rceil \epsilon$.
Hence there exists a path $P_{ab}$ in $\bigcup_{z \in E(Q)}Q_z \subseteq G$ between $a$ and $b$ such that $V(P_{ab}) \subseteq N_{(G,\phi)[X_{T_e}]}^{\leq 3\ell}[X_e]$, and the length in $(G,\phi)$ of any edge of $P_{ab}$ is at most $\lceil \frac{3\ell}{\epsilon} \rceil \epsilon$.
Combining this with Claim 2, the definition of $\P_e$ implies that $Y_a \cup Y_b$ is contained in some member of $\P_{e,\lceil \frac{3\ell}{\epsilon} \rceil}$.
So $v_Y$ and $v_{Y'}$ are contained in the same component $C$ of $H_e-X_e$.
Let $Y^*$ be the unique member of $\P_{e,\max I_e}$ such that $v_{Y^*} \in V(C)$.
Then the distance between $v_{Y^*}$ and any vertex in $V(C)$ in $(G_0,\phi_0)$ is at most $$(\max I_e) \cdot \frac{\epsilon}{8(\theta+\frac{\lambda}{\ell})} \leq \lceil \frac{3\ell+\lambda}{\epsilon} \rceil \cdot \frac{\epsilon\ell}{8(\theta\ell+\lambda)} \leq \frac{\ell}{2}.$$
So $\dist_{(G_0,\phi_0)}(v_Y,v_{Y'}) \leq \ell$.
Therefore, $v_Yv_{Y'} \in E((G_0,\phi_0)^{\ell})$.

This proves $E(M') \subseteq E((G_0,\phi_0)^\ell)$, and the claim follows.
$\Box$

\medskip

\noindent{\bf Claim 9:} For any $e \in B$ and $x,y \in X_{T_e}-X_e$, if there exists a path $P$ in $G$ between $x,y$ such that $V(P) \subseteq N_{(G,\phi)[X_{T_e}]}^{\leq k}[X_e]$ for some real number $k$, then there exists a path in $G$ between $x$ and $y$ of length in $(G,\phi)$ at most $(2k+2\mu+\ell')\theta$.

\noindent{\bf Proof of Claim 9:}
Let $e \in B$.
Let $x,y \in X_{T_e}-X_e$ such that there exists a path $P$ in $G$ between $x,y$ such that $V(P) \subseteq N_{(G,\phi)[X_{T_e}]}^{\leq k}[X_e]$ for some real number $k$.
If $\lvert X_e \rvert \leq \theta$, then let $A_e=X_e$.
So by the assumption of this lemma, there exists $A_e \subseteq X_e$ with $|A_e| \leq \theta$ such that $X_e \subseteq N_{(G,\phi)}^{\leq \mu}[A_e]$.
Hence $V(P) \subseteq N_{(G,\phi)[X_{T_e}]}^{\leq k}[X_e] \subseteq N_{(G,\phi)}^{\leq k+\mu}[A_e]$.
So for every $z \in V(P)$, there exists $v_z \in A_e$ such that $N_{(G,\phi)}^{\leq k+\mu}[\{v_z\}]$, and hence there exists $P_z$ in $G$ between $z$ and $v_z$ of length in $(G,\phi)$ at most $k+\mu$.

So for any two distinct vertices $z_1$ and $z_2$ in $P$ with $v_{z_1}=v_{z_2}$, the subpath of $P$ between $z_1$ and $z_2$ can be replaced by $P_{z_1} \cup P_{z_2}$ which has length in $(G,\phi)$ at most $2k+2\mu$.
Since $\lvert A_e \rvert \leq \theta$ and $\phi(z) \leq \ell'$ for every $z \in E(P) \subseteq E(G)$, there exists a path in $G$ between $x$ and $y$ of length in $(G,\phi)$ at most $\theta \cdot (2k+2\mu) + \ell'\theta = (2k+2\mu+\ell')\theta$.
$\Box$

\medskip

\noindent{\bf Claim 10:} For any $e \in B$ and $x \in V(M) \cap (Z_2-Z_1) \cap X_{T_e}$, there exists $y \in V(M) \cap Z_1 \cap X_{T_e}-X_e$ such that there exists a path in $G$ between $x$ and $y$ of length in $(G,\phi)$ at most $(6\ell+2\mu+\ell')\theta$. 

\noindent{\bf Proof of Claim 10:}
Let $e \in B$ and $x \in V(M) \cap (Z_2-Z_1) \cap X_{T_e}$.
By Claim 5, $V(M) \cap Z_1 \cap X_{T_e}-X_e \neq \emptyset$.
So there exists a path $P$ in $M$ between $x$ and $V(M) \cap Z_1 \cap X_{T_e}-X_e$ internally disjoint from $V(M) \cap Z_1 \cap X_{T_e}-X_e$.
Since $M \subseteq (G,\phi)^\ell[V(G)]$, by Claim 6, there exists a path $P'$ in $G[X_{T_e} \cap Z_3]$ between $x$ and $V(M) \cap Z_1 \cap X_{T_e}-X_e$.
Let $y$ be the end of $P'$ in $V(M) \cap Z_1 \cap X_{T_e}-X_e$.
By Claim 9, there exists a path in $G$ between $x$ and $y$ of length in $(G,\phi)$ at most $(6\ell+2\mu+\ell')\theta$. 
$\Box$

\medskip

By Claim 8, $M'$ is contained in a $c_0$-monochromatic component in $(G_0,\phi_0)^\ell[V(G_0)]-R$.
So the weak diameter of $M'$ in $(G_0,\phi_0)^\ell$ is at most $\nu$.

For every vertex $x \in V(M) \cap ((V(G_0) \cap V(G)) \cup \bigcup_{e \in B}(Z_1 \cap X_{T_e}))$, if $x \in V(G_0) \cap V(G)$, then define $h(x)=x$; if $x \in Z_1-(V(G_0) \cap V(G))$, define $h(x)=v_{Y_x}$.
Note that $h(x) \in V(M')$ for every $x \in V(M) \cap ((V(G_0) \cap V(G)) \cup \bigcup_{e \in B}(Z_1 \cap X_{T_e}))$.

\medskip

\noindent{\bf Claim 11:} For any vertices $x,y \in V(M) \cap ((V(G_0) \cap V(G)) \cup \bigcup_{e \in B}(Z_1 \cap X_{T_e}))$, there exists a path in $G$ between $x$ and $y$ with length in $(G,\phi)$ at most $(4(\nu+3)(3\theta+1)(\theta+\frac{\lambda}{\ell})+2)\ell$. 

\noindent{\bf Proof of Claim 11:}
Let $x,y \in V(M) \cap ((V(G_0) \cap V(G)) \cup \bigcup_{e \in B}(Z_1 \cap X_{T_e}))$.
Note that $h(x)$ and $h(y)$ are defined and are contained in $V(M')$.
Since the weak diameter of $M'$ in $(G_0,\phi_0)^\ell$ is at most $\nu$, there exists a path $P_0$ in $(G_0,\phi_0)^\ell$ between $h(x)$ and $h(y)$ with $\leng_{(G_0,\phi_0)^\ell}(P_0) \leq \nu$.
So there exists a path $P_0'$ in $G_0$ between $h(x)$ and $h(y)$ with $\leng_{(G_0,\phi_0)}(P_0') \leq \ell \nu$.

If $\{h(x),h(y)\} \subseteq V(G_0) \cap V(G)$, then by Statement 2 of Lemma \ref{con_qu_iso}, there exists a path $\widehat{P_0'}$ in $G$ between $h(x)=x$ and $h(y)=y$ of length in $(G,\phi)$ at most $4(3\theta+1)(\theta+\frac{\lambda}{\ell})\ell \nu$, so we are done.
Hence we may assume that $\{h(x),h(y)\} \not \subseteq V(G_0) \cap V(G)$.

For each $u \in \{x,y\}$, define $h'(u)$ and $h''(u)$ as follows.
	\begin{itemize}
		\item If $h(u) \in V(G_0) \cap V(G)$, then $h(u)=u$, and we let $h'(u)=u$ and $h''(u)=u$.
		\item Otherwise, there exists $e_u \in U_E'$ such that $u \in (X_{T_{e_u}}-X_{e_u}) \cap N_{(G,\phi)[X_{T_{e_u}}]}^{\leq \ell}[Y_u]$, so there exists a path $P_u$ in $G[X_{T_e}]$ from $u$ to $Y_u$ with length in $(G,\phi)$ at most $\ell$, and  
			\begin{itemize}
				\item we let $h'(u)$ be the end of $P_u$ in $Y_u$, and
				\item if $X_{e_u} \cap V(P_0')=\emptyset$, then let $h''(u)=h'(u)$,  
				\item otherwise, let $h''(u)$ be the vertex in $X_{e_u} \cap V(P_0')$ such that the subpath of $P_0'$ between $h''(u)$ and $h(u)$ is contained in $G_0[V(H_{e_u})]$. 
			\end{itemize}
	\end{itemize}
Note that if $h(u) \not \in V(G_0) \cap V(G)$, then there exists $Y^*_u \in \P_{e_u,\lceil \frac{3\ell}{\epsilon} \rceil}$ such that $Y_u \cup \{h''(u)\} \subseteq Y^*_u$, so there exists a path $Q_u$ in $G_0[V(H_{e_u})]$ between $h'(u)$ and $h''(u)$ internally disjoint from $X_{e_u}$ of length in $(G_0,\phi_0)$ at most $\ell$ by Lemma \ref{hierachy_distance}. 

Hence for each $u \in \{x,y\}$, 
	\begin{itemize}
		\item there exists a path in $G$ from $u$ to $h'(u)$ of length in $(G,\phi)$ at most $\ell$, and
		\item by Statement 2 in Lemma \ref{con_qu_iso}, there exists a path in $G$ from $h'(u)$ to $h''(u)$ of length in $(G,\phi)$ at most $4(3\theta+1)(\theta+\frac{\lambda}{\ell})\ell$, 
	\end{itemize}
so there exists a path $P^*_u$ in $G$ from $u$ to $h''(u)$ with $\leng_{(G,\phi)}(P^*_u) \leq (4(3\theta+1)(\theta+\frac{\lambda}{\ell})+1)\ell$.

We first assume that either $h(x) \in V(G_0) \cap V(G)$ or $V(P_0') \cap X_{e_x} \neq \emptyset$.
So $h''(x) \in V(P_0') \cap V(G_0) \cap V(G)$.
Since $\{h(x),h(y)\} \not \subseteq V(G_0) \cap V(G)$, we know $V(P_0') \cap X_{e_y} \neq \emptyset$.
So $h''(y) \in V(P_0') \cap V(G_0) \cap V(G)$.
Note that the subpath $P'$ of $P_0'$ between $h''(x)$ and $h''(y)$ is a path in $G_0$ of length in $(G_0,\phi_0)$ at most $\ell \nu$.
Since $h''(x), h''(y) \in V(G_0) \cap V(G)$, by Lemma \ref{con_qu_iso}, there exists a path $\widehat{P'}$ in $G$ between $h''(x)$ and $h''(y)$ of length in $(G,\phi)$ at most $4(3\theta+1)(\theta+\frac{\lambda}{\ell})\ell \nu$.
Therefore, $P^*_x \cup \widehat{P'} \cup P^*_y$ is a walk in $G$ from $x$ to $y$ of length in $(G,\phi)$ at most $2(4(3\theta+1)(\theta+\frac{\lambda}{\ell})+1)\ell + 4(3\theta+1)(\theta+\frac{\lambda}{\ell})\ell \nu$.

So we may assume that $h(x) \not\in V(G_0) \cap V(G)$ and $V(P_0') \cap X_{e_x}=\emptyset$.
Hence $h(y) \not \in V(G_0) \cap V(G)$, $e_x=e_y$ and $V(P_0') \cap V(G)=\emptyset$. 
So $Y_x \cup Y_y$ is contained in a member of $\P_{i,\lceil \frac{3\ell}{\epsilon} \rvert}$.
Hence there exists a path $W_0$ in $G_0[V(H_{e_x})]$ between $h'(x)$ and $h'(y)$ of length $(G_0,\phi_0)$ at most $\ell$.
By Lemma \ref{con_qu_iso}, there exists a path $\widehat{W_0}$ in $G$ between $h'(x)$ and $h'(y)$ of length in $(G,\phi)$ at most $4(3\theta+1)(\theta+\frac{\lambda}{\ell})\ell$.
Therefore, $P_x \cup \widehat{W_0} \cup P_x$ contains a path in $G$ between $x$ and $y$ of length in $(G,\phi)$ at most $2\ell+4(3\theta+1)(\theta+\frac{\lambda}{\ell})\ell = (4(3\theta+1)(\theta+\frac{\lambda}{\ell})+2)\ell$.
$\Box$

\medskip

\noindent{\bf Claim 12:} The weak diameter of $M$ in $(G,\phi)$ is at most $\frac{\ell \nu^*}{2}$. 

\noindent{\bf Proof of Claim 12:}
Suppose to the contrary that there exist vertices $x,y \in V(M)$ such that $\dist_{(G,\phi)}(x,y)>\ell \nu^*/2$.
For each $u \in \{x,y\}$, if $u \in (V(G_0) \cap V(G)) \cup \bigcup_{e \in B}(Z_1 \cap X_{T_e})$, then let $q_u=u$; otherwise, $u \in (Z_2-Z_1) \cap X_{T_{e_u}}$ by Claim 6, so by Claim 10, there exists $q_u \in V(M) \cap Z_1 \cap X_{T_{e_u}}-X_{e_u}$ such that $\dist_{(G,\phi)}(u,q_u) \leq (6\ell+2\mu+\ell')\theta$. 
By Claim 11, $\dist_{(G,\phi)}(q_x,q_y) \leq (4(\nu+3)(3\theta+1)(\theta+\frac{\lambda}{\ell})+2)\ell$. 
Therefore, $$\dist_{(G,\phi)}(x,y) \leq (12\ell+4\mu+2\ell')\theta + (4(\nu+3)(3\theta+1)(\theta+\frac{\lambda}{\ell})+2)\ell \leq \frac{\ell \nu^*}{2},$$ a contradiction. 
$\Box$

\medskip

Hence the weak diameter of $M$ in $(G,\phi)$ is at most $\frac{\ell \nu^*}{2}$ by Claim 12. 
Then by Lemma \ref{dist_original_power}, the weak diameter of $M$ in $(G,\phi)^\ell$ is at most $\lceil \frac{2}{\ell} \cdot \frac{\ell \nu^*}{2} \rceil = \lceil \nu^* \rceil = \nu^*$, a contradiction.
This proves the lemma.
\end{pf}

\subsection{Construction with bounded adhesion} \label{sec:bdd_adhesion}

The goal of this subsection is to prove Lemma \ref{weighted_tree_extension_clean}, which roughly states that if a weighted graph $(G,\phi)$ has a tree-decomposition with bounded adhesion such that every weighted graph that can be obtained from the subgraph induced by a bag by adding ``gadgets'' has a nice coloring, then $(G,\phi)$ has a nice coloring.
It would allow us to reduce the coloring problem on $(G,\phi)$ essentially to the coloring problem on its bags.
In order to prove Lemma \ref{weighted_tree_extension_clean}, we prove a more technical version (Lemma \ref{weighted_tree_extension}) for the sake of having a stronger inductive hypothesis.

Let $S$ be a subset of ${\mathbb R}^+$.
We say that a class of weighted graphs is \defn{$S$-bounded} if every weighted graph in this class is $S$-bounded.
Let $\F$ and $\F'$ be $S$-bounded classes of weighted graphs.
Let $\eta,\theta,\lambda$ be nonnegative integers with $\eta \leq \theta$.
We say that a weighted graph $(G,\phi)$ is \defn{$(\eta,\theta,\F,\F',\lambda)$-constructible} if there exists a rooted tree-decomposition $(T,\X)$ of $G$ of adhesion at most $\theta$ such that the following hold.
	\begin{itemize}
		\item[(EC1)] For every edge $tt' \in E(T)$, if $\lvert X_t \cap X_{t'} \rvert > \eta$, then one end of $tt'$ has no child, say $t'$, and $\lvert X_{t'}-X_t \rvert \leq \lambda$,
		\item[(EC2-EC6)] For every $t \in V(T)$,
			\begin{itemize}
				\item[(EC2)] if $t$ is the root of $T$, then $\lvert X_t \rvert \leq \theta$, 
				\item[(EC3)] if $\eta>0$ and $t$ is the root of $T$, then $X_t \neq \emptyset$,
				\item[(EC4)] if $t$ has a child in $T$, then $(G,\phi)[X_t] \in \F$,
				\item[(EC5)] if $t$ has no child in $T$, then $(G,\phi)[X_t] \in \F^{+\lambda,S}$, and
				\item[(EC6)] $\F'$ contains every weighted graph that is obtained from $(G,\phi)[X_t]$ by for each child $t'$ of $t$, adding at most $\lambda$ new vertices and some new edges such that each new edge is between two new vertices or between a new vertex and $X_t \cap X_{t'}$, and the weight of each new edge is in $S$.
			\end{itemize}
	\end{itemize}
In this case, we call $(T,\X)$ an \defn{$(\eta,\theta,\F,\F',\lambda)$-construction} of $(G,\phi)$.

A class $\F$ of weighted graphs is \defn{hereditary} if for every $(G,\phi) \in \F$ and every $X \subseteq V(G)$, $(G,\phi)[X] \in \F$.
Note that if $\F$ is hereditary and $S$-bounded for some $S \subseteq {\mathbb R}^+$, then so is $\F^{+\lambda,S}$ for any nonnegative integer $\lambda$.

Recall that if $\ell \in {\mathbb R}^+$ and $(G,\phi)$ is a weighted graph such that the image of $\phi$ is contained in $(0,\ell]$, then $(G,\phi,\ell)$ is obtained from $(G,\phi)$ by duplicating edges and defining weights of the copied edges as their original, and $G$ is a spanning subgraph of $(G,\phi)^\ell$.

\begin{lemma} \label{weighted_tree_extension}
For any positive real numbers $\ell,\nu$, any positive integer $m$ with $m \geq 2$, any nonnegative integer $\theta$, and any $(m,\ell,\nu)$-nice hereditary $(0,\ell]$-bounded classes $\F,\F'$ of weighted graphs, there exists a function $f^*: {\mathbb N} \cup \{0\} \rightarrow {\mathbb R}^+$ such that the following holds. 
Let $\eta$ be a nonnegative integer with $\eta \leq \theta$.
Let $(G,\phi)$ be an $(\eta,\theta,\F,\F',\theta^2)$-constructible $(0,\ell]$-bounded weighted graph with an $(\eta,\theta,\F,\F',\theta^2)$-construction $(T,\X)$.
Denote $\X$ by $(X_t: t \in V(T))$.
Let $t^*$ be the root of $T$.
Let $Z \subseteq N_{(G,\phi)}^{\leq 3\ell}[X_{t^*}]$.
If $c: Z \rightarrow [m]$ is a function, then $c$ can be extended to an $m$-coloring of $(G,\phi)^\ell$ with weak diameter in $(G,\phi)^\ell$ at most $f^*(\eta)$.
\end{lemma}

\begin{pf}
Let $\ell,\nu$ be positive real numbers, $m$ be a positive integer with $m \geq 2$, $\theta$ be a nonnegative integer, and $\F,\F'$ be $(m,\ell,\nu)$-nice hereditary $(0,\ell]$-bounded classes of weighted graphs.
By Lemma \ref{apex_extension}, there exists an integer $\nu_{\F^{+}}$ such that $\F^{+\theta^2,(0,\ell]}$ is $(m,\ell,\nu_{\F^+})$-nice.
Note that $\F \subseteq \F^{+\theta^2,(0,\ell]}$, so we may assume that $\nu_{\F^+} \geq \nu$.
We define the following.
	\begin{itemize}
		\item Let $f_1: {\mathbb R}^+ \rightarrow {\mathbb R}^+$ be the function such that for every $x \in {\mathbb R}^+$, $f_1(x)$ is the number $\nu^*$ given by Lemma \ref{deleting_centered_set} by taking $(k,r,\ell,\nu)=(\theta,3\ell,\ell,x)$. 
		\item Let $\nu'_\theta$ be the number $\nu$ given by Lemma \ref{all_centered} by taking $(k,r,\ell)=(\theta,3\ell,\ell)$. 
		\item Let $f_2: {\mathbb R}^+ \rightarrow {\mathbb R}^+$ be the function such that for every $x \in {\mathbb R}^+$, $f_2(x)$ is the number $\nu^*$ given by Lemma \ref{con_color} by taking $(\ell,\ell',\nu,\mu,\lambda,m,\theta)=(\ell,\ell,x,1,0,m,\max\{\theta,1\})$. 
        \item Let $\nu''_\theta$ be the number $\nu$ given by Lemma \ref{all_centered} by taking $(k,r,\ell)=(\theta+\theta^2,0,\ell)$. 
		\item Define $f^*: ({\mathbb N} \cup \{0\}) \rightarrow {\mathbb R}^+$ to be the function such that $f^*(0)=\nu_{\F^+} +\nu'_\theta+\nu''_\theta+f_1(\nu)$, and for every $x \in {\mathbb N}$, $f^*(x)= f_2(f_1(f^*(x-1)))$. 
	\end{itemize}
Note that $f^*$ is an increasing function by Lemmas \ref{deleting_centered_set} and \ref{con_color}.

Let $\eta,(G,\phi),(T,\X),t^*,Z,c$ be the ones as defined in the lemma.
Suppose to the contrary that $c$ cannot be extended to an $m$-coloring of $(G,\phi)^\ell$ with weak diameter in $(G,\phi)^\ell$ at most $f^*(\eta)$, and subject to this, the tuple $(\eta,\lvert V(T) \rvert+\lvert V(G)-Z \rvert+\lvert V(G) \rvert)$ is lexicographically minimal.

Since $(T,\X)$ is an $(\eta,\theta,\F,\F',\theta^2)$-construction, $\lvert X_{t^*} \rvert \leq \theta$ by (EC2).
So $Z$ is $(\theta,3\ell)$-centered.
If $Z=V(G)$, then $c$ can be extended to an $m$-coloring of $(G,\phi)^\ell$ with weak diameter in $(G,\phi)^\ell$ at most $\nu'_\theta \leq f^*(\eta)$ by Lemma \ref{all_centered} (with $R=\emptyset$), a contradiction. 

So $Z \neq V(G)$.

\medskip

\noindent{\bf Claim 1:} $\eta \geq 1$.

\noindent{\bf Proof of Claim 1:}
Suppose to the contrary that $\eta=0$.
Let $W=\{tt' \in E(T): X_t \cap X_{t'}=\emptyset\}$.
Note that for every component $C$ of $G$, we have $V(C) \subseteq X_{C'}$ for some component $C'$ of $T-W$.
Since every component of $(G,\phi)^\ell$ is $(C,\phi|_{E(C)})^\ell$ for some component $C$ of $G$, there exists a component $C$ of $G$ such that $c|_{Z \cap V(C)}$ cannot be extended to an $m$-coloring of $(C,\phi|_{E(C)})^\ell$ with weak diameter in $(C,\phi|_{E(C)})^\ell$ at most $f^*(\eta)$.

Let $T_C$ be the component of $T-W$ with $V(C) \subseteq X_{T_C}$.
Let $t_C$ be the root of $T_C$.
Since $(T,\X)$ is a $(0,\theta,\F,\F',\theta^2)$-construction, (EC1) implies that $T_C$ is a star, and $(G,\phi)[X_{T_C}]$ is a $(0,\ell]$-bounded weighted graph obtained from $(G,\phi)[X_{t_C}]$ by for each child $t'$ of $t_C$ in $T_C$, adding at most $\theta^2$ new vertices and some new edges such that each new edge is between two new vertices or between a new vertex and $X_{t_C} \cap X_{t'}$.
So $(G,\phi)[X_{T_C}] \in \F'$ by (EC6).
Let $(G_C,\phi_C) = (G,\phi)[X_{T_C}]-(Z \cap X_{T_C})$.
Since $\F'$ is hereditary, $(G_C,\phi_C) \in \F'$.
Since $\F'$ is $(m,\ell,\nu)$-nice, $(G_C,\phi_C)^\ell$ is $m$-colorable with weak diameter in $(G_C,\phi_C)^\ell$ at most $\nu \leq \nu_{\F^+} \leq f^*(\eta)$.

Since $\eta=0$ and $Z \subseteq N_{(G,\phi)}^{\leq 3\ell}[X_{t^*}]$, if $t_C \neq t^*$, then $Z \cap V(C) \subseteq Z \cap X_{T_C}=\emptyset$, so $c|_{Z \cap V(C)}$ can be extended to an $m$-coloring of $(C,\phi|_{E(C)})^\ell$ with weak diameter in $(C,\phi|_{E(C)})^\ell$ at most $f^*(\eta)$, a contradiction.
Hence $t_C=t^*$.
So $Z \cap X_{T_C} = Z$ is $(\theta,3\ell)$-centered in $(G,\phi)[X_{T_C}]$.
By Lemma \ref{deleting_centered_set}, $c|_{Z \cap V(C)}$ can be extended to an $m$-coloring of $((G,\phi)[X_{T_C}])^\ell$ with weak diameter in $((G,\phi)[X_{T_C}])^\ell$ at most $f_1(\nu) \leq f^*(\eta)$.
Since $(C,\phi|_{E(C)})^\ell$ is a component of $((G,\phi)[X_{T_C}])^\ell$, $(C,\phi|_{E(C)})^\ell$ is $m$-colorable with weak diameter in $(C,\phi|_{E(C)})^\ell$ at most $f^*(\eta)$, a contradiction.
$\Box$

\medskip

\noindent{\bf Claim 2:} $G$ is connected.

\noindent{\bf Proof of Claim 2:}
Suppose to the contrary that $G$ is not connected.
Since every component of $(G,\phi)^\ell$ is $(C,\phi|_{E(C)})^\ell$ for some component $C$ of $G$, there exists a component $C$ of $G$ such that $c|_{Z \cap V(C)}$ cannot be extended to an $m$-coloring of $(C,\phi|_{E(C)})^\ell$ with weak diameter in $(C,\phi|_{E(C)})^\ell$ at most $f^*(\eta)$.
Let $T_C$ be the subtree of $T$ induced by $\{t \in V(T): X_t \cap V(C) \neq \emptyset\}$.
Let $\X_C = (X_t \cap V(C): t \in V(T_C))$.

If $t^* \in V(T_C)$, then $t^*$ is the root of $T_C$ and $X_{t^*} \cap V(C) \neq \emptyset$, so $(T_C,\X_C)$ is an $(\eta,\theta,\F,\F',\theta^2)$-construction of $(C,\phi|_{E(C)})$ with $(\eta,\lvert V(T_C) \rvert+\lvert V(C)-(Z \cap V(C)) \rvert+\lvert V(C) \rvert)$ lexicographically smaller than $(\eta,\lvert V(T) \rvert+\lvert V(G)-Z \rvert+\lvert V(G) \rvert)$, so $c|_{Z \cap V(C)}$ can be extended to an $m$-coloring of $(C,\phi|_{E(C)})^\ell$ with weak diameter in $(C,\phi|_{E(C)})^\ell$ at most $f^*(\eta)$, a contradiction.
So $t^* \not \in V(T_C)$.
Since $Z \subseteq N_{(G,\phi)}^{\leq 3\ell}[X_{t^*}]$, $Z \cap V(C)=\emptyset$.
Let $T_C'$ be the rooted tree obtained from $T_C$ by adding a node $t^*_C$ adjacent to the root of $T_C$, where $t^*_C$ is the root of $T_C'$.
Let the bag at $t^*_C$ be the set consisting of a single vertex in the intersection of $V(C)$ and the bag of the root of $T_C$.
Then since $\eta \geq 1$ by Claim 1, we obtain an $(\eta,\theta,\F,\F',\theta^2)$-construction of $(C,\phi|_{E(C)})$ with underlying tree $T_C'$.
Since $t^* \not \in V(T_C)$, $(\eta,\lvert V(T_C') \rvert+\lvert V(C)-(Z \cap V(C)) \rvert+\lvert V(C) \rvert)$ is lexicographically smaller than $(\eta,\lvert V(T) \rvert+\lvert V(G)-Z \rvert+\lvert V(G) \rvert)$.
Since $Z \cap V(C)=\emptyset$, the minimality implies that $c|_{Z \cap V(C)}$ can be extended to an $m$-coloring of $(C,\phi|_{E(C)})^\ell$ with weak diameter in $(C,\phi|_{E(C)})^\ell$ at most $f^*(\eta)$, a contradiction.
$\Box$

\medskip

\noindent{\bf Claim 3:} $Z=N_{(G,\phi)}^{\leq 3\ell}[X_{t^*}]$ and $Z-X_{t^*} \neq \emptyset$.

\noindent{\bf Proof of Claim 3:}
Suppose to the contrary that there exists $v \in N_{(G,\phi)}^{\leq 3\ell}[X_{t^*}]-Z$.
Let $Z' = Z \cup \{v\}$.
Let $c': Z' \rightarrow [m]$ be the function obtained from $c$ by further defining $c'(v)=m$.
Then the minimality of $(\eta,\lvert V(T) \rvert+\lvert V(G)-Z \rvert+\lvert V(G) \rvert)$ implies that $c'$ (and hence $c$) can be extended to an $m$-coloring of $(G,\phi)^\ell$ with weak diameter in $(G,\phi)^\ell$ at most $f^*(\eta)$, a contradiction.

So $Z=N_{(G,\phi)}^{\leq 3\ell}[X_{t^*}]$.
Suppose to the contrary that $Z \subseteq X_{t^*}$.
Then $N_{(G,\phi)}^{\leq 3\ell}[X_{t^*}]= Z \subseteq X_{t^*}$.
By Claim 1 and (EC3) in the definition of $(\eta,\theta,\F,\F',\theta^2)$-constructions, $X_{t^*} \neq \emptyset$.
Since $G$ is connected by Claim 2, $V(G)=X_{t^*}= Z$, a contradiction.
$\Box$

\medskip

Note that Claim 3 implies that $Z \neq \emptyset$ and $X_{t^*} \neq \emptyset$, and Claim 1 implies that $\theta \geq \eta \geq 1$.

For each $v \in X_{t^*}$, let $T_v$ be the subgraph of $T$ induced by $\{t \in V(T): N_{(G,\phi)}^{\leq 3\ell}[\{v\}] \cap X_t \neq \emptyset\}$.
Since $(T,\X)$ is a tree-decomposition, $T_v$ is a subtree of $T$ containing $t^*$.
So $\bigcup_{v \in X_{t^*}}T_v$ is a subtree of $T$ containing $t^*$.

Let $U_E = \{e \in E(T):$ exactly one end of $e$ is in $\bigcup_{v \in X_{t^*}}V(T_v)\}$.
Note that edges in $U_E$ are pairwise incomparable.
For each $e \in E(T)$, define $X_e$ to be the intersection of the bags of the ends of $e$; define $T_e$ to be the component of $T-e$ disjoint from $t^*$.
Since the adhesion of $(T,\X)$ is at most $\theta$, $\lvert X_e \rvert \leq \theta$ for every $e \in E(T)$.
Let $\epsilon = \min_{e \in E(G)}\phi(e)$.

Define $(G_0,\phi_0)$ to be the $(U_E,U_E,\ell,\theta,0)$-condensation of $(G,\phi,T,\X)$.
For each $e \in U_E$, let $(H_e,\phi_{H_e},I_e,B_e)$ be the $(\ell,\epsilon,\theta,0)$-hierarchy mentioned in the definition of the condensation.

Let $T_0 = \bigcup_{v \in X_{t^*}}T_v$.
Let $Z_0 = N_{(G,\phi)}^{\leq 3\ell}[X_{t^*}]$.
So $Z_0=Z$ by Claim 3; for every $t \in V(T_0)$, $X_t \cap Z_0 \neq \emptyset$; for every $t \in V(T)-V(T_0)$, $X_t \cap Z_0=\emptyset$.
Note that $Z =Z_0 \subseteq V(G[X_{T_0}]) \subseteq V(G_0)$.

\medskip

\noindent{\bf Claim 4:} There exists an $m$-coloring $c_0$ of $(G_0,\phi_0)^\ell$ with weak diameter in $(G_0,\phi_0)^\ell$ at most $f_1(f^*(\eta-1))$ such that $c_0(v)=c(v)$ for every $v \in Z$. 

\noindent{\bf Proof of Claim 4:}
Define $T'$ to be the rooted tree obtained from $T_0$ by the following operation: for each $e \in U_E$, adding a node $t_e$ adjacent to the end of $e$ in $V(T_0)$.
For each $t' \in V(T_0)$, define $X'_{t'}=X_t$; for each $t' \in V(T')-V(T_0)$, we know $t'=t_e$ for some $e \in U_E$, and we define $X'_{t'}=V(H_e)$. 
Let $\X'=(X'_t: t \in V(T'))$.

Clearly, $(T',\X')$ is a tree-decomposition of $G_0$ of adhesion at most $\max_{e \in U_E}\{\theta,\lvert X_e \rvert\} = \theta$.
For each $tt' \in E(T')$, say $t'$ is a child of $t$, if $tt' \in E(T_0)$, then $X'_t=X_t$, $X'_{t'}=X_{t'}$, $t$ has a child in both $T$ and $T'$, and $t'$ has a child in $T'$ if and only if $t'$ has a child in $T$; if $tt' \not \in E(T_0)$, then $t \in V(T_0)$ and $t' \not \in V(T_0)$, and $\lvert X'_{t'}-X'_t \rvert = \lvert V(H_e) \rvert - \lvert B_e \rvert \leq \theta^2$ by Lemma \ref{hierachy_distance}.
Hence for every $tt' \in E(T')$, if $\lvert X'_t \cap X'_{t'} \rvert > \eta$, then one end of $tt'$, say $t'$, has no child, and $\lvert X'_{t'}-X'_t \rvert \leq \theta^2$.
So $(T',\X')$ satisfies (EC1) for being a $(\eta,\theta,\F,\F',\theta^2)$-construction of $(G_0,\phi_0)$.

Furthermore, $t^* \in V(T_0) \subseteq V(T')$ and $X'_{t^*}=X_{t^*}$, so $\lvert X'_{t^*} \rvert = \lvert X_{t^*} \rvert \leq \theta$ by (EC2).
Since $\eta \geq 1$ by Claim 1, $X'_{t^*} = X_{t^*} \neq \emptyset$ by (EC3).
So $(T',\X')$ satisfies (EC2) and (EC3) for being a $(\eta,\theta,\F,\F',\theta^2)$-construction of $(G_0,\phi_0)$.

In addition, for every $t \in V(T')$, if $t$ has a child in $T'$, then $t \in V(T_0) \subseteq V(T)$ has a child in $T$, so $(G_0,\phi_0)[X'_t]=(G,\phi)[X_t] \in \F$; if $t$ has no child in $T'$, then either $t \in V(T)$ has no child in $T$, or $t \in V(T')-V(T)$, so $(G_0,\phi_0)[X'_t] \in \F^{+\theta^2,(0,\ell]}$ by (EC5).
Hence $(T',\X')$ satisfies (EC4) and (EC5) for being a $(\eta,\theta,\F,\F',\theta^2)$-construction of $(G_0,\phi_0)$.

Now we show that $(T',\X')$ satisfies (EC6) for being a $(\eta,\theta,\F,\F',\theta^2)$-construction of $(G_0,\phi_0)$.
If $t \in V(T_0) \subseteq V(T)$, then $(T',\X')$ satisfies (EC6) since $(G_0,\phi_0)[X'_t]=(G,\phi)[X_t]$.
If $t \in V(T')-V(T_0)$, then $t$ has no child, so it suffices to show $(G_0,\phi_0)[X'_t] \in \F'$; since $t \in V(T')-V(T_0)$, we know that $(G_0,\phi_0)[X'_t]$ can be obtained from an induced subgraph of $(G,\phi)[X_{p_t}]$, where $p_t$ is the parent of $t$, by adding at most $\theta^2$ new vertices and some new edges such that each new edge is between two new vertices or between a new vertex and $X_{p_t} \cap X_{t}$ and hence belongs to $\F'$.
Hence $(T',\X')$ satisfies (EC6) for being a $(\eta,\theta,\F,\F',\theta^2)$-construction of $(G_0,\phi_0)$.

In summary, $(T',\X')$ is an $(\eta,\theta,\F,\F',\theta^2)$-construction of $(G_0,\phi_0)$.

For every $t \in V(T')$, let $X''_t = X'_t-Z_0$.
Let $\X''=(X''_t: t\in V(T'))$.
So $(T',\X'')$ is a tree-decomposition of $G_0-Z_0$ of adhesion at most $\theta$.
Since $Z_0=N_{(G,\phi)}^{\leq 3\ell}[X_{t^*}]$ and $X_t \cap Z_0 \neq \emptyset$ for every $t \in V(T_0)$, we know $X_e \cap Z_0 \neq \emptyset$ for every $e \in E(T_0)$.
Note that $X''_{t^*} = X_{t^*}-Z_0=\emptyset$.

Now we define $T'''$ and $\X'''$.
	\begin{itemize}
		\item If $\eta-1=0$, then let $T'''=T'$ and $\X'''=\X''$; 
		\item otherwise, $\eta-1 \geq 1$ by Claim 1, and we 
			\begin{itemize}
				\item let $t_0$ be a node of $T'$ with $X''_{t_0} \neq \emptyset$ closest to $t^*$, let $v_0$ be a vertex in $X''_{t_0}$, let $T'''$ be the rooted tree obtained from $T'$ by adding a new node $t_0^*$ adjacent to $t^*$, where $t_0^*$ is the root of $T'''$, and 
				\item let $\X'''=(X'''_t: t \in V(T'''))$, where 
					\begin{itemize}
						\item $X'''_{t_0^*}=\{v_0\}$, 
						\item $X'''_t=X''_t \cup \{v_0\}$ if $t \neq t_0^*$ and $t$ is in the path in $T'$ between $t^*$ and $t_0$, and 
						\item $X'''_{t}=X''_t$ otherwise.
					\end{itemize}
			\end{itemize}
	\end{itemize}

Then $(T''',\X''')$ is a tree-decomposition of $G_0-Z_0$.
Since $(T',\X')$ is an $(\eta,\theta,\F,\F',\theta^2)$-construction of $G_0$, and $G$ is connected, and $\F,\F'$ and $\F^{+\theta^2,(0,\ell]}$ are hereditary, $(T''',\X''')$ is an $(\eta-1,\theta,\F,\F',\theta^2)$-construction of $(G_0,\phi_0)-Z_0$.

By the minimality of $\eta$, there exists an $m$-coloring $c_0'$ of $((G_0,\phi_0)-Z_0)^\ell$ with weak diameter in $((G_0,\phi_0)-Z_0)^\ell$ at most $f^*(\eta-1)$.
Let $c_0=c \cup c_0'$.
By Lemma \ref{deleting_centered_set}, $c_0$ is an $m$-coloring of $(G_0,\phi_0)^\ell$ with weak diameter in $(G_0,\phi_0)^\ell$ at most $f_1(f^*(\eta-1))$ with $c_0(v)=c(v)$ for every $v \in Z$.
$\Box$

\medskip

Note that the coloring $c_0$ can be obtained by extending $c$.
By Claim 4 and Lemma \ref{con_color} (taking $R=\emptyset$), $c_0$ (and hence $c$) can be extended to an $m$-coloring $c_3$ of $(G,\phi)^\ell[X_{T_0} \cup \bigcup_{e \in U_E} N^{\leq 3\ell}_{(G,\phi)[X_{T_e}]}[X_e]]$ with weak diameter in $(G,\phi)^\ell$ at most $f_2(f_1(f^*(\eta-1)))$ such that 
	\begin{itemize}
		\item[(i)] for every $m$-coloring $c'$ of $(G,\phi)^\ell$ that can be obtained by extending $c_3$, and for every $c'$-monochromatic component $M$ in $(G,\phi)^\ell$ intersecting $X_{T_0} \cup \bigcup_{e \in U_E}N_{(G,\phi)[X_{T_e}]}^{\leq \ell}[X_e]$, $M$ has weak diameter in $(G,\phi)^\ell$ at most $f_2(f_1(f^*(\eta-1))) \leq f^*(\eta)$.
	\end{itemize}

Let $(G_1,\phi_1) = (G,\phi)[\bigcup_{e \in U_E}X_{T_e}]$.
By the definition of $T_0$ and $U_E$, $Z \cap V(G_1)=\emptyset$. 
Let $S_1 = \bigcup_{e \in U_E}X_e$.
For $i \in [3]$, let $Z_i = N_{(G_1,\phi_1)}^{\leq i\ell}[S_1]$. 

For every $e \in U_E$, let $(G_e,\phi_e) = (G,\phi)[X_{T_e}]$ and $Z_e = Z_3 \cap V(G_e)$, and let $c_e: Z_e \rightarrow [m]$ such that $c_e(v)=c_3(v)$ for every $v \in Z_e$.
Note that for every $e \in U_E$, $Z_e \subseteq N_{(G_e,\phi_e)}^{\leq 3\ell}[X_e]$.

\medskip

\noindent{\bf Claim 5:} For every $e \in U_E$, the coloring $c_e$ can be extended to an $m$-coloring $c_e'$ of $(G_e,\phi_e)^\ell$ with weak diameter in $(G_e,\phi_e)^\ell$ at most $f^*(\eta)$.

\noindent{\bf Proof of Claim 5:}
Let $e \in U_E$.
If $\lvert X_e \rvert > \eta$, then since $(T,\X)$ is an $(\eta, \theta,\F,\F',\theta^2)$-construction of $G$, (EC1) implies that $\lvert V(G_e) \rvert \leq \lvert X_e \rvert+\theta^2 \leq \theta+\theta^2$, so $V(G_e)$ is $(\theta+\theta^2,0)$-centered in $(G_e,\phi_e)$, and hence by Lemma \ref{all_centered}, $c_e$ can be extended to an $m$-coloring of $(G_e,\phi_e)^\ell$ with weak diameter in $(G_e,\phi_e)^\ell$ at most $\nu''_\theta \leq f^*(\eta)$. 

So we may assume $\lvert X_e \rvert \leq \eta$.
If $X_e = \emptyset$, then $V(G_e)=\emptyset$ since $G$ is connected and $X_{t^*} \neq \emptyset$ by Claims 1 and 2, so we are done.

So we may assume $X_e \neq \emptyset$.
Define $Q_e$ to be the tree obtained from $T_e$ by adding a node $r_e$ adjacent to the end of $e$ in $V(T_e)$.
We make $Q_e$ a rooted tree by assigning $r_e$ to be the root of $Q_e$.
Let $W_{r_e} = X_e$; for every $t \in V(T_e)$, let $W_t=X_t$.
Let $\W=(W_t: t \in V(Q_e))$.
Then $(Q_e,\W)$ is a rooted tree-decomposition of $G_e$ of adhesion at most $\theta$ such that $\lvert W_{r_e} \rvert = \lvert X_e \rvert \leq \eta$.
So if $tt' \in E(Q_e)$ with $\lvert W_t \cap W_{t'} \rvert > \eta$, then $tt' \in E(T_e)$, so $W_t=X_t$ and $W_{t'}=X_{t'}$.
Since $(G_e,\phi_e)[W_{r_e}]=(G,\phi)[X_e]$ and $\F$ and $\F'$ are hereditary, $(Q_e,\W)$ is an $(\eta,\theta,\F,\F',\theta^2)$-construction of $G_e$.

Note that $\lvert V(Q_e) \rvert \leq \lvert V(T) \rvert$, and equality holds only when $t^*$ is an end of $e$. 
If $t^*$ is an end of $e$, then since $X_e \neq \emptyset$, by Claim 3, $X_t \cap Z \supseteq X_t \cap X_{t^*} \neq \emptyset$, where $t$ is the end of $e$ other than $t^*$, so $t \in V(T_0)$, a contradiction.
Hence $\lvert V(Q_e) \rvert < \lvert V(T) \rvert$.

Recall that $Z_e \subseteq N_{(G_e,\phi_e)}^{\leq 3\ell}[X_e] = N_{(G_e,\phi_e)}^{\leq 3\ell}[W_{r_e}]$.
So by the minimality of $(\eta,\lvert V(T) \rvert+\lvert V(G)-Z \rvert+\lvert V(G) \rvert)$, $c_e$ can be extended to an $m$-coloring of $(G_e,\phi_e)^\ell$ with weak diameter in $(G_e,\phi_e)^\ell$ at most $f^*(\eta)$.
$\Box$

\medskip

For every $e \in U_E$, let $c_e'$ be the $m$-coloring mentioned in Claim 5.

Define $c^* = c_3|_{V(G) \cap V(G_0)} \cup \bigcup_{e \in U_E}c'_e$.
Note that $c^*$ is well-defined by the definition of $c_e$, and $c^*$ is an $m$-coloring of $(G,\phi)^\ell$ that can be obtained by extending $c_3$ and hence by extending $c$. 
So there exists a $c^*$-monochromatic component $M$ in $(G,\phi)^\ell$ with weak diameter in $(G,\phi)^\ell$ greater than $f^*(\eta)$.

Since $c^*$ can be obtained by extending $c_3$, $V(M) \cap (X_{T_0} \cup \bigcup_{e \in U_E}N_{(G,\phi)[X_{T_e}]}^{\leq \ell}[X_e])=\emptyset$ by (i).
Since $M \subseteq (G,\phi)^\ell$, there exists $e \in U_E$ such that $V(M) \subseteq X_{T_e}-Z_1$, and $M$ is a $c_e'$-monochromatic component in $(G_e,\phi_e)^\ell$.
By Claim 5, $M$ has weak diameter in $(G_e,\phi_e)^\ell$ at most $f^*(\eta)$.
But $(G_e,\phi_e)^\ell \subseteq (G,\phi)^\ell$, so $M$ has weak diameter in $(G,\phi)^\ell$ at most $f^*(\eta)$, a contradiction.
This proves the lemma.
\end{pf}

\begin{lemma} \label{weighted_tree_extension_clean}
For any positive real numbers $\ell,\nu$, any positive integers $m,\theta$ with $m \geq 2$, every set $I \subseteq (0,\ell]$, and any $(m,\ell,\nu)$-nice hereditary $I$-bounded classes $\F,\F'$ of weighted graphs, there exists a real number $\nu^*$ such that the following holds.
Let $\F^*$ be an $I$-bounded class of weighted graphs such that for every $(G,\phi) \in \F^*$, there exists a tree-decomposition $(T,\X)$ of $G$ of adhesion at most $\theta$, where $\X=(X_t: t \in V(T))$, such that for every $t \in V(T)$, 
	\begin{itemize}
		\item $(G,\phi)[X_t] \in \F$, and 
		\item $\F'$ contains every $I$-bounded weighted graph that can be obtained from $(G,\phi)[X_t]$ by, for each neighbor $t'$ of $t$ in $T$, adding at most $\theta^2$ new vertices and new edges such that every new edge is either between two new vertices or between a new vertex and $X_t \cap X_{t'}$.
	\end{itemize}
Then $\F^*$ is $(m,\ell,\nu^*)$-nice.
\end{lemma}

\begin{pf}
Let $\ell,\nu$ be positive real numbers, $m,\theta$ be positive integers with $m \geq 2$, $I$ be a subset of $(0,\ell]$, and $\F,\F'$ be $(m,\ell,\nu)$-nice hereditary $I$-bounded classes of weighted graphs.
Let $\F''$ be the union of $\F'$ and the class consisting of all $I$-bounded weighted graphs on at most $2\theta^2+\theta$ vertices.
Since $\F'$ is hereditary, $\F'' \subseteq (\F')^{+(2\theta^2+\theta),(0,\ell]}$ is hereditary.
By Lemma \ref{apex_extension}, $\F''$ is an $(m,\ell,\nu_1)$-nice hereditary $I$-bounded class of weighted graphs, where $\nu_1$ is the real number $\nu^*$ given by Lemma \ref{apex_extension} by taking $(\ell,\nu,n)=(\ell,\nu,2\theta^2+\theta)$.
Define $\nu^*$ to be the integer $f^*(\theta)$, where $f^*$ is the function given by Lemma \ref{weighted_tree_extension} by taking $(\ell,\nu,m,\theta,\F,\F')=(\ell,\nu_1,m,\theta,\F,\F'')$. 
Note that $\nu^* \geq \nu$.

Let $\F^*$ be a class of weighted graphs as stated in the lemma.
Let $(G,\phi) \in \F^*$.
To show that $\F^*$ is $(m,\ell,\nu^*)$-nice, it suffices to show that $(G,\phi)^\ell$ is $m$-colorable with weak diameter in $(G,\phi)^\ell$ at most $\nu^*$.

Since $(G,\phi) \in \F^*$, there exists a tree-decomposition $(T,\X)$ of $G$ of adhesion at most $\theta$, where $\X=(X_t: t \in V(T))$, such that for every $t \in V(T)$, 
	\begin{itemize}
		\item[(i)] $(G,\phi)[X_t] \in \F$, and 
		\item[(ii)] $\F'$ contains every $I$-bounded weighted graph that can be obtained from $(G,\phi)[X_t]$ by, for each neighbor $t'$ of $t$ in $T$, adding at most $\theta^2$ new vertices and new edges such that every new edge is either between two new vertices or between a new vertex and $X_t \cap X_{t'}$.
	\end{itemize}
We may assume that subject to (i) and (ii), $\lvert V(T) \rvert$ is as small as possible.

We shall modify $(T,\X)$ to make it a $(\theta,\theta,\F,\F'',\theta^2)$-construction of $(G,\phi)$ and then apply Lemma \ref{weighted_tree_extension}.

We first show that we may assume that $G$ is connected.
Note that $(G,\phi)^\ell$ is $m$-colorable with weak diameter in $(G,\phi)^\ell$ at most $\nu^*$ if and only if for every component $C$ of $G$, $(C,\phi|_{E(C)})^\ell$ is $m$-colorable with weak diameter in $(C,\phi|_{E(C)})^\ell$ at most $\nu^*$.
In addition, since $\F$ and $\F'$ are hereditary, for every component $C$ of $G$, $(T,(X_t \cap V(C): t \in V(T))$ is a tree-decomposition of adhesion at most $\theta$ satisfying (i) and (ii) (where $(G,\phi)$ is replaced by $(C,\phi|_{E(C)})$).

Hence we may assume that $G$ is connected.
Since $G$ is connected, by the minimality of $|V(T)|$, for every $e \in E(T)$, we know $X_{t_e} \cap X_{t_e'} \neq \emptyset$, where $t_e$ and $t_e'$ are the ends of $e$.

We first assume $\lvert V(T) \rvert=1$.
Then $(G,\phi)=(G,\phi)[X_t] \in \F$, where $t$ is the unique vertex of $T$.
Since $\F$ is $(m,\ell,\nu)$-nice, $(G,\phi)^\ell$ is $m$-colorable with weak diameter in $(G,\phi)^\ell$ at most $\nu \leq \nu^*$.

Hence we may assume $\lvert V(T) \rvert \geq 2$.
Let $e_0 \in E(T)$.
Let $T'$ be the tree obtained from $T$ by subdividing $e_0$ once, and we call the new vertex $t^*$.
Let $X'_{t^*}=X_{t_{e_0}} \cap X_{t_{e_0}'}$, where $t_{e_0}$ and $t_{e_0}'$ are the ends of $e_0$.
Note that $X'_{t^*} \neq \emptyset$.
For every $t \in V(T')-\{t^*\}$, let $X'_t = X_t$.
Let $\X'=(X'_t: t \in V(T'))$.
Then $(T',\X')$ is a tree-decomposition of $G$ of adhesion at most $\theta$.
We make $T'$ a rooted tree by defining the root to be $t^*$.

Since $\F$ is hereditary, $(G,\phi)[X'_{t^*}] \in \F$, and hence $(G,\phi)[X'_{t}] \in \F$ for every $t \in V(T')$.
Note that for every $t \in V(T')-\{t^*\}$, if there exists a neighbor $t'$ of $t$ in $T'$, then there exists a neighbor $t''$ of $t$ in $T$ such that $X'_{t} \cap X'_{t'} = X_t \cap X_{t''}$.
So for every $t \in V(T')-\{t^*\}$, $\F'' \supseteq \F'$ contains every $I$-bounded weighted graph that can be obtained from $(G,\phi)[X'_t]$ by for each neighbor $t'$ of $t$ in $T'$, adding at most $\theta^2$ new vertices and new edges such that every new edge is either between two new vertices or between a new vertex and $X'_t \cap X'_{t'}$.
Since $t^*$ has degree 2 in $T'$ and $\lvert X'_{t^*} \rvert \leq \theta$, every $I$-bounded weighted graph that can be obtained from $(G,\phi)[X'_{t^*}]$ by, for each neighbor $t'$ of $t^*$ in $T'$, adding at most $\theta^2$ new vertices and new edges such that every new edge is either between two new vertices or between a new vertex and $X'_{t^*} \cap X'_{t'}$ has at most $2\theta^2+\theta$ vertices and hence belongs to $\F''$.
Therefore, $(T',\X')$ is a $(\theta,\theta,\F,\F'',\theta^2)$-construction of $(G,\phi)$.
By taking $Z=\emptyset$ in Lemma \ref{weighted_tree_extension}, $(G,\phi)^\ell$ is $m$-colorable with weak diameter in $(G,\phi)^\ell$ at most $\nu^*$.
This proves the lemma.
\end{pf}

\subsection{Construction with weak control} \label{sec:weak_control}

The goal of this subsection is to prove Lemma \ref{strong_weighted_tree_extension_control_clean}, which states that every weighted graph that has a tree-decomposition such that every bag is $(k,r)$-centered (for some $k$ and $r$) has a 2-coloring with small weak diameter.
It will be used to deal with tree-decompositions with arbitrarily large adhesion and is a counterpart of Lemma \ref{weighted_tree_extension_clean}.
The strategy for proving Lemma \ref{strong_weighted_tree_extension_control_clean} is similar to the one for proving Lemma \ref{weighted_tree_extension_clean}, but we do not afford to delete vertices because it would make the bags become non-$(k,r)$-centered. 
We shall leave the vertices that we tried to delete uncolored instead, but there are many technical details that require close attention when actually applying this strategy. 
To prove Lemma \ref{strong_weighted_tree_extension_control_clean}, we prove a more technical version (Lemma \ref{strong_weighted_tree_extension_control}) that allows us to apply a stronger inductive hypothesis.

\subsubsection{Intuition for constructions}
In this section, we briefly explain the intuition of the $(\eta,\theta,\mu,\ell,\ell')$-constructions that will be formally defined in Section \ref{subsubsec:weak_control_construction_def} and will be used to prove Lemma \ref{strong_weighted_tree_extension_control} in Section \ref{subsubsec:weak_control_main_lemma}.
Recall that our goal is to prove an analog of Lemma \ref{weighted_tree_extension}.
But now the adhesion of the tree-decomposition $(T,\X)$ is unbounded, and we only assume that the bags and the adhesion sets are $(k,r)$-centered.

We frequently delete vertices in the proof of Lemma \ref{weighted_tree_extension}.
But this operation is not allowed when proving Lemma \ref{strong_weighted_tree_extension_control}, as it would destroy the property of being $(k,r)$-centered.
The set $R$ in the definition of an $(\eta,\theta,\mu,\ell,\ell')$-construction $(R,T,\X)$ is used to deal with the issue about deleting vertices.
We will move vertices into the set $R$ to make them look like deleted without actually deleting anything.

Moreover, the reason for deleting vertices in the proof of Lemma \ref{weighted_tree_extension} is to reduce the adhesion of the tree-decomposition.
Similarly, in the proof of Lemma \ref{strong_weighted_tree_extension_control}, we would like to show that the subsets consisting of the vertices in the bags or the adhesion sets not covered by $R$ become $(k-1,r)$-centered after we move vertices into $R$.
The formal way that we deal with it is to track the centers of those centered sets.
In particular, $D_e,S_e,O_e$ are subsets of the adhesion set $X_e$ for a tree edge $e$, and we require that $|D_e \cup S_e \cup O_e|$ is small and that $X_e$ is not far from $D_e \cup S_e \cup O_e$.
Roughly, $D_e$ is a set so that $X_e-R$ is not far from $D_e$, and $S_e$ is a set so that $X_e \cap R$ is not far from $S_e$.
However, for some technical reason, it is hard to modify $D_e$ and $S_e$ when we move vertices into $R$.
Instead, we maintain an extra set $O_e$ so that $X_e-R$ is not far from $D_e \cup O_e$, and $X_e \cap R$ is not far from $S_e \cup O_e$.
It is the motivation for the conditions (C4)-(C6) in the definition of an $(\eta,\theta,\mu,\ell,\ell')$-construction; (C1)-(C3) are the analogous conditions for the root bag, like the condition (EC2) used in the past.
Note that $|D_e \cup O_e|$ should be considered as the ``fake adhesion'', and we will do induction on it.
In particular, the conditions (C9) and (C11) are the analogs of (EC1) and (EC3), respectively, under our usage of fake adhesion.

On the other hand, $O_e$ is too powerful in the sense that it is used to deal with both $X_e-R$ and $X_e \cap R$, so we should limit its use.
Conditions (C7) and (C8) show this limitation.
For some technical reason to make the proof work, we should also keep track the ``boundary'' of $R$, which is the motivation of (C10).

\subsubsection{Formal definition of constructions} \label{subsubsec:weak_control_construction_def}

Let $\eta,\theta$ be nonnegative integers with $\eta \leq \theta$.
Let $\mu$ be a nonnegative real number.
Let $\ell,\ell'$ be positive real numbers.
Let $a=6\ell+2\mu$.
We say that a $(0,\ell']$-bounded weighted graph $(G,\phi)$ is \defn{$(\eta,\theta,\mu,\ell,\ell')$-constructible} if there exist $R \subseteq V(G)$ and a rooted tree-decomposition $(T,\X)$ of $G$, where $\X=(X_t: t \in V(T))$ and the root of $T$ is denoted by $t^*$, such that the following hold:
\begin{itemize}
	\item[(C1)] There exist $D_{t^*} \subseteq X_{t^*}-R$, $S_{t^*} \subseteq X_{t^*} \cap R$, and $O_{t^*} \subseteq X_{t^*}$ such that $\lvert D_{t^*} \cup S_{t^*} \cup O_{t^*} \rvert \leq \theta$. 
	\item[(C2)] $X_{t^*} -R \subseteq N_{(G,\phi)}^{\leq \mu}[D_{t^*} \cup O_{t^*}]$. 
	\item[(C3)] $R \cap X_{t^*} \subseteq N_{(G,\phi)}^{\leq \mu}[S_{t^*} \cup O_{t^*}]$. 
	\item[(C4-C9)] For every edge $e=tt' \in E(T)$, (we define $X_e=X_t \cap X_{t'}$)
		\begin{itemize}
			\item[(C4)] there exist $D_e \subseteq X_e-R$, $S_e \subseteq X_e \cap R$, and $O_e \subseteq X_e$ such that $\lvert D_e \cup S_e \cup O_e \rvert \leq \theta$, 
			\item[(C5)] $X_e-R \subseteq N_{(G,\phi)}^{\leq \mu}[D_e \cup O_e]$, 
			\item[(C6)] $R \cap X_e \subseteq N_{(G,\phi)}^{\leq \mu}[S_e \cup O_e]$, 
			\item[(C7)] $(X_e \cap R-N_{(G,\phi)}^{\leq \mu}[O_e]) \cap N_{(G,\phi)}^{\leq 3\ell+\mu}[(D_{t^*} \cup O_{t^*})-R] = \emptyset$, 
			\item[(C8)] $(X_e \cap R-N_{(G,\phi)}^{\leq \mu}[O_e]) \cap N_{(G,\phi)}^{\leq a}[V(G)-(R \cup D_{t^*} \cup O_{t^*})] = \emptyset$,  
			\item[(C9)] if $\lvert D_e \cup O_e \rvert > \eta$, then one end of $e$ has no child, say $t'$, and $X_{t'} \subseteq N_{(G,\phi)}^{\leq \ell}[X_t \cap X_{t'}]$.
		\end{itemize}
    \item[(C10)] There exist $R_- \subseteq V(G)$ with $|R_-| \leq (3\theta+3)^{\theta-\eta}$ and $t_- \in V(T)$ such that for every $rv \in E(G)$ with $r \in R$ and $v \in V(G)-R$, 
        \begin{itemize}
            \item $v \in N_{(G,\phi)}^{\leq \ell}[X_{t^*}] \cup N_{(G,\phi)}^{\leq 2\ell'}[X_{t_-}] \cup N_{(G,\phi)}^{\leq \rho+2\ell'}[R_-]$, and
            \item if $N_{(G,\phi)}^{\leq 3\ell+\mu}[(D_{t^*} \cup O_{t^*})-R]=\emptyset$, then $v \in N_{(G,\phi)}^{\leq 2\ell'}[X_{t_-}] \cup N_{(G,\phi)}^{\leq \rho+2\ell'}[R_-]$,
        \end{itemize}
                where $\rho = a+7\ell+3\mu$.  
    \item[(C11)] If $\eta \geq 1$, then $X_{t^*} \neq \emptyset$ and $X_e \neq \emptyset$ for every $e \in E(T)$.
\end{itemize}
In this case, we call $(R,T,\X)$ an \defn{$(\eta,\theta,\mu,\ell,\ell')$-construction} of $(G,\phi)$.

\subsubsection{Main technical lemma} \label{subsubsec:weak_control_main_lemma}

We prove the main technical lemma (Lemma \ref{strong_weighted_tree_extension_control}) in this section.
It is an analog of Lemma \ref{weighted_tree_extension}, and we will prove it by using a similar strategy.
That is, we will do induction on the ``fake adhesion'' mentioned before, and we will move vertices into $R$ to reduce the fake adhesion.

We sketch its proof before formally stating and proving this lemma.
We first prove the case that the fake adhesion is 0 in Claims 1 and 2, which builds the induction base.
Then we show that we can assume that the precolored set $Z$ is ``full'' in Claims 3-5.
The ``fullness'' of $Z$ allows us to define the ``central part'' $X_{T_0}$ and the condensation $(G_0,\phi_0)$ by seeing what bags intersect the precolored set $Z$.
Claim 6 shows a relationship between the distance in $(G,\phi)$ and in $(G_0,\phi_0)$, which helps us analyze the distance later.
Next, we define a tree-decomposition $(T',\X')$ of the condensation, and then we move vertices to create a new set $R'$.
We show that $(R',T',\X')$ is a construction with smaller fake adhesion than $(R,T,\X)$ in Claim 9 after the preparation work Claims 7 and 8.
Then the inductive hypothesis provides a good coloring $c_0$ of the condensation. 
A tool that we developed earlier allows us to extend $c_0$ by further coloring all vertices close to the central part to obtain a good coloring $c'$ on those vertices.
Note that the vertices in each peripheral part already colored by $c'$ form the precolored set for the peripheral part.
Then, in Claim 10, we apply induction to each peripheral part to show that we are done unless some peripheral part $X_{T_e}$ cannot be further colored.
The remaining goal is to color $X_{T_e}$ nicely to obtain a contradiction.
We move all vertices in $G$ but not in $X_{T_e}$ into $R$ (with some other modification) to obtain a construction $(R^*,T^*,\X^*)$ of $X_{T_e}$ in Claims 11-13.
In Claim 14, we show that we can apply induction to $(R^*,T^*,\X^*)$ unless we are in a ``degenerate case''.
Finally, we deal with the degenerate case by constructing a new construction $(R^1,T,\X)$ in Claim 15 and show that the induction hypothesis can be applied to finish the proof.

\begin{lemma} \label{strong_weighted_tree_extension_control}
For any positive real numbers $\ell,\ell',\mu$ with $\ell' \geq \ell$, any positive integer $m$ with $m \geq 2$, and any nonnegative integer $\theta$, there exists a function $f^*: {\mathbb N} \cup \{0\} \rightarrow {\mathbb R}^+$ such that the following holds.

Let $\eta$ be a nonnegative integer with $\eta \leq \theta$.
Let $(G,\phi)$ be an $(\eta,\theta,\mu,\ell,\ell')$-constructible $(0,\ell']$-bounded weighted graph with an $(\eta,\theta,\mu,\ell,\ell')$-construction $(R,T,\X)$ such that $(G,\phi)$ is $(T,\X,\ell)$-bounded.
Denote $\X$ by $(X_t: t \in V(T))$. 
Let $t^*$ be the root of $T$.
Let $Z \subseteq N_{(G,\phi)}^{\leq 3\ell+\mu}[(D_{t^*} \cup O_{t^*})-R]$, where $D_{t^*}$ and $O_{t^*}$ are the subsets of $X_{t^*}$ mentioned in the definition of an $(\eta,\theta,\mu,\ell,\ell')$-construction. 
Assume that for every $x \in V(T)$, there exists $A_x \subseteq X_x$ with $\lvert A_x \rvert \leq \theta$ such that $X_x \subseteq N_{(G,\phi)}^{\leq 3\ell+\mu}[A_x]$.

If $c: Z \rightarrow [m]$ is a function, then $c|_{Z-R}$ can be extended to an $m$-coloring of $(G,\phi)^\ell[V(G)]-R$ with weak diameter in $(G,\phi)^\ell$ at most $f^*(\eta)$.
\end{lemma}

\begin{pf}
Let $\ell,\ell',\mu$ be positive real numbers with $\ell' \geq \ell$, $m$ be a positive integer with $m \geq 2$, and $\theta$ be a nonnegative integer. 
We define the following.
	\begin{itemize}
        \item Let $a=6\ell+2\mu$ and let $\rho=a+7\ell+3\mu$. 
        Note that $a$ and $\rho$ are the real numbers mentioned in the definition of an $(\eta,\theta,\mu,\ell,\ell')$-construction.
		\item Let $\nu_1$ be the number $\nu$ given by Lemma \ref{all_centered} by taking $(k,r,\ell)=((3\theta+3)^\theta+\theta,6\ell'+\mu+\rho,\ell)$. 
		\item Let $f_2: {\mathbb R}^+ \rightarrow {\mathbb R^+}$ be the function such that for every $x \in {\mathbb R}^+$, $f_2(x)$ is the number $\nu^*$ given by Lemma \ref{patching_centered_set} by taking $(k,r,\ell,\nu)=(\theta,3\ell+3\mu+a,\ell,x)$. 
        \item Let $\tau = 2\rho+16\ell'+2\mu$.
        \item Let $\ell_0' = \ell' + 3\ell+\tau$.
		\item Let $f_3: {\mathbb R}^+ \rightarrow {\mathbb R^+}$ be the function such that for every $x \in {\mathbb R}^+$, $f_3(x)$ is the number $\nu^*$ given by Lemma \ref{con_color} by taking $(\ell,\ell',\nu,\mu,\lambda,m,\theta)=(\ell,\ell',x,\mu,\tau,m,\max\{\theta,1\})$.  
		\item Define $f^*: ({\mathbb N} \cup \{0\}) \rightarrow {\mathbb R}^+$ to be the function such that $f^*(0)=\nu_1$, and for every $x \in {\mathbb N}$, $f^*(x)= f_3(f_2(f^*(x-1)))$. 
	\end{itemize}
Note that $f^*$ is an increasing function.

Let $\eta,(G,\phi),(R,T,\X),t^*,Z,c$ be the ones as defined in the lemma.
Let $R_-,D_{t^*},S_{t^*},O_{t^*},D_e, \allowbreak S_e,O_e$ (for every $e \in E(T)$) be the sets and $t_-$ be the node of $T$ stated in the definition of an $(\eta,\theta,\mu,\ell,\ell')$-construction. 
Suppose to the contrary that $c|_{Z-R}$ cannot be extended to an $m$-coloring of $(G,\phi)^\ell[V(G)]-R$ with weak diameter in $(G,\phi)^\ell$ at most $f^*(\eta)$, and subject to this, the tuple $(\eta,|V(G) -(Z \cup R)|+|V(G)-(R \cup D_{t^*} \cup O_{t^*})|, |V(G)-R|, |V(G)-Z|)$ is lexicographically minimal among all tuples $(\eta,(G,\phi),(R,T,\X),t^*,Z,c)$ satisfying the condition of this lemma.  

For every $e=tt' \in E(T)$, let $X_e=X_t \cap X_{t'}$.

\medskip

\noindent{\bf Claim 1:} For every $e \in E(T)$ with $D_e \cup O_e=\emptyset$, if $X_{T_1}-R \neq \emptyset$ and $X_{T_2}-R \neq \emptyset$, where $T_1$ and $T_2$ are the components of $T-e$, then there exists no path in $G$ of length in $(G,\phi)$ at most $\ell$ between $V(G)-R$ and $X_e$.

\noindent{\bf Proof of Claim 1:}
Let $e$ be an edge of $T$ satisfying the condition of this claim. 
Suppose to the contrary that there exists a path $P$ in $G$ with $\leng_{(G,\phi)}(P) \leq \ell$ between $V(G)-R$ and $X_e$.
Since $D_e \cup O_e=\emptyset$, (C5) implies $X_e-R \subseteq N_{(G,\phi)}^{\leq \mu}[D_e \cup O_e]=\emptyset$, so $X_e=X_e \cap R = X_e \cap R-N_{(G,\phi)}^{\leq \mu}[O_e]$.
Let $v$ be a vertex in $V(P) \cap X_e$, and let $u$ be the end of $P$ in $V(G)-R$.
So $v \in X_e \cap N_{(G,\phi)}^{\leq \ell}[\{u\}] = (X_e \cap R-N_{(G,\phi)}^{\leq \mu}[O_e]) \cap N_{(G,\phi)}^{\leq \ell}[\{u\}]$.
Since $u \not \in R$, we know $u \in (V(G)-(R \cup D_{t^*} \cup O_{t^*})) \cup ((D_{t^*} \cup O_{t^*})-R)$.
Since $a\geq 3\ell+\mu>\ell$, by (C7) and (C8), 
	\begin{align*}
		v \in & (X_e \cap R-N_{(G,\phi)}^{\leq \mu}[O_e]) \cap N_{(G,\phi)}^{\leq \ell}[\{u\}] \\
		\subseteq & (X_e \cap R-N_{(G,\phi)}^{\leq \mu}[O_e]) \cap N_{(G,\phi)}^{\leq \ell}[(V(G)-(R \cup D_{t^*} \cup O_{t^*})) \cup ((D_{t^*} \cup O_{t^*})-R)] = \emptyset,
	\end{align*}
a contradiction. 
$\Box$

\medskip

\noindent{\bf Claim 2:} $\theta \geq \eta \geq 1$.

\noindent{\bf Proof of Claim 2:}
Since $\theta \geq \eta$ by assumption, it suffices to show $\eta \geq 1$.
Suppose to the contrary that $\eta=0$.
Let $W=\{e \in E(T): D_e \cup O_e=\emptyset\}$.
By (C5), for every $e \in W$, we have $X_e-R \subseteq N_{(G,\phi)}^{\leq \mu}[D_e \cup O_e]=\emptyset$.
So for every component $C$ of $G-R$, we know $V(C) \subseteq X_{C'}$ for some component $C'$ of $T-W$.
Let $T_1,T_2$ be the components of $T-e$. 

Suppose that there exist $e \in W$ and an edge in $(G,\phi)^\ell$ between $X_{T_1}-R$ and $X_{T_2}-R$. 
Then there exists a path $P$ in $(G,\phi)$ between a vertex $v_1 \in X_{T_1}-R$ and a vertex $v_2 \in X_{T_2}-R$ of length in $(G,\phi)$ at most $\ell$.
Since $e \in W$, we know $X_e-R=\emptyset$, so $v_1,v_2 \not \in X_e$.
So there exists a subpath $P'$ of $P$ between $v_1 \in X_{T_1}-R$ and $X_e$ of length in $(G,\phi)$ at most $\ell$, contradicting Claim 1.

Hence for every $e \in W$, there exists no edge in $(G,\phi)^\ell$ between $X_{T_1}-R$ and $X_{T_2}-R$. 
Therefore, if for every component $C$ of $T-W$, $c|_{Z \cap X_C-R}$ can be extended to an $m$-coloring of $(G,\phi)^\ell[X_C]-R$ with weak diameter in $(G,\phi)^\ell$ at most $f^*(\eta)$, then $c|_{Z-R}$ can be extended to an $m$-coloring of $(G,\phi)^\ell[V(G)]-R$ with weak diameter in $(G,\phi)^\ell$ at most $f^*(\eta)$.
Recall that we assume that $c|_{Z-R}$ cannot be extended to an $m$-coloring of $(G,\phi)^\ell[V(G)]-R$ with weak diameter in $(G,\phi)^\ell$ at most $f^*(\eta)$.
So there exists a component $C$ of $T-W$ such that $c|_{Z \cap X_C-R}$ cannot be extended to an $m$-coloring of $(G,\phi)^\ell[X_C]-R$ with weak diameter in $(G,\phi)^\ell$ at most $f^*(\eta)$.

Since $(R,T,\X)$ is an $(0,\theta,\mu,\ell,\ell')$-construction, $C$ is a star by (C9) and the definition of $W$.
Let $t_C$ be the root of $C$.
By (C9), $X_C \subseteq N_{(G,\phi)}^{\leq \ell}[X_{t_C}]$. 
Since there exists $A_{t_C} \subseteq X_{t_C}$ with $\lvert A_{t_C} \rvert \leq \theta$ such that $X_{t_C} \subseteq N_{(G,\phi)}^{3\ell+\mu}[A_{t_C}]$, we know $X_C \subseteq N_{(G,\phi)}^{\leq 4\ell+\mu}[A_{t_C}]$. 
So any $m$-coloring of $(G,\phi)^\ell[X_C-R]$ has weak diameter in $(G,\phi)^\ell$ at most $\nu_1 \leq f^*(\eta)$ by Lemma \ref{all_centered} (taking $R=R \cup (V(G)-X_C)$).
Hence $c|_{Z \cap X_C-R}$ can be extended to an $m$-coloring of $(G,\phi)^\ell[X_C]-R$ with weak diameter in $(G,\phi)^\ell$ at most $f^*(\eta)$, a contradiction.
$\Box$

\medskip

Next, we study the set $Z$ in Claims 3 and 5.

\medskip

\noindent{\bf Claim 3:} $Z = N_{(G,\phi)}^{\leq 3\ell+\mu}[(D_{t^*} \cup O_{t^*})-R]$. 

\noindent{\bf Proof of Claim 3:}
Suppose to the contrary that there exists $v \in N_{(G,\phi)}^{\leq 3\ell+\mu}[(D_{t^*} \cup O_{t^*})-R]-Z$.
Let $Z' = Z \cup \{v\}$.
Let $c': Z' \rightarrow [m]$ be the function obtained from $c$ by further defining $c'(v)=m$.
Since $Z' \subseteq N_{(G,\phi)}^{\leq 3\ell+\mu}[(D_{t^*} \cup O_{t^*})-R]$, the minimality of $(\eta,|V(G)-(Z \cup R)|+|V(G)-(R \cup D_{t^*} \cup O_{t^*})|,|V(G)-R|,|V(G)-Z|)$ implies that $c'|_{Z'-R}$ (and hence $c|_{Z-R}$) can be extended to an $m$-coloring of $(G,\phi)^\ell[V(G)]-R$ with weak diameter in $(G,\phi)^\ell$ at most $f^*(\eta)$, a contradiction.
$\Box$

\medskip

\noindent{\bf Claim 4:} If $Z = \emptyset$, then there exists $t_0 \in V(T)$ with $X_{t_0}-R \neq \emptyset$ such that either $t_0=t^*$, or $t_0$ has a child in $T$, or $X_{t_0} \not \subseteq N_{(G,\phi)}^{\leq \ell}[X_{e_{t_0}}]$, where $e_{t_0}$ is the edge of $T$ between $t_0$ and its parent. 

\noindent{\bf Proof of Claim 4:}
Assume $Z=\emptyset$.
Suppose to the contrary that for every $t \in V(T)$ with $X_t-R \neq \emptyset$, $t$ has the parent $p_t$ in $T$, $t$ has no child in $T$ and the edge $e_t$ of $T$ between $t$ and $p_t$ satisfies $X_{t} \subseteq N_{(G,\phi)}^{\leq \ell}[X_{e_{t}}]$.
So for every $t \in V(T)$ with $X_t-R \neq \emptyset$, we have $X_{p_t} \subseteq R$.
Hence for every $v \in X_t-R$, there exists a path $P_v$ in $G$ from $v$ to $X_{e_t} \subseteq X_{p_t} \subseteq R$ with $\leng_{(G,\phi)}(P_v) \leq \ell$.
So there exists $r_vq_v \in E(P_v)$ with $r_v \in R$ and $q_v \in V(G)-R$.
By Claim 3, $N_{(G,\phi)}^{\leq 3\ell+\mu}[(D_{t^*} \cup O_{t^*})-R]=Z=\emptyset$.
So by (C10), $q_v \in N_{(G,\phi)}^{\leq 2\ell'}[X_{t_-}] \cup  N_{(G,\phi)}^{\leq \rho+2\ell'}[R_-]$.
Hence $v \in N_{(G,\phi)}^{\leq \leng_{(G,\phi)}(P_v)}[\{q_v\}] \subseteq N_{(G,\phi)}^{\leq 3\ell'}[X_{t_-}] \cup  N_{(G,\phi)}^{\leq \rho+3\ell'}[R_-]$. 
So $$V(G)-R \subseteq N_{(G,\phi)}^{\leq 3\ell'}[X_{t_-}] \cup  N_{(G,\phi)}^{\leq \rho+3\ell'}[R_-] \subseteq N_{(G,\phi)}^{\leq 6\ell'+\mu}[A_{t_-}] \cup N_{(G,\phi)}^{\leq \rho+3\ell'}[R_-].$$ 
Note that $|A_{t_-}|+|R_-| \leq \theta+(3\theta+3)^\theta$, so $V(G)-R$ is $((3\theta+3)^\theta+\theta,6\ell'+\mu+\rho)$-centered in $(G,\phi)$. So by Lemma \ref{all_centered}, $c|_{Z-R}=c|_\emptyset$ can be extended to an $m$-coloring of $(G,\phi)^\ell[V(G)]-R$ with weak diameter in $(G,\phi)^\ell$ at most $\nu_1 \leq f^*(\eta)$, a contradiction.
$\Box$

\medskip

\noindent{\bf Claim 5:} $Z \neq \emptyset$.

\noindent{\bf Proof of Claim 5:}
Suppose to the contrary that $Z=\emptyset$.

By Claim 4, there exists $t_0 \in V(T)$ with $X_{t_0}-R \neq \emptyset$ such that either $t_0=t^*$, or $t_0$ has a child in $T$, or $X_{t_0} \not \subseteq N_{(G,\phi)}^{\leq \ell}[X_{e_{t_0}}]$, where $e_{t_0}$ is the edge of $T$ between $t_0$ and its parent.
Let $v \in X_{t_0}-R$. 

Since $Z=\emptyset$, by Claim 3, $(D_{t^*} \cup O_{t^*})-R=\emptyset$. 
Since $v \not \in R$ and $(D_{t^*} \cup O_{t^*})-R=\emptyset$, $v \not \in D_{t^*} \cup O_{t^*}$.

Let $T'$ be the rooted tree obtained from $T$ by adding a new node $t_1$ adjacent to $t_0$, where $t_1$ is the root of $T'$.
Let $X'_{t_1}=\{v\}$.
For every $t \in V(T)$, let $X'_t=X_t$.
Let $\X'=(X_t': t\in V(T'))$.
Then $(T',\X')$ is a tree-decomposition of $G$.
By the choice of $t_0$, $(G,\phi)$ is $(T',\X',\ell)$-bounded.

We shall prove that $(R,T',\X')$ is an $(\eta,\theta,\mu,\ell,\ell')$-construction of $(G,\phi)$.

For any $x \in \{t_1,t_0t_1\}$, let $D_x=\{v\}$, $O_x=\{v\}$ and $S_x=\emptyset$.
So for $x \in \{t_1,t_0t_1\}$, we have $X'_x=\{v\}$, $D_x=\{v\} \subseteq X'_x-R$, $S_x=\emptyset \subseteq X'_x \cap R$, $O_x=\{v\} \subseteq X'_x$, $|D_x \cup S_x \cup O_x|=1 \leq \theta$ (where $1 \leq \theta$ follows from Claim 2), $X'_x-R = \{v\} \subseteq N_{(G,\phi)}^{\leq \mu}[D_x \cup O_x]$, and $R \cap X'_x = \emptyset$.
Since $(R,T,\X)$ satisfies (C4)-(C6), $(R,T',\X')$ satisfies (C1)-(C6).

Now we show that $(R,T',\X')$ satisfies (C7) and (C8).
Note that (C7) and (C8) hold for the edge $t_0t_1$ since $X'_{t_0t_1} \cap R = \emptyset$.
Since $3\ell+\mu \leq a$, and $(D_{t_1} \cup O_{t_1})-R = \{v\} \subseteq V(G)-(R \cup D_{t^*} \cup O_{t^*})$, we know that for every $e \in E(T')-\{t_0t_1\} = E(T)$, 
	\begin{align*}
		& (X'_e \cap R -N_{(G,\phi)}^{\leq \mu}[O_e]) \cap N_{(G,\phi)}^{\leq 3\ell+\mu}[(D_{t_1} \cup O_{t_1})-R] \\
		\subseteq & (X_e \cap R -N_{(G,\phi)}^{\leq \mu}[O_e]) \cap N_{(G,\phi)}^{\leq a}[V(G)-(R \cup D_{t^*} \cup O_{t^*})] =\emptyset
	\end{align*}
by (C8).
So $(R,T',\X')$ satisfies (C7).
Since $(D_{t^*} \cup O_{t^*})-R=\emptyset$, we know $V(G)-(R \cup D_{t_1} \cup O_{t_1}) \subseteq V(G)-R = V(G)-(R \cup D_{t^*} \cup O_{t^*})$.
So for every $e \in E(T')-\{t_0t_1\}=E(T)$, 
	\begin{align*}
		& (X'_e \cap R -N_{(G,\phi)}^{\leq \mu}[O_e]) \cap N_{(G,\phi)}^{\leq a}[V(G)-(R \cup D_{t_1} \cup O_{t_1})] \\
		\subseteq & (X_e \cap R -N_{(G,\phi)}^{\leq \mu}[O_e]) \cap N_{(G,\phi)}^{\leq a}[V(G)-(R \cup D_{t^*} \cup O_{t^*})] =\emptyset
	\end{align*}
by (C8).
Hence $(R,T',\X')$ satisfies (C8).

Now we show that $(R,T',\X')$ satisfies (C9).
Let $e \in E(T')$ with $\lvert D_e \cup O_e \rvert >\eta$.
Since $\eta \geq 1$ by Claim 2, $|D_{t_0t_1} \cup O_{t_0t_1}|=|\{v\}|=1 \leq \eta$.
So $e \in E(T)$.
The condition (C9) for $(R,T,\X)$ implies that one end of $e$ has no child in $T$.
Moreover, the condition (C9) and the choice of $t_0$ imply that either $t_0=t^*$, or $t_0$ has a child in $T$, or $e$ is not the edge of $T$ between $t_0$ and its parent.
Hence the parent end of $e$ in $T'$ is the parent end of $e$ in $T$.
So $(R,T',\X')$ satisfies (C9).

Now we show that $(R,T',\X')$ satisfies (C10).
Note that $t_- \in V(T)$, so $X'_{t_-}=X_{t_-}$.
For every $ry \in E(G)$ with $r \in R$ and $y \in V(G)-R$, since $(R,T,\X)$ satisfies (C10) and $N_{(G,\phi)}^{\leq 3\ell+\mu}[(D_{t^*} \cup O_{t^*})-R]=Z=\emptyset$ by Claim 3, we know $y \in N_{(G,\phi)}^{\leq 2\ell'}[X'_{t_-}] \cup N_{(G,\phi)}^{\leq \rho+2\ell'}[R_-] \subseteq N_{(G,\phi)}^{\leq \ell}[X'_{t_1}] \cup N_{(G,\phi)}^{\leq 2\ell'}[X'_{t_-}] \cup N_{(G,\phi)}^{\leq \rho+2\ell'}[R_-]$.
Hence $(R,T',\X')$ satisfies (C10).

Moreover, $X'_{t_1}=X'_{t_0t_1}=\{v\}$.
So $(R,T',\X')$ satisfies (C11) since $(R,T,\X)$ satisfies (C11).
Hence $(R,T',\X')$ is an $(\eta,\theta,\mu,\ell,\ell')$-construction of $(G,\phi)$.

Let $Z' = \{v\}$.
So $Z' \subseteq N_{(G,\phi)}^{\leq 3\ell+\mu}[(D_{t_1} \cup O_{t_1})-R]$.
Let $c'$ be the $m$-coloring obtained from $c$ by further defining $c'(v)=m$.
Since $D_{t^*} \cup O_{t^*}-R=\emptyset=Z$, we know that $(\eta,\lvert V(G)-(Z' \cup R) \rvert + \lvert V(G)-(R \cup D_{t_1} \cup O_{t_1}) \rvert)$ is lexicographically smaller than $(\eta,\lvert V(G)-(Z \cup R) \rvert + \lvert V(G)-(R \cup D_{t^*} \cup O_{t^*}) \rvert)$.
So $c'|_{Z'-R}$ (and hence $c|_{Z-R}$) can be extended to an $m$-coloring of $(G,\phi)^{\ell}[V(G)]-R$ with weak diameter in $(G,\phi)^\ell$ at most $f^*(\eta)$, a contradiction.
$\Box$

\medskip

Let $$T_0=T[\{t \in V(T): X_t \cap Z \neq \emptyset\}].$$
By (C1) and Claims 3 and 5, $T_0$ is a subtree of $T$ containing $t^*$.
By Claim 3, for every $e \in E(T)$, $X_e \cap Z \neq \emptyset$ if and only if $e \in E(T_0)$.
We define the following:
	\begin{itemize}
		\item Let $U_E$ be the set of edges of $T$ with exactly one end in $T_0$.
		\item Let $(G_0,\phi_0)$ be the $(U_E,\emptyset,\ell,\theta,\tau)$-condensation of $(G,\phi,T,\X)$.
	\end{itemize}
So $V(G_0) \subseteq V(G)$.

\medskip

\noindent{\bf Claim 6:} For any $W,W' \subseteq V(G_0)$ and for every real number $x$ with $0 \leq x \leq 3\ell+\tau$, $$N_{(G,\phi)}^{\leq x}[W] \cap W' \subseteq N_{(G_0,\phi_0)}^{\leq x}[W] \cap W'.$$

\noindent{\bf Proof of Claim 6:}
Since $V(G_0) \subseteq V(G)$ and $(G_0,\phi_0)$ is the $(U_E,\emptyset,\ell,\theta,\tau)$-condensation of $(G,\phi,T,\X)$, by Statement 1 of Lemma \ref{con_qu_iso}, we know that for any $W,W' \subseteq V(G_0)$ and for every real number $x$ with $0 \leq x \leq 3\ell+\tau$, $N_{(G,\phi)}^{\leq x}[W] \cap W' \subseteq N_{(G_0,\phi_0)}^{\leq x}[W] \cap W'$.
$\Box$

\medskip

For every $e \in E(T_0)$, let $$R_e = \{v \in D_e \cup O_e: N_{(G,\phi)}^{\leq \mu}[\{v\}] \cap Z \cap X_e \neq \emptyset\}.$$

\medskip

\noindent{\bf Claim 7:} For every $e \in E(T_0)$, $R_e \neq \emptyset$. 

\noindent{\bf Proof of Claim 7:}
Let $e \in E(T_0)$.
By (C7), $X_e \cap R \cap N_{(G,\phi)}^{\leq 3\ell+\mu}[(D_{t^*} \cup O_{t^*})-R] \subseteq X_e \cap N_{(G,\phi)}^{\leq \mu}[O_e]$.
By (C5), $X_e - R \subseteq N_{(G,\phi)}^{\leq \mu}[D_e \cup O_e]$.
Since $e \in E(T_0)$, we know $Z \cap X_e \neq \emptyset$.
By Claim 3, $Z = N_{(G,\phi)}^{\leq 3\ell+\mu}[(D_{t^*} \cup O_{t^*})-R]$.
So 
	\begin{align*}
		\emptyset \neq X_e \cap Z = & X_e \cap Z \cap N_{(G,\phi)}^{\leq 3\ell+\mu}[(D_{t^*} \cup O_{t^*})-R] \\
		= & Z \cap ((X_e \cap R) \cup (X_e-R)) \cap N_{(G,\phi)}^{\leq 3\ell+\mu}[(D_{t^*} \cup O_{t^*})-R] \\
		\subseteq & Z \cap X_e \cap (N_{(G,\phi)}^{\leq \mu}[O_e] \cup N_{(G,\phi)}^{\leq \mu}[D_e \cup O_e]) \\
		= & Z \cap X_e \cap N_{(G,\phi)}^{\leq \mu}[D_e \cup O_e].
	\end{align*}
Hence $R_e \neq \emptyset$.
$\Box$

\medskip

Our next goal is to construct a tuple $(R',T',\X')$ and show that it is an $(\eta-1,\theta,\mu,\ell,\ell_0')$-construction of $(G_0,\phi_0)$.
We define the following:
	\begin{itemize}
		\item Let $R_{t^*}=(D_{t^*} \cup O_{t^*})-R$. (By Claims 3 and 5, $R_{t^*} \neq \emptyset$.)
		\item Define $R' = (R \cap V(G_0)) \cup N_{(G_0,\phi_0)}^{\leq a+\mu}[\bigcup_{e \in E(T_0) \cup \{t^*\}}R_e]$.
		\item For every $e \in U_E$, let $T_e$ be the component of $T-e$ disjoint from $t^*$.
		\item Let $T'$ be the rooted tree obtained from $T_0$ by, for each $e \in U_E$, adding a new node $t_e$ adjacent to the end of $e$ in $T_0$.
		\item For every $t \in V(T_0)$, let $X'_t = X_t$.
		\item For every $t \in V(T')-V(T_0)$, we know $t=t_e$ for some $e \in U_E$, and we let $X'_t = N_{(G,\phi)[X_{T_e}]}^{\leq \ell}[X_e]$. 
		\item Let $\X'=(X'_t: t \in V(T'))$.
	\end{itemize}
Clearly, $(T',\X')$ is a tree-decomposition of $G_0$.

Recall that $\ell_0' = \ell' + 3\ell+\tau$.
Since $(G_0,\phi_0)$ is the $(U_E,\emptyset,\ell,\theta,\tau)$-condensation of $(G,\phi,T,\X)$, we know that $(G_0,\phi_0)$ is a $(0,\max\{\ell',3\ell+\tau\}]$-bounded weighted graph and hence is a $(0,\ell_0']$-bounded weighted graph.
By the definition of $(T',\X')$, we know that $(G_0,\phi_0)$ is $(T',\X',\ell)$-bounded.

We shall prove that $(R',T',\X')$ is an $(\eta-1,\theta,\mu,\ell,\ell_0')$-construction of $(G_0,\phi_0)$ in Claim 9.
The following claim is a preparation.

\medskip

\noindent{\bf Claim 8:} There exist $R'_- \subseteq V(G_0)$ with $|R'_-| \leq (3\theta+3)^{\theta-(\eta-1)}$ and $t_-' \in V(T')$ such that for every $rv \in E(G_0)$ with $r \in R'$ and $v \in V(G_0)-R'$, we have $v \in N_{(G_0,\phi_0)}^{\leq 2\ell_0'}[X'_{t'_-}] \cup N_{(G_0,\phi_0)}^{\leq \rho+2\ell_0'}[R'_-]$.

\noindent{\bf Proof of Claim 8:}
If $t_- \in V(T_0)$, then define $t'_-=t_-$ and $A'_{t'_-}=A_{t_-}$; otherwise, there exists $e_- \in U_E$ such that $t_- \in V(T_{e_-})$, and we define $t_-'=t_{e_-}$ and $A'_{t'_-}=D_{e_-} \cup O_{e_-} \cup S_{e_-}$.
Define $$R'_- = (R_- \cap V(G_0)) \cup \bigcup_{e \in U_E, R_- \cap X_{T_e}-X_e \neq \emptyset}(D_e \cup S_e \cup O_e) \cup (D_{t^*} \cup O_{t^*} \cup S_{t^*}) \cup A'_{t_-'}.$$ 
Note that for distinct $e_1,e_2 \in U_E$, we know that $R_- \cap X_{T_{e_1}}-X_{e_1}$ and $R_- \cap X_{T_{e_2}}-X_{e_2}$ are disjoint.
So there are at most $|R_-|$ edges $e$ in $U_E$ such that $R_- \cap X_{T_e}-X_e \neq \emptyset$.
Moreover, $|D_e \cup S_e \cup O_e| \leq \theta$ for each $e \in U_E \cup \{t^*\}$ by (C1) and (C4). 
So by (C10), 
	\begin{align*}
		|R'_-| \leq & |R_-|+ |R_-| \cdot \theta +|D_{t^*} \cup O_{t^*} \cup S_{t^*}|+|A'_{t_-'}| \\
		\leq & |R_-|(\theta+1)+2\theta \\
		\leq & (3\theta+3)^{\theta-\eta}(\theta+1)+2\theta \leq 3(3\theta+3)^{\theta-\eta}(\theta+1) = (3\theta+3)^{\theta-(\eta-1)}.
	\end{align*}

Note that $X'_{t^*}=X_{t^*} \subseteq N_{(G,\phi)}^{\leq \mu}[D_{t^*} \cup O_{t^*} \cup S_{t^*}] \subseteq N_{(G,\phi)}^{\leq \mu}[R'_-]$ by (C2) and (C3).
Since $0 \leq \mu \leq \tau$ and $X_{t^*} \cup R'_- \subseteq V(G_0)$, we know $X'_{t^*} = N_{(G,\phi)}^{\leq \mu}[R'_-] \cap X_{t^*} \subseteq N_{(G_0,\phi_0)}^{\leq \mu}[R'_-]$ by Claim 6.

Suppose that this claim does not hold.
So there exists $rv \in E(G_0)$ with $r \in R'$ and $v \in V(G_0)-R'$ such that $v \not \in N_{(G_0,\phi_0)}^{\leq 2\ell_0'}[X'_{t'_-}] \cup N_{(G_0,\phi_0)}^{\leq \rho+2\ell_0'}[R'_-]$.

Suppose to the contrary that $r \in R'-(R \cap V(G_0))$.
Then $r \in N_{(G_0,\phi_0)}^{\leq a+\mu}[\bigcup_{e \in E(T_0) \cup \{t^*\}}R_e]$.
Recall that $R_{t^*} = (D_{t^*} \cup O_{t^*})-R$.
Note that for $e \in E(T_0)$, we have $R_e \subseteq N_{(G,\phi)}^{\leq \mu}[Z] \subseteq N_{(G,\phi)}^{\leq 3\ell+2\mu}[(D_{t^*} \cup O_{t^*})-R]$.
So $\bigcup_{e \in E(T_0) \cup \{t^*\}}R_e \subseteq N_{(G,\phi)}^{\leq 3\ell+2\mu}[(D_{t^*} \cup O_{t^*})-R] \cap V(G_0)$.
Since $0 \leq 3\ell+2\mu \leq \tau$ and $(D_{t^*} \cup O_{t^*})-R \subseteq V(G_0)$, by Claim 6, $$\bigcup_{e \in E(T_0) \cup \{t^*\}}R_e \subseteq N_{(G,\phi)}^{\leq 3\ell+2\mu}[(D_{t^*} \cup O_{t^*})-R] \cap V(G_0) \subseteq N_{(G_0,\phi_0)}^{\leq 3\ell+2\mu}[(D_{t^*} \cup O_{t^*})-R] \subseteq N_{(G_0,\phi_0)}^{\leq 3\ell+2\mu}[R'_-].$$
Hence $r \in N_{(G_0,\phi_0)}^{\leq a+\mu}[\bigcup_{e \in E(T_0) \cup \{t^*\}}R_e] \subseteq N_{(G_0,\phi_0)}^{\leq a+3\ell+3\mu}[R'_-] \subseteq N_{(G_0,\phi_0)}^{\leq \rho}[R'_-]$.
Since $(G_0,\phi_0)$ is $(0,\ell_0']$-bounded, $v \in N_{(G_0,\phi_0)}^{\leq \ell_0'}[\{r\}] \subseteq N_{(G_0,\phi_0)}^{\leq \rho+\ell_0'}[R'_-]$, a contradiction.

So $r \in R \cap V(G_0)$.
Note that $v \in V(G_0)-R' \subseteq V(G)-R$.
If $rv \in E(G)$, then since $(R,T,\X)$ satisfies (C10), we have $v \in N_{(G,\phi)}^{\leq \ell}[X_{t^*}] \cup N_{(G,\phi)}^{\leq 2\ell'}[X_{t_-}] \cup N_{(G,\phi)}^{\leq \rho+2\ell'}[R_-]$, and we let $y_v=v$.
If $rv \not \in E(G)$, then by the definition of $G_0$, there exist $e_1 \in U_E$ such that $\{r,v\} \subseteq N_{(G,\phi)[X_{T_{e_1}}]}^{\leq \ell}[X_{e_1}]-X_{e_1}$ and $\dist_{(G,\phi)[X_{T_{e_1}}]}(r,v) \leq 3\ell+\tau$, so there exists a path $P$ in $(G,\phi)[X_{T_{e_1}}]$ between $r \in R$ and $v \in V(G)-R$ with length in $(G,\phi)[X_{T_{e_1}}]$ at most $3\ell+\tau$, and hence there exists $r'v' \in E(P)$ such that $r' \in R$ and $v' \in V(G)-R$, implying that $v' \in N_{(G,\phi)}^{\leq \ell}[X_{t^*}] \cup N_{(G,\phi)}^{\leq 2\ell'}[X_{t_-}] \cup N_{(G,\phi)}^{\leq \rho+2\ell'}[R_-]$ by (C10); note that there exists $y_v \in \{r,v\}$ such that $\dist_{(G,\phi)}(y_v,v') \leq \frac{1}{2} \cdot \leng_{(G,\phi)}(P) \leq (3\ell+\tau)/2$, so $y_v \in N_{(G,\phi)}^{\leq (3\ell+\tau)/2}[\{v'\}] \subseteq N_{(G,\phi)}^{\leq (3\ell+\tau)/2+\ell}[X_{t^*}] \cup N_{(G,\phi)}^{\leq (3\ell+\tau)/2+2\ell'}[X_{t_-}] \cup N_{(G,\phi)}^{\leq (3\ell+\tau)/2+\rho+2\ell'}[R_-]$.

Hence in either case, $$y_v \in \{r,v\} \cap (N_{(G,\phi)}^{\leq (3\ell+\tau)/2+\ell}[X_{t^*}] \cup N_{(G,\phi)}^{\leq (3\ell+\tau)/2+2\ell'}[X_{t_-}] \cup N_{(G,\phi)}^{\leq (3\ell+\tau)/2+\rho+2\ell'}[R_-]).$$
Since $(G_0,\phi_0)$ is $(0,\ell_0']$-bounded and $v \not \in N_{(G_0,\phi_0)}^{\leq 2\ell_0'}[X'_{t'_-}] \cup N_{(G_0,\phi_0)}^{\leq \rho+2\ell_0'}[R'_-]$, we know $$y_v \not \in N_{(G_0,\phi_0)}^{\leq \ell_0'}[X'_{t'_-}] \cup N_{(G_0,\phi_0)}^{\leq \rho+\ell_0'}[R'_-].$$

Recall $X_{t^*}=X'_{t^*} \subseteq N_{(G_0,\phi_0)}^{\leq \mu}[R'_-]$.
If $y_v \in N_{(G,\phi)}^{\leq (3\ell+\tau)/2+\ell}[X_{t^*}]$, then since $(3\ell+\tau)/2+\ell \leq \tau$, by Claim 6, $y_v \in N_{(G_0,\phi_0)}^{\leq (3\ell+\tau)/2+\ell}[X_{t^*}] \subseteq N_{(G_0,\phi_0)}^{\leq (3\ell+\tau)/2+\ell+\mu}[R'_-] \subseteq N_{(G_0,\phi_0)}^{\leq \ell_0'}[R'_-]$, a contradiction.

Hence $$y_v \in (N_{(G,\phi)}^{\leq (3\ell+\tau)/2+2\ell'}[X_{t_-}] \cup N_{(G,\phi)}^{\leq (3\ell+\tau)/2+\rho+2\ell'}[R_-]) - (N_{(G_0,\phi_0)}^{\leq \ell_0'}[X'_{t'_-}] \cup N_{(G_0,\phi_0)}^{\leq \rho+\ell_0'}[R'_-]).$$

Suppose $y_v \in N_{(G,\phi)}^{\leq (3\ell+\tau)/2+2\ell'}[X_{t_-}]$.
If $t_- \in V(T_0)$, then $X_{t_-} \subseteq N_{(G,\phi)}^{\leq 3\ell+\mu}[A'_{t_-'}] \subseteq N_{(G,\phi)}^{\leq 3\ell+\mu}[R'_-]$, so 
	\begin{align*}
		y_v \in & N_{(G,\phi)}^{\leq (3\ell+\tau)/2+2\ell'}[X_{t_-}] \cap V(G_0) \\
		\subseteq & N_{(G,\phi)}^{\leq (3\ell+\tau)/2+2\ell'+3\ell+\mu}[R'_-] \cap V(G_0) \\
		\subseteq & N_{(G_0,\phi_0)}^{\leq (3\ell+\tau)/2+2\ell'+3\ell+\mu}[R'_-] \subseteq N_{(G_0,\phi_0)}^{\leq \rho+\ell_0'}[R'_-],
	\end{align*}
where the second last inclusion follows from Claim 6, a contradiction.
So $t_- \not \in V(T_0)$ and $e_-$ is defined.
If $y_v \in X_{T_{e_-}}$, then $y_v \in N_{(G,\phi)}^{\leq \ell}[X_{e_-}]$ by the definition of a condensation, but 
	\begin{align*}
		y_v \in N_{(G,\phi)}^{\leq \ell}[X_{e_-}] \cap V(G_0) \subseteq & N_{(G,\phi)}^{\leq \ell+\mu}[A'_{t'_-}] \cap V(G_0) & \text{\ \ (by (C5) and (C6))} \\
		\subseteq & N_{(G,\phi)}^{\leq \ell+\mu}[R'_-] \cap V(G_0) & \text{\ \ (by the definition of $R'_-$)} \\
		\subseteq & N_{(G_0,\phi_0)}^{\leq \ell+\mu}[R'_-] & \text{\ \ (by Claim 6)} \\
		\subseteq & N_{(G_0,\phi_0)}^{\leq \ell_0'}[R'_-], &
	\end{align*}
a contradiction.
So $y_v \in V(G_0)-X_{T_{e_-}}$.
Since $y_v \in N_{(G,\phi)}^{\leq (3\ell+\tau)/2+2\ell'}[X_{t_-}]$, there exists a path $P_v$ in $G$ with $\leng_{(G,\phi)}(P_v) \leq (3\ell+\tau)/2+2\ell'$ from $y_v$ to $X_{t_-}$.
Since $y_v \in V(G_0)-X_{T_{e_-}}$, $P_v$ intersects $X_{e_-} \subseteq N_{(G,\phi)}^{\leq \mu}[R'_-]$.
So 
	\begin{align*}
		y_v \in N_{(G,\phi)}^{\leq \mu+\leng_{(G,\phi)}(P_v)}[R'_-] \cap V(G_0) \subseteq & N_{(G,\phi)}^{\leq \mu+(3\ell+\tau)/2+2\ell'}[R'_-] \cap V(G_0) \\
		\subseteq & N_{(G_0,\phi_0)}^{\leq \mu+(3\ell+\tau)/2+2\ell'}[R'_-] \subseteq N_{(G_0,\phi_0)}^{\leq \ell_0'}[R'_-]
	\end{align*}
by Claim 6, a contradiction.

Hence $y_v \in N_{(G,\phi)}^{\leq (3\ell+\tau)/2+\rho+2\ell'}[R_-]$.
So there exists a path $Q_v$ in $G$ between $y_v$ and $R_-$ internally disjoint from $R_-$ with $\leng_{(G,\phi)}(Q_v) \leq (3\ell+\tau)/2+\rho+2\ell'$.
Let $q$ be the vertex in $V(Q_v) \cap R_-$.
If $q \in V(G_0)$, then $q \in R'_-$, so $$y_v \in N_{(G,\phi)}^{\leq (3\ell+\tau)/2+\rho+2\ell'}[\{q\}] \cap V(G_0) \subseteq N_{(G_0,\phi_0)}^{\leq (3\ell+\tau)/2+\rho+2\ell'}[\{q\}] \subseteq N_{(G_0,\phi_0)}^{\leq (3\ell+\tau)/2+\rho+2\ell'}[R'_-] \subseteq N_{(G_0,\phi_0)}^{\leq \rho+\ell_0'}[R'_-]$$ by Claim 6, a contradiction.
So $q \not \in V(G_0)$.
Hence $q \in X_{T_{e_q}}-X_{e_q}$ for some $e_q \in U_E$.
Since $q \in R_- \cap X_{T_{e_q}}-X_{e_q}$, we know $D_{e_q} \cup S_{e_q} \cup O_{e_q} \subseteq R'_-$.
So $X_{e_q} \subseteq N_{(G,\phi)}^{\leq \mu}[R'_-]$ by (C5) and (C6).
If $y_v \in X_{T_{e_q}}$, then since $y_v \in V(G_0)$, $$y_v \in N_{(G,\phi)}^{\leq \ell}[X_{e_q}] \cap V(G_0) \subseteq N_{(G,\phi)}^{\leq \ell+\mu}[R'_-] \cap V(G_0) \subseteq N_{(G_0,\phi_0)}^{\leq \ell+\mu}[R'_-] \subseteq N_{(G_0,\phi_0)}^{\leq \ell_0'}[R'_-]$$ by Claim 6, a contradiction.
So $y_v \not \in X_{T_{e_q}}$.
Hence $V(Q_v) \cap X_{e_q} \neq \emptyset$.
So 
	\begin{align*}
		y_v \in N_{(G,\phi)}^{\leq \leng_{(G,\phi)}(Q_v)}[X_{e_q}] \cap V(G_0) \subseteq & N_{(G,\phi)}^{\leq (3\ell+\tau)/2+\rho+2\ell'}[X_{e_q}] \cap V(G_0) \\
		\subseteq & N_{(G,\phi)}^{\leq (3\ell+\tau)/2+\rho+2\ell'+\mu}[R'_-] \cap V(G_0) \\
		\subseteq & N_{(G_0,\phi_0)}^{\leq (3\ell+\tau)/2+\rho+2\ell'+\mu}[R'_-] \subseteq N_{(G_0,\phi_0)}^{\leq \rho+\ell_0'}[R'_-]
	\end{align*}
by Claim 6, a contradiction.
$\Box$

\medskip

For every $e \in U_E \cup (E(T')-E(T_0))$, let $R_e = \emptyset$.
For every $e \in E(T') \cup \{t^*\}$, 
	\begin{itemize}
		\item let $D'_e = D_e-R'$, 
		\item let $O'_e = (O_e \cup \{v \in D_e: N_{(G_0,\phi_0)}^{\leq \mu}[\{v\}] \cap R' \neq \emptyset\})-R_e$, and 
		\item let $S'_e = S_e \cup R_e$, 
	\end{itemize}
so $D_e \cup O_e \supseteq D_e' \cup O_e' \supseteq (D_e \cup O_e)-R_e$.

\medskip

\noindent{\bf Claim 9:} $(R',T',\X')$ is an $(\eta-1,\theta,\mu,\ell,\ell_0')$-construction of $(G_0,\phi_0)$.

\noindent{\bf Proof of Claim 9:}
Note that we may treat every edge $e \in U_E$ as an edge of $T'$.
So $X'_e = X_e$ for every $e \in E(T')$.

By (C1) and (C4) for $(R,T,\X)$, for every $e \in E(T') \cup \{t^*\}$, we know $D'_e \subseteq D_e-R' \subseteq X_e'-R'$, and $O'_e \subseteq O_e \cup D_e \subseteq X'_e$, and $S'_e = S_e \cup R_e \subseteq (X_e \cap R) \cup (X_e \cap R') = X'_e \cap R'$.
Moreover, for every $e \in E(T') \cup \{t^*\}$, since $R_e \subseteq D_e \cup O_e$, so $D'_e \cup O'_e \cup S'_e \subseteq (D_e-R') \cup (O_e \cup D_e) \cup (S_e \cup R_e) \subseteq D_e \cup O_e \cup S_e$, and hence $\lvert D'_e \cup O'_e \cup S'_e \rvert \leq \lvert D_e \cup O_e \cup S_e \rvert \leq \theta$ by (C1) and (C4) for $(R,T,\X)$.
So $(R',T',\X')$ satisfies (C1) and (C4).

By Claim 8, $(R',T',\X')$ satisfies (C10).
In addition, $X'_{e}=X_{e}$ for every $e \in E(T') \cup \{t^*\}$, so $(R',T',\X')$ satisfies (C11).

Suppose to the contrary that $(R',T',\X')$ is not an $(\eta-1,\theta,\mu,\ell,\ell_0')$-construction of $(G_0,\phi_0)$.
So there exists $e \in E(T') \cup \{t^*\}$ witnessing that $(R',T',\X')$ is not an $(\eta-1,\theta,\mu,\ell,\ell_0')$-construction of $(G_0,\phi_0)$ by violating (C2), (C3) or one of (C5)-(C9).

Suppose that (C2) or (C5) is violated.
Then there exists $u \in (X'_e-R')-N_{(G_0,\phi_0)}^{\leq \mu}[D'_e \cup O'_e]$.
Since $u \in X'_e-R' \subseteq X_e-R$, the conditions (C2) and (C5) for $(R,T,\X)$ implies that there exists $z_u \in D_e \cup O_e$ such that $u \in N_{(G,\phi)}^{\leq \mu}[\{z_u\}] \cap X_e \subseteq N_{(G_0,\phi_0)}^{\leq \mu}[\{z_u\}]$, where the last inclusion follows from Claim 6. 
Since $u \not \in R'$, we know $z_u \not \in R_e$.
Since $u \not \in N_{(G_0,\phi_0)}^{\leq \mu}[D'_e \cup O'_e]$, we have $z_u \not \in D'_e \cup O'_e$.
Hence $z_u \in (D_e \cup O_e)-(D_e' \cup O_e' \cup R_e) = \emptyset$, a contradiction.

Hence $X_e'-R' \subseteq N_{(G_0,\phi_0)}^{\leq \mu}[D'_e \cup O'_e]$ and (C2) and (C5) are satisfied for $e$.

Now we show that (C3) and (C6) are satisfied for $e$.
Let $v \in X'_e \cap R'$.
Since $v \in X'_e=X_e \subseteq N_{(G,\phi)}^{\leq \mu}[D_e \cup O_e \cup S_e]$ by (C2), (C3), (C5) and (C6), we know that there exists $x_v \in D_e \cup O_e \cup S_e$ such that $v \in N_{(G,\phi)}^{\leq \mu}[\{x_v\}]$.
By Claim 6, $v \in N_{(G_0,\phi_0)}^{\leq \mu}[\{x_v\}]$.
If $x_v \in D_e$, then since $v \in R'$, we have $x_v \in O_e' \cup R_e \subseteq O_e' \cup S_e'$, so $v \in N_{(G_0,\phi_0)}^{\leq \mu}[O_e' \cup S_e']$.
If $x_v \not \in D_e$, then $x_v \in O_e \cup S_e$, so $v \in N_{(G_0,\phi_0)}^{\leq \mu}[\{x_v\}] \subseteq N_{(G_0,\phi_0)}^{\leq \mu}[O_e \cup S_e] \subseteq N_{(G_0,\phi_0)}^{\leq \mu}[O_e' \cup S_e']$. 
Hence $X_e' \cap R' \subseteq N_{(G_0,\phi_0)}^{\leq \mu}[O_e' \cup S_e']$.
So (C3) and (C6) holds for $e$.

Therefore, $e \in E(T')$ and one of (C7)-(C9) is violated by $e$.
Recall that we can treat $e$ as an edge of $T$ and $X'_e=X_e$.

Suppose that (C9) is violated by $e$.
So $\lvert D'_e \cup O'_e \rvert > \eta-1$, and either both ends of $e$ have children in $T'$, or one end $t'$ of $e$ has no child and $X'_{t'} \not \subseteq N_{(G_0,\phi_0)}^{\leq \ell}[X'_e]$. 
If $e \in E(T')-E(T_0)$, then $e \in U_E$ and $t'=t_e$, so $$X'_{t'} = N_{(G,\phi)[X_{T_e}]}^{\leq \ell}[X_e] \cap V(G_0) \subseteq N_{(G,\phi)}^{\leq \ell}[X_e] \cap V(G_0) \subseteq N_{(G_0,\phi_0)}^{\leq \ell}[X_e] = N_{(G_0,\phi_0)}^{\leq \ell}[X'_e]$$ by Claim 6, a contradiction.
Hence $e \in E(T_0)$.
Note that $D_e \cup O_e \supseteq D'_e \cup O'_e$ and $R_e \subseteq D_e \cup O_e$ and $(D'_e \cup O'_e) \cap R_e=\emptyset$.
By Claim 7, $R_e \neq \emptyset$, so $\lvert D_e \cup O_e \rvert \geq |D_e' \cup O_e'|+1 > \eta$.
Since $(R,T,\X)$ satisfies (C9), some end $t'$ of $e$ has no child, and $X'_{t'}=X_{t'} \subseteq N_{(G,\phi)}^{\leq \ell}[X_e] = N_{(G,\phi)}^{\leq \ell}[X'_e] \cap V(G_0) \subseteq N_{(G_0,\phi_0)}^{\leq \ell}[X'_e]$ by Claim 6, a contradiction.

Hence (C7) or (C8) is violated by $e$.
Since $$(D'_{t^*} \cup O'_{t^*})-R' \subseteq (D_{t^*} \cup O_{t^*})-R' \subseteq R_{t^*}-R' =\emptyset,$$ we know $(X'_e \cap R' -N_{(G_0,\phi_0)}^{\leq \mu}[O_e']) \cap N_{(G_0,\phi_0)}^{\leq 3\ell+\mu}[(D'_{t^*} \cup O'_{t^*})-R'] = \emptyset$.

Hence (C8) is violated by $e$.
Therefore, there exists $$u \in (X'_e \cap R' -N_{(G_0,\phi_0)}^{\leq \mu}[O_e']) \cap N_{(G_0,\phi_0)}^{\leq a}[V(G_0)-(R' \cup D'_{t^*} \cup O'_{t^*})].$$
So there exists $x_u \in V(G_0)-(R' \cup D'_{t^*} \cup O'_{t^*})$ such that $\dist_{(G_0,\phi_0)}(x_u,u) \leq a$. 
Since $\{u,x_u\} \subseteq V(G_0) \subseteq V(G)$, by Statement 2(b) in Lemma \ref{con_qu_iso}, $\dist_{(G,\phi)}(u,x_u) \leq \dist_{(G_0,\phi_0)}(u,x_u) \leq a$.
Moreover, $$V(G_0)-(R' \cup D'_{t^*} \cup O'_{t^*}) \subseteq V(G)-(R \cup D_{t^*} \cup O_{t^*}).$$
So $x_u \in V(G_0)-(R' \cup D'_{t^*} \cup O'_{t^*}) \subseteq V(G)-(R \cup D_{t^*} \cup O_{t^*})$.
Hence $$u \in N_{(G,\phi)}^{\leq a}[V(G)-(R \cup D_{t^*} \cup O_{t^*})].$$

Since $u \in X_e$, there exists $y_u \in D_e \cup O_e \cup S_e$ such that $u \in N_{(G,\phi)}^{\leq \mu}[\{y_u\}] \cap X_e \subseteq N_{(G_0,\phi_0)}^{\leq \mu}[\{y_u\}]$ by (C5), (C6) and Claim 6.
Since $u \in X'_e \cap R' -N_{(G_0,\phi_0)}^{\leq \mu}[O_e']$, we know $y_u \not \in O_e'$.
If $y_u \in R_e$, then $e \in E(T_0)$ (for otherwise $R_e=\emptyset$), so $x_u \in N_{(G_0,\phi_0)}^{\leq a}[\{u\}] \subseteq N_{(G_0,\phi_0)}^{\leq a+\mu}[\{y_u\}] \subseteq N_{(G_0,\phi_0)}^{\leq a+\mu}[R_e] \subseteq R'$, a contradiction.
If $y_u \in D_e-R_e$, then $N_{(G_0,\phi_0)}^{\leq \mu}[\{y_u\}] \cap R' \supseteq N_{(G_0,\phi_0)}^{\leq \mu}[\{y_u\}] \cap \{u\} \neq \emptyset$, so $y_u \in O_e'$, a contradiction.
Hence $y_u \not \in R_e \cup D_e$.
So $y_u \in (D_e \cup O_e \cup S_e) - (O_e' \cup R_e \cup D_e) \subseteq S_e-(O_e \cup D_e)$.
Note that it is true for any vertex $y_u$ satisfying $y_u \in D_e \cup O_e \cup S_e$ and $u \in N_{(G,\phi)}^{\leq \mu}[\{y_u\}]$.
Hence $u \not \in N_{(G,\phi)}^{\leq \mu}[D_e \cup O_e]$.
So $u \in X_e \cap R$ by (C5).

If $u \in N_{(G_0,\phi_0)}^{\leq \mu}[R_e]$, then $e \in E(T_0)$ (for otherwise $R_e=\emptyset$), so $x_u \in N_{(G_0,\phi_0)}^{\leq a}[\{u\}] \subseteq N_{(G_0,\phi_0)}^{\leq a+\mu}[R_e] \subseteq R'$, a contradiction.
Hence $u \not \in N_{(G_0,\phi_0)}^{\leq \mu}[R_e]$.
Since $O_e \subseteq O_e' \cup R_e$, by Claim 6, $N_{(G,\phi)}^{\leq \mu}[O_e] \cap X_e \subseteq N_{(G_0,\phi_0)}^{\leq \mu}[O'_e \cup R_e] \cap X_e$, so $u \in X_e-N_{(G,\phi)}^{\leq \mu}[O_e]$.
Recall that we proved $u \in X_e \cap R$ in the previous paragraph.
So $u \in (X_e \cap R)-N_{(G,\phi)}^{\leq \mu}[O_e]$.
Since $(R,T,\X)$ satisfies (C8), $u \not \in N_{(G,\phi)}^{\leq a}[V(G)-(R \cup D_{t^*} \cup O_{t^*})]$.
Recall that we showed $u \in N_{(G,\phi)}^{\leq a}[V(G)-(R \cup D_{t^*} \cup O_{t^*})]$, a contradiction. 
This proves the claim.
$\Box$

\medskip

By the assumption of this lemma, for every $t \in V(T_0)$, there exists $A_t \subseteq X_t=X'_t$ with $\lvert A_t \rvert \leq \theta$ such that $X'_t = X_t \subseteq N_{(G,\phi)}^{\leq 3\ell+\mu}[A_t] \cap V(G_0) \subseteq N_{(G_0,\phi_0)}^{\leq 3\ell+\mu}[A_x]$ by Claim 6. 
For every $t \in V(T')-V(T_0)$, $t=t_e$ for some $e \in U_E$, so by Claim 9, $D'_e \cup O'_e \cup S'_e$ is a subset of $X'_e \subseteq X'_t$ with size at most $\theta$ such that $X'_t \subseteq N_{(G_0,\phi_0)}^{\leq \ell}[X'_e] \subseteq N_{(G_0,\phi_0)}^{\leq 3\ell+\mu}[D'_e \cup O_e' \cup S_e']$ by the definition of the condensation.
Recall that $(G_0,\phi_0)$ is a $(0,\ell_0']$-bounded weighted graph and is $(T',\X',\ell)$-bounded.
Hence by Claim 9 and the minimality of $\eta$ (by taking $Z=\emptyset$ and $c$ to be the function with the empty domain), $(G_0,\phi_0)^\ell[V(G_0)]-R'$ is $m$-colorable with weak diameter in $(G_0,\phi_0)^\ell$ at most $f^*(\eta-1)$.

Let $c_0$ be an $m$-coloring of $(G_0,\phi_0)^\ell[V(G_0)]-R'$ with weak diameter in $(G_0,\phi_0)^\ell$ at most $f^*(\eta-1)$.

Note that for every $e \in E(T_0)$, $R_e \subseteq N_{(G,\phi)}^{\leq \mu}[Z] \subseteq N_{(G,\phi)}^{\leq \mu+3\ell+\mu}[(D_{t^*} \cup O_{t^*})-R]$.
Since $R_{t^*} = (D_{t^*} \cup O_{t^*})-R$, we know $$\bigcup_{e \in E(T_0) \cup \{t^*\}}R_e \subseteq N_{(G,\phi)}^{\leq 3\ell+2\mu}[(D_{t^*} \cup O_{t^*})-R] \cap V(G_0) \subseteq N_{(G_0,\phi_0)}^{\leq 3\ell+2\mu}[(D_{t^*} \cup O_{t^*})-R]$$ by Claim 6.
Hence $$R'-R \subseteq R'-(R \cap V(G_0)) \subseteq N_{(G_0,\phi_0)}^{\leq a+\mu}[\bigcup_{e \in E(T_0) \cup \{t^*\}}R_e] \subseteq N_{(G_0,\phi_0)}^{\leq 3\ell+3\mu+a}[(D_{t^*} \cup O_{t^*})-R].$$

By Claims 3 and 6, $Z=N_{(G,\phi)}^{\leq 3\ell+\mu}[(D_{t^*} \cup O_{t^*})-R] \cap V(G_0)=N_{(G,\phi)}^{\leq 3\ell+\mu}[R_{t^*}] \cap V(G_0) \subseteq N_{(G_0,\phi_0)}^{\leq 3\ell+\mu}[R_{t^*}]$.
Since $a \geq 3\ell+\mu$, we know $Z \subseteq N_{(G_0,\phi_0)}^{\leq 3\ell+\mu}[R_{t^*}] \subseteq R'$ by the definition of $R'$.
So $Z-R \subseteq R'-R$ and $c|_{Z-R}$ can be extended to an $m$-coloring $c_{R'-R}$ of $R'-R$.
By Lemma \ref{patching_centered_set} (taking $(k,r,\ell,\nu,(G,\phi),S,Z,m,R,c_Z,c)=(\theta,3\ell+3\mu+a,\ell,f^*(\eta-1),(G_0,\phi_0),(D_{t^*} \cup O_{t^*})-R,R'-R,m,R \cap V(G_0),c_{R'-R},c_0)$), $c|_{Z-R}$ can be extended to the $m$-coloring $c_0 \cup c_{R'-R}$ of $(G_0,\phi_0)^\ell[V(G_0)]-R$ with weak diameter in $(G_0,\phi_0)^\ell$ at most $f_2(f^*(\eta-1))$.

By Lemma \ref{con_color} (taking $(\ell,\ell',\nu,\mu,\lambda,m,\theta)=(\ell,\ell',f_2(f^*(\eta-1)),\mu,\tau,m,\theta))$, $c_0 \cup c_{R'-R}$ (and hence $c|_{Z-R}$) can be extended to an $m$-coloring $c'$ of $(G,\phi)^\ell[X_{T_0} \cup \bigcup_{e \in U_E} N_{(G,\phi)[X_{T_e}]}^{\leq 3\ell}[X_e]]-R$ with weak diameter in $(G,\phi)^\ell$ at most $f_3(f_2(f^*(\eta-1))) = f^*(\eta)$ such that for every $m$-coloring $c''$ of $(G,\phi)^\ell[V(G)]-R$ that can be obtained by extending $c'$, every $c''$-monochromatic component in $(G,\phi)^\ell[V(G)]-R$ intersecting $X_{T_0} \cup \bigcup_{e \in U_E} N_{(G,\phi)[X_{T_e}]}^{\leq \ell}[X_e]$ has weak diameter in $(G,\phi)^\ell$ at most $f^*(\eta)$.

For every $e \in U_E$, let $Z_e =  N_{(G,\phi)[X_{T_e}]}^{\leq 3\ell}[X_e]$.

\medskip

\noindent{\bf Claim 10:} There exists $e \in U_E$ such that $c'|_{Z_e-R}$ cannot be extended to an $m$-coloring of $(G,\phi)^\ell[X_{T_e}]-R$ with weak diameter in $(G,\phi)^\ell$ at most $f^*(\eta)$.

\noindent{\bf Proof of Claim 10:}
Suppose to the contrary that for every $e \in U_E$. $c'|_{Z_e-R}$ can be extended to an $m$-coloring $c_e$ of $(G,\phi)^\ell[X_{T_e}]-R$ with weak diameter in $(G,\phi)^\ell$ at most $f^*(\eta)$.
Let $c^* = c' \cup \bigcup_{e \in U_E}c_e$.
Then $c^*$ is a well-defined $m$-coloring of $(G,\phi)^\ell[V(G)]-R$ that can be obtained by extending $c'$ (and hence by extending $c|_{Z-R}$).
So there exists a $c^*$-monochromatic component $M$ in $(G,\phi)^\ell[V(G)]-R$ such that $M$ has weak diameter in $(G,\phi)^\ell$ greater than $f^*(\eta)$.
By the property of $c'$, we know $V(M) \cap (X_{T_0} \cup \bigcup_{e \in U_E} N_{(G,\phi)[X_{T_e}]}^{\leq \ell}[X_e]) = \emptyset$.
Since $(G,\phi)$ is $(T,\X,\ell)$-bounded and $M \subseteq (G,\phi)^\ell[V(G)]$, there exists $e \in U_E$ such that $M$ is a $c_e$-monochromatic component in $(G,\phi)^\ell[X_{T_e}]-R$.
Then the weak diameter in $(G,\phi)^\ell$ of $M$ is at most $f^*(\eta)$, a contradiction.
$\Box$

\medskip

By Claim 10, there exists $e \in U_E$ such that $c'|_{Z_e-R}$ cannot be extended to an $m$-coloring of $(G,\phi)^\ell[X_{T_e}]-R$ with weak diameter in $(G,\phi)^\ell$ at most $f^*(\eta)$.
In particular, $X_{T_e}-Z_e \neq \emptyset$, so $\lvert D_e \cup O_e \rvert \leq \eta$ by (C9).

Our next goal is to construct an $(\eta,\theta,\mu,\ell,\ell')$-construction $(R^*,T^*,\X^*)$ of $(G,\phi)$.
We will first provide its definition and then prove that it indeed gives an $(\eta,\theta,\mu,\ell,\ell')$-construction in Claim 13.
This construction $(R^*,T^*,\X^*)$ will be used in Claim 14 to give us information about $Z$ and $R$.

\medskip

\noindent{\bf Claim 11:} There exists $W \subseteq D_e \cup O_e \cup S_e$ such that 
	\begin{itemize}
		\item $Z_e-R \subseteq N_{(G,\phi)}^{\leq 3\ell+\mu}[W]$, and
		\item for every $e' \in E(T)$, $(X_{e'} \cap R-N_{(G,\phi)}^{\leq \mu}[O_{e'}]) \cap N_{(G,\phi)}^{\leq 3\ell+\mu}[W] = \emptyset$.
	\end{itemize}

\noindent{\bf Proof of Claim 11:}
Let $W = \{v \in D_e \cup O_e \cup S_e: (Z_e-R) \cap N_{(G,\phi)}^{3\ell+\mu}[\{v\}] \neq \emptyset\}$.
Since $Z_e-R \subseteq N_{(G,\phi)[X_{T_e}]}^{\leq 3\ell}[X_e]-R \subseteq N_{(G,\phi)}^{\leq 3\ell+\mu}[D_e \cup O_e \cup S_e]-R$ by (C5) and (C6), $Z_e-R \subseteq N_{(G,\phi)}^{\leq 3\ell+\mu}[W]$.

Suppose to the contrary that there exists $e' \in E(T)$ such that $(X_{e'} \cap R-N_{(G,\phi)}^{\leq \mu}[O_{e'}]) \cap N_{(G,\phi)}^{\leq 3\ell+\mu}[W] \neq \emptyset$.
So there exist $x \in X_{e'} \cap R-N_{(G,\phi)}^{\leq \mu}[O_{e'}]$ and $y \in W$ such that $\dist_{(G,\phi)}(x,y) \leq 3\ell+\mu$.

Since $y \in W$, there exists $z \in Z_e-R$ such that $\dist_{(G,\phi)}(y,z) \leq 3\ell+\mu$.
So $x \in (X_{e'} \cap R-N_{(G,\phi)}^{\leq \mu}[O_{e'}]) \cap N_{(G,\phi)}^{\leq 6\ell+2\mu}[\{z\}]$.
Since $a \geq 6\ell+2\mu$, we know $z \not \in V(G)-(R \cup D_{t^*} \cup O_{t^*})$ by (C8).
That is, $z \in R \cup D_{t^*} \cup O_{t^*}$.
Since $z \in Z_e-R$, we have $z \in ((D_{t^*} \cup O_{t^*})-R) \cap X_{T_e} \subseteq Z \cap X_{T_e}$ by Claim 3.
But $Z \cap X_{T_e} = \emptyset$ by the definition of $U_E$, a contradiction.
$\Box$

\medskip

Let $T^*$ be the rooted tree obtained from $T$ by subdividing $e$ once to create a new node $q^*$ and making $q^*$ be the root of $T^*$.
Let $X^*_{q^*}=X_e$.
For every $t \in V(T)$, let $X^*_t = X_t$.
Let $\X^*=(X^*_t: t \in V(T^*))$.
Then $(T^*,\X^*)$ is a tree-decomposition of $G$ such that $(G,\phi)$ is $(T^*,\X^*,\ell)$-bounded.

For every edge $tt'$ of $T^*$, let $X^*_{tt'}=X^*_t \cap X^*_{t'}$.
Let $e_0$ be the edge of $T^*$ between $q^*$ and $V(T_e)$.
Let $e_1$ be the edge of $T^*$ between $q^*$ and $V(T)-V(T_e)$.

Let $W$ be the set given by Claim 11.
Let $$R^* = ((X_{T_e} \cap (R \cup D_{t^*} \cup O_{t^*})-W) \cup (V(G)-X_{T_e}).$$

\medskip

We define the following.
	\begin{itemize}
		\item Let $D^*_{q^*}=D_e-R^*$, $O^*_{q^*} = D_e \cup O_e \cup W$, and $S^*_{q^*}=S_e-W$.
		\item For every $z \in \{e_0,e_1\}$, let $D^*_z=D_e-R^*$, $O^*_z = D_e \cup O_e$, and $S^*_z=S_e-W$. 
		\item For every $e' \in E(T_e)$, let $D^*_{e'}=D_{e'}-R^*$, $O^*_{e'} = D_{e'} \cup O_{e'}$, and $S^*_{e'}=S_{e'}-W$. 
		\item For every $e' \in E(T)-E(T_e)$, let $D^*_{e'}=D_{e'}-R^*$, $O^*_{e'}=D_{e'} \cup O_{e'}$ and $S^*_{e'}=(D_{e'} \cup O_{e'} \cup S_{e'})-((W \cap R) \cup (X_e-R))$. 
	\end{itemize}

\medskip

\noindent{\bf Claim 12:} $(R^*,T^*,\X^*)$ satisfies (C9) and (C11) for being an $(\eta,\theta,\mu,\ell,\ell')$-construction of $(G,\phi)$.

\noindent{\bf Proof of Claim 12:}
Recall that $\lvert D_e \cup O_e \rvert \leq \eta$.
So for every $z \in \{e_0,e_1\}$, $\lvert D^*_z \cup O^*_z \rvert = \lvert D_e \cup O_e \rvert \leq \eta$.
For every $z \in E(T^*)-\{e_0,e_1\}$, $\lvert D^*_z \cup O^*_z \rvert \leq \lvert D_z \cup O_z \rvert$.
So for every $z \in E(T^*)$, if $\lvert D^*_z \cup O^*_z \rvert >\eta$, then $z \in E(T)$ and $\lvert D_z \cup O_z \rvert>\eta$.
Since (C9) holds for $(R,T,\X)$, (C9) also holds for $(R^*,T^*,\X^*)$.

By Claim 2, $\eta \geq 1$.
Since (C11) holds for $(R,T,\X)$, $X_{e'} \neq \emptyset$ for every $e' \in E(T)$. 
So $X^*_{q^*}=X_e=X_{e_0}=X_{e_1} \neq \emptyset$, and $X^*_{e'}=X_{e'} \neq \emptyset$ for every $e' \in E(T^*)-\{e_0,e_1\}=E(T)$.
Hence (C11) also holds for $(R^*,T^*,\X^*)$.
$\Box$

\medskip

\noindent{\bf Claim 13:} $(R^*,T^*,\X^*)$ is an $(\eta,\theta,\mu,\ell,\ell')$-construction of $(G,\phi)$.

\noindent{\bf Proof of Claim 13:}
We first show that (C1) holds for $(R^*,T^*,\X^*)$.
By (C4) for $(R,T,\X)$, we know that $D^*_{q^*} = D_e-R^* \subseteq (X_e-R)-R^* \subseteq X^*_{q^*}-R^*$, and $O^*_{q^*} \subseteq X_e = X^*_{q^*}$, and $S^*_{q^*} = S_e-W \subseteq R \cap X_e-W \subseteq X^*_{q^*} \cap R^*$, and $D^*_{q^*} \cup O^*_{q^*} \cup S^*_{q^*} \subseteq D_e \cup O_e \cup S_e$.
So (C1) holds for $(R^*,T^*,\X^*)$.

Now we show that (C4) holds for $(R^*,T^*,\X^*)$.
For every $z \in \{e_0,e_1\}$, (C4) for $(R,T,\X)$ implies that $D^*_z=D_e-R^* \subseteq (X_e-R)-R^* \subseteq X^*_z-R^*$ and $O^*_z \subseteq X_e = X^*_z$ and $S^*_z = S_e-W \subseteq R \cap X_e-W \subseteq X^*_z \cap R^*$ and $D^*_z \cup O^*_z \cup S^*_z \subseteq D_e \cup O_e \cup S_e$.
For every $e' \in E(T_e)$, we know $D^*_{e'}=D_{e'}-R^* \subseteq (X_{e'}-R)-R^* \subseteq X^*_{e'}-R^*$, and $O^*_{e'} \subseteq X_{e'} = X^*_{e'}$, and $S^*_{e'} = S_{e'}-W \subseteq R \cap X_{e'}-W \subseteq R \cap X_{T_e} \cap X^*_{e'}-W \subseteq X^*_{e'} \cap R^*$, and $D^*_{e'} \cup O^*_{e'} \cup S^*_{e'} \subseteq D_{e'} \cup O_{e'} \cup S_{e'}$.
For every $z \in E(T)-E(T_e)$, we know $D^*_z \subseteq X^*_z-R^*$, and $O^*_z=D_z \cup O_z \subseteq X^*_z$, and 
	\begin{align*}
		S^*_z = & (D_{z} \cup O_{z} \cup S_{z})-((W \cap R) \cup (X_e-R)) \\
		\subseteq & X^*_z-((W \cap R) \cup (X_e-R)) \\
		= & X^*_z \cap V(G) -((W \cap R) \cup (X_e-R)) \\
		\subseteq & X^*_z \cap \big((V(G)-X_{T_e}) \cup (X_z \cap X_{T_e}-((W \cap R) \cup (X_e-R)))\big) \\
		\subseteq & X^*_z \cap \big((V(G)-X_{T_e}) \cup (X_e-((W \cap R) \cup (X_e-R)))\big) \\
		\subseteq & X^*_z \cap \big((V(G)-X_{T_e}) \cup (R \cap X_{T_e}-W)\big) \subseteq X^*_z \cap R^*,
	\end{align*}
and $D^*_z \cup O^*_z \cup S^*_z \subseteq D_z \cup O_z \cup S_z$.
So (C4) holds for $(R^*,T^*,\X^*)$.

Now we show that (C2), (C3), (C5) and (C6) hold for $(R^*,T^*,\X^*)$.

By Claim 11, for every $e' \in E(T)$, $X_{e'} \cap R \cap N_{(G,\phi)}^{\leq 3\ell+\mu}[W] \subseteq N_{(G,\phi)}^{\leq \mu}[O_{e'}]$, so we have 
	\begin{itemize}
		\item[(i)] $X_{e'} \cap W \cap R \subseteq X_{e'} \cap R \cap N_{(G,\phi)}^{\leq 3\ell+\mu}[W] \subseteq N_{(G,\phi)}^{\leq \mu}[O_{e'}]$.
	\end{itemize}

Hence for every $z \in \{q^*,e_0,e_1\}$, 
	\begin{align*}
		X^*_z-R^* \subseteq & (X_e-R) \cup (X_e \cap W \cap R) \\
		\subseteq & N_{(G,\phi)}^{\leq \mu}[D_e \cup O_e] \cup N_{(G,\phi)}^{\leq \mu}[O_e] \subseteq N_{(G,\phi)}^{\leq \mu}[D_e \cup O_e] \subseteq N_{(G,\phi)}^{\leq \mu}[D^*_z \cup O^*_z]
	\end{align*}
by (C5) and (i); and by (C5), (C6) and (i), 
	\begin{align*}
		X^*_z \cap R^* \subseteq & (X_e \cap R) \cup (X_e \cap (D_{t^*} \cup O_{t^*})-R) \\
		\subseteq & (X_e \cap R \cap N_{(G,\phi)}^{\leq \mu}[O_e \cup S_e]) \cup N_{(G,\phi)}^{\leq \mu}[D_e \cup O_e] \\
		\subseteq & (N_{(G,\phi)}^{\leq \mu}[O_e \cup (S_e-W)] \cup (X_e \cap R \cap N_{(G,\phi)}^{\leq \mu}[W])) \cup N_{(G,\phi)}^{\leq \mu}[O^*_z] \\
		\subseteq & N_{(G,\phi)}^{\leq \mu}[O_e \cup (S_e-W)] \cup N_{(G,\phi)}^{\leq \mu}[O^*_z] \subseteq N_{(G,\phi)}^{\leq \mu}[O^*_z \cup S^*_z].
	\end{align*}
So (C2) and (C3) hold for $(R^*,T^*,\X^*)$.

Similarly, for every $z \in E(T_e)$, $X^*_{z}-R^* \subseteq (X_{z}-R) \cup (X_{z} \cap W \cap R) \subseteq N_{(G,\phi)}^{\leq \mu}[D_{z} \cup O_z] \subseteq N_{(G,\phi)}^{\leq \mu}[D^*_{z} \cup O^*_z]$ by (C5) and (i); and by (C5) and (C6) and (i), 
	\begin{align*}
		X^*_z \cap R^* \subseteq & (X_z \cap R) \cup (X_z \cap (D_{t^*} \cup O_{t^*})-R) \\
		= & (X_z \cap R \cap N_{(G,\phi)}^{\leq \mu}[O_z \cup S_z]) \cup N_{(G,\phi)}^{\leq \mu}[D_z \cup O_z] \\
		\subseteq & N_{(G,\phi)}^{\leq \mu}[O_z \cup (S_z-W)] \cup (X_z \cap R \cap N_{(G,\phi)}^{\leq \mu}[W]) \cup N_{(G,\phi)}^{\leq \mu}[O^*_z] \\
		\subseteq & N_{(G,\phi)}^{\leq \mu}[O_z \cup (S_z-W)] \cup N_{(G,\phi)}^{\leq \mu}[O^*_z] \subseteq N_{(G,\phi)}^{\leq \mu}[O^*_z \cup S^*_z].
	\end{align*}

For every $z \in E(T)-E(T_e)$, since $X_z \cap X_{T_e} \subseteq X_e$, we know $$X^*_z-R^* \subseteq (X_z-R) \cup (X_z \cap X_e \cap W \cap R) \subseteq N_{(G,\phi)}^{\leq \mu}[D_z \cup O_z] \cup N_{(G,\phi)}^{\leq \mu}[O_z] \subseteq N_{(G,\phi)}^{\leq \mu}[D^*_z \cup O^*_z]$$ by (C5) and (i); and by (C5), (C6) and (i), 
	\begin{align*}
		& X^*_z \cap R^* \\
		\subseteq & (X_z - R) \cup (X_z \cap R) \\
		\subseteq & N_{(G,\phi)}^{\leq \mu}[D_z \cup O_z] \cup (X_z \cap R \cap N_{(G,\phi)}^{\leq \mu}[O_z \cup S_z]) \\
		\subseteq & N_{(G,\phi)}^{\leq \mu}[O^*_z] \cup (X_z \cap R \cap N_{(G,\phi)}^{\leq \mu}[O^*_z \cup S^*_z \cup (W \cap R) \cup ((D_z \cup O_z \cup S_z) \cap X_e-R)]) \\
		\subseteq & N_{(G,\phi)}^{\leq \mu}[O^*_z \cup S^*_z] \cup (X_z \cap R \cap N_{(G,\phi)}^{\leq \mu}[W \cap R]) \cup (X_z \cap R \cap N_{(G,\phi)}^{\leq \mu}[(D_z \cup O_z \cup S_z) \cap X_e-R]) \\
		\subseteq & N_{(G,\phi)}^{\leq \mu}[O^*_z \cup S^*_z] \cup N_{(G,\phi)}^{\leq \mu}[O_z] \cup (N_{(G,\phi)}^{\leq \mu}[D_z \cup O_z] \cup N_{(G,\phi)}^{\leq \mu}[S_z-R])) \\
		= & N_{(G,\phi)}^{\leq \mu}[O^*_z \cup S^*_z] \cup N_{(G,\phi)}^{\leq \mu}[O_z] \cup N_{(G,\phi)}^{\leq \mu}[D_z \cup O_z] \subseteq N_{(G,\phi)}^{\leq \mu}[O^*_z \cup S^*_z].
	\end{align*}
Therefore, (C5) and (C6) hold for $(R^*,T^*,\X^*)$.

So it suffices to show that (C7)-(C11) hold for $(R^*,T^*,\X^*)$.

Now we show that (C7) and (C8) hold for $(R^*,T^*,\X^*)$.
Suppose to the contrary that there exist $z \in E(T^*)$ and $$u \in (X^*_z \cap R^*-N_{(G,\phi)}^{\leq \mu}[O^*_z]) \cap (N_{(G,\phi)}^{\leq 3\ell+\mu}[(D^*_{q^*} \cup O^*_{q^*})-R^*] \cup N_{(G,\phi)}^{\leq a}[V(G)-(R^* \cup D^*_{q^*} \cup O^*_{q^*})]).$$
For every $e' \in \{e_0,e_1\}$, we let $X_{e'}=X_e$, $D_{e'}=D_e$, $O_{e'}=O_e$ and $S_{e'}=S_e$.
So $u \in X_z \cap R^*-N_{(G,\phi)}^{\leq \mu}[O^*_z]$.

Since $O^*_z \supseteq D_z \cup O_z$, we know $u \in X_z-N_{(G,\phi)}^{\leq \mu}[O^*_z] \subseteq X_z-N_{(G,\phi)}^{\leq \mu}[D_z \cup O_z] \subseteq X_z-(X_z-R) \subseteq R$ by (C5).
Hence $u \in R$.
That is, $$u \in X_z \cap R-N_{(G,\phi)}^{\leq \mu}[O^*_z] \subseteq X_z \cap R-N_{(G,\phi)}^{\leq \mu}[O_z].$$

Suppose that $u \in N_{(G,\phi)}^{\leq 3\ell+\mu}[(D^*_{q^*} \cup O^*_{q^*})-R^*]$.
Hence there exists $v \in (D^*_{q^*} \cup O^*_{q^*})-R^*$ such that $\dist_{(G,\phi)}(u,v) \leq 3\ell+\mu$.
So $v \in (D^*_{q^*} \cup O^*_{q^*})-R^* = (D_e \cup O_e \cup W)-R^* \subseteq W \cup ((D_e \cup O_e)-R)$.
If $v \in W$, then $u \in (X_z \cap R-N_{(G,\phi)}^{\leq \mu}[O_z]) \cap N_{(G,\phi)}^{\leq 3\ell+\mu}[W] =\emptyset$ by Claim 11, a contradiction.
So $v \not \in W$.
Hence $v \in (D_e \cup O_e)-(R \cup W) \subseteq V(G)-R \subseteq ((D_{t^*} \cup O_{t^*})-R) \cup (V(G)-(R \cup D_{t^*} \cup O_{t^*}))$.
Then $u \in (X_z \cap R-N_{(G,\phi)}^{\leq \mu}[O_z]) \cap N_{(G,\phi)}^{\leq 3\ell+\mu}[((D_{t^*} \cup O_{t^*})-R) \cup (V(G)-(R \cup D_{t^*} \cup O_{t^*}))] =\emptyset$ by (C7) and (C8) since $a \geq 3\ell+\mu$, a contradiction.

Hence $u \in N_{(G,\phi)}^{\leq a}[V(G)-(R^* \cup D^*_{q^*} \cup O^*_{q^*})]$.
So there exists $y \in V(G)-(R^* \cup D^*_{q^*} \cup O^*_{q^*})$ such that $\dist_{(G,\phi)}(u,y) \leq a$.
Note that $y \in V(G)-(R^* \cup D^*_{q^*} \cup O^*_{q^*}) \subseteq V(G)-(R \cup D_e \cup O_e \cup W \cup D_{t^*} \cup O_{t^*}) \subseteq V(G)-(R \cup D_{t^*} \cup O_{t^*})$.
Then $u \in (X_z \cap R-N_{(G,\phi)}^{\leq \mu}[O_z]) \cap N_{(G,\phi)}^{\leq a}[V(G)-(R \cup D_{t^*} \cup O_{t^*})] =\emptyset$ by (C8), a contradiction.

Therefore, (C7) and (C8) hold for $(R^*,T^*,\X^*)$.
By Claim 12, it suffices to show that (C10) holds for $(R^*,T^*,\X^*)$.

Let $R^*_-=R_-$.
If $t_- \in V(T_e)$, then let $t^*_-=t_-$; otherwise, $t_- \in V(T)-V(T_e)$ and we let $t^*_-=q^*$.
So $|R^*_-|=|R_-| \leq (3\theta+3)^{\theta-\eta}$ (by (C10)) and $t^*_- \in V(T^*)$.
Let $rv \in E(G)$ with $r \in R^*$ and $v \in V(G)-R^*$.
So $v \in X_{T_e}$.

It suffices to show $v \in N_{(G,\phi)}^{\leq \ell}[X^*_{q^*}] \cup N_{(G,\phi)}^{\leq 2\ell'}[X^*_{t^*_-}] \cup N_{(G,\phi)}^{\leq \rho+2\ell'}[R^*_-]$, and if $N_{(G,\phi)}^{\leq 3\ell+\mu}[D^*_{q^*} \cup O^*_{q^*}-R^*]=\emptyset$, then $v \in N_{(G,\phi)}^{\leq 2\ell'}[X^*_{t^*_-}] \cup N_{(G,\phi)}^{\leq \rho+2\ell'}[R^*_-]$.

We first assume $v \in X_{T_e}-N_{(G,\phi)}^{\leq \ell}[X_e]$.
Since $(G,\phi)$ is $(T,\X,\ell)$-bounded and $rv \in E(G)$, $r \in X_{T_e}-X_e$, so $r \in X_{T_e}-(D_{t^*} \cup O_{t^*})$.
This implies $r \in R$ since $r \in R^*$.
Since $W \subseteq X_e$ and $v \in X_{T_e}-N_{(G,\phi)}^{\leq \ell}[X_e]$, we know $v \not \in W$.
So $v \in X_{T_e}-(R^* \cup W)$, and hence $v \not \in R$.
Hence by (C10) for $(R,T,\X)$, $v \in N_{(G,\phi)}^{\leq \ell}[X_{t^*}] \cup N_{(G,\phi)}^{\leq 2\ell'}[X_{t_-}] \cup N_{(G,\phi)}^{\leq \rho+2\ell'}[R_-]$.
Since $v \in X_{T_e}-N_{(G,\phi)}^{\leq \ell}[X_e]$, $v \not \in N_{(G,\phi)}^{\leq \ell}[X_{t^*}]$.
Since $R^*_-=R_-$, if $v \in N_{(G,\phi)}^{\leq \rho+2\ell'}[R_-]$, then $v \in N_{(G,\phi)}^{\leq \rho+2\ell'}[R^*_-]$ and we are done.
So we may assume $v \in N_{(G,\phi)}^{\leq 2\ell'}[X_{t_-}]$.
If $t^*_-=t_-$, then $v \in N_{(G,\phi)}^{\leq 2\ell'}[X_{t_-}] = N_{(G,\phi)}^{\leq 2\ell'}[X^*_{t_-}]$ and we are done.
So we may assume $t_- \in V(T)-V(T_e)$ and $t^*_-=q^*$.
Since $v \in X_{T_e} \cap N_{(G,\phi)}^{\leq 2\ell'}[X_{t_-}]$ and $t_- \in V(T)-V(T_e)$, we have $v \in N_{(G,\phi)}^{\leq 2\ell'}[X_e] = N_{(G,\phi)}^{\leq 2\ell'}[X^*_{q^*}] = N_{(G,\phi)}^{\leq 2\ell'}[X^*_{t^*_-}]$, so we are done. 

Hence we may assume $v \not \in X_{T_e}-N_{(G,\phi)}^{\leq \ell}[X_e]$.
That is, $v \in X_{T_e} \cap N_{(G,\phi)}^{\leq \ell}[X_e]$.
In particular, $v \in X_{T_e} \cap N_{(G,\phi)}^{\leq \ell}[X_e] \subseteq N_{(G,\phi)}^{\leq \ell}[X^*_{q^*}]$.
So we may assume $N_{(G,\phi)}^{\leq 3\ell+\mu}[D^*_{q^*} \cup O^*_{q^*}-R^*]=\emptyset$, for otherwise we are done.
In particular, $D^*_{q^*} \cup O^*_{q^*}-R^* = \emptyset$.
If $W \neq \emptyset$, then $D^*_{q^*} \cup O^*_{q^*}-R^* \supseteq W-R^* = X_e \cap W-R^* = W \neq \emptyset$ (since $W \subseteq X_e$ and $R^* \cap X_e \cap W=\emptyset$), a contradiction.
So $W=\emptyset$.
Hence $R^* \cap X_{T_e} = X_{T_e} \cap (R \cup D_{t^*} \cup O_{t^*})$.
By Claim 11, $Z_e \subseteq R$, so $X_{T_e} \cap (D_{t^*} \cup O_{t^*}) \subseteq X_e \subseteq Z_e \subseteq R$.
So $R^* \cap X_{T_e} = X_{T_e} \cap R$.
But $v \in X_{T_e} \cap N_{(G,\phi)}^{\leq \ell}[X_e] \subseteq Z_e \subseteq X_{T_e} \cap R \subseteq R^*$, a contradiction.
This proves the claim.
$\Box$

\medskip

\noindent{\bf Claim 14:} $Z_e \subseteq R$ and $R^* \subseteq Z \cup R$.

\noindent{\bf Proof of Claim 14:}
By Claim 13, $(R^*,T^*,\X^*)$ is an $(\eta,\theta,\mu,\ell,\ell')$-construction of $(G,\phi)$.
Note that $W \subseteq D^*_{q^*} \cup O^*_{q^*}$ by definition; since $W \subseteq X_e$, we know $W \cap R^* = \emptyset$.
So by Claim 11, $Z_e-R^* \subseteq (Z_e-R) \cup W \subseteq N_{(G,\phi)}^{\leq 3\ell+\mu}[W] \subseteq N_{(G,\phi)}^{\leq 3\ell+\mu}[D^*_{q^*} \cup O^*_{q^*}-R^*]$.

Suppose that $\lvert V(G)-(Z_e \cup R^*) \rvert + \lvert V(G)-(R^* \cup D^*_{q^*} \cup O^*_{q^*}) \rvert < \lvert V(G)-(Z \cup R) \rvert + \lvert V(G)-(R \cup D_{t^*} \cup O_{t^*}) \rvert$.
Then $(\eta,\lvert V(G)-(Z_e \cup R^*) \rvert + \lvert V(G)-(R^* \cup D^*_{q^*} \cup O^*_{q^*}) \rvert)$ is lexicographically smaller than $(\eta,\lvert V(G)-(Z \cup R) \rvert + \lvert V(G)-(R \cup D_{t^*} \cup O_{t^*}) \rvert)$.
It implies that $c'|_{Z_e-R^*}$ can be extended to an $m$-coloring $c''$ of $(G,\phi)^\ell[V(G)]-R^*$ with weak diameter in $(G,\phi)^\ell$ at most $f^*(\eta)$.
Let $c'''=c'' \cup c'|_{X_{T_e} \cap (D_{t^*} \cup O_{t^*})-R}$.
Since $X_{T_e}-R \subseteq (X_{T_e}-R^*) \cup (X_{T_e} \cap (D_{t^*} \cup O_{t^*}-R))$, we know that $c'''|_{X_{T_e}-R}$ is an $m$-coloring of $(G,\phi)^\ell[X_{T_e}]-R$ that can be obtained by extending $c'|_{Z_e-R}$ such that every $c'''$-monochromatic component in $(G,\phi)^\ell[X_{T_e}]-R$ with weak diameter in $(G,\phi)^\ell$ greater than $f^*(\eta)$ intersects $X_{T_e} \cap (D_{t^*} \cup O_{t^*})-R$.
By Claim 10, $c'|_{Z_e-R}$ cannot be extended to an $m$-coloring of $(G,\phi)^\ell[X_{T_e}]-R$ with weak diameter in $(G,\phi)^\ell$ at most $f^*(\eta)$, so there exists a $c'''$-monochromatic component $M$ in $(G,\phi)^\ell[X_{T_e}]-R$ with weak diameter in $(G,\phi)^\ell$ greater than $f^*(\eta)$.
So $M$ intersects $X_{T_e} \cap (D_{t^*} \cup O_{t^*})-R$.
Since $X_{T_e} \cap (D_{t^*} \cup O_{t^*}) \subseteq X_e \subseteq X_{T_0}$ and $c''' \cup c'$ can be obtained by extending $c'$, $M$ has weak diameter in $(G,\phi)^\ell$ at most $f^*(\eta)$ by the property of $c'$, a contradiction.

So $\lvert V(G)-(Z_e \cup R^*) \rvert + \lvert V(G)-(R^* \cup D^*_{q^*} \cup O^*_{q^*}) \rvert \geq \lvert V(G)-(Z \cup R) \rvert + \lvert V(G)-(R \cup D_{t^*} \cup O_{t^*}) \rvert$.

Since $Z_e \cup R^* \supseteq Z \cup R$, we know $\lvert V(G)-(Z_e \cup R^*) \rvert \leq \lvert V(G)-(Z \cup R) \rvert$.
Since $R^* \cup D^*_{q^*} \cup O^*_{q^*} \supseteq R \cup D_{t^*} \cup O_{t^*}$, we know $\lvert V(G)-(R^* \cup D_{q^*} \cup O_{q^*}) \rvert \leq \lvert V(G)-(R \cup D_{t^*} \cup O_{t^*}) \rvert$.
So $Z_e \cup R^* = Z \cup R$. 
It implies $R^* \subseteq Z \cup R$.
Moreover, it implies $Z_e \subseteq R$ since $Z_e \subseteq X_{T_e}$ and $Z \cap X_{T_e} = \emptyset$.
$\Box$.

\medskip

Finally, we shall use the information given by Claim 14 to finish the proof.

Let $R^1=Z \cup R$.
Note that for every $z \in \{t^*\} \cup (E(T)-E(T_e))$, we have $X_{z} \subseteq (V(G)-X_{T_e}) \cup X_e \subseteq R^* \cup Z_e \subseteq Z \cup R = R^1$ by Claim 14.
We define the following:
	\begin{itemize}
		\item Let $D^1_{t^*}=\emptyset$, $O^1_{t^*}=D_{t^*} \cup O_{t^*}$ and $S^1_{t^*}=S_{t^*}$.
		\item For every $e' \in E(T)-E(T_e)$, let $D^1_{e'}=\emptyset$, $O^1_{e'}=D_{e'} \cup O_{e'}$, and $S^1_{e'}=S_{e'}$.
		\item For every $e' \in E(T_e)$, let $D^1_{e'}=D_{e'}$, $O^1_{e'}=O_{e'}$, and $S^1_{e'}=S_{e'}$.
	\end{itemize}

\medskip

\noindent{\bf Claim 15:} $(R^1,T,\X)$ is an $(\eta,\theta,\mu,\ell,\ell')$-construction of $(G,\phi)$.

\noindent{\bf Proof of Claim 15:}
Clearly, (C1) and (C4) hold for $(R^1,T,\X)$.

Now we show that (C2), (C3), (C5) and (C6) hold for $(R^1,T,\X)$.
For every $z \in \{t^*\} \cup (E(T)-E(T_e))$, recall that we know $X_{z} \subseteq R^1$, so we have $X_{z}-R^1 = \emptyset$ and $X_{z} \cap R^1 = X_{z} \subseteq N_{(G,\phi)}^{\leq \mu}[D_{z} \cup O_{z} \cup S_{z}] = N_{(G,\phi)}^{\leq \mu}[O^1_{z} \cup S^1_z]$ by (C2), (C3), (C5) and (C6).
For every $z \in E(T_e)$, we have $X_z-R \subseteq N_{(G,\phi)}^{\leq \mu}[D_z \cup O_z] = N_{(G,\phi)}^{\leq \mu}[D^1_z \cup O^1_z]$ by (C5), and we have $X_z \cap R \subseteq N_{(G,\phi)}^{\leq \mu}[S_z \cup O_z] = N_{(G,\phi)}^{\leq \mu}[S^1_z \cup O^1_z]$ by (C6).
So (C2), (C3), (C5) and (C6) hold for $(R^1,T,\X)$.

Now we show that (C7) and (C8) hold for $(R^1,T,\X)$.
Note that $$D^1_{t^*} \cup O^1_{t^*} \subseteq X_{t^*} \subseteq (V(G)-X_{T_e}) \cup X_e \subseteq R^* \cup Z_e \subseteq Z \cup R = R^1$$ by Claim 14.
It implies that $(D^1_{t^*} \cup O^1_{t^*})-R^1 = \emptyset$, so (C7) holds for $(R^1,T,\X)$.
It also implies that $R^1 \supseteq R \cup X_{t^*} \supseteq R \cup D_{t^*} \cup O_{t^*}$.
So $V(G)-(R^1 \cup D^1_{t^*} \cup O^1_{t^*}) \subseteq V(G)-(R \cup D_{t^*} \cup O_{t^*})$.
For every $z \in E(T_e)$, we know $X_z \cap R^1 = X_z \cap (R \cup Z) = X_z \cap R$ (since $Z \subseteq X_{T_0}-X_e$), so 
	\begin{align*}
		& (X_z \cap R^1 - N_{(G,\phi)}^{\leq \mu}[O^1_z]) \cap N_{(G,\phi)}^{\leq a}[V(G)-(R^1 \cup D^1_{t^*} \cup O^1_{t^*})] \\
		\subseteq & (X_z \cap R - N_{(G,\phi)}^{\leq \mu}[O_z]) \cap N_{(G,\phi)}^{\leq a}[V(G)-(R \cup D_{t^*} \cup O_{t^*})] = \emptyset
	\end{align*}
by (C8) for $(R,T,\X)$. 
For every $z \in E(T)-E(T_e)$, recall that we know $X_z \subseteq R^1$, so $X_z \cap R^1 = X_z$, and hence
	\begin{align*}
		& (X_z \cap R^1 - N_{(G,\phi)}^{\leq \mu}[O^1_z]) \cap N_{(G,\phi)}^{\leq a}[V(G)-(R^1 \cup D^1_{t^*} \cup O^1_{t^*})] \\
		\subseteq & (X_z - N_{(G,\phi)}^{\leq \mu}[D_z \cup O_z]) \cap N_{(G,\phi)}^{\leq a}[V(G)-(R \cup D_{t^*} \cup O_{t^*})] \\
		\subseteq & (X_z \cap R - N_{(G,\phi)}^{\leq \mu}[D_z \cup O_z]) \cap N_{(G,\phi)}^{\leq a}[V(G)-(R \cup D_{t^*} \cup O_{t^*})] = \emptyset
	\end{align*}
by (C5) and (C8) for $(R,T,\X)$. 
Hence (C8) holds for $(R^1,T,\X)$.

Now we show that (C9) and (C11) hold for $(R^1,T,\X)$.
Since $D^1_{e'} \cup O^1_{e'} = D_{e'} \cup O_{e'}$ for every $e \in E(T)$, (C9) holds for $(R^1,T,\X)$ by the condition (C9) for $(R,T,\X)$.
In addition, (C11) holds for $(R^1,T,\X)$ by the condition (C11) for $(R,T,\X)$.

So it suffices to show that (C10) hold for $(R^1,T,\X)$.
Let $rv \in E(G)$ with $r \in R^1$ and $v \in V(G)-R^1$.
Since $R^1 \supseteq Z_e \cup R^*$ by Claim 14, $v \in X_{T_e}-N_{(G,\phi)}^{\leq 3\ell}[X_e]$.
Since $(G,\phi)$ is $(T,\X,\ell)$-bounded, $r \in X_{T_e}-X_e$.
Since $Z \subseteq X_{T_0}-X_e$, we know $r \not \in Z$, so $r \in R^1-Z \subseteq R$.
Since $v \not \in R^1 = R \cup Z$, we know $v \not \in R$.
By (C10) for $(R,T,\X)$, $v \in N_{(G,\phi)}^{\leq \ell}[X_{t^*}] \cup N_{(G,\phi)}^{\leq 2\ell'}[X_{t_-}] \cup N_{(G,\phi)}^{\leq \rho+2\ell'}[R_-]$.
Since $v \in X_{T_e}-N_{(G,\phi)}^{\leq 3\ell}[X_e]$, we know $v \not \in N_{(G,\phi)}^{\leq \ell}[X_{t^*}]$.
So $v \in N_{(G,\phi)}^{\leq 2\ell'}[X_{t_-}] \cup N_{(G,\phi)}^{\leq \rho+2\ell'}[R_-]$.
Hence (C10) holds for $(R^1,T,\X)$.
$\Box$

\medskip

By Claim 15, $(R^1,T,\X)$ is an $(\eta,\theta,\mu,\ell,\ell')$-construction of $(G,\phi)$ such that $(G,\phi)$ is $(T,\X,\ell)$-bounded.
Let $Z^1=\emptyset$.
Note that $Z^1 \cup R^1 = R \cup Z$ and $R^1 \cup D^1_{t^*} \cup O^1_{t^*} \supseteq R \cup D_{t^*} \cup O_{t^*}$.
By Claims 3 and 5, $D_{t^*} \cup O_{t^*}-R \neq \emptyset$, so $Z-R \supseteq D_{t^*} \cup O_{t^*}-R \neq \emptyset$.
Hence $|V(G)-R^1|<|V(G)-R|$.
So $(\eta,|V(G) -(Z^1 \cup R^1)|+|V(G)-(R^1 \cup D^1_{t^*} \cup O^1_{t^*})|, |V(G)-R^1|)$ is lexicographically smaller than $(\eta,|V(G) -(Z \cup R)|+|V(G)-(R \cup D_{t^*} \cup O_{t^*})|, |V(G)-R|)$.
Hence the coloring with empty domain can be extended to an $m$-coloring $c^1$ of $(G,\phi)^\ell[V(G)]-R^1 = (G,\phi)^\ell[V(G)]-(R \cup Z)$ with weak diameter in $(G,\phi)^\ell$ at most $f^*(\eta)$.
Since $Z \cap X_{T_e}=\emptyset$, we know that $c^1|_{X_{T_e}-R}$ is an $m$-coloring of $(G,\phi)^\ell[X_{T_e}]-R$ with weak diameter in $(G,\phi)^\ell$ at most $f^*(\eta)$.
By Claim 14, $Z_e-R=\emptyset$, so $c'|_{Z_e-R}$ is the coloring with empty domain, and hence Claim 10 implies that there exists no $m$-coloring of $(G,\phi)^\ell[X_{T_e}]-R$ with weak diameter in $(G,\phi)^\ell$ at most $f^*(\eta)$, a contradiction.
This proves the lemma.
\end{pf}

\begin{lemma} \label{strong_weighted_tree_extension_control_clean}
For any positive real numbers $\ell,\mu$ and every positive integer $\theta$, there exists a real number $\nu^*$ such that if $(G,\phi)$ is a $(0,\ell]$-bounded weighted graph with a tree-decomposition $(T,(X_t: t \in V(T)))$ of $G$ such that for every $t \in V(T)$, there exists $A_t \subseteq X_t$ with $\lvert A_t \rvert \leq \theta$ and $X_t \subseteq N_{(G,\phi)}^{\leq \mu}[A_t]$, then there exists a 2-coloring of $(G,\phi)^\ell$ with weak diameter in $(G,\phi)^\ell$ at most $\nu^*$.
\end{lemma}

\begin{pf}
Let $\ell$ and $\mu$ be positive real numbers.
Let $\theta$ be a positive integer.
Define $\nu^*=f^*(\theta)$, where $f^*$ is the function given by Lemma \ref{strong_weighted_tree_extension_control} by taking $(\ell,\ell',\mu,m,\theta)=(\ell,\ell,2\mu,2,\theta)$.

Let $(G,\phi)$ be a $(0,\ell]$-bounded weighted graph with a tree-decomposition $(T,\X=(X_t: t \in V(T)))$ of $G$ such that for every $t \in V(T)$, there exists $A_t \subseteq X_t$ with $\lvert A_t \rvert \leq \theta$ and $X_t \subseteq N_{(G,\phi)}^{\leq \mu}[A_t]$. 
To prove this lemma, it suffices to prove the case that $G$ is connected.
So we may assume that $G$ is connected and $X_t \neq \emptyset$ for every $t \in V(T)$.
For each $tt' \in E(T)$, let $X_{tt'}=X_t \cap X_{t'}$.
Since $G$ is connected and $X_t \neq \emptyset$ for every $t \in V(T)$, we know $X_e \neq \emptyset$ for every $e \in E(T)$.

For every $e \in E(T)$, there exists $t_e \in V(T)$ incident with $e$, so $X_e \subseteq X_{t_e} \subseteq N_{(G,\phi)}^{\leq \mu}[A_{t_e}]$; let $A'_{t_e} = \{v \in A_{t_e}: X_e \cap N_{(G,\phi)}^{\leq \mu}[\{v\}] \neq \emptyset\}$, so $X_e \subseteq N_{(G,\phi)}^{\leq \mu}[A'_{t_e}]$; hence for each $e \in E(T)$ and $v \in A'_{t_e}$, there exists $x_v \in X_e \cap N_{(G,\phi)}^{\leq \mu}[\{v\}] \neq \emptyset$; let $A_e = \{x_v: v \in A'_{t_e}\}$.
Hence for each $e \in E(T)$, $X_e \subseteq N_{(G,\phi)}^{\leq \mu}[A'_{t_e}] \subseteq N_{(G,\phi)}^{\leq 2\mu}[A_e]$ and $\lvert A_e \rvert \leq \lvert A'_{t_e} \rvert \leq \lvert A_{t_e} \rvert \leq \theta$.

Let $t^*$ be a node of $t$.
For every $z \in \{t^*\} \cup E(T)$, let $D_z=\emptyset$, $O_z=A_z$, and $S_z=\emptyset$.
Then $(\emptyset,T,\X)$ is a $(\theta,\theta,2\mu,\ell,\ell)$-construction of $(G,\phi)$ such that $(G,\phi)$ is $(T,\X,\ell)$-bounded.
By taking $Z=R=\emptyset$ and applying Lemma \ref{strong_weighted_tree_extension_control}, there exists a 2-coloring of $(G,\phi)^\ell[V(G)]=(G,\phi)^\ell$ with weak diameter in $(G,\phi)^\ell$ at most $\nu^*$.
\end{pf}

\section{Bounded tree-width and bounded layered tree-width classes} \label{sec:app_ad}

Now we are ready to prove Theorem \ref{tw_AN_intro} which immediately implies Theorem \ref{minor_planar_AN_intro} by the Grid Minor Theorem \cite{rs_V}.
The following theorem is equivalent to Theorem \ref{tw_AN_intro} by Lemma \ref{scaling_closed_equiv}. 

\begin{theorem} \label{tw_ad}
Let $w$ be a positive integer.
Let $\F$ be the class of weighted graphs whose underlying graphs have tree-width at most $w$. 
Then $\ad(\F)=1$.
\end{theorem}

\begin{pf}
Since every path has tree-width at most 1, and there exists no integer $\nu$ such that every path is 1-colorable with weak diameter at most $\nu$, $\ad(\F) \geq 1$.
So it suffices to show that $\ad(\F) \leq 1$.
We may assume $w \geq 2$. 

We shall prove $\ad(\F) \leq 1$ by applying Theorem \ref{weighted_ad_wd_equiv} to reduce the problem to a coloring problem, and applying Lemma \ref{weighted_tree_extension_clean}.
In order to apply Lemma \ref{weighted_tree_extension_clean}, we shall define some auxiliary graph classes.

For every positive real number $\ell$, let $\F_{\ell,1}$ be the class of $(0,\ell]$-bounded weighted graphs on at most $w+1$ vertices, and let $\F_{\ell,2}$ be the class of $(0,\ell]$-bounded weighted graphs that have a $(w+1)^2$-vertex-cover of size at most $w+1$.
Note that $\F_{\ell,1}$ and $\F_{\ell,2}$ are hereditary classes, and $\F_{\ell,1} \subseteq \F_{\ell,2}$.
For every positive real number $\ell$, by Lemma \ref{weighted_vc_color}, there exists a positive integer $\nu_\ell$ such that $\F_{\ell,2}$ is $(2,\ell,\nu_\ell)$-nice.

For every positive real number $\ell$, let $\F_\ell$ be the class of $(0,\ell]$-bounded weighted graphs $(G,\phi)$ in which $G$ has tree-width at most $w$.
Since for any positive real number $\ell$ and $(G,\phi) \in \F_\ell$, there exists a tree-decomposition $(T,\X)$ of $G$ of adhesion at most $w+1$, where $\X=(X_t: t \in V(T))$, such that for every $t \in V(T)$, $(G,\phi)[X_t] \in \F_{\ell,1}$, and $\F_{\ell,2}$ contains every $(0,\ell]$-bounded weighted graph that can be obtained from $(G,\phi)[X_t]$ by, for each neighbor $t'$ of $t$ in $T$, adding at most $(w+1)^2$ new vertices and new edges such that every new edge is between two new vertices or between a new vertex and $X_t \cap X_{t'}$.
Hence for every positive real number $\ell$, by Lemma \ref{weighted_tree_extension_clean}, there exists a positive integer $\nu^*_\ell$ such that $\F_\ell$ is $(2,\ell,\nu^*_\ell)$-nice.

Let $\F_0$ be the class of graphs of tree-width at most $w$.
Since $w \geq 2$, $\F_0$ is closed under taking subdivision and duplication.
Define $f: {\mathbb R}^+ \rightarrow {\mathbb R}^+$, to be the function such that for every $x \in {\mathbb R}^+$, $f(x)=\nu^*_x$.
Then for any $\ell \in {\mathbb R}^+$, $G \in \F_0$ and $\phi:E(G) \rightarrow (0,\ell]$, we have $(G,\phi) \in \F_\ell$, so $(G,\phi)^\ell$ is 2-colorable with weak diameter in $(G,\phi)^\ell$ at most $\nu^*_\ell = f(\ell)$.
By Theorem \ref{weighted_ad_wd_equiv}, $\ad(\F) \leq 2-1=1$.
\end{pf}

\bigskip

Let $m$ be a nonnegative integer.
Let $X,Y$ be metric spaces.
We denote the metrics of $X$ and $Y$ by $d_X$ and $d_Y$, respectively.
Let $f: X \rightarrow Y$ be a function.
For a real number $r$, a subset $S$ of $X$ is \defn{$(\infty,r)$-bounded with respect to $f$} if $d_Y(f(x_1),f(x_2))\leq r$ for any $x_1,x_2 \in S$.

A \defn{real projection} of a weighted graph $(G,\phi)$ is a function $f: V(G) \rightarrow {\mathbb R}^+$ such that for every $x,y \in V(G)$, $\lvert f(x)-f(y) \rvert \leq \dist_{(G,\phi)}(x,y)$.
Let $\L = (\L_i)_{i \in {\mathbb N}}$ be a sequence of class of weighted graphs such that $\L_i \subseteq \L_{i+1}$ for every $i \in {\mathbb N}$.
A class $\F$ of weighted graphs is \defn{$\L$-layerable} if there exists a function $f: {\mathbb R}^+ \rightarrow {\mathbb N}$ such that every $(G,\phi) \in \F$ has a real projection $g$ such that for every $\ell >0$ and every $(\infty,\ell)$-bounded set $A$ in $(G,\phi)$ with respect to $g$, $(G,\phi)[A] \in \L_{f(\ell)}$.

\begin{theorem}[{\cite[Theorem 4.3]{bbeglps_merged}}] \label{layerable_exten}
Let $n$ be a nonnegative integer.
Let $\L=(\L_i)_{i \in {\mathbb N}}$ be a sequence of class of weighted graphs such that $\ad(\L_i) \leq n$ for every $i \in {\mathbb N}$.
If $\F$ is a $\L$-layerable class of weighted graphs, then $\ad(\F) \leq n+1$.
\end{theorem}

Now we are ready to prove Theorem \ref{layered_tw_AN_intro}. 
The following is a restatement.

\begin{theorem} \label{layered_tw_ad_lower}
Let $\epsilon_0$ be a positive real number.
Let $w$ be a positive integer.
Let $\F$ be an $[\epsilon_0,\infty)$-bounded class of weighted graphs such that for every $(G,\phi) \in \F$, the layered tree-width of $G$ is at most $w$.
Then $\ad(\F) \leq 2$.
\end{theorem}

\begin{pf}
For every $i \in {\mathbb N}$, let $\L_i$ be the class of weighted graphs such that for every $(G,\phi) \in \L_i$, the tree-width of $G$ is at most $i$.
Let $f: {\mathbb R}^+ \rightarrow {\mathbb N}$ such that for every $x \in {\mathbb R}^+$, $f(x) = (2\lceil \frac{x}{\epsilon_0} \rceil +1)w$.

Let $(G,\phi) \in \F$.
Since $(G,\phi) \in \F$, there exist a tree-decomposition $(T,\X)$ of $G$ and a layering $(V_1,V_2,...)$ of $G$, where $\X =(X_t: t\in V(T))$, such that for every $t \in V(T)$ and $i \in {\mathbb N}$, $\lvert X_t \cap V_i \rvert \leq w$.
Note that for every $v \in V(G)$, there exists $i_v \in {\mathbb N}$ such that $v \in V_{i_v}$.
Let $g: V(G) \rightarrow {\mathbb R}^+$ such that for every $v \in V(G)$, $g(v)=i_v\epsilon_0$.
Note that for any $x,y \in V(G)$, every path in $G$ between $x$ and $y$ contains at least $\lvert i_x-i_y \rvert$ edges, so $\dist_{(G,\phi)}(x,y) \geq \epsilon_0 |i_x-i_y| = |g(x)-g(y)|$.
So $g$ is a real projection of $(G,\phi)$.

Let $\ell$ be a positive real number.
Let $A$ be an $(\infty,\ell)$-bounded set in $(G,\phi)$ with respect to $g$.
Let $a \in A$.
Since $A$ is $(\infty,\ell)$-bounded, $\lvert i_x-i_a \rvert = \frac{1}{\epsilon_0} \cdot \lvert g(x)-g(a) \rvert \leq \frac{\ell}{\epsilon_0}$ for every $x \in A$.
So $A \subseteq \bigcup_{j=i_a- \lceil \frac{\ell}{\epsilon_0} \rceil}^{i_a+\lceil \frac{\ell}{\epsilon_0} \rceil}V_j$.
Hence $G[A]$ is a graph of tree-width at most $(2\lceil \frac{\ell}{\epsilon_0} \rceil +1)w$.
So $(G,\phi)[A] \in \L_{(2\lceil \frac{\ell}{\epsilon_0} \rceil +1)w} =\L_{f(\ell)}$.

Hence $\F$ is $\L$-layerable.
Since $\ad(\L_i) \leq 1$ for every $i \in {\mathbb N}$ by Theorem \ref{tw_ad}, $\ad(\F) \leq 2$ by Theorem \ref{layerable_exten}.
\end{pf}

\section{Near embeddings and geodesic bags} \label{sec:near_embedding_geo}

A \defn{rooted spanning tree} of a graph $G$ is a rooted tree $T$ whose underlying graph is a subgraph of $G$ such that $V(T)=V(G)$.
Note that every rooted tree is a directed graph, so every subgraph of this rooted tree is also a directed graph.
But for convenience, for every subgraph $H$ of a rooted spanning tree $T$ of $G$, we treat $H$ a subgraph of $G$ by taking the underlying graph.

The following is implicitly proved in \cite[Theorem 12]{dmw}.

\begin{theorem}[{\cite[Theorem 12]{dmw}}] \label{planar_tree_path_bags}
Let $G$ be a connected planar graph.
Let $T_0$ be a rooted spanning tree.
Let $r$ be the root of $T_0$.
Then there exists a tree-decomposition $(T,\X=(X_t: t \in V(T))$ of $G$ such that for every $t \in V(T)$, $X_t$ is a union of the vertex-sets of at most 3 directed paths in $T_0$.
\end{theorem}

A rooted spanning tree $T$ of a graph $G$ is \defn{BFS} if for every positive integer $i$ and every vertex $v$ of $T$ whose distance in $G$ from the root of $T$ is $i$, the parent of $v$ in $T$ is a vertex whose distance in $G$ from the root of $T$ is $i-1$.

The following is implicitly proved in \cite[Lemma 21]{djmmuw}.

\begin{lemma}[{\cite[Lemma 21]{djmmuw}}] \label{genus_cutting}
Let $G$ be a connected graph of Euler genus $g$.
For every BFS rooted spanning tree $T$ of $G$ rooted at a vertex $r$, there exists $Z \subseteq V(G)$ such that the following hold.
	\begin{enumerate}
		\item $Z$ is a union of the vertex-sets of at most $2g$ directed paths in $T$ each containing $r$.
		\item $G-Z$ is planar.
		\item There exists a connected planar graph $G^+$ containing $G-Z$ as a subgraph, and there exists a BFS rooted spanning tree $T^+$ of $G^+$ such that for every directed path $P$ in $T^+$, the subgraph of $P$ induced by $V(P) \cap V(G)-Z$ is a directed path in $T$.
	\end{enumerate}
\end{lemma}

Let $g$ be a nonnegative integer.
Let $\Sigma$ be a surface of Euler genus at most $g$.
Let $p$ be a nonnegative integer.
A graph $G$ is \defn{$p$-nearly embeddable in $\Sigma$} if the following hold.
	\begin{itemize}
		\item There exist edge-disjoint graphs $G_0,G_1,...,G_p$ such that $G=\bigcup_{i=0}^pG_i$.
		\item For each $i \in [p]$, there exist a cyclic ordering $\Omega_i$ of $V(G_0) \cap V(G_i)$.
		\item There exist pairwise disjoint closed disks $\Delta_1,\Delta_2,...,\Delta_p$ in $\Sigma$ such that $G_0$ can be drawn in $\Sigma$ in such a way that the drawing only intersects $\bigcup_{i=1}^p\Delta_i$ at vertices of $G_0$, no vertex of $G_0$ is contained in the interior of $\bigcup_{i=1}^p\Delta_i$, and for each $i \in [p]$, the set of vertices of $G_0$ contained in the boundary of $\Delta_i$ is $V(G_0) \cap V(G_i)$, and the ordering of the vertices in $V(G_0) \cap V(G_i)$ appearing in the boundary of $\Delta_i$ equals $\Omega_i$.
		\item For each $i \in [p]$, if we denote $\Omega_i$ by $(v_1,v_2,...,v_{k_i})$ for some integer $k_i$, then there exists a path decomposition $(Q_i,\X_i)$ of $G_i$ such that $\lvert V(Q_i) \rvert =k_i$, and for each $j \in [k_i]$, the bag at the $j$-th vertex of $Q_i$ contains $v_j$ and has size at most $p$.
	\end{itemize}
We call $(G_0,G_1,...,G_p,\Delta_1,...,\Delta_p,((Q_i,\X_i): i \in [p]))$ a \defn{witness} of $G$.

For a tree-decomposition $(T,\X=(X_t: t \in V(T)))$ of a graph $G$ and for $t \in V(T)$, the \defn{torso} at $t$ is the graph obtained from $G[X_t]$ by, for each neighbor $t'$ of $t$ in $T$, adding edges such that $X_t \cap X_{t'}$ is a clique.

The following is a slightly stronger form of \cite[Theorem 1.3]{rs_XVI} which is actually implicitly proved in \cite{rs_XVI}.

\begin{theorem}[{\cite[Theorem 1.3]{rs_XVI}}] \label{rs_structure}
For every graph $H$, there exists a positive integer $p$ such that the following holds.
Let $G$ be an $H$-minor free graph.
Then $G$ admits a tree-decomposition $(T,\X=(X_t: t \in V(T)))$ of adhesion at most $p$ such that for every $t \in V(T)$, if $G_t$ is the torso at $t$, then there exists $Z_t \subseteq V(G_t)$ with $\lvert Z_t \rvert \leq p$ such that $G_t-Z_t$ is $p$-nearly embeddable in a surface $\Sigma$ in which $H$ cannot be drawn with a witness $(G_{t,0},G_{t,1},...,G_{t,p},\Delta_{t,1},...,\Delta_{t,p},((Q_{t,i},\X_{t,i}): i \in [p]))$ such that for every neighbor $t'$ of $t$ in $T$, $X_t \cap X_{t'}-Z_t$ is either a clique in $G_{t,0}$ of size at most three or contained in a bag of $(Q_{t,i},\X_{t,i})$ for some $i \in [p]$.
\end{theorem}

Let $(G,\phi)$ be a weighted graph.
Let $v \in V(G)$.
We say that a path $P$ in $G$ is a \defn{$v$-geodesic in $(G,\phi)$} if $v$ is an end of $P$ and the length in $(G,\phi)$ of $P$ equals the distance in $(G,\phi)$ between the ends of $P$.
Note that every subpath of a $v$-geodesic containing $v$ is also a $v$-geodesic.

The proof of the following lemma uses ideas in \cite[Lemma 24]{djmmuw}.

\begin{lemma} \label{geodesic_bags}
Let $p$ be a positive integer.
Let $g$ be a nonnegative integer.
Let $\Sigma$ be a surface of Euler genus at most $g$.
Let $\lambda$ be a positive real number. 
Let $(G,\phi)$ be a weighted graph with $G$ connected.
If there exists a tree-decomposition $(T,\X=(X_t: t \in V(T)))$ of $G$ such that 
	\begin{enumerate}
		\item $T$ is a star centered at its root $t^*$,
		\item the torso at $t^*$ is $p$-nearly embeddable in $\Sigma$ with a witness $(G_0,G_1,...,G_p,\Delta_1,...,\Delta_p,((Q_i,\X_i): i \in [p]))$, 
		\item for each neighbor $t'$ of $t^*$, $\lvert X_{t'} \rvert \leq p$ and $X_{t^*} \cap X_{t'}$ is either a clique in $G_{0}$ of size at most three or contained in a bag of $(Q_{i},\X_{i})$ for some $i \in [p]$, and
		\item for each neighbor $t'$ of $t^*$, $X_{t'} \subseteq N_{(G,\phi)}^{\leq \lambda}[X_{t^*} \cap X_{t'}]$,
	\end{enumerate}
then there exist a vertex $v^* \in V(G)$ and a tree-decomposition $(T^*,\X^*=(X^*_t: t \in V(T^*)))$ of $G$ such that for every $t \in V(T^*)$, there exist $i_t \in [(4g+8p+7)(2p+1)] \cup \{0\}$ and paths $P_{t,1},P_{t,2},...,P_{t,i_t}$ in $G$ such that 
	\begin{enumerate}
		\item for each $i \in [i_t]$, $P_{t,i}$ is a subpath of a $v^*$-geodesic in $(G,\phi)$ with $V(P_{t,i}) \subseteq X^*_t$, and
		\item $X^*_t \subseteq N_{(G,\phi)}^{\leq \lambda}[\bigcup_{i \in [i_t]}V(P_{t,i})]$.
	\end{enumerate}
\end{lemma}

\begin{pf}
Let $v^* \in V(G_0) \cap V(G_1)$.
Let $T_1$ be a spanning tree of $G$ such that for every $v \in V(G)$, the unique path in $T_1$ between $v^*$ and $v$ is a $v^*$-geodesic in $(G,\phi)$ between $v^*$ and $v$.
Note that such $T_1$ exists and can be constructed greedily.
We treat $T_1$ a rooted tree rooted at $v^*$.

Let $G'$ be the simple graph obtained from $G_0$ by  
	\begin{itemize}
		\item adding new edges such that for each $i \in [p]$, there exists a cycle whose vertex-set is $V(G_0) \cap V(G_i)$ and passing through them in an order consistent with $\Omega_i$, and
		\item for each neighbor $t$ of $t^*$ in $T$ in which $X_t \cap X_{t^*}$ is a clique in $G_0$ of size 3, if there exists a vertex $v \in X_t-X_{t^*}$ and three paths in $T_1$ from $v$ to $X_t \cap X_{t^*}$ intersecting only at $v$, then adding a new vertex $v_t$ and a new edge between $v_t$ and each vertex in $X_t \cap X_{t^*}$.
	\end{itemize}
Clearly, $G'$ can be drawn in $\Sigma$ with no edge-crossing.
Let $T_2$ be the subgraph of $G'$ obtained from the subgraph of $T_1$ induced by $V(G') \cap V(T_1)$ by, for each neighbor $t$ of $t^*$ in $T$ with $X_t \cap X_{t^*}$ a clique in $G_0$, doing the following operations: 
	\begin{itemize}
		\item If $\lvert X_t \cap X_{t^*} \rvert=2$ and there exists a path in $T_1[X_t]$ between the two vertices in $X_t \cap X_{t^*}$, then adding an edge between the two vertices $X_t \cap X_{t^*}$ into $T_2$.
		\item If $\lvert X_t \cap X_{t^*} \rvert=3$, then
			\begin{itemize}
				\item if there exists a vertex $v \in X_t-X_{t^*}$ and three paths in $T_1$ from $v$ to $X_t \cap X_{t^*}$ intersecting only at $v$, then adding $v_t$ and the three edges in $G'$ between $v_t$ and $X_t \cap X_{t^*}$ into $T_2$,
				\item otherwise, for each path in $T_1[X_t]$ between two distinct vertices $x,y$ in $X_t \cap X_{t^*}$ internally disjoint from $X_t \cap X_{t^*}$, adding $xy$ into $T_2$.
			\end{itemize}
	\end{itemize}
Note that $V(T_2)=V(G')$, $E(T_2) \subseteq E(G')$ and $T_2$ is a forest in $G$.
Since $v^* \in V(G_0)$, each edge in $E(T_2)-E(T_1)$ corresponds to a directed path in $T_1$, so we may assign a direction on each edge in $E(T_2)-E(T_1)$ consistent with the direction of the corresponding directed path in $T_1$ such that $T_2$ is a directed graph with no directed cycle.

Let $G''$ be the graph obtained from $G'$ by adding an edge $v^*v$ for each $v \in V(G_0) \cap (\bigcup_{i=1}^pV(G_i)-\{v^*\})$.
Let $\Sigma'$ be a surface that is obtained from $\Sigma$ by adding $p$ handles.
Note that $G''$ can be drawn in $\Sigma'$ with no edge-crossing.

Let $T_3$ be the graph obtained from $T_2$ by adding an edge $v^*v$ for each $v \in V(G_0) \cap (\bigcup_{i=1}^pV(G_i)-\{v^*\})$.
So $T_3$ is a connected spanning subgraph of $G''$.
Let $T_4$ be a breadth-first-search tree of $T_3$ rooted at $v^*$.
Note that every vertex in $V(G_0) \cap (\bigcup_{i=1}^pV(G_i)-\{v^*\})$ is a neighbor of $v^*$ in $T_4$, and every edge of $T_4$ is either an edge of $T_2$ or equal to $v^*v$ for some $v \in V(G_0) \cap (\bigcup_{i=1}^pV(G_i)-\{v^*\})$.

Let $G'''$ be the graph obtained from $G''$ by for each edge $uv$ in $E(G'')-E(T_4)$, subdividing $uv$ $\lvert V(G) \rvert$ times to create a path $P_{uv}$ between $u$ and $v$ with $\lvert V(G) \rvert+1$ edges.
Note that for each edge $uv \in E(G'')-E(T_4)$, there exist two disjoint subpaths $P_{uv,u}$ and $P_{uv,v}$ of $P_{uv}$, where $P_{uv,u}$ contains $u$ but not $v$ and $P_{uv,v}$ contains $v$ but not $u$, such that the tree, denoted by $T_5$, obtained from $T_4$ by adding $P_{uv,u} \cup P_{uv,v}$ for each edge $uv \in E(G'')-E(T_4)$ is a BFS rooted spanning tree of $G'''$ rooted at $v^*$.
Hence $T_5$ is a BFS rooted spanning tree of $G'''$ rooted at $v^*$ such that $E(T_4) \subseteq E(T_5)$.

Since $G'''$ is a subdivision of $G''$, $G'''$ can be drawn in $\Sigma'$ with no edge-crossing.
Since $\Sigma'$ is a surface of Euler genus at most $g+2p$, by Lemma \ref{genus_cutting}, there exists $Z \subseteq V(G''')$ such that
	\begin{itemize}
		\item[(i)] $Z$ is a union of the vertex-sets of at most $2g+4p$ directed paths in $T_5$ each containing $v^*$,
		\item[(ii)] $G'''-Z$ is planar, and
		\item[(iii)] there exists a connected planar graph $G^+$ containing $G'''-Z$ as a subgraph, and there exists a BFS rooted spanning tree $T^+$ of $G^+$ such that for every directed path $P$ in $T^+$, the subgraph of $P$ induced by $V(P) \cap V(G''')-Z$ is a directed path in $T_5$.
	\end{itemize}

Since $G^+$ is connected and planar, by Theorem \ref{planar_tree_path_bags}, there exists a tree-decomposition $(T^*,\X^1=(X^1_t: t \in V(T^*)))$ of $G^+$ such that for every $t \in V(T^*)$, $X^1_t$ is a union of the vertex-sets of at most 3 directed paths in $T^+$.
For each $t \in V(T^*)$, let $X^2_t = (X^1_t \cap V(G''')) \cup Z$.
Let $\X^2 = (X^2_t: t \in V(T^*))$.
Clearly, $(T^*,\X^2)$ is a tree-decomposition of $G'''$.
For each $t \in V(T^*)$, since $X^1_t$ is a union of the vertex-sets of at most 3 directed paths in $T^+$, $X^2_t-Z$ is a union of the vertex-sets of at most 3 directed paths in $T_5$ by (iii), so $X^2_t$ is a union of the vertex-sets of at most $2g+4p+3$ directed paths in $T_5$ by (i).

For each $t \in V(T^*)$, let $X^3_t$ be the set obtained from $X^2_t$ by, for each edge $uv \in E(G'')-E(T_4)$ and $x \in \{u,v\}$ with $V(P_{uv,x}) \cap X^2_t \neq \emptyset$, deleting $V(P_{uv,x})-\{x\}$ and adding $x$.
Let $\X^3=(X^3_t: t \in V(T^*))$.
Note that $G''$ can be obtained from $G'''$ by, for each $uv \in E(G'')-E(T_4)$ and $x \in \{u,v\}$, contracting $P_{uv,x}$ into $x$.
So $(T^*,\X^3)$ is the tree-decomposition of $G''$ obtained by the corresponding contraction.
Since $T_5$ is obtained from $T_4$ by adding $P_{uv,u} \cup P_{uv,v}$ for each edge $uv \in E(G'')-E(T_4)$, we know for every $t \in V(T^*)$, $X^3_t$ is a union of the vertex-sets of at most $2g+4p+3$ directed paths in $T_4$.

Recall that every edge of $T_4$ is either an edge of $T_2$ or equal to $v^*v$ for some $v \in V(G_0) \cap (\bigcup_{i=1}^pV(G_i)-\{v^*\})$.
So for every directed path $P$ in $T_4$, the underlying graph of $P-v^*$ is a path in $T_2$ and hence its vertex-set is a union of at most two directed paths in $T_2$.
Hence for every $t \in V(T^*)$, $X^3_t$ is a union of the vertex-sets of at most $1+2(2g+4p+3)=4g+8p+7$ directed paths in $T_2$.

For every $t \in V(T^*)$, let $X^4_t$ be the set obtained from $X^3_t$ by for each $v \in X^3_t \cap V(G_0) \cap \bigcup_{i=1}^pV(G_i)$, adding the bag in $(Q_{i_v},\X_{i_v})$ indexed by $v$, where $i_v$ is the unique element of $[p]$ with $v \in V(G_{i_v})$.
Recall that every vertex in $V(G_0) \cap (\bigcup_{i=1}^pV(G_i)-\{v^*\})$ is a neighbor of $v^*$ in $T_4$.
So every directed path in $T_4$ contains at most 2 vertices in $V(G_0) \cap \bigcup_{i=1}^pV(G_i)$.
Hence for every $t \in V(T^*)$, $\lvert X^4_t-X^3_t \rvert \leq 2p \cdot (2g+4p+3)=(4g+8p+6)p$, so 
	\begin{itemize}
		\item[(iv)] $X^4_t$ is a union of a subset of $V(G)$ with size at most $(4g+8p+6)p$ and the vertex-sets of at most $4g+8p+7$ directed paths in $T_2$.
	\end{itemize}
Let $\X^4 = (X^4_t: t \in V(T^*))$.

Let $G_{t^*}$ be the graph obtained from the torso of $(T,\X)$ at $t^*$ by adding $v_t$ and the three edges of $G''$ incident with $v_t$ for each neighbor $t$ of $t^*$ in $T$ for which $v_t$ is defined.

\medskip

\noindent{\bf Claim 1:} $(T^*,\X^4)$ is a tree-decomposition of $G_{t^*}$. 

\noindent{\bf Proof of Claim 1:}
Clearly, $\bigcup_{t \in V(T^*)}X^4_t \supseteq X_{t^*}=V(G_{t^*})$.

Note that every edge of $G_{t^*}$ is either an edge of $G''$ or an edge contained in a bag of $(Q_i,\X_i)$ for some $i \in [p]$. 
So for every $e \in E(G_{t^*})$, if $e$ is an edge of $G''$, then there exists $t \in V(T^*)$ with $e \subseteq X^3_t \subseteq X^4_t$ since $(T^*,\X^3)$ is a tree-decomposition of $G''$; otherwise, there exists $v \in V(G_0) \cap \bigcup_{i \in [p]}V(G_i)$ such that $e$ is contained in the bag of $(Q_{i_v},\X_{i_v})$ indexed by $v$, so there exists $t \in V(T^*)$ such that $e \subseteq X^4_t$.

Note that every vertex of $G_{t^*}$ is either a vertex of $G''-\bigcup_{i \in [p]}V(G_i)$ or a vertex of $\bigcup_{i \in [p]}G_i$.
In addition, for every $t \in V(T^*)$, $X^4_t-X^3_t \subseteq \bigcup_{i \in [p]}V(G_i)$.
So for every $x \in V(G_{t^*})$, if $x \in V(G'')-\bigcup_{i \in [p]}V(G_i)$, then the set $\{t \in V(T^*): x \in X^4_t\}=\{t \in V(T^*): x \in X^3_t\}$ induces a connected subgraph of $T^*$ since $(T^*,\X^3)$ is a tree-decomposition of $G''$; otherwise, the set $\{t \in V(T^*): x \in X^4_t\}$ equals the set $$\{t \in V(T^*): x \text{ is contained in a bag of $(Q_{i_v},\X_{i_v})$ indexed by $v$ for some } v \in X^3_t \cap V(G_0) \cap \bigcup_{i \in [p]}V(G_i)\},$$ and this set induces a connected subgraph of $T^*$ since for each $i \in [p]$, $(Q_i,\X_i)$ is a path-decomposition and $G''$ contains a cycle consisting of the vertices in $V(G_i) \cap V(G_0)$ and passing through them in an order consistent with $\Omega_i$.
Therefore, $(T^*,\X^4)$ is a tree-decomposition of $G_{t^*}$.
$\Box$

\medskip

By the definition of $T_2$, for each edge $e$ in $E(T_2)-E(T_1)$, there exists a directed path $P_e$ in $T_1$ connecting the ends of $e$ (for the case that $e$ is incident with $v_t$ for some neighbor $t$ of $T^*$ in $T$, we treat $v_t$ as the vertex $v$ in the definition of $v_t$).

For every $t \in V(T^*)$, define $X^*_t$ to be the set obtained from $X^4_t \cap V(G)$ by, for each neighbor $t'$ of $t^*$ in $T$, doing the following operations:
	\begin{itemize}
		\item If $X^4_t \supseteq X_{t'} \cap X_{t^*}$ and $v_{t'}$ is undefined, then adding $X_{t'}$ into $X^*_t$.
		\item If $v_{t'}$ is defined and $X^4_t \supseteq (X_{t'} \cap X_{t^*}) \cup \{v_{t'}\}$, then adding $X_{t'}$ into $X^*_t$.
		\item If $X^4_t$ contains both ends of some edge $e$ in $E(T_2)-E(T_1)$, then adding $V(P_e)$ into $X^*_t$. 
	\end{itemize}
Let $\X^*=(X^*_t: t \in V(T^*))$.

\medskip

\noindent{\bf Claim 2:} $(T^*,\X^*)$ is a tree-decomposition of $G$. 

\noindent{\bf Proof of Claim 2:}
For every neighbor $t'$ of $T$ of $t^*$, since $X_{t^*} \cap X_{t'}$ is either a clique of $G_0$ or contained in a bag of $(Q_i,\X_i)$ for some $i \in [p]$, we know $\bigcup_{t \in V(T^*)}X^*_t \supseteq V(G)$ by Claim 1.
Similarly, for every edge $e$ of $G$, there exists $t \in V(T^*)$ such that $X^*_t$ contains the both ends of $e$.

Let $x \in V(G)$.
Let $S = \{t \in V(T^*): x \in X^*_t\}$.
Clearly, if $x \not \in V(P_e)$ for some $e \in E(T_2)-E(T_1)$, then $S$ induces a connected subgraph of $T^*$.
So we may assume that $x \in V(P_e)$ for some $e \in E(T_2)-E(T_1)$.

By the definition of $T_2$, for every edge $e \in E(T_2)-E(T_1)$, either 
	\begin{itemize}
		\item $e \subseteq X_t \cap X_{t^*}$ for some neighbor $t$ of $t^*$ in $T$, and $v_t$ is undefined for all such $t$, or 
		\item $e \subseteq (X_t \cap X_{t^*}) \cup \{v_t\}$ is incident with $v_t$ for a unique neighbor $t$ of $t^*$ in $T$, and $v_t$ is defined.
	\end{itemize}
For the former, the set $\{z \in V(T^*): e \subseteq X^4_z\}$ induces a connected subgraph of $T^*$ containing all nodes $t'$ of $T^*$ with $X^4_{t'} \supseteq X_t \cap X_{t^*}$; for the latter, the set $\{z \in V(T^*): e \subseteq X^4_z\}$ induces a connected subgraph of $T^*$ containing all nodes $t'$ of $T^*$ with $X^4_{t'} \supseteq (X_t \cap X_{t^*}) \cup \{v_t\}$.
Hence $S$ induces a connected subgraph of $T^*$.
This shows that $(T^*,\X^*)$ is a tree-decomposition of $G$.
$\Box$

\medskip

Note that for each directed path $P$ in $T_2$, by replacing each edge $e$ in $E(P)-E(T_1) \subseteq E(T_2)-E(T_1)$ by $P_e$, we obtain a directed path $\overline{P}$ in $T_1$.
By the definition of $\X^*$, 
	\begin{itemize}
		\item[(v)] for any $t \in V(T^*)$ and directed path $P$ in $T_2$ with $V(P) \subseteq X^4_t$, we know $V(\overline{P}) \subseteq X^*_t$.
	\end{itemize}

By (iv), for every $t \in V(T^*)$, $X^4_t$ is a union of a subset of $V(G)$ with size at most $(4g+8p+6)p$ and the vertex-sets of at most $4g+8p+7$ directed paths in $T_2$.
So for every $t \in V(T^*)$, $X^4_t$ is a union of a subset $W_t$ of $V(G)-V(G')=V(G)-V(T_2)$ with size at most $(4g+8p+6)p$ and the vertex-sets of $j_t$ directed paths $R_{t,1},R_{t,2},...,R_{t,j_t}$ in $T_2$ for some $j_t \in [4g+8p+7+(4g+8p+6)p] \cup \{0\}$.
For each $t \in V(T^*)$ and $i \in [j_t]$, let $P_{t,i}=\overline{R_{t,i}}$, so $V(P_{t,i}) \subseteq X^*_t$ by (v).

Note that for every $t \in V(T^*)$, if an edge $e$ of $T_2$ satisfies $e \subseteq X^4_t$, then both ends of $e$ are contained in $\bigcup_{i\in [j_t]}R_{t,i}$.
Since $T_2$ is a rooted tree, for each $t \in V(T^*)$, we may assume that if an edge $e$ of $T_2$ satisfies $e \subseteq X^4_t$, then $e \subseteq \bigcup_{i\in [j_t]}R_{t,i}$.
So by the definition of $\X^*$, for each $t \in V(T^*)$ and $x \in X^*_t-X^4_t$, either $x \in X_{t'}-X_{t^*}$ and $X^4_t \supseteq X_{t'} \cap X_{t^*}$ for some neighbor $t'$ of $t^*$ in $T$, or $x \in \bigcup_{i \in [j_t]}V(P_{t,i})$.
Hence for every $t \in V(T^*)$, $X^*_t$ is the union of $W_t \cup \bigcup_{i \in [j_t]}V(P_{t,i})$ and a subset of $\bigcup_{t' \in S^*}X_{t'}$, where $S^* = \{t' \in V(T): t't^* \in E(T),X_{t'} \cap X_{t^*} \subseteq X^4_t\}$.
By Assumption 4 of this lemma, $$\bigcup_{t' \in S^*}X_{t'} \subseteq \bigcup_{t' \in S^*}N_{(G,\phi)}^{\leq \lambda}[X_{t'} \cap X_{t^*}] \subseteq N_{(G,\phi)}^{\leq \lambda}[X^4_t \cap V(G)] \subseteq N_{(G,\phi)}^{\leq \lambda}[W_t \cup \bigcup_{i \in [j_t]}V(P_{t,i})].$$
So for every $t \in V(T^*)$, $X^*_t \supseteq W_t \cup \bigcup_{i \in [j_t]}V(P_{t,i})$ and $X^*_t \subseteq N_{(G,\phi)}^{\leq \lambda}[W_t \cup \bigcup_{i \in [j_t]}V(P_{t,i})]$.

For every $t \in V(T^*)$, let $i_t = j_t+ \lvert W_t \rvert$, and for each $i \in [i_t]-[j_t]$, let $P_{t,i}$ be a path consisting of one vertex of $W_t$.
Therefore, $X^*_t \supseteq \bigcup_{i \in [i_t]}V(P_{t,i})$ and $X^*_t \subseteq N_{(G,\phi)}^{\leq \lambda}[\bigcup_{i \in [i_t]}V(P_{t,i})]$.
Note that for each $t \in V(T^*)$ and $i \in [i_t]$, $P_{t,i}$ is a directed path in $T_1$, so it is a subpath of a $v^*$-geodesic in $(G,\phi)$.
In addition, $i_t \leq j_t+\lvert W_t \rvert \leq (4g+8p+7+(4g+8p+6)p) + (4g+8p+6)p \leq (4g+8p+7)(2p+1)$.
This proves the lemma.
\end{pf}

\section{Asymptotic dimension of weighted minor-closed families} \label{sec:weighted_minor}

Let $m$ be a nonnegative integer.
Let $X,Y$ be metric spaces.
We denote the metrics of $X$ and $Y$ by $d_X$ and $d_Y$, respectively.
Let $f: X \rightarrow Y$ be a function.
Recall that for a real number $r$, a subset $S$ of $X$ is $(\infty,r)$-bounded with respect to $f$ if $d_Y(f(x_1),f(x_2))\leq r$ for any $x_1,x_2 \in S$.
An \defn{$m$-dimensional control function of $f$} is a function $D_f: {\mathbb R}^+ \times {\mathbb R}^+ \rightarrow {\mathbb R}^+$ such that for any positive real numbers $r_X$ and $r_Y$ and any $(\infty,r_Y)$-bounded set $S$ in $X$ with respect to $f$, there exist collections $\X_1,\X_2,...,\X_{m+1}$ of subsets of $X$ such that 
	\begin{itemize}
		\item $\bigcup_{i=1}^{m+1}\bigcup_{W \in \X_i}W =S$,
		\item for each $i \in [m+1]$, any distinct members of $W_1,W_2$ of $\X_i$, and any elements $w_1 \in W_1$ and $w_2 \in W_2$, we have $d_X(w_1,w_2)>r_X$, and
		\item for each $i \in [m+1]$, $W \in \X_i$, and any elements $w_1,w_2 \in W$, we have $d_X(w_1,w_2) \leq D_f(r_X,r_Y)$.
	\end{itemize}
We say that $f$ is \defn{large-scale uniform} if there exists a function $c_f: {\mathbb R}^+ \rightarrow {\mathbb R}^+$ such that $d_Y(f(x_1),f(x_2)) \allowbreak \leq c_f(d_X(x_1,x_2))$ for any $x_1,x_2 \in X$.

The following is an immediate corollary of a combination of results of \cite{bdlm}.

\begin{theorem}[\cite{bdlm}] \label{large_scale_metric_projection}
Let $m,n$ be nonnegative integers.
Let $X,Y$ be metric spaces with $\ad(Y) \leq n$.
Let $f: X \rightarrow Y$ be a large-scale uniform function.
If there exists an $m$-dimensional control function of $f$, then $\ad(X) \leq m+n$.
\end{theorem}

\begin{pf}
By \cite[Proposition 4.7]{bdlm}, since there exists an $m$-dimensional control function of $f$, there exists an $(m,m+n+1)$-control function of $f$.
By \cite[Theorem 4.9]{bdlm}, $\ad(X) \leq m+n$.
(We omit the definition of an  $(m,m+n+1)$-control function of $f$.)
\end{pf}

\begin{lemma} \label{ad_geo_control}
Let $q$ be a positive integer.
Let $\lambda$ be a positive real number. 
Let $\F$ be the class of weighted graphs such that for every $(G,\phi) \in \F$, there exist a vertex $v^* \in V(G)$ and a tree-decomposition $(T,\X=(X_t: t \in V(T)))$ of $G$ such that for every $t \in V(T)$, there exist $i_t \in [q] \cup \{0\}$ and paths $P_{t,1},P_{t,2},...,P_{t,i_t}$ in $G$ such that 
	\begin{enumerate}
		\item for each $i \in [i_t]$, $P_{t,i}$ is a subpath of a $v^*$-geodesic in $(G,\phi)$ with $V(P_{t,i}) \subseteq X_t$, and
		\item $X_t \subseteq N_{(G,\phi)}^{\leq \lambda}[\bigcup_{i \in [i_t]}V(P_{t,i})]$.
	\end{enumerate}
Then $\ad(\F) \leq 2$. 
\end{lemma}

\begin{pf}
Let $q$ be a positive integer.
Let $\lambda$ be a positive real number.
Let $g: {\mathbb R}^+ \times {\mathbb R}^+ \rightarrow {\mathbb R}^+$ be the function such that for any $x,y \in {\mathbb R}^+$, $g(x,y)$ equals the number $\nu^*$ given by Lemma \ref{strong_weighted_tree_extension_control_clean} by taking $(\ell,\mu,\theta)=(x,y,q+3)$.

Let $\F$ be the class of weighted graphs mentioned in the statement of the lemma.
Let $(G,\phi) \in \F$.
So there exist a vertex $v^* \in V(G)$ and a tree-decomposition $(T,\X=(X_t: t \in V(T)))$ of $G$ such that for every $t \in V(T)$, there exist $i_t \in [q] \cup \{0\}$ and paths $P_{t,1},P_{t,2},...,P_{t,i_t}$ in $G$ such that 
	\begin{itemize}
		\item[(i)] for each $i \in [i_t]$, $P_{t,i}$ is a subpath of a $v^*$-geodesic in $(G,\phi)$ with $V(P_{t,i}) \subseteq X_t$, and
		\item[(ii)] $X_t \subseteq N_{(G,\phi)}^{\leq \lambda}[\bigcup_{i \in [i_t]}V(P_{t,i})]$.
	\end{itemize}

Let $X$ be the metric space whose underlying space is $V(G)$, and the metric $d_X$ of $X$ is the distance function in $(G,\phi)$.
Let $d_Y: {\mathbb R} \times {\mathbb R} \rightarrow {\mathbb R}$ be the function such that $d_Y(a,b)=\lvert a-b \rvert$ for any $a,b \in {\mathbb R}$.
Let $Y$ be the metric space $({\mathbb R},d_Y)$.
Let $f: X \rightarrow Y$ be the function such that $f(x)=d_X(x,v^*)$ for every $x \in X$.

Define $D_f: {\mathbb R}^+ \times {\mathbb R}^+ \rightarrow {\mathbb R}^+$ to be the function such that $D_f(a,b)=a \cdot g(a,2\lambda+a+b)$ for any $a,b \in {\mathbb R}^+$. 

Since $d_Y(f(x_1),f(x_2))=\lvert d_X(x_1,v^*)-d_X(x_2,v^*) \rvert \leq d_X(x_1,x_2)$ for any $x_1,x_2 \in X$, we know that $f$ is a large-scale uniform function.
Since $\ad(Y) \leq 1$, if $D_f$ is a $1$-dimensional control function of $f$, then $\ad(X) \leq 1+1=2$ by Theorem \ref{large_scale_metric_projection}; note that the $2$-dimensional control function for $X$ only depends on $D_f$, so $\ad(\F) \leq 2$.

Hence to prove this lemma, it suffices to show that $D_f$ is a $1$-dimensional control function of $f$

Let $r_X$ and $r_Y$ be positive real numbers.
Let $S$ be an $(\infty,r_Y)$-bounded set in $X$ with respect to $f$.
So $\lvert d_X(a,v^*)-d_X(b,v^*) \rvert = \lvert f(a)-f(b) \rvert \leq r_Y$ for any $a,b \in S$.
Hence there exists a nonnegative real number $k$ such that $k \leq d_X(a,v^*) \leq k+r_Y$ for any $a \in S$.

To prove this lemma, it suffices to show that there exists a desired cover for $S$.
So we may assume that $S$ is a maximal $(\infty,r_Y)$-bounded set in $X$ with respect to $f$.
Hence $S = \{v \in X: k \leq d_X(v,v^*) \leq k+r_Y\}$.

For any $t \in V(T)$ and $i \in [i_t]$, since $P_{t,i}$ is a subpath of a $v^*$-geodesic in $(G,\phi)$ by (i), we know that $P_{t,i}[S \cap V(P_{t,i})]$ is a subpath of $P_{t,i}$, and we let $v_{t,i}$ be the vertex in $P_{t,i}[V(P_{t,i}) \cap S]$ closest to $v^*$ in $(G,\phi)$; so $V(P_{t,i}) \cap S \subseteq N_{(G,\phi)[V(P_{t,i}) \cap S]}^{\leq r_Y}[\{v_{t,i}\}]$, and hence $$\bigcup_{i \in [i_t]}V(P_{t,i}) \cap S \subseteq N_{(G,\phi)[\bigcup_{i \in [i_t]}V(P_{t,i}) \cap S]}^{\leq r_Y}[\{v_{t,i}: i \in [i_t]\}].$$

Let $S'' = \{v \in X: k-r_X \leq d_X(v,v^*) \leq k+r_Y+r_X\}$.

\medskip

\noindent{\bf Claim 1:} For any $t \in V(T)$ and $z \in S'' \cap X_t$, there exists $i_z \in [i_t]$ such that $$\dist_{(G,\phi)}(v_{t,i_z},z) \leq 2\lambda+r_X+r_Y.$$ 

\noindent{\bf Proof of Claim 1:}
Let $t \in V(T)$ and $z \in S'' \cap X_t$.
Since $X_t \subseteq N_{(G,\phi)}^{\leq \lambda}[\bigcup_{i \in [i_t]}V(P_{t,i})]$ by (ii), there exist $i_z \in [i_t]$ and $v_z \in V(P_{t,i_z})$ such that $z \in N_{(G,\phi)}^{\leq \lambda}[\{v_z\}]$.
Since $z \in S''$, we know $k-r_X \leq \dist_{(G,\phi)}(z,v^*) \leq \dist_{(G,\phi)}(z,v_z)+\dist_{(G,\phi)}(v_z,v^*) \leq \lambda + \dist_{(G,\phi)}(v_z,v^*)$.
So $\dist_{(G,\phi)}(v_z,v^*) \geq k-r_X-\lambda$.
Similarly, $\dist_{(G,\phi)}(v_z,v^*) \leq \dist_{(G,\phi)}(v_z,z) + \dist_{(G,\phi)}(z,v^*) \leq \lambda + k+r_Y + r_X$.
Hence $v_z$ is a vertex in $V(P_{t,i_z})$ such that $k-r_X - \lambda \leq \dist_{(G,\phi)}(v_z,v^*) \leq \lambda + k +r_Y+r_X$.
Since $P_{t,i_z}$ is a subpath of a $v^*$-geodesic and $k \leq \dist_{(G,\phi)}(v_{t,i_z},v^*) \leq k+r_Y$, we have $\dist_{(G,\phi)}(v_{t,i_z},v_z) \leq \lambda+r_X+r_Y$.
Hence $\dist_{(G,\phi)}(v_{t,i_z},z) \leq \dist_{(G,\phi)}(v_{t,i_z},v_z) + \dist_{(G,\phi)}(v_z,z) \leq 2\lambda+r_X+r_Y$.
$\Box$

\medskip

Let $(G',\phi')$ be the weighted graph obtained from $(G,\phi)$ by, 
	\begin{itemize}
		\item for any $t \in V(T)$ and any two distinct vertices $a,b$ in $X_t$, adding an edge $ab$, (note that we create a parallel edge if $a$ and $b$ are already adjacent in $G$),
		\item for every $e \in E(G)$, defining $\phi'(e)=\phi(e)$, and
		\item for every edge $xy \in E(G')-E(G)$, defining $\phi'(xy) = \dist_{(G,\phi)}(x,y)$.
	\end{itemize}

Note that for every $e \in E(G')-E(G)$, both ends of $e$ are contained in $X_t$ for some $t \in V(T)$.
Hence $(T,\X)$ is a tree-decomposition of $G'$.

\medskip

\noindent{\bf Claim 2:} For any $x,y \in V(G)=V(G')$, $\dist_{(G,\phi)}(x,y) = \dist_{(G',\phi')}(x,y)$.

\noindent{\bf Proof of Claim 2:}
Since $(G,\phi)$ is a subgraph of $(G',\phi')$, $\dist_{(G,\phi)}(x,y) \geq \dist_{(G',\phi')}(x,y)$.
Let $P$ be a path in $(G',\phi')$ between $x$ and $y$ with $\leng_{(G',\phi')}(P) = \dist_{(G',\phi')}(x,y)$.
By the definition of $(G',\phi')$, for every $e \in E(G')$, there exists a path $P_e$ in $G$ between the ends of $e$ such that $\leng_{(G,\phi)}(P_e) \leq \phi'(e)$.
So by replacing each edge $e$ of $P$ by the path $P_e$ in $G$, we obtain a walk in $G$ between $x$ and $y$ with length in $(G,\phi)$ at most $\leng_{(G',\phi')}(P) = \dist_{(G',\phi')}(x,y)$.
Hence $\dist_{(G,\phi)}(x,y) \leq \dist_{(G',\phi')}(x,y)$.
$\Box$

\medskip

Let $(G'',\phi'') = (G',\phi')[S'']$.
For every $t \in V(T)$, let $X''_t = X_t \cap S''$. 
Let $\X'' = (X_t'': t \in V(T))$.
Then $(T,\X'')$ is a tree-decomposition of $(G'',\phi'')$.

\medskip

\noindent{\bf Claim 3:} For any $t \in V(T)$ and distinct vertices $x,y \in X''_t=X_t \cap S''$, we have $\dist_{(G,\phi)}(x,y) = \dist_{(G'',\phi'')}(x,y)$. 

\noindent{\bf Proof of Claim 3:}
Since $x,y \in X_t$, we know $xy \in E(G'') \subseteq E(G')$ and $\dist_{(G'',\phi'')}(x,y) \leq \phi'(xy) = \dist_{(G,\phi)}(x,y) = \dist_{(G',\phi')}(x,y)$, where the last equality follows from Claim 2.
Since $(G'',\phi'')$ is a subgraph of $(G',\phi')$, $\dist_{(G'',\phi'')}(x,y) \geq \dist_{(G',\phi')}(x,y)$.
So $\dist_{(G'',\phi'')}(x,y) = \dist_{(G',\phi')}(x,y) = \dist_{(G,\phi)}(x,y)$.
$\Box$

\medskip

For any $t \in V(T)$ and $z \in X_t''$, since $\{z,v_{t,i_z}\} \subseteq X_t''$, Claims 1 and 3 imply that 
$$z \in N_{(G'',\phi'')}^{\leq \dist_{(G,\phi)}(z,v_{t,i_z})}[\{v_{t,i_z}\}] \subseteq N_{(G'',\phi'')}^{\leq 2\lambda+r_X+r_Y}[\{v_{t,i}: i \in [i_t]\}].$$
For every $t \in V(T)$, let $A_t = \{v_{t,i}: i \in [i_t]\}$, so $A_t \subseteq X_t \cap S \subseteq X_t''$ (by (i)) with $|A_t| \leq i_t \leq q$ and $X_t'' \subseteq N_{(G'',\phi'')}^{\leq 2\lambda+r_X+r_Y}[A_t]$.

Let $(G_1,\phi_1)=(G'',\phi'',r_X)$.
Since $(G_1,\phi_1)$ is obtained from $(G'',\phi'')$ by duplicating edges and subdividing edges, there exists a tree-decomposition $(T^1,\X^1=(X^1_t: t \in V(T^1)))$ of $(G_1,\phi_1)$ such that $T$ is a subtree of $T^1$, and 
	\begin{itemize}
		\item for every $t \in V(T)$, $X^1_t=X_t'$, and
		\item for every $t \in V(T^1)-V(T)$, $|X^1_t| \leq 3$.
	\end{itemize}
For every $t \in V(T^1)-V(T)$, let $A_t = X^1_t$, so $|A_t| \leq 3$.
Therefore, for every $t \in V(T^1)$, we have $A_t \subseteq X^1_t$ with $|A_t| \leq q+3$ and $X^1_t \subseteq N_{(G_1,\phi_1)}^{\leq 2\lambda+r_X+r_Y}[A_t]$ by Lemma \ref{length_subdiv}.

Hence $(G_1,\phi_1)$ is a $(0,r_X]$-bounded weighted graph that has a tree-decomposition $(T^1,\X^1)$ such that for every $t \in V(T^1)$, $X^1_t \subseteq N_{(G_1,\phi_1)}^{\leq 2\lambda+r_X+r_Y}[A_t]$ for some $A_t \subseteq X^1_t$ with $|A_t| \leq q+3$.
By Lemma \ref{strong_weighted_tree_extension_control_clean}, there exists a 2-coloring $c$ of $(G_1,\phi_1)^{r_X}$ with weak diameter in $(G_1,\phi_1)^{r_X}$ at most $g(r_X,2\lambda+r_X+r_Y)$.

For $i \in [2]$, let $\X_i = \{V(M) \cap S: M$ is a $c$-monochromatic component in $(G_1,\phi_1)^{r_X}$ with color $i\}$.
Clearly, $\bigcup_{i=1}^2\bigcup_{W \in \X_i}W = V(G_1) \cap S = S$.

\medskip

\noindent{\bf Claim 4:} For any $i \in [2]$, any distinct $W_1,W_2 \in \X_i$, any elements $w_1 \in W_1$ and $w_2 \in W_2$, we have $\dist_{(G,\phi)}(w_1,w_2) > r_X$.

\noindent{\bf Proof of Claim 4:}
Since $W_1$ and $W_2$ are contained in different $c$-monochromatic components in $(G_1,\phi_1)^{r_X}$ with the same color, we have $\dist_{(G_1,\phi_1)}(w_1,w_2)>r_X$.
Since $\{w_1,w_2\} \subseteq S \subseteq V(G'')$, $\dist_{(G'',\phi'')}(w_1,w_2) = \dist_{(G_1,\phi_1)}(w_1,w_2)>r_X$ by Lemma \ref{length_subdiv}.

Suppose to the contrary that $\dist_{(G,\phi)}(w_1,w_2) \leq r_X$.
Then there exists a path $P$ in $G$ between $w_1$ and $w_2$ with $\leng_{(G,\phi)}(P) \leq r_X$.
Let $v \in V(P)$.
Then $\dist_{(G,\phi)}(v,w_1) \leq \leng_{(G,\phi)}(P) \leq r_X$.
Since $w_1 \in S$, we know $\dist_{(G,\phi)}(v,v^*) \leq \dist_{(G,\phi)}(v,w_1)+\dist_{(G,\phi)}(w_1,v^*) \leq r_X+k+r_Y$ and $\dist_{(G,\phi)}(v,v^*) \geq \dist_{(G,\phi)}(v^*,w_1)-\dist_{(G,\phi)}(w_1,v) \geq k-r_X$.
So $v \in S''$.
Since $v$ is an arbitrary vertex of $P$, we know $V(P) \subseteq S''$, so $P \subseteq G[S''] \subseteq G'[S''] \subseteq G''$ and $\leng_{(G,\phi)}(P) = \leng_{(G'',\phi'')}(P)$.
Hence $\dist_{(G'',\phi'')}(w_1,w_2) \leq \leng_{(G'',\phi'')}(P) = \leng_{(G,\phi)}(P) \leq r_X$, a contradiction.
$\Box$

\medskip

By Claim 4, to show that $D_f$ is a $1$-dimensional control function of $f$, it suffices to show that for any $i \in [2]$ and $W \in \X_i$, $W$ has weak diameter in $(G,\phi)$ at most $D_f(r_X,r_Y)$.

Let $i \in [2]$ and $W \in \X_i$.
Since $W$ is contained in a $c$-monochromatic component in $(G_1,\phi_1)^{r_X}$, the weak diameter in $(G_1,\phi_1)^{r_X}$ of $W$ is at most $g(r_X,2\lambda+r_X+r_Y)$.
So the weak diameter of $W$ in $(G_1,\phi_1)$ is at most $r_X \cdot g(r_X,2\lambda+r_X+r_Y)=D_f(r_X,r_Y)$.
Since $V(W) \subseteq S \subseteq V(G'')$, by Lemma \ref{length_subdiv}, the weak diameter of $W$ in $(G'',\phi'')$ is at most $D_f(r_X,r_Y)$.
Since $(G'',\phi'')$ is a subgraph of $(G',\phi')$, the weak diameter in $(G',\phi')$ of $W$ is at most $D_f(r_X,r_Y)$.
By Claim 2, the weak diameter in $(G,\phi)$ of $W$ is at most $D_f(r_X,r_Y)$.
This proves the lemma.
\end{pf}

\begin{lemma} \label{wd_geo_bags}
For any positive integer $q$ and positive real numbers $\lambda,\ell$, there exists a positive real number $\nu^*$ such that the following holds.
Let $\F$ be the class of $(0,\ell]$-bounded weighted graphs such that for every $(G,\phi) \in \F$, there exist a vertex $v^* \in V(G)$ and a tree-decomposition $(T,\X=(X_t: t \in V(T)))$ of $G$ such that for every $t \in V(T)$, there exist $i_t \in [q] \cup \{0\}$ and paths $P_{t,1},P_{t,2},...,P_{t,i_t}$ in $G$ such that 
	\begin{enumerate}
		\item for each $i \in [i_t]$, $P_{t,i}$ is a subpath of a $v^*$-geodesic in $(G,\phi)$ with $V(P_{t,i}) \subseteq X_t$, and
		\item $X_t \subseteq N_{(G,\phi)}^{\leq \lambda}[\bigcup_{i \in [i_t]}V(P_{t,i})]$.
	\end{enumerate}
Then for every $(0,\ell]$-bounded weighted graph $(G,\phi) \in \F$, $(G,\phi)^\ell$ is $3$-colorable with weak diameter in $(G,\phi)^\ell$ at most $\nu^*$. 
\end{lemma}

\begin{pf}
Let $q$ be a positive integer.
Let $\lambda$ and $\ell$ be positive real numbers.
Let $\F$ be the class mentioned in the statement of this lemma.
Note that $\F$ only depends on $q,\ell,\lambda$.

Let $\F_1$ be the class of $(0,\ell]$-bounded weighted graphs such that for every $(G,\phi) \in \F_1$, there exists $(G_0,\phi_0) \in \F$ such that $(G,\phi)$ can be obtained from $(G_0,\phi_0)$ by repeatedly applying any of the following operations. 
	\begin{itemize}
		\item Duplicating one edge, and defining the weight of the copied edge to be equal to the weight of its original.
		\item Subdividing an edge $e$, and defining the weight of each of the two resulting edges to be half of the weight of $e$.
	\end{itemize}

\medskip

\noindent{\bf Claim 1:} For every $(G,\phi) \in \F_1$, there exist a vertex $v^* \in V(G)$ and a tree-decomposition $(T,\X=(X_t: t \in V(T)))$ of $G$ such that for every $t \in V(T)$, there exist $i_t \in [q] \cup \{0\}$ and paths $P_{t,1},P_{t,2},...,P_{t,i_t}$ in $G$ such that 
	\begin{itemize}
		\item for each $i \in [i_t]$, $P_{t,i}$ is a subpath of a $v^*$-geodesic in $(G,\phi)$ with $V(P_{t,i}) \subseteq X_t$, and
		\item $X_t \subseteq N_{(G,\phi)}^{\leq \lambda+\ell}[\bigcup_{i \in [i_t]}V(P_{t,i})]$.
	\end{itemize}

\noindent{\bf Proof of Claim 1:}
Let $(G,\phi) \in \F_1$.
So there exists $(G_0,\phi_0) \in \F$ such that $(G,\phi)$ is obtained from $(G_0,\phi_0)$ by repeatedly applying the operations. 
Since $(G_0,\phi_0) \in \F$, there exist a vertex $v^* \in V(G_0)$ and a tree-decomposition $(T,\X^0=(X^0_t: t \in V(T)))$ of $G$ such that for every $t \in V(T)$, there exist $i_t \in [q] \cup \{0\}$ and paths $P^0_{t,1},P^0_{t,2},...,P^0_{t,i_t}$ in $G_0$ such that 
	\begin{itemize}
		\item for each $i \in [i_t]$, $P^0_{t,i}$ is a subpath of a $v^*$-geodesic in $(G_0,\phi_0)$ with $V(P^0_{t,i}) \subseteq X^0_t$, and
		\item $X^0_t \subseteq N_{(G_0,\phi_0)}^{\leq \lambda}[\bigcup_{i \in [i_t]}V(P^0_{t,i})]$.
	\end{itemize}

Note that $v^* \in V(G_0) \subseteq V(G)$, and for every $e \in E(G_0)$, there exists a subgraph $H_e$ of $G$ consisting of the vertices and edges derived from $e$ in the process of generating $(G,\phi)$ from $(G_0,\phi_0)$.
Note that it is possible that $H_e$ consists of the edge $e$.
Observe that $\bigcup_{e \in E(G_0)}H_e=G$.

For each $t \in V(T)$, let $X_t = X^0_t \cup \bigcup_{e \in E(G_0), e \subseteq X^0_t}V(H_e)$.
Let $\X=(X_t: t \in V(T))$.
Then $(T,\X)$ is a tree-decomposition of $G$.

Note that for each $e \in E(G_0)$, there exists a path $P_e$ in $H_e$ connecting the ends of $e$ such that the length in $(G,\phi)$ of $P_e$ equals $\phi_0(e)$.
For each $t \in V(T)$ and $i \in [i_t]$, let $P_{t,i}$ be the path in $G$ obtained from $P^0_{t,i}$ by replacing each edge $e$ of $P^0_{t,i}$ by $P_e$.
Since the ends of each $P_{t,i}$ are in $V(G_0)$, each $P_{t,i}$ is a subpath of a $v^*$-geodesic in $(G,\phi)$ contained in $X_t$.
In addition, since $V(P^0_{t,i}) \subseteq X^0_t$, $V(P_{t,i}) \subseteq X_t$.

In addition, for every $e=uv \in E(G_0)$, $V(H_e) \subseteq N_{(G,\phi)}^{\leq \phi_0(e)}[\{u,v\}] \subseteq N_{(G,\phi)}^{\leq \ell}[\{u,v\}]$ since $(G_0,\phi_0)$ is $(0,\ell]$-bounded.
So for every $t \in V(T)$, $X_t \subseteq N_{(G,\phi)}^{\leq \ell}[X^0_t]$.
Since for every $t \in V(T)$, $X^0_t \subseteq N_{(G_0,\phi_0)}^{\leq \lambda}[\bigcup_{i \in [i_t]}V(P^0_{t,i})] \subseteq N_{(G,\phi)}^{\leq \lambda}[\bigcup_{i \in [i_t]}V(P^0_{t,i})]$, we have $X_t \subseteq N_{(G,\phi)}^{\lambda+\ell}[\bigcup_{i \in [i_t]}V(P^0_{t,i})] \subseteq N_{(G,\phi)}^{\lambda+\ell}[\bigcup_{i \in [i_t]}V(P_{t,i})]$.
$\Box$

\medskip

By Claim 1 and Lemma \ref{ad_geo_control}, $\ad(\F_1) \leq 2$.
Let $\F_0 = \{G: (G,\phi) \in \F_1\}$.
So $\F_0$ is closed under duplication and taking subdivision.
By Lemma \ref{weighted_ad_wd_converse}, since every member of $\F_1$ is $(0,\ell]$-bounded, there exists a positive real number $\nu^*$ such that for every $(G,\phi) \in \F_1$, $(G,\phi)^\ell$ is 3-colorable with weak diameter in $(G,\phi)^\ell$ at most $\nu^*$.
Note that $\nu^*$ only depends on $\F_0$ and $\F_1$ which only depend on $q,\lambda,\ell$.
So $\nu^*$ only depends on $q,\lambda,\ell$.
Since $\F \subseteq \F_1$, the lemma follows.
\end{pf}

\begin{lemma} \label{hereditary_near_embedding}
For any positive integer $p$, nonnegative integer $g$, and positive real number $\ell$, there exists a positive real number $\nu^*$ such that if $\F$ is the class of $(0,\ell]$-bounded weighted graphs such that for every $(H,\phi_H) \in \F$, there exist a $(0,\ell]$-bounded weighted graph $(G,\phi)$ and a tree-decomposition $(T,\X)$ of $G$ such that 
	\begin{enumerate}
		\item[(i)] $(G,\phi)$ is obtained from $(H,\phi_H)$ by attaching leaves,
		\item[(ii)] $T$ is a star centered at its root $t^*$,
		\item[(iii)] the torso at $t^*$ is $p$-nearly embeddable in a surface $\Sigma_G$ of Euler genus at most $g$ with a witness $(G_0,G_1,...,G_p,\Delta_1,...,\Delta_p,((Q_i,\X_i): i \in [p]))$, 
		\item[(iv)] for each neighbor $t'$ of $t^*$, $X_{t^*} \cap X_{t'}$ is either a clique in $G_{0}$ of size at most three or contained in a bag of $(Q_{i},\X_{i})$ for some $i \in [p]$,
		\item[(v)] for each neighbor $t'$ of $t^*$, $\lvert X_{t'}-X_{t^*} \rvert \leq p$, and 
		\item[(vi)] $V(G)-V(H) \subseteq V(G_0) \cap \bigcup_{i \in [p]}V(G_i)$,
	\end{enumerate}
then $\F^{+p,(0,\ell]}$ is hereditary and $(3,\ell,\nu^*)$-nice.
\end{lemma}

\begin{pf}
Let $p$ be a positive integer.
Let $g$ be a nonnegative integer.
Let $\ell$ be a positive real number.
Let $\nu$ to be the number $\nu^*$ given by Lemma \ref{wd_geo_bags} by taking $(q,\lambda,\ell)=((4g+16p+31)(4p+7),\ell p,\ell)$.
Define $\nu^*$ to be the number $\nu^*$ given by Lemma \ref{apex_extension} by taking $(\ell,\nu,n)=(\ell,\nu,p)$.

\medskip

\noindent{\bf Claim 1:} $\F$ is hereditary.

\noindent{\bf Proof of Claim 1:}
Let $(H,\phi_H) \in \F$.
Let $S \subseteq V(H)$.

Let $(G,\phi)$ be a $(0,\ell]$-bounded weighted graph and $(T,\X)$ be a tree-decomposition of $G$ witnessing the membership of $(H,\phi_H)$ in $\F$.
Let $L=V(G)-V(H)$.
So $L$ is a set of leaves in $G$.
Let $t^*$ be the root of $T$.
Let $(G_0,G_1,...,G_p,\Delta_1,...,\Delta_p,((Q_i,\X_i): i \in [p]))$ be a witness certifying that the torso at $t^*$ is $p$-nearly embeddable in a surface of Euler genus at most $g$.
For every $v \in V(G_0) \cap \bigcup_{i \in [p]}V(G_i)$, let $W_v$ be the bag of $(Q_{i_v},\X_{i_v})$ indexed by $v$, where $i_v$ is the member of $[p]$ with $v \in V(G_{i_v})$, and let $Y_v = W_v \cup \bigcup_{tt^* \in E(T),X_t \cap X_{t^*} \subseteq W_v} X_{t}$. 

Let $S_0 = S - (V(G_0) \cap \bigcup_{i \in [p]}V(G_i))$. 
Let $\X^1=(X_t-S_0: t \in V(T))$.
Clearly, $(T,\X^1)$ is a tree-decomposition of $G-S_0$ satisfying (ii)-(v).

Construct the graph $G'$ from $G-S_0$ by 
	\begin{itemize}
		\item[(OP1)] for each $i \in [p]$ and $v \in V(G_0) \cap V(G_i) \cap S$ with $Y_v \not \subseteq L \cup S$, adding a new vertex $z_v$ drawn at the place where $v$ was drawn and adding an edge between $z_v$ and a vertex in $Y_v-(S \cup L)$, and
		\item[(OP2)] for each $i \in [p]$ and $v \in V(G_0) \cap V(G_i) \cap L-S$,
			\begin{itemize}
				\item if $Y_v \subseteq L \cup S$, then deleting $v$ from $G$,
				\item otherwise, there exists a vertex $v' \in Y_v-(S \cup L)$, and if the neighbor of $v$ in $G$ is in $S$, then adding an edge $vv'$,
			\end{itemize}
	\end{itemize}
For each $e \in E(G')-E(G)$, let $\phi'(e)=\ell$; for each $e \in E(G') \cap E(G)$, let $\phi'(e)=\phi(e)$.
Then $(G',\phi')$ can be obtained from $(H,\phi_H)-S$ by attaching leaves such that $V(G')-(V(H)-S) \subseteq V(G_0) \cap \bigcup_{i \in [p]}V(G_i)$. 
And the bags of $(T,\X^1)$ can be easily revised to become of tree-decomposition $(T,\X^2)$ of $G'$ satisfying (ii)-(v) by removing vertices in $V(G)-V(G')$ and replacing each vertex $v$ in (OP1) by $z_v$.
So $(G',\phi')$ witnesses that $(H,\phi_H)-S \in \F$.
Hence $\F$ is hereditary.
$\Box$

\medskip

By Claim 1, $\F$ is hereditary.
So $\F^{+p,(0,\ell]}$ is hereditary.

Now we show that $\F^{+p,(0,\ell]}$ is $(3,\ell,\nu^*)$-nice.
By Lemma \ref{apex_extension}, it suffices to show that $\F$ is $(3,\ell,\nu)$-nice.

Let $(H,\phi_H) \in \F$.
It suffices to prove that $(H,\phi_H)^\ell$ is $3$-colorable with weak diameter in $(H,\phi_H)^\ell$ at most $\nu$.
Let $(G,\phi)$ be the weighted graph witnessing the membership of $(H,\phi_H)$ in $\F$.
For every vertex $v$ in $V(G)-V(H)$, since $v$ is a leaf in $G$, the neighbors of $v$ in $(G,\phi)^\ell$ are neighbors in $(G,\phi)^\ell$ of the unique neighbor of $v$ in $G$.
Hence it suffices to show that $(G,\phi)^\ell$ is 3-colorable with weak diameter in $(G,\phi)^\ell$ at most $\nu$.
Since we can color vertices component by component, we may assume that $G$ is connected.

Let $(T,\X)$ be the rooted tree-decomposition of $G$ satisfying (ii)-(v).
Let $t^*$ be the root of $T$.
By (iii)-(v), for every neighbor $t'$ of $t^*$ in $T$, $|X_{t'}| = |X_{t'}-X_t| + |X_{t'} \cap X_t| \leq p+\max\{3,p\} \leq 2p+3$.
Since $G$ is connected and $(G,\phi)$ is $(0,\ell]$-bounded, for every neighbor $t'$ of $t^*$ in $T$, $X_{t'} \subseteq N_{(G,\phi)}^{\leq \ell p}[X_{t'} \cap X_{t^*}]$ by (v).
By Lemma \ref{geodesic_bags}, there exist a vertex $v^* \in V(G)$ and a tree-decomposition $(T^*,\X^*=(X^*_t: t \in V(T^*)))$ of $G$ such that for every $t \in V(T^*)$, there exist $i_t \in [(4g+8(2p+3)+7)(2(2p+3)+1)] \cup \{0\}$ and paths $P_{t,1},P_{t,2},...,P_{t,i_t}$ in $G$ such that 
	\begin{itemize}
		\item for each $i \in [i_t]$, $P_{t,i}$ is a subpath of a $v^*$-geodesic in $(G,\phi)$ with $V(P_{t,i}) \subseteq X^*_t$, and
		\item $X^*_t \subseteq N_{(G,\phi)}^{\leq \ell p}[\bigcup_{i \in [i_t]}V(P_{t,i})]$.
	\end{itemize} 
By Lemma \ref{wd_geo_bags}, $(G,\phi)^\ell$ is 3-colorable with weak diameter in $(G,\phi)^\ell$ at most $\nu$.
This proves the lemma.
\end{pf}

\bigskip

Now we are ready to prove Theorem \ref{minor_AN_intro}, which immediately follows from Lemma \ref{scaling_closed_equiv} and the following theorem.

\begin{theorem}
Let $H$ be a graph.
If $\F$ is the class of weighted graphs whose underlying graphs have no $H$-minor, then $\ad(\F) \leq 2$.
\end{theorem}

\begin{pf}
By subdividing edges of $H$, we may assume that $H$ is simple.
By adding isolated vertices into $H$, we may assume that $\lvert V(H) \rvert \geq 5$.
By replacing $H$ by $K_{\lvert V(H) \rvert}$, we may assume that $H$ is a complete graph on at least five vertices.
So duplicating edges and subdividing edges of an $H$-minor free graph keep it $H$-minor free.
By Lemma \ref{weight_wd_ad_2}, it suffices to show that for every positive real number $\ell$, there exists a positive real number $\nu^*_\ell$ such that for every $(0,\ell]$-bounded $(G,\phi) \in \F$, $(G,\phi)^\ell$ is $3$-colorable with weak diameter in $(G,\phi)^\ell$ at most $\nu^*_\ell$.

Let $\ell$ be a positive real number.
We define the following.
	\begin{itemize}
		\item Let $p$ be the positive integer $p$ given by Theorem \ref{rs_structure} by taking $H=H$.
		\item Let $g$ be the maximum integer such that $H$ cannot be drawn in a surface of Euler genus at most $g$.
			(Note that $g$ is nonnegative since $H$ is a complete  graph on at least five vertices.)
		\item Let $\nu_1$ be the number $\nu^*$ given by Lemma \ref{hereditary_near_embedding} by taking $(p,g,\ell)=(p^2,g,\ell)$.
		\item Let $\F_1$ be the class of $(0,\ell]$-bounded graphs such that for every $(H_0,\phi_{H_0}) \in \F_1$, there exist a $(0,\ell]$-bounded graph $(G',\phi')$ and a tree-decomposition $(T',\X')$ of $G'$ such that 
			\begin{enumerate}
				\item[(i)] $(G',\phi')$ is obtained from $(H_0,\phi_{H_0})$ by attaching leaves,
				\item[(ii)] $T'$ is a star centered at its root $t^*$,
				\item[(iii)] the torso at $t^*$ is $p^2$-nearly embeddable in a surface of Euler genus at most $g$ with a witness $(G'_0,G'_1,...,G'_{p^2},\Delta'_1,...,\Delta'_{p^2},((Q'_i,\X'_i): i \in [p^2]))$, 
				\item[(iv)] for each neighbor $t'$ of $t^*$, $X'_{t^*} \cap X'_{t'}$ is either a clique in $G'_{0}$ of size at most three or contained in a bag of $(Q'_{i},\X'_{i})$ for some $i \in [p^2]$,
				\item[(v)] for each neighbor $t'$ of $t^*$, $\lvert X'_{t'}-X'_{t^*} \rvert \leq p^2$, and 
				\item[(vi)] $V(G')-V(H_0) \subseteq V(G'_0) \cap \bigcup_{i \in [p^2]}V(G'_i)$.
			\end{enumerate}
		\item Let $\F'=\F_1^{+p^2,(0,\ell]}$.
			(Note that by Lemma \ref{hereditary_near_embedding}, $\F'$ is hereditary and $(3,\ell,\nu_1)$-nice.)
		\item Define $\nu^*_\ell$ to be the number $\nu^*$ given by Lemma \ref{weighted_tree_extension_clean} by taking $(\ell,\nu,m,\theta,I,\F,\F') = (\ell,\nu_1,3,p, \allowbreak (0,\ell],\F',\F')$.
	\end{itemize}

Let $(G,\phi)$ be a member of $\F$ such that $(G,\phi)$ is $(0,\ell]$-bounded.
By Theorem \ref{rs_structure}, there exists a tree-decomposition $(T,\X=(X_t: t \in V(T)))$ of adhesion at most $p$ such that for every $t \in V(T)$, if $G_t$ is the torso at $t$, then there exists $Z_t \subseteq V(G_t)$ with $\lvert Z_t \rvert \leq p$ such that $G_t-Z_t$ is $p$-nearly embeddable in a surface $\Sigma_t$ of Euler genus at most $g$ with a witness $(G_{t,0},G_{t,1},...,G_{t,p},\Delta_1,...,\Delta_p,((Q_{t,i},\X_{t,i}): i \in [p]))$ such that for every neighbor $t'$ of $t$ in $T$, $X_t \cap X_{t'}-Z_t$ is either a clique in $G_{t,0}$ of size at most three or contained in a bag of $(Q_{t,i},\X_{t,i})$ for some $i \in [p]$.

For every $t \in V(T)$, let $T^t$ be the maximal star centered at $t$ contained in $T$, and let $X^t_x = X_x \cap X_t$ for each $x \in V(T^t)$; then $(T^t,(X^t_x:x \in V(T^t)))$ is a tree-decomposition of $(G,\phi)[X_t]$ such that the torso at $t$ in it is the same as the torso at $t$ in $(T,\X)$.
Hence for every $t \in V(T)$, we know $(G,\phi)[X_t] \in \F_1^{+p,(0,\ell]} \subseteq \F'$, as we may take $(H_0,\phi_{H_0})=(G,\phi)[X_t]-Z_t$, $(G',\phi')=(G,\phi)[X_t]-Z_t$ and $(T',\X')=(T^t,(X^t_x-Z_t:x \in V(T^t)))$.
In addition, for every $t \in V(T)$, $\F'$ contains every $(0,\ell]$-bounded weighted graph that can be obtained from $(G,\phi)[X_t]$ by for each neighbor $t'$ of $t$ in $T$, adding at most $p^2$ new vertices and new edges such that every new edge is either between two new vertices or between a new vertex and $X_t \cap X_{t'}$.
By Lemma \ref{weighted_tree_extension_clean}, $(G,\phi)^\ell$ is 3-colorable with weak diameter in $(G,\phi)^\ell$ at most $\nu^*_\ell$.
This proves the theorem.
\end{pf}

\section{Concluding remarks} \label{sec:concluding_remarks}

In this paper we determine the Assouad-Nagata dimension of minor-closed families of weighted graphs.
In particular, weighted graphs of bounded Euler genus have Assouad-Nagata dimension at most 2.
The asymptotic dimension version of this particular result was proved in an unpublished draft \cite{bbegps} via fat minors and was strengthened to Assouad-Nagata dimension in its journal version \cite{bbeglps_merged} due to the discovery of Lemma \ref{scaling_closed_equiv}.
This paper gives a different proof for this bounded genus result without considering fat minors.

The notion of fat minors was implicitly introduced by Fujiwara and Papasoglu \cite{fp} when considering the asymptotic dimension of geodesic metric spaces homeomorphic to ${\mathbb R}^2$ and explicitly defined by Bonamy et al.\ \cite{bbeglps_merged,bbegps} when considering the asymptotic dimension of graphs of bounded Euler genus mentioned above.
For an integer $k$, we say that a graph $H$ is a \defn{$k$-fat minor} of another graph $G$ if there exist a set $\{H_v: v \in V(H)\}$ of pairwise disjoint connected subgraphs of $G$ and a set $\{P_e: e \in E(H)\}$ of pairwise disjoint paths in $G$ internally disjoint from $\bigcup_{v \in V(H)}V(H_v)$ such that for every $uv \in E(H)$, $P_{uv}$ is between $V(H_u)$ and $V(H_v)$, and the distance in $G$ between any two distinct members in $\{H_v, P_e: v \in V(H),e \in E(H)\}$ is at least $k$ unless one of them is $H_v$ and the other is $P_e$ for some $e \in E(H)$ incident with $v \in V(H)$.
Note that every graph containing $H$ as a $k$-fat minor also contains $H$ as a minor.
Fat minors enjoy nice properties from geometric viewpoints, and we refer readers to \cite{gp} for details.
In particular, Georgakopoulos and Papasoglu \cite{gp} proposed the following conjecture.

\begin{conjecture}[\cite{gp}] \label{conj_quasi}
For every (finite) graph $H$, there exists a function $f$ such that for every graph or length space $G$, $G$ does not contain $H$ as a $k$-fat minor if and only if $G$ is $f(k)$-quasi-isometric to an $H$-minor free graph.
\end{conjecture}

Fujiwara and Papasoglu \cite{fp_quasi_cacti} proved that Conjecture \ref{conj_quasi} holds when $H=K_4^-$, where $K_4^-$ is the graph obtained from $K_4$ by deleting an edge.
Georgakopoulos and Papasoglu \cite{gp} extended this result to every subdivision of $K_4^-$.
Since Assouad-Nagata dimension is an invariant under quasi-isometry, Theorem \ref{minor_planar_AN_intro} immediately implies the following.

\begin{corollary}
For every integer $k$, if $H$ is a (finite) graph that is a subdivision of $K_4^-$, then the Assouad-Nagata dimension of the class of (unweighted) graphs with no $k$-fat $H$-minor is 1.
\end{corollary}

Moreover, if a weaker version of Conjecture \ref{conj_quasi} that ensures that every graph with no $k$-fat $H$-minor is $f(k)$-quasi-isometric to an $H'$-minor free graph for some graph $H'$ only depending on $H$ and $k$ satisfying that $H'$ is planar when $H$ is planar holds, then Theorems \ref{minor_AN_intro} and \ref{minor_planar_AN_intro} imply the following conjecture.

\begin{conjecture}
Let $H$ be a (finite) graph.
Let $k$ be an integer.
Let $\F$ be the class of (unweighted) graphs with no $k$-fat $H$-minor.
    \begin{enumerate}
        \item Then the Assouad-Nagata dimension of $\F$ is at most 2.
        \item If $H$ is planar, then the Assouad-Nagata dimension of $\F$ is at most 1.
    \end{enumerate}
\end{conjecture}

Note that the above conjecture follows from Conjecture \ref{conj_quasi} and is stronger than another conjecture proposed by Georgakopoulos and Papasoglu in the same paper \cite{gp}.

\begin{conjecture}[\cite{gp}]
For any positive integers $t$ and $k$, the Assouad-Nagata dimension of the class of (unweighted) graphs that do not contain $K_{t+1}$ as a $k$-fat minor is at most $t-1$.
\end{conjecture}

\bigskip

\bigskip

\noindent{\bf Acknowledgement:}
The author thanks Marthe Bonamy, Nicolas Bousquet, Louis Esperet, Carla Groenland, Fran\c{c}ois Pirot and Alex Scott for comments about the first draft of this paper.
He thanks Agelos Georgakopoulos for bringing \cite{gp} and David Wood for bringing \cite{d} to his attention.
He also thanks them and Panos Papasoglu for discussions.
He thanks anonymous reviewers for their diligence, careful reading and suggestions.

\printindex

\end{document}